\documentclass [12pt]{article}
\usepackage{amssymb}
\usepackage[T2A]{fontenc}

\begin{document}

\Large{

\title{
\Large{
      NUMERICAL METHODS  FOR SOLUTION OF SINGULAR
      INTEGRAL EQUATIONS
      }
\author{\normalsize{I. V. Boykov}\\
 Penza State  University,\\ Krasnay Str.,40, Penza,\\
 440026, Russia\\
 ( boikov@pnzgu.ru)}
\date{\today }}
\maketitle

Abstract \dotfill

\marginpar{4}

Preface \dotfill

\marginpar{4}

Introduction \dotfill

\marginpar{10}

1. Classes of Functions \dotfill

\marginpar{10} 

2. Designations and Auxiliary Statements \dotfill

\marginpar{12}

2.1.  Designations of Optimal Algorithms \dotfill

\marginpar{12}

2.2. Elements of Functional Analysis \dotfill

\marginpar{14}

2.3. Elements of Approximation Theory \dotfill

\marginpar{19}

2.4. Inverse Theorems \dotfill

\marginpar{20}

Chapter 1 \dotfill

\marginpar{22}

Approximate Solution of Singular Integral Equations \dotfill

\marginpar{22}

1. An Smoothness of  Solutions of  Singular Integral Equations  \dotfill

\marginpar{22}

1.1. The Integral Operators on the Smooth Functions \dotfill

\marginpar{23}

1.1.1. Fundamental Statements \dotfill

\marginpar{23}

1.2. On Smoothness of  Solutions of  Singular Integral Equations
on Closed Contours \dotfill

\marginpar{24}

 1.2.1. Fundamental Statements \dotfill
 
 \marginpar{24}

2. Approximate Solution of Linear Singular Integral Equations on the
   Closed     Circuits of Integration (Basis in Holder Spaces) \dotfill
   
   \marginpar{24}
   
 2.1. Methods of Collocations and Mechanical Quadratures \dotfill
 
 \marginpar{26}
 
 2.1.1. Fundamental Statements  \dotfill
 
 \marginpar{26}
 
 2.1.2. Proofs of Theorems \dotfill
 
 \marginpar{31}

3. Approximate Solution of Singular Integral Equations on 
Closed Paths of Integration (Basis in Space  $ L_2 $)} \dotfill

\marginpar{52}

3.1. Basis Statements \dotfill

\marginpar{53}

3.2. Proofs of Theorems \dotfill

\marginpar{54}

4.  Approximate Solution of Singular Integral Equations with
   Discontinuous  Coefficients and on Open   Circuits of Integration \dotfill
   
   \marginpar{60}
   
    4.1. Fundamental Statements \dotfill
    
    \marginpar{60}

5. Singular Integral Equations with Constant Coefficients \dotfill

\marginpar{65}

 5.1. Index $ \xi=0$ \dotfill
 
 \marginpar{66} 
 
  5.2. Index $ \xi=1 $ \dotfill
  
  \marginpar{68}

6. Approximate Methods for Solution of Nonlinear Singular Integral Equations \dotfill

\marginpar{70}

 6.1. Projective Methods for Solution of Nonlinear Equations on Closed
 Paths of Integration \dotfill
 
 \marginpar{70}
 
  6.2. Projection Methods for Solving Nonlinear
  Singular Integral Equations
  on Open Contours of Integration \dotfill
  
  \marginpar{86}

   7. Spline-Collocation Method \dotfill
   
   \marginpar{88}
   
    7.1. Linear Singular Integral Equations \dotfill
    
    \marginpar{88}
    
    7.2. Nonlinear Equations \dotfill
    
    \marginpar{95}

8. Singular Integral Equations in Exceptional Cases \dotfill

\marginpar{100}

 8.1. Equations on Closed Circles \dotfill
 
 \marginpar{101}
 
  8.2. Equations on Segments \dotfill
  
  \marginpar{106}
  
   8.3. Approximate Solution of the Equation (8.1) \dotfill
   
   \marginpar{110}
   
    8.4. Approximate Solution of the Equation (8.2) \dotfill
    
    \marginpar{110}
    
 8.5. Approximate Solution of the Equation (8.3) \dotfill
 
 \marginpar{112}

9. Approximate Solution of Singular Integro-Differential Equa-\\tions \dotfill

\marginpar{113}

9.1. Linear Equations \dotfill

\marginpar{113}

 9.1.1. Approximate Solution of the Boundary Value Problem (9.1), (9.2) \dotfill
 
 \marginpar{114}
 
	9.1.2. Approximation  Solution of the Boundary Value Problem (9.3), (9.2) \dotfill
	
	\marginpar{122}
	
 9.3. Approximate Solution of Nonlinear
 Singular Integro-Differen-\\tial 
 Equations on Closed
 Contours of Integration \dotfill
 
 \marginpar{126}	
	 
 9.4. Approximate Solution of Linear
 Singular Integro-Differential Equations with Discontinuous
 Coefficients and on Open Contours of Integration \dotfill
 
 \marginpar{130}
 
 9.5. Approximate Solution of Nonlinear
 Singular Integro-Differen-\\tial Equations on the Open
 Contour of Integration \dotfill
 
 \marginpar{136}

Chapter 2 \dotfill

\marginpar{138}

Approximate Solution of Multi-Dimensional Singular Integral Equations \dotfill

\marginpar{8}

  1. Bisingular Integral Equations \dotfill
  
  \marginpar{139}

2. Riemann Boundary Value Task \dotfill

\marginpar{144}

3. Approximate Solution  of  Multi-Dimensional Singular Integral Equations \dotfill

\marginpar{147}

3.1.  Approximate \  Solution \  of \ Multi-Dimensional \ Singular
Integral Equations on Holder Classes of Functions \dotfill

\marginpar{148}

3.2. Approximate Solution of Linear Multi-Dimensional Singular Integral Equations
on Sobolev Classes of Functions \dotfill

\marginpar{153}

 3.3. Parallel Method for Solution of Multi-Dimensional
 Singular Integral Equations \dotfill
 
 \marginpar{156}
 
  Appendix 1 \dotfill
  
  \marginpar{157}
  
  1. Stability of  Solutions of Numerical Schemes \dotfill
  
  \marginpar{157}
  
  1.1. Stability of Solutions of Bisingular Integral Equation \dotfill
  
  \marginpar{158}
  
  1.2.  Stability of Solutions of Multi-Dimensional Singular Integral Equations \dotfill
  
  \marginpar{161}
  
   1.3. Stability \  of \  Solutions \  of \  Nonlinear Singular\  Integral Equations \dotfill
   
   \marginpar{162}

References \dotfill

\marginpar{165}

{\bf Abstract}

This paper is devoted to an overview of the authors own works for \ numerical  \ solution\  of \ singular integral equations (SIE), polysingular integral
equations and multi-dimensional singular integral \  equations \  of the\  second \ 
kind. \ Considered iterative - projective methods and parallel methods for
solution of
      singular integral\  equations,\  polysingular integral
equations \  and multi- dimensional singular integral equations.
The paper is the second part of the overview of the authors own works devoted to numerical methods for calculation singular and hypersingular integrals \cite{Boy23} and to approximate methods for solution of singular integral equations.


\begin{center}
{\bf Preface}
\end{center}

Hilbert and Poincare created the theory of singular integral equations early in the 20th century, which has undergone an intense growth during the last 100 years.
 The theory is associated with numerous applications of singular  integral equations, as well as with Riemann boundary value problem.  
 The Riemann boundary value problem and singular integral equations are widely used as basic techniques of mathematical modeling in physics 
(quantum field theory \cite{Bogol}, theory of close and long-distance inter-\\action \cite{Br},\  soliton theory \cite{Fad},  theory of elasticity and
 thermo-\\elasticity \cite{Ioa},
 aerodynamics \cite{Lif}, and electrodynamics \cite{Lif1}, etc.

 A closed-form solution of singular  integral equations is only possible in exceptional cases.  A comprehensive presentation and an extensive 
literature survey associated with all methods of solution of singular integral equations of the first and second kinds are found in \cite{Boy25}, \cite{Dzi},
 \cite{Dzi1}, \cite{Dzi2},   \cite{Ell}, \cite{Ell1}, \cite{Goh}, \cite{Gol2},  \cite{Ivan1},  \cite{Lif},  \cite{Mich},  \cite{Pr3},  \cite{Siz},
 \cite{Sr} \cite{Tsa}, \cite{Zh}. 

Various numerical methods for solution of singular integral
equations have been the topic of a great many of papers, most of which can be
divided in two fields.
The first field is devoted to numerical methods for solution of singular
integral equations of the first kind. The second field is devoted to numerical
methods for solution of singular integral equations of the second kind.

The early works belong to the first field.
Singular integral equations are the main tool for simulation aerodynamics
tasks. M.A. Lavrent'ev paper an theory of aerofoil \cite{Lavr} was the first work
devoted to numerical methods for solution of the singular integral equations
of the first kind. Later  numerical methods for solution of the first kind singular integral
equations was devoted many papers based on different methods. A. I. Kalandia
proposed the method of mechanical quadrature  \cite{Kal}, F. Erdogan and G.D. Gupta 
proposed the Gauss-Chybyshev method  \cite{Erd}, N.I. Ioakimidis and P.S. Theocaris
proposed the Lobatto-Chebyshev method  \cite{Ioa1}, 
S. Krenk proposed Gauss-Jacobi,
Lobatto-Jacobi methods  \cite{Kr}. Discrete wortexes method for solution of the singular
integral equations of the first kind was investigated by S.M. Belotserkovsky
and I.K. Lifanov  \cite{Bel}, I.K. Lifanov  \cite{Lif} , N.F. Vorob'ev  \cite{Vor} .
A.V. Dzhishkariani was
printed the detailed analyses of numerical methods for solution of the
singular integral equations of the first kind \cite{Dzi}, \cite{Dzi1}.

The numerical methods for solution of the singular integral equations of the
second kind  we can divide on two
parts: direct
and indirect methods.
For different cases of singular integral equations \ \ \  Du \  Jinyuan \  \ \cite{Du}, \  \cite{Du1}, \  \cite{Du2}, \ \  A. \   Gerasolus \ \ \cite{Ger},  J.G. Graham \cite{Grah}, N.L. Gori \cite{Gori},
D. Elliot \cite{Ell}, \cite{Ell1}, \cite{Ell2}, \cite{Ell3} , B.I. Musaev \cite{Musa}, J.G. Sanikidze \cite{San},  I. Ioakimidis \cite{Ioa}, \cite{Ioa1}, \  M.A. Sheshko \ \cite{She1}, \ 
R.P. Srivastav and E Jen \cite{Sr},   G. Tsamasphyras and P.S. Theocaris \cite{Tsa}
constructed various
numerical schemes of indirect methods and proved its convergence.

We shall not touch the part which are devoted to indirect methods
for solution of singular integral equations. One of these methods consists in
transformation of singular integral equations 
\[
a(t)x(t) + \frac{b(t)}{\pi i}\int\limits^1_{-1}
\frac{x(\tau)}{\tau-t}d\tau + \int\limits^1_{-1}
h(t,\tau)x(\tau)d\tau = f(t)
\]
to the equivalent Fredholm integral equations
   \[
x(t) + \int\limits^1_{-1} h^*(t,\tau)x(\tau)d\tau =f^*(t).
\]
Numerical methods are applied to the last equations.

One can become acquainted with this method by papers\\
G. Tsamasphyras and P.S. Theocaris \cite{Tsa}, M.A. Sheshko \cite{She}.

For the first time direct methods for solution of SIE as
$$
a(t)x(t) + \frac{b(t)}{\pi i}\int\limits_{\gamma}
\frac{x(\tau)d\tau}{\tau-t}+ \int\limits_{\gamma}
h(t,\tau)x(\tau)d\tau = f(t),
\eqno (1.1)
$$
where $\gamma= \{z, z\in C, |z|=1\},$
was considered by V.V. Ivanov \cite{Ivan1}, \ 
\cite{Ivan}. He \  proved the convergence of collocation
method, moments method, Galerkin-Petrov method for the equation (1.1).

For solution of SIE as
$$
a(t)x(t) + b(t)\int\limits^1_{-1}
\frac{x(\tau)d\tau}{\tau-t} + \int\limits^1_{-1}
h(t,\tau)x(\tau) = f(t)
\eqno (1.2)
$$
V.V. Ivanov proposed a method for transformation  the equation
(1.2) to the
equation (1.1).

For solution of the equation (1.1) V.V. Ivanov used the method which is based
on connection between singular integral equations and Riemann boundary value
problem. According this method the singular integral equation is transformed
to the Riemann boundary value problem. Simultaneously, \ numerical \ scheme for solution of
the singular integral equation is transformed to numerical scheme for
approximate solution of Riemann boundary value problem. Using the Kantorovich common
theory of approxi-\\mate methods of analysis \cite{Kan} V.V. Ivanov proved existence and
unique \ of approximate \ solution of Riemann boundary value problem and, as
corollary, existence and unique of approximate solution of a singular integral
equation. In this way he investigated collocation method, moments method,
Bubnov - Galerkin method,   least square method for one-dimensional singular
integral equations. Convergence of these methods are given in the spaces $W$
and $L_2.$ Historical this method was the first common method for solution of
singular integral equations as (1.1).

Development of this method make up the first direction in numerical methods
for solution of second kind singular integral equations.

In the frame of this direction was received following results.

B.G. Gabdulhaev \cite{Gab}  have proved  convergence of the mechanical quadrature method for
the equation (1.1) in the space $H_{\beta}.$

I.V. Boykov \cite{Boy1},  \cite{Boy5},  \cite{Boy7},  \cite{Boy25} offered modifications of collocation  and mechanical
quadrature methods for solution of linear singular integral equations 
$$
a(t)x(t) +\frac{1}{\pi i}\int\limits_\gamma
\frac{h(t,\tau)x(\tau)}{\tau-t}d\tau = f(t)
\eqno (1.3)
$$
and nonlinear singular integral equations 
$$
a(t,x(t)) +\frac{1}{\pi i}\int\limits_\gamma
\frac{h(t,\tau,x(\tau))}{\tau-t}d\tau = f(t).
\eqno (1.4)
$$
Convergence of these methods are given in $H_\beta$ and $L_2.$

 I.V. Boykov and I.I.  Zhechev \cite{Boy32}, \cite{Boy33}, \cite{Boy34}
proved convergence of numerical methods  for solution of 
singular integro - differential integral equations
\[
\sum\limits^m_{k=0} \left[a_k(t)x^{(k)}(t) + \frac{b_k(t)}{\pi i}
\int\limits_\gamma \frac{x^{(k)}(\tau)}{\tau-t}d\tau +
\frac{1}{2\pi i}\int\limits_\gamma
h_k(t,\tau)x^{(k)}(\tau)d\tau \right] =
\]
$$
= f(t)
\eqno (1.5)
$$
with boundary dates
$$
\int\limits_\gamma x(\tau)\tau^{-k-1}d\tau=0,
\eqno (1.6)
$$
$k=0,1,\ldots,m-1.$

Approximate solution of nonlinear singular
integro - differential equations as
\[
a(t,x(t),x'(t),\ldots,x^{(m)}(t)) +
\]
\[
\frac{1}{\pi i}\int\limits_\gamma
\frac{h(t,\tau,x(\tau),x'(\tau),\ldots,x^m(\tau))}{\tau-t}d\tau
= f(t)
\]
under boundary dates (1.6) with method of mechanical quadrature was investigated by I.V. Boykov and I.I.  Zhechev \cite{Boy35}, \cite{Boy36}.

Collocation method and method of mechanical quadrature for system of linear
and nonlinear singular integral equations was investigated by I.I.  Zhechev
\cite{Zh}.

Projective methods for solution of singular integral equations as (1.2) in
weight spaces was considered by I.V. Boykov \cite{Boy8}, \cite{Boy14}.

Some results, received for solution of one-dimensional singular integral equations, was
diffused by I.V. Boykov \cite{Boy16}, \cite{Boy24}, \cite{Boy25} to bisingular integral equations as
\[
a(t_1,t_2)x(t_1,t_2) + \frac{b(t_1,t_2)}{\pi i}\int\limits_{\gamma_1}
\frac{x(\tau_1,t_2)}{\tau_1-t_1}d\tau_1 +
\]
\[
+ \frac{c(t_1,t_2)}{\pi i}\int\limits_{\gamma_2}
\frac{x(t_1,\tau_2)}{\tau_2-t_2}d\tau_2 +
\frac{d(t_1,t_2)}{\pi^2 }\int\limits_{\gamma_1} \int\limits_{\gamma_2}
\frac{x(\tau_1,\tau_2)d\tau_1 d\tau_2}{(\tau_1-t_1)(\tau_2-t_2)} +
\]
\[
+ \int\limits_{\gamma_1} \int\limits_{\gamma_2}
h(t_1,t_2,\tau_1,\tau_2)x(\tau_1,\tau_2)d\tau_1 d\tau_2=f(t_1,t_2),
\]
where $\gamma_i = \{z_i, z_i \in C, |z_i|=1, i=1,2\}$.

Some results, which are received in the first direction, are printed in the
books \cite{Boy16}, \cite{Boy25}.

The second direction in approximate methods for solution of singular integral
equation is connected with common projective methods for solution of integral
equations with convolution.

 I.Tz. Gohberg and I.A. Feldman \cite{Goh} offered common
method for investigation \ of \  integral\   equations \ with convolution. The method
is based on the theory of rings.

Later this method was diffused to singular integral equations of non normal type
and polysingular integral equations. Basic results in this direction was
printed in the books \cite{Pr}, \cite{Mich}, \cite{Pr3}.

The third direction in approximate methods for solution of singular integral
equations is spline - collocation methods. 

Apparently, the first   works in this \
direction \ was the papers of S. Pressdorf \cite{Pr1},
S. Pressdorf and G. Shmidt \cite{Pr2}. In these papers
verification of spline-projective method was based on study of the
spectrum of the linear system of algebraic equations, which are approximated
the given singular integral equation.

During the last twenty five years spline-collocation methods was diffused on
many types of singular and polysingular integral equations.

Basic results of this direction was printed in the works \cite{Gori}, \cite{Grah}, \cite{Mich}, \cite{Pr3}.

For solution of aerodynamics problems  discrete vortex method was offered by
S.M. Belotserkovsky \cite{Bel}. This method was diffused to many types of
singular and polysingular integral equations   which are used as basic techniques of mathematical modeling in aerodynamics \cite{Lif1}.

Recently the wavelet collocation method and moment method, based on wavelet
functions, begin the rapid development \cite{Bey}, \cite{Blat}.

Apparently approximate methods for solution of singular and  polysingular integral equations,
based on theory of wavelets, will form the fifth direction in numerical
methods for solution of singular integral equations.

For solution of the equation (1.2) was offered some   special methods which we
can not associate to concrete direction.
So direct methods for solution of SIE as (1.2) given by D. Elliot in the seria of
important papers \cite{Ell} - \cite{Ell3}.

    In this paper we do not give review of all directions. We only give
    the short review of some results in the first direction of numerical methods
for solution of singular integral equations, polysingular integral equations
and multi-dimensional singular integral equations of the second kind which was received
author and his disciples.

{\it Note.} The author's name in various publications was translated from Russian to English as Boikov or Boykov.

The author's work in numerical methods for solution of singular integral equations was supported by Russian Foundation for Basic Research (grant 16-01-00594).

\begin{center}
 {\bf  Introduction}
\vskip 10 pt

 {\bf 1. Classes of Functions}
\vskip 10 pt
\end{center}

In this section we will list several classes of functions, which will be used later.

Let $\gamma$ be the unit circle: $\gamma = \{z:|z|=1\}.$

To measure the continuity of a function $f \in C[a,b]$ we proceed as follows \cite{Lor}. We consider the first difference with step $h$ 
$$
\Delta_h(f(x)) = f(x+h) - f(x)
$$
of the function $f$ and put
$$
\omega(f,\delta) = \omega(\delta) = \max\limits_{x,h  (|h| \leq \delta)}|f(x+h) - f(x)|.
$$
The function $\omega(\delta),$ called the modulus of continuity of $f,$ is defined for $0 \leq \delta \leq b-a.$

{\bf Definition 1.1.} A function $f,$ defined on $\Delta = [a,b]$ or $\Delta = \gamma,$ satisfies a Lipschits condition with constant $M$ 
and exponent $\alpha,$ or belongs to the class $H_\alpha(M),$ $M \geq 0,$ $0 < \alpha \leq 1,$ if 
$$
|f(x') - f(x'')| \leq M|x'-x''|^\alpha, \, x', x'' \in \Delta.
$$

{\bf Definition 1.2.} A function $f,$ defined on $\Delta = [a,b]$ or $\Delta = \gamma,$ satisfies a Zigmund condition with constant $M$,  or belongs to the class $Z(M),$ $M \geq 0,$  if 
$$
|f(x') - f(x'')| \leq M|x'-x''||\ln |x'-x''||, \, x', x'' \in \Delta.
$$

{\bf Definition 1.3.} The class $W^r(M,\Delta),$ $r=1,2,\ldots,$ $\Delta = [a,b]$ or $\Delta = \gamma,$ consists of all functions $f \in C(\Delta),$ 
which have an absolutely continuous derivative $f^{(r-1)}(x)$ and piecewise derivative $f^{(r)}(x)$ with $|f^{(r)}(x)| \le M.$
  
{\bf Definition 1.4 } \cite{Lor}. \ Let  $r=0,1,\ldots,$  $M_i \geq 0,$  $i=0,1,\ldots,\\ r+1,$ let $\omega$ be a modulus of continuity and let $\Delta=[a,b]$ 
or $\Delta= \gamma.$ Then $W^r_\omega = W^r_\omega(M_0,\ldots, M_{r+1}; \Delta)$ is the set of all functions $f \in C(\Delta),$ which have continuous
 derivatives $f, f',\ldots, f^{(r)}$ on $\Delta,$ satisfying
$$
|f^{(i)}(x)| \leq M_i, \, x \in \Delta, \, i=0,1,\ldots,r, \, \omega(f^{(r)}, \delta) \leq M_{r+1}\omega(\delta).
$$
They write $W^r H_\alpha$ if $\omega(\delta) = \delta^\alpha,$ $0 < \alpha \leq 1.$

Let us consider functions $f(x_1,\ldots,x_l)$ of $l$ variables on $\Delta,$ where $\Delta$ is either an $l$-dimensional parallelepiped (that is, the product of $l$
 intervals $a_k \leq x_k \leq b_k,$ $k=1,2,\ldots,l)$ or an $l$-dimensional torus $\gamma^{(l)},$ the product of $l$ circles $\gamma.$ The modulus of continuity 
$\omega(f, \delta)$ are defined as the maximum of 
$$|f(y_1,\ldots,y_l) - f(x_1,\ldots,x_l)| \ {\rm for \ }
|y_k-x_k| \leq \delta, k=1,2,\ldots,l.
\eqno (1.1)
$$

Sometimes the partial modulus of continuity $\omega^{(k)}(f,\delta)$ of $f$ are used. The $ kth$ modulus among them is the maximum of (1.1) when the increment is
 only with respect to the $kth$ coordinate: $x_i = y_i,$ $i \ne k,$ $|y_k - x_k| \leq \delta.$

In \cite{Lor} was defined classes $W^r_{l \omega} = W^r_{l \omega}(M_0,\ldots, M_{r+1}; \Delta)$ of functions $f(x_1,\ldots,x_l)$ of $l$ variables
 on $l$-dimensional set $\Delta,$ which is either a parallelepiped or a torus. 

{\bf Definition 1.5} \cite{Lor}. A function $f$ belongs to $W^r_{l \omega}$ if and only if all its partial derivatives $D^jf$ of 
order $j=0,1,\ldots,r$ exist and continuous and satisfy the following conditions:
 For each partial derivative $D^j f$ of order $j,$ $\|D^j f\| \leq M_j,$ $j=0,1,\ldots,r,$ and in addition for each derivative of order 
$r,$ $\omega(D^r f, \delta) \leq M_{r+1} \omega(\delta).$

They write $W^r_l H_\alpha$ if $\omega(\delta) = \delta^\alpha,$ $0 < \alpha \leq 1.$

If coefficients $A$ and $M$ are not essential we use designations
$H_{\alpha},$ \ \  $W^r H_{\alpha},$ \ \ $W^{rr} H_{\alpha \alpha}$ \ \ instead of \
$H_{\alpha}(A),$\ \ $W^r H_{\alpha}(A,M),$ \\ $W^{rr} H_{\alpha \alpha}(A,M)$
respectively.

\begin{center}
{\bf 2. Designations and Auxiliary Statements}
\end{center}

{\bf 2.1.  Designations of Optimal Algorithms}

In this paper
we will use definitions of optimal algorithms for solution of problems of mathematical physics, given N.S. Bakhvalov \cite{Bakh}.
These definitions we will use in the study of optimal algorithms for solving singular integral equations.

Let $H=\{h \}$ be a class of smooth functions. Let $F=\{f \}$
be a class of smooth functions. Let $C=\{c \}$ be a class of smooth
functions.

Let $\Psi$ be a class of vector-functionals, defined on $H.$
Let $\Psi^{*}$ be a class of vector-functionals, defined on $C.$
Let $\Psi^{**}$ be a class of vector-functionals, defined on $F.$

Let $M$ be a set of Markov algorithms.\\
Let $R=R(K,h,c,f,A,\{\psi_v\}_1^{N^2},
\{\psi^*_v\}_1^{N},\{\psi^{**}_v\}_1^{N},t)$
be a result of numerical solution of the equation
$$
Kx \equiv c(t)x(t) + \frac{1}{\pi i}\int\limits_\gamma
\frac{h(t,\tau)x(\tau)d\tau}{\tau-t} = f(t)
\eqno (2.1)
$$
with a algorith $A \in M.$

Algorithm $A$ uses  $N^2$ functionals $\psi_v(h),$
$v=1,2,\ldots,N^2,$ $\psi_v \in \Psi,$  $N$ functionals
$\psi_v^{**}(f),$ $v=1,2,\ldots,N,$ $\psi_v^{**} \in \Psi^{**},$ and
$N$
functionals $\psi_v^{*}(c),$ $v=1,2,\ldots,N,$ $\psi_v^{*} \in \Psi^{*}.$

Let us
introduce the following designations:

\[
v(K,h,c,f,A,\{\psi_v\}^{N^2}_1,\{\psi^*_v\}^{N}_1 ,
\{\psi^{**}_v\}^{N}_1   )=\rho(x^*,R),
\]
\[
v(K,H,C,F,A,\{\psi_v\}^{N^2}_1,\{\psi^*_v\}^{N}_1,\{\psi^{**}_v\}^{N}_1  )=
\]
\[
=\sup_{c\in C,f \in F, h\in H} v(K,h,c,f,A,
\{\psi_v\}^{N^2}_1, \{\psi^*_v\}^N_1),\{\psi^{**}_v\}^{N}_1),
\]
\[
v(K,H,C,F,M,\{\psi_v\}^{N^2}_1 ,\{\psi^*_v\}^N_1,\{\psi^{**}_v\}^{N}_1   )=
\]
\[
=\inf_{A \in M} v(K,H,C,F,A,\{\psi_v\}^{N^2}_1,\{\psi^*_v\}^N_1),\{\psi^{**}_v\}^{N}_1,
\]
\[
v_N(K,H,C,F,M,\Psi ,\Psi^* ,\Psi^{**}  )=
\]
\[
=\inf_{\{\psi_v\}^{N^2}_1 \in \Psi,  \{\psi^*_v\}^{N}_1 \in \Psi^*
,\{\psi^{**}_v\}^{N}_1 \in \Psi^{**}  }
v(K,H,C,F,M,\{\psi_v\}^{N^2}_1  ,\{\psi^*_v\}^{N}_1
,\{\psi^{**}_v\}^{N}_1    ),
\]
\[
v(K,H,C,F,\{\psi_v\}^{N^2}_1 ,\{\psi^*_v\}^{N}_1 ,\{\psi^{**}_v\}^{N}_1  )=
\]
\[
=\inf_{A}v(K,H,C,F,A,\{\psi_v\}^{N^2}_1 ,\{\psi^*_v\}^{N}_1
,\{\psi^{**}_v\}^{N}_1    ),
\]
\[
v_N(K,H,C.F,\Psi ,\Psi^*  ,\Psi^{**} )=
\]
\[
=\inf_{\{\psi_v\}^{N^2}_1 \in \Psi,\{\psi^*_v\}^{N}_1 \in \Psi^*
,\{\psi^{**}_v\}^{N}_1 \in \Psi^{**}}
v(K,H,C,F,\{\psi_v\}^{N^2}_1,\{\psi^*_v\}^N_1),\{\psi^{**}_v\}^N_1),
\]
\[
v_N(K,H,C,F)=
\]
\[
=\inf_{\{\psi_v\}^{N^2}_1,
\{\psi^*_v\}^N_1,\{\psi^{**}_v\}^N_1}
v(K,H,C,F,\{\psi_v\}^{N^2}_1,\{\psi^*_v\}^N_1),\{\psi^{**}_v\}^N_1).
\]

Here $\rho(x^*,R)$ is the distance between
exact solution $x^*$  of the equation $Kx=f$ and approximate solution $R$ of the equation (2.1).

Functional\ \ $v_N(K,H,C,F,\Psi,\Psi^*,\Psi^{**})$ \ is\  infimum of values\\
$v_N(K,H,C,F,\{\psi_v(h)\}_1^{N^2}$, $\{\psi_v^*(c)\}_1^{N}$, $
\{\psi_v^{**}(f)\}_1^{N})$ on functionals $\Psi,\Psi^*,\Psi^{**}.$

Functional $v_N(K,H,C,F)$ is infimum of values
\[
v_N(K,H,C,F,\Psi,\Psi^*,\Psi^{**})
\] using all kind of
functional classes $\Psi,\Psi^*,\Psi^{**}.$

Algorithm  $A,$  using  functionals
$\{\bar \psi_v\}_1^{N^2} \in \Psi,$ $\{\bar \psi_v^*\}_1^{N} \in \Psi^*,$
$\{\bar \psi_v^{**}\} \in \Psi^{**},$ called optimal, asymptotically optimal,
optimal with respect to order if
$$
\frac{v(K,H,C,F,A,\{\overline\psi_v\}_1^{N^2},\{\overline\psi_v^*\}_1^N
,\{\overline\psi_v^{**}\}_1^N)}
{v_N(K,H,C,F,\Psi,\Psi^*,\Psi^{**})}=1,\sim 1,\asymp 1. $$

{\bf 2.2. Elements of Functional Analysis
	}

Let us remind some statements of functional analysis.

Let $X$ is a normed space. Let $K$ is a linear bounded operator from a normed
space $X$ to a normed space $Y.$ This fact we will write as $K \in [X,Y].$

{\bf Banach Theorem} \cite{Lus}. Let $B$ be a Banach space. Let $A \in [B,B]$ with a norm
$\|A\| = q < 1.$ In this case the operator $K=I+A$ has the linear inverse
operator $K^{-1}$ with the norm $\|K^{-1}\| \le 1/(1-q).$

{\bf Generalized Banach Theorem} \cite{Lus}. Let $A$ and $B$ are linear bounded operators from
Banach space $X$ to Banach space $Y.$ Let the operator $A$ has the inverse
operator $A^{-1} \in [Y,X].$ Let $q=\|A-B\|.$ If $q\|A^{-1}\| < 1,$ then the
operator $B$ has the linear bounded operator $B^{-1} \in [Y,X]$ and the
inequality
$$
\|A^{-1} - B^{-1}\| \le \frac{q\|A^{-1}\|^2}{1-\|A^{-1}\|q}
$$
is valid.

{\it Note.} From the proof of the Generalized Banach Theorem follows the
inequality
$$
\|B^{-1}\| \le \|A^{-1}\|/(1-q).
$$

We will need in Kantorovich theory for approximate methods of analysis \cite{Kan}.

Let $X$ be a Banach space and let $X_n$ be a closed subspace of $X.$

Let us consider two equations:
$$
Kx \equiv x +Hx =f, \quad K \in [X,X]
\eqno (2.2)
$$
and
$$
K_nx_n \equiv x_n +H_nx_n = f_n, \quad K_n \in [X_n,X_n].
\eqno (2.3)
$$

The equation (2.2) is a given equation. The equation (2.3)  is a
 approximate  equation.

Let the operator $K$ has the inverse bounded operator $K^{-1}.$ Let $P_n$ is a
projector from $X$ on $X_n.$ Let the operators $K$ and $K_n$ are connected by
following conditions:\\
I.  For each $x \in X$ exists such $x_n \in X_n,$ that
\[
\|Hx - x_n\| \le \eta_1(n)\|x\|;
\]
II. For each $x_n \in X_n$
\[
\|P_n Hx - H_n x_n\| \le \eta_2(n)\|x_n\|;
\]
III. For each $x \in X$ exists such $x_n \in X_n$ such that
\[
\|x-x_n\| \le \eta_3(n,x)\|x\|.
\]
{\it Note.} In the conditions I and II constants $\eta_1$,  $\eta_2$ are
independent from $x$ or $x_n.$

Under conditions I-III the following statements is valid.

{\bf Theorem 2.1} \cite{Kan}. Let the operator $K$ is inverse. Let the conditions I and II
are valid. If
\[
q=[\eta_2(n) + \|I-P_n\|\eta_2(n)]\|K^{-1}\| <1,
\]
then the operator $K_n$ has the inverse operator $K_n^{-1}$ with the norm
\[
\|K_n^{-1}\| \le \|K^{-1}\|/(1-q).
\]

{\bf Theorem 2.2} \cite{Kan}. Let the conditions of the Theorem 1 are fulfilled. Let the
condition III is valid. The equation (2.3) has a unique solution $x_n^*$ and have
place the inequality
\[
\|x^* - x^*_n\| \le 2\eta_2(n) \|K_n^{-1}\| +(\eta_1(n) + \eta_3(x^*)\|K\|
(1+ \|K_n^{-1} P_n K\|),
\]
where $x^*$ is a unique solution of the equation (2.2).

During the paper we will use the following well known Hadamard Theorem.

{\bf Hadamard Theorem} \cite{Gant}. Let $|c_{jj}|>\sum\limits^n_{k=0, k \ne j}|c_{jk}|$
for $j=1,2,\ldots,n.$ Then the system $Cx=b,$ where $C=\{c_{ij}\},$ $i,j=1,2,\ldots,n,$
$x=(x_1,\ldots,x_n)^{T},$ $b=(b_1,\ldots,b_n)^{T},$ has a unique solution.

The convergence of numerical methods for solution nonlinear singular integral equations   follows from convergence theorems for
Newton $-$ Kantorovich's method in Banach spaces. 
Let us recall convergence criteria for the Newton $-$ Kantorovich method.

Let $X,Y$ be Banach spaces.
Consider the equation
$$
Kx=0, 
\eqno (2.4)
$$
where $K$ is a nonlinear operator  from $X$ into  $Y$.
Let the operator $K$ has the Frechet derivative in a neighborhood of  initial point
$x_0$. Assume that  $[K'(x)]^{-1}$, or in a general case, the right inverse operator
$[K'(x_0)]^{-1}_r$ exists.

We will seek for a solution of the equation (2.4)  in the form of the following iteration processes -- the basic
$$
x_{n+1}=x_n - [K'(x_n)]^{-1}K(x_n)  
\eqno (2.5)
$$
and the modified 
$$
x_{n+1}=x_n - [K'(x_0)]^{-1}K(x_n).
\eqno (2.6)
$$

If an operator $Kx$ has the second order  Frechet 
derivative, the following assertion  is used to justifiy the 
iteration processes (2.5)  and (2.6).

	{\bf Theorem 2.3}  \cite{Kan}. Let \  the \  operator \  $K$ is defined in $\Omega$ \ \ \  $(\|x-x_0\|<R) $ and has  continuous
	second order derivative in $\Omega_0$ \ \ \  $(\|x-x_0\|\le r, r< R)$. Let\\
		i) a countinuous linear operator  $\Gamma_0=[K{}'(x_0)]^{-1}$ exists; 
		\\
		ii) $\|\Gamma_0K(x_0)\|\leq \eta ;$\\
		iii) $\|\Gamma_0K{}''(x)\|\leq b_0 \ \ \   (x\in \Omega_0)$.\\
	Then, if 
	$
	h=b_0\eta\leq 1/2
	$
	and
	$
	r\geq r_0=(1-\sqrt{1-2h})\eta/h,
	$
	the Newton methods (basic and modified) converge to a solution $x^*$ of
	the equation (2.4). Moreover,
	$
	\|x^*-x_0 \|\leq r_0.
	$\\
	If for $h<1/2$\ \ \ 
	$
	r<r_1=(1+\sqrt{1-2h})\eta/h,
	$
	and for $h=1/2$\ \ \ 
	$
	r\leq r_1,
	$
	then a solution $x^*$ is unique in the ball $\Omega_0$.
	
	The convergence speed for the basic method is characterized by
	$$
	\|x^*-x_n\| \leq \frac{1}{2^n}(2h)^{2^{n}}\frac{\eta}{h} \ \ \ (n=0,1,...),
	$$
	the convergence speed for the modified method for $h<1/2$ is characterized by 
	$$
	\|x^*-x{}'_n\|\leq \frac{\eta}{h}(1-\sqrt{1-2h})^{n+1} \ \ \ (n=0,1,...).
	$$

If the  Frechet derivative of $Kx$ satisfies the Lipschitz
condition, convergence of the modified Newton - Kantorovich 
method follows from the following

{\bf Theorem 2.4} \cite{Kras}.  Let  an operator $K$ is defined and 
	Frechet differentiable on a ball 
	$\Omega(x_0,R) \ \ \  (\|x-x_0\|<R)$, and its derivative $K'(x)$ satisfies
	on $\Omega(x_0,R)$ the Lipschitz condition $\|K'(x)-K'(y)\| \le 
	L\|x-y\|$. Let a linear bounded operator 
	$\Gamma_0=[K'(x_0)]^{-1}$ exists, and $\| \Gamma_0 \| \le b_0,
	\|$ $\Gamma_0 K(x_0) \| \le \eta_0$. Let 
	$h_0=b_0 L \eta_0 <1/2, r_0 = (1-\sqrt{1-2 h_0}) \eta_0 /h_0 
	\le R$. Then successive approximations  (2.5) converge to a solution 
	$x^{*} \in \Omega $ of the equation $K(x)=0$.

In the more common case we have the following statement.

{\bf Theorem 2.5} \cite{Boy25}. Let $X$ and $Y$ be Banach spaces. Suppose the following
	conditions are fulfilled:\\
		i) $\Bigl\| K(x_0) \Bigr\| \equiv \eta_0;$\\
	ii) The operator $K$ has the Frechet derivative in  neighbourhood of  initial point
		$x_0$;\\
	iii)		 There exists a right reverse operator $[K'(x_0)]_r^{-1}$ with the norm
	$$
		\Bigl\| [K'(x_0)]_r^{-1} \Bigr\| = B_0; 
	\eqno (2.7) 
		$$
	iiii)	$$			
		\Bigl\| K'(x_1)-K'(x_2) \Bigr\| \le q/(B_0(1+q)) 
	\eqno (2.8)
	$$
	in the sphere $S\{x:\| x-x_0 \| \le \frac{B_0 \eta_0}{1-q} \} \quad (q<1)$.
	
	Then,  equation (2.4) has a solution $x^*$ in  $S$. 
	A sequence
$$
	x_{n+1}=x_n-[K'(x_n)]^{-1}_r Kx_n
\eqno (2.9)
$$
	converges to $x^*$. The estimate
	$ \| x^*-x_n \| \le q^n \eta_0 B_0/(1-q)$ is valid.
	
	{\it Note.} Let a solution $x^*$ of   (2.4) enters 
	the domain $S_0= S \cap \Delta  R(x)$.  Here  $\Delta
	R(x)$ is a domain in which $\| [K'(x)]_r^{-1} \| \le B_0$, and $ \| K'(x_1)-K'(x_2) \Bigr\| \le q/(B_0(1+q)),  x_1, x_2 \in \Delta  R(x) $ are valid.  Then, 
	$x^*$ is a unique solution in $S_0$.

Let us solve the equation  (2.4) with the modified Newton -- Kantorovich method.

{\bf Theorem 2.6} \cite{Boy25}. Let  $X$ and $Y$ be Banach spaces. Suppose the following conditions
	are fulfilled:\\
		i) $\Bigl\| K(x_0) \Bigr\| \equiv \eta_0;$\\
	ii) The operator $K$ has the Frechet  derivative in a neighborhood of a initial point 
		$x_0$, and there exists the right inverse operator $[K'(x_0)]_r^{-1}$ with the norm
		$\Bigl\| [K'(x_0)]_r^{-1} \Bigr\| = B_0;$\\ 
	iii) The condition $ \Bigl\| K'(x_1)-K'(x_2) \Bigr\| \le q/B_0$ is fulfilled 
		in the sphere  $S\{x:\| x-x_0 \| \le \frac{B_0 \eta_0}{1-q} \} \quad (q<1)$.

	Then,  equation (2.4) has a  solution $x^*$ in $S$.
	A sequence $x_{n+1}=x_n-[K'(x_0)]^{-1}_r Kx_n$ converges to $x^*$.
	The estimate  $ \| x^*-x_n \| \le q^n \eta_0 B_0/(1-q)$ is valid. 
	A solution $x^*$ is unique in $S \cap (\Delta R(x))$.

Let us preseent the proof of the Theorem 2.6.
First recall the Rakovchik Lemma.

	{\bf Lemma 2.1} \cite{Rak}.
	Let $X$ and $Y$ be Banach spaces, $A$
	and $B$ be
	linear bounded operators mapping $X$ into $Y$. 
	If $A$ has a bounded right inverse operator   
	$A_r^{-1}$ and  if an operator $B$ satisfies 
	$\|A-B\| \|A_r^{-1}\| <1$, it has a bounded right inverse 
	operator $B_r^{-1}$,  $\|B_r^{-1}\| \le \|A_r^{-1}\| /
	(1-\|A-B \|\|A_r^{-1}\|)$.

{\bf   Proof of the Theorem 2.6.} First we prove that a uniformly bounded right inverse operator $[K'(x_n)]_r^{-1}$ exists  in the domain $S$.
Indeed, from Rakovchik's lemma it follows 
$$
\| [K'(x_n)]_r^{-1} \Bigr\|_Y \le (1+q)B_0.  
\eqno (2.9) 
$$

Show that all the approximations obtained by the iterative process (2.6) are in $S$.

Clearly,
$ \| x_1-x_0 \|_X = \Bigl\| [K'(x_0)]_r^{-1}K(x_0)
\Bigr\|_X \le B_0 \eta_0 <  B_0\eta_0/(1-q)$.
So. $x_1 \in S$.

Let $x_m \in  S$ for  $m  \le  n$ have already been proved. 
Since  $x_n-x_{n-1}= [K'(x_{n-1})]_r^{-1}K(x_{n-1})$,  then by Lemma  we have
\[
\| x_{n+1}-x_n \|_X \le
\]
\[ \le \Bigl\|  [K'(x_n)]_r^{-1}[K(x_n)-K(x_{n-1})-K'(x_{n-1})
(x_n-x_{n-1})] \Bigr\|_X \le
\]
$$
\le q \| x_n-x_{n-1} \|_X. 
\eqno (2.10)
$$

Therefore
$ \| x_{n+1}-x_0 \|_X \le \sum\limits_{k=0}^n q^k B_0 \eta_0$,
i.e. $ x_{n+1} \in S$.

From  (2.10)  it follows that     $\{x_n\}$
is a fundamental sequence. Hence  
$x^*=\lim x_n$ exists.   Since $K(x_n)=-K'(x_n)(x_{n+1}-x_n)$, then $K(x^*)=0$. 
It follows from (2.10)  that 
\[
 \|  x^*-x_n  \|_X  \le  \sum_{k=n}^\infty  q^k  B_0
\eta_0 \le \frac{q^n B_0 \eta_0}{1-q}.
\]

The theorem is proved.

{\bf 2.3. Elements of Approximation Theory}

We  give some well known results from  theory of approximation functions of real
variable.

{\bf Theorem 2.7} \cite{Nat}.
If a trigonometric polynomial $T_n(x)$ of order $n$ satisfies on $[0,2\pi]$  the inequality
$|T_n(x)| \le M,$ then
$$
|T'{}_n(x)| \le nM, \quad x \in [0,2\pi].
\eqno (2.11)
$$

Let $P_n(x), \ -1\le x \le1,$ is a  algebraic polynomial of degree $n.$

{\bf Theorem 2.8}  \cite{Nat}. Let $P_n(x) \le M$ for $-1 \le x \le 1.$
Then
$$
|P'_n(x)| \le \frac{Mn}{\sqrt{1-x^2}}, \quad -1<x<1.
\eqno (2.12)
$$

{\bf Theorem 2.9} (A.A. Markov)  \cite{Nat}. Let $P_n(x), \ -1\le x \le1.$ If $P_n(x) \le M$ on $[a,b],$ then
\[
|P'_n(x)| \le \frac{2Mn^2}{b-a}, \quad a \le x \le b.
\]

{\bf Theorem 2.10} (S.M.  Nikolskii) \cite{Nik}.
Let $T_n(x)$ is $n$-order trigonometric polynomial on $[0,2\pi].$ Then
\[
\|T_n(x)\|_C \le n^{1/p}\|T_n\|_{L_p} \quad (1 \le p \le \infty).
\]

The nth degree of the best approximation of a function $f \in \tilde C [0,2\pi]$  by
trigonometrical polynomial $T_n$ of degree $n$ is defined by
\[
\tilde E_n(f) =\min\limits_{T_n}\max\limits_{x} |f(x) - T_n(x)|.
\]

In this section we will give estimations of $
\tilde E_n(f).$

{\bf Theorem 2.11}  \cite{Nat}, \cite{Lor}.
There exists a constant $M$ such that, for each $f \in \tilde C [0,2\pi]$
\[
\tilde E_n(f) \le M \omega(f; \frac{1}{n}), \quad n=1,2,\ldots .
\]
{\it Note.}  N.P. Korneichuk \cite{Korn} proved that $M=1.$

{\bf Theorem 2.12} \cite{Nat}. For each $r=1,2,\ldots$ there is a constant $M_r$ with the
property that if $f \in \tilde C [0,2\pi]$ has a continuous derivative
$f^{(r)}(x),$ then
\[
\tilde E_n(f) \le \frac{M_r}{n^r}\omega(f^{(r)};\frac{1}{n}), \quad n=1,2,\ldots .
\]

{\bf Theorem 2.13} \cite{Tim}. If the function $f(x)$ has continuous derivatives
$f'(x),\ldots,f^{(r)}(x)$ on $[-1,1],$ then is a sequence of polynomials
$P_n(x)$ for which
\[
|f(x)-P_n(x)| \le 
\]
\[
\le M_r\left[\frac{1}{n}\left(\sqrt{1-x^2} + \frac{1}{n}\right)\right]^r
\omega \left(f^{(r)}; \frac{1}{n}\left(\sqrt{1-x^2} +
\frac{1}{n}\right)\right),
\]
where the constant $M_r$ depend only upon $r.$

\begin{center}
{\bf 2.4. Inverse Theorems}
\end{center}

In this section we will give some statements that are derived smoothness \ 
properties \ of a function \ $f \in \tilde C [0,2\pi]$ from the hypothesis that
the numbers $\tilde E_n(f)$ approach zero with a given rapidity. Theorems of this
kind were first obtained by Bernstein. His \  results \  have been developed \ 
further by Zigmund, Gaier, Stechkin and other. (For detail see \cite{Nat}, \cite{Lor}.)

{\bf Theorem 2.14}  \cite{Nat}, \cite{Lor}.
A function $f \in \tilde C [0,2\pi]$ belongs to the class $H_{\alpha}$
$(0 < \alpha < 1),$ if and only if
\[
\tilde E_n(f) = O(n^{-\alpha}).
\]

{\bf Theorem 2.15}  \cite{Nat}, \cite{Lor}. A
function $f \in \tilde C [0,2\pi]$ belongs to the class $Z$ if and only if
\[
\tilde E_n(f) = O(n^{-1}).
\]

{\bf Theorem 2.16} \cite{Nat}, \cite{Lor}. The  $2\pi$-periodic continuous function $f$ belongs to the
class $W^r H_\alpha$ $(0 < \alpha <1)$ if and only if
\[
\tilde E_n(f) = O\left(\frac{1}{n^{r+\alpha}}\right).
\]

Here we give the designations which are used in the paper.

Let $f(x) \in \tilde C[0,2\pi ]$ be a periodic function with $2\pi $ period.
Let $ s_k=2k\pi /(2n+1), k=0,1,\dots,2n.$ Then interpolating
polynomial $ P_n(x) $ can be written in the form
\[
P_n(x)=\sum_{k=0}^{2n}f(s_k)\psi_k (x),
\]
where
\[
\psi_k(x)=\frac{1}{2n+1}\frac{\sin\frac{2n+1}{2}
(x-s_k)}{\sin\frac{1}{2}(x-s_k)}. 
\]

Let $ \gamma=\{z:|z|=1\}.$ For a function $f(t),t\in \gamma $, the interpolating
polynomial  can be written in the form
\[
P_n(f)=\sum_{k=0}^{2n}f(t_k)\psi_k(s),
\]
where \  $ t_k=e^{is_k},$ $ s_k=2k\pi /(2n+1)$, $k=0,1,\dots,2n,$ $ s\in [0,2\pi ].$

In this \  paper \ 
we will\  use the Holder \  space of functions
$x(t) \in H_{\beta}$ $(0<\beta \le 1)$ with the  norm
\[
\|x(t)\|=M(x)+H(x,\beta)=\max\limits_{t \in \gamma} |x(t)|+
\sup\limits_{t_1,t_2 \in \gamma,t_1 \ne t_2} \frac{|x(t_1)-x(t_2)|}
{|t_1-t_2|^\beta}
\]
and its subspace $X_n$, consisting of polynomials  $x_n(t)=\sum^n_{k=-n}\alpha_k t^k.$

\vskip 25 pt
\begin{center}
{\bf Chapter 1}

\vskip 25 pt

{\bf Approximate Solution of Singular Integral Equations}
\vskip 15 pt

{\bf 1. An Smoothness of  Solutions of Singular Integral Equations}

\end{center}

Consider the connection between a smoothness of coefficients
of
singular integral equations and a smoothness of their solutions.

\newpage

{\bf 1.1. The Integral Operators on the Smooth Functions}

{\bf 1.1.1. Fundamental Statements}

{\bf Lemma 1.1} \cite{Boy16}, \cite{Boy25}.
If $ h(t,\tau) \in H_{\alpha\alpha},0<\alpha \le 1 $, then the operator
\[
 Kx\int\limits_0^{2\pi}(h(t,\tau) - h(t,t))\mid \cot \frac {\tau-t}{2}
\mid ^{\eta}x(\tau) d\tau, 
\]
$ 0\le\eta<1, $ transforms every function $x(t) \in C[0,2\pi]$ to the
function belonging to the class $ H_{\alpha} (0 < \alpha \le 1). $

{\bf Lemma 1.2} \cite{Boy16}, \cite{Boy25}.
Let  $ x \in H_{\alpha} (0 < \alpha \le 1). $
Then the function $ v(t), $
\[
v(t)=\int\limits_0^{2\pi} x(\tau) \mid \cot \frac
{\tau-t}{2} \mid^{\eta} d\tau, 
\]
belongs to the class $ Z $ if
$ \alpha=\eta $ or belongs to the class $ H_\gamma $ if
$ \alpha \ne \eta, $ moreover $ \gamma=1 $ if $ \alpha>\eta $ and
$ \gamma=\alpha+1-\eta $ if $ \alpha<\eta. $

{\bf Lemma 1.3} \cite{Boy16}, \cite{Boy25}.
Let the integral equation
\[
 Hx \equiv x(t)+ \int\limits_0^{2\pi} h(t, \tau) \mid \cot \frac {\tau-t}
{2} \mid^{\eta} x(\tau)d\tau=f(t) , \hskip   10 pt
0 \le \eta<1, 
\]
where $ h(t, \tau) \in H_{\alpha \alpha}, f(t) \in H_{\alpha} (0 < \alpha \le 1), $
has a unique solution $ x^*(t). $
Then $ x^*(t) \in H_{\alpha}. $

{\bf Lemma 1.4}  \cite{Boy16}, \cite{Boy25}.
Let the integral equation $ Hx=f, $ where $ h(t, \tau)
\in W^{rr}H_{\alpha\alpha}, f(t) \in  W^{r}H_{\alpha}, $
has a unique solution $ x^*(t). $ Then $ x^*(t) \in W^{r}H_{\alpha}. $

\newpage

\begin{center}
{\bf 1.2. On Smoothness of  Solutions of the Singular Integral Equations
on Closed Contours}
\end{center}

{\bf 1.2.1. Fundamental Statements}

Let us consider the singular integral equation
$$
Kx \equiv a(t) x(t)+b(t)(S_{\gamma}x)(t)+U_{\gamma}(h(t,\tau)
\mid \tau-t \mid^{-\eta} x(\tau)) =f(t) \eqno (1.1)
$$
where
$\gamma=\{z:\mid z \mid=1\},$
\[
 S_{\gamma}(x) \equiv \frac{1}{\pi i} \int\limits_{\gamma} 
 \frac{x(\tau) d\tau}{\tau-t}, U_{\gamma}(hx)=\frac{1}{2\pi i}
\int\limits_{\gamma} h(t,\tau) x(\tau) d\tau, 0 \le \eta<1. 
\]

Let us put the characteristic operator
$ K^0x=ax+bS_{\gamma}x $ in accordance to  the operator $ K. $

{\bf Theorem 1.1}  \cite{Boy16}, \cite{Boy25}.
Let  $ a(t), b(t), f(t) \in H_{\alpha}, h(h,\tau) \in H_{\alpha \alpha}$ $(0<\alpha<1). $
Suppose that the index of the operator
$ K^0 $ is equal to zero, and the operator $ K $ is continuously
invertible in the space $ H_{\beta} (0<\beta\le \alpha). $
Then a unique solution $ x^*(t) $ of the equation (1.1) belongs to
the class $ H_{\alpha}. $

{\it Note.}
If the conditions of Theorem 1.1 are fulfilled, then the
operator $ K $  operating from $ H_{\alpha} $ into $ H_{\alpha} $ is
continuously invertible.

{\bf Theorem 1.2}  \cite{Boy16}, \cite{Boy25}.
Let  $ a(t), b(t), f(t) \in W^{r}H_{\alpha},
h(t, \tau) \in W^{rr}H_{\alpha\alpha}, $ $0<\alpha<1.$
Suppose that the index of the operator
$ K^0 $ is equal to zero, and the operator $ K $ is continuously
invertible in the space $ H_{\beta} (0<\beta\le \alpha). $
Then the solution $ x^*(t) $ of the equation (1.1) belongs to
the space $W^{r}H_{\alpha}. $

\begin{center}
{\bf 2. Approximate Solution of Linear Singular Integral Equations on
     Closed Contours of Integration\\
    (Basis in Holder Space)}
\end{center}

In this section we investigate approximate methods for solution of
singular integral equations
\[
Kx \equiv a(t)x(t)+b(t)S_{\gamma}(x)+U_{\gamma}(h(t,\tau){\mid
\tau-t \mid}^{-\eta} x(\tau))=
\]
$$
=f(t), \eqno (2.1)
$$
$$
Lx \equiv a(t)x(t)+S_{\gamma}(h(t, \tau)x(\tau))=f(t). \eqno (2.2)
$$
Here $a(t), b(t), f(t) \in H_\alpha,$ $h(t,\tau) \in H_{\alpha,\alpha}, 0<\alpha \le 1.$

General form of one-dimensional singular integral equations is
\[
a(t)x(t) + \int\limits_{\gamma}
\frac{k(\tau,t)}{\tau-t}d\tau =f(t).
\]

Let us put
\[
b(t)=k(t,t), \ \ \ \frac{h(t,\tau)}{|\tau-t|^{\eta}} =
\frac{k(t,\tau)-k(t,t)}{\tau-t}.
\]
In result we receive the singular integral equation in the form  (2.1).

Here we use the following designations
\[
S_{\gamma}x=\frac{1}{\pi i} \int\limits_{\gamma}\frac{x(\tau)d\tau}
{\tau-t},\hskip 10 pt  U_{\gamma}(h(t,\tau) x(\tau))=\frac{1}{2\pi i}
\int\limits_{\gamma}h(t,\tau)x(\tau)d\tau,
\]
where
$ \gamma $ is the unit circle with the center in origin of coordinates:
$\gamma=\{z:\mid z \mid=1\}.$

In this section we will consider collocation method and method of mechanical quadrature.

Proofs of these methods we give in Holder space $X$ of functions
$x(t) \in H_{\beta}$ $(0<\beta \le \alpha < 1)$ with the  norm
\[\|x(t)\|=\max\limits_{t \in \gamma} |x(t)|+
\sup\limits_{t_1,t_2 \in \gamma,t_1 \ne t_2} \frac{|x(t_1)-x(t_2)|}
{|t_1-t_2|^\beta}
\]
and its subspace $X_n$, consisting of $n$-order polynomials  $x_n(t)=\sum^n_{k=-n}\alpha_k t^k.$

Also we will need in the space $W(\gamma)$ introduced in the book \cite{Ivan1}.

Let us associate with a function $\varphi(z)$ functions $\varphi^+(t)$ and $\varphi^-(t)$, which
are analytical inside and outside of the contour
$\gamma$ respectively and connected with
$\varphi(t)$ by Sohotzky - Plemel formulas
\[
\varphi^+(t) - \varphi^-(t) = \varphi(t), \quad
\varphi^+(t) - \varphi^-(t) = S_{\gamma}\varphi(t).
\]
Well known \cite{Gakh} that Sohotzky - Plemel formulas have place if the
function $\varphi(t)$ belongs to common enough classes of functions.

The space $W(\gamma)$ consists of functions $\varphi(t),$ $t \in \gamma,$
for which functions $\varphi^+(t)$ and $\varphi^-(t)$ are continuous.
The norm in the space $W(\gamma)$ is introduced with the formula
\[
\|\varphi\|_{W(\gamma)} = \max\limits_{t \in \gamma}|\varphi^+(t)| +
\max\limits_{t \in \gamma}|\varphi^-(t)|.
\]

\begin{center}
{\bf 2.1. Methods of Collocations and Mechanical Quadratures}
\end{center}

{\bf 2.1.1. Fundamental Statements}

Many works was devoted to numerical solution of the equations (2.1) with $\eta
= 0.$ First of them was the V.V. Ivanov's paper  \cite{Ivan}, in which the collocation
method was used for solution of the equation (2.1) with $\eta =0.$

An approximate solution of equation (2.1) is sought in the form
of the polynomial
$$
 x_n(t)=\sum_{k=-n}^n \alpha_k t^k. \eqno (2.3)
$$

According the collocation method, coefficients $\alpha_k,$
$k=\overline{-n,n},$ are defined from the system of linear
algebraic equations, that in the operator form is written as
$$
P_n[a(t)x_n(t) + b(t)S_\gamma(x_n)
+ U_\gamma(h(t,\tau)x_n(\tau)] = P_n[f],
\eqno (2.4)
$$
where operator $P_n$ is the \  projector \ of \ interpolation   onto the set of polynomials, \ 
constructed \ on the knots $t_k=\exp(is_k),$ \  $s_k=2k\pi/(2n+1),$
$k=0,1,\ldots,2n.$

{\bf Theorem 2.1} \cite{Ivan},  \cite{Ivan1}.  Let the functions $a$, $b$, $f \in H_\alpha,$
$h(t,\tau) \in H_{\alpha \alpha} \quad (0< \alpha \leq 1),$ $\eta = 0$
and the operator $K$ is continuously invertible in the space $W(\gamma).$
Then for large enough $n$ the system (2.4) has a unique solution $x^*_n(t)$
and the estimate 
\[
\|x^*(t) - x^*_n(t)\|_{W(\gamma)} \leq
A \ln n\|x^*(t) - T_n(x^*)\|_{W(\gamma)}
\]
 is valid.
Here $x^*(t)$ is a solution of the equation (2.1),
$T_n(t) = \sum\limits^n_{k=-n}\beta_k t^k$ is the polynomial of the best
approximation to $x^*(t)$ in the metric of the space $W(\gamma).$

{\it Note.} The similar results is valid when the index of the operator $K$ do not
equal to zero.

Later the similar results was received in the Holder space.

{\bf Theorem 2.2} \cite{Gab}.  Let the functions $a(t),b(t),f(t) \in H_\alpha,$
$h(t,\tau) \in H_{\alpha \alpha},$ $\eta = 0$
and the operator $K$ is continuously invertible in the Holder space $X =
H_\beta$ $(0<\beta<\alpha).$ Then for $n$
such that $q=An^{-\alpha+\beta}\ln n <1,
$ the system (2.4) has a unique solution $x^*_n(t)$ and the estimation
${\|x^*-x^*_n\|}_{\beta} \leq An^{-(\alpha-\beta)}\ln n,$ where $x^*(t)$ is
a solution of the equation (2.1), is valid.

Apparently the first works devoted to numerical methods for solution of the
equations (2.1) (for $ 0< \eta <1 $) and (2.2) was papers \cite{Boy1}, \cite{Boy5}. In these papers was
investigated the collocation method and mechanical quadrature method for
solution of the equations (2.1) and (2.2).

Let the index of the operator $K$ is equal to zero.

An approximate solution of the equation (2.1) we look for in the form of the
polynomial (2.3).

The coefficients ${\alpha_k},$ $k=-n,\ldots,-1,0,1,\ldots,n,$
are defined from the system of linear
algebraic  equations 
$$ K_n x_n \equiv P_n [a(t) x_n(t)+b(t) S_{\gamma} (x_n)+U_{\gamma}
(P_{n}^{\tau}[h(t,\tau) d(t,\tau) x_n(\tau)])]= $$
$$
= P_n[f(t)], \eqno (2.5)
$$
where
\[
 d(t,\tau)=
\left \{ \begin{array}{ccc}
\mid \tau-t \mid^{-\eta}, & \mid \sigma-s \mid \ge 2\pi/(2n+1),\\
\mid e^{i2\pi/(2n+1)}-1 \mid^{-\eta}, &
 \mid \sigma-s \mid < 2\pi/(2n+1),\\
\end{array} \right.
\]
$\tau=e^{i\sigma}, t=e^{is}.$

{\bf Theorem 2.3} \cite{Boy1}, \cite{Boy5}, \cite{Boy7} .
Let the functions
$ a,b, f \in H_{\alpha},$ $ h \in H_{\alpha,\alpha}$ $(0<\alpha \le 1) $ and the operator $ K $ is
continuously invertible in the Holder space
$ X=H_{\beta} (0<\beta< \min(\alpha, 1-\eta). $
Then for $ n $ such that $ q=An^{-\xi} \ln n<1\  
(\xi=\rm {min}(\alpha-\beta, 1-\eta-\beta, \beta)), $ the system of equations 
(2.5) has a unique solution $ x^*_n $ and the estimation 
$ \Vert x^*-x^*_n \Vert_\beta <An^{-\xi} \ln n $ is valid.
Here $ x^* $ is a unique solution of the equation (2.1).

The estimate of error obtained in preceding Theorem depends on the
constant $\eta.$ Let us build a calculating scheme the estimation
of the error of which depends only on a smoothness of the functions $a,b,h,f.$

     For simplicity  we will take $G(t)=(a(t)-b(t))/(a(t)+b(t)), $ $ a(t)+b(t)\equiv 1.$

With each function $ x(t) \in H_{\beta} $ we associate
functions $ x^+(t) $ and $ x^-(t) $ which are analytical inside
and outside $ \gamma $ respectively and connected with $ x(t) $
by Sohotzky - Plemel formulas $x^+(t)-x^-(t)=x(t)$,
$ S_{\gamma}x=x^+(t)+x^-(t). $

     The approximate solution of the equation (2.1) we will look for in
the form of the polynomial (2.3), the coefficients of which are defined
from the system of  equations, that representable in the operator form
by the expression
$$
K_n x_n \equiv P_n[a(t) x_n(t) + b(t)
S_{\gamma} (x_n(\tau)) +
$$
$$
\left.+\sum_{k=0}^{2n}h(t,t_k)
 x_n(t_k)
\int\limits_{t'_k}^ {t'_{k+1}}{\mid \tau-t \mid}^{-\eta}
d\tau\right]=P_n[f(t)], \eqno (2.6)
$$
where
$ t'_{k+1}=e^{is'_{k+1}}, $ $ s'_{k+1}=(2k+1) \pi/(2n+1),$
$k=0,1,\ldots,2n.$

{\bf Theorem 2.4} \cite{Boy1},  \cite{Boy16},  \cite{Boy25}.
Let the conditions of Theorem 2.3 are fulfilled.
Then among all possible approximate methods for solution of the
equation (2.1), using $ n $ values of the functions $ a,b,f $ and
$ n^2 $ values of the function $ h, $ the method, being described by
the calculating scheme (2.6),  is optimal with respect to
order on the class $ H_{\alpha}. $ The error of this method is equal to
$\Vert x^* -x^*_n \Vert_\beta = 0((n^{-\alpha+\beta} +n^{-(1-\gamma)+\beta})\ln n), $ where $x^* $ and
$ x^*_n  $ are solutions of the equations (2.1) and (2.6) respectively.

      Now assume that  $ a(t), b(t), f(t) \in W^r, h(t,
\tau) \in W^{r,r}. $

     The approximate solution of the equation (2.1) we will seek in
the form of the polynomial (2.3), the coefficients $ {\alpha_k} $ of
which are defined from the system of  linear algebraic equations
$$
K^{(1)}_nx_n \equiv P_n[a(t)x_n(t) + b(t)S_{\gamma}(x_n(\tau)) +
$$
$$
+ U_{\gamma}[P_{n}^{\tau}[h(t,\tau)]x_n(\tau)|\tau-t|^{-\eta}]] =
P_n[f(t)].
\eqno (2.7)
$$

The next Theorem is a generalization of some optimal with respect \  to order \ 
algorithms \  for solution of singular integral equations, printed in \cite{Boy1},  \cite{Boy16},  \cite{Boy25}.

{\bf Theorem 2.5}. Let the operator $ K $ is continuously invertible
in the Holder spaces $ X =H_\beta (0<\beta <\alpha, \  0< \alpha <1). $
 Let
$$
 a,b,f\in W^rH_\alpha (M), h\in W^{r,r}H_{\alpha,\alpha}(M),
 r=0,1,\dots.
\eqno (2.8)
 $$
 Among all
possible algorithms for  solution of the equation (2.1), using $ n $
values of the functions $ a,b,f $ and $ n^2 $ values of the function
$ h $, the calculating scheme (2.3), (2.7) is optimal with respect to order. The
error of this calculating scheme is equal to $ \Vert x^*-\widetilde
x_n^* \Vert_\beta \le An^{-(r+\alpha-\beta)}\ln n, $ where $ x^* $ and $ x_n^* $ are
solutions of the equations (2.1), (2.7) respectively.

     The approximate solution of the equation (2.2) is sought in the form
of the polynomial (2.3), the coefficients $ {\alpha_k} $ of which are
defined from the system of  linear algebraic equations
$$
L_n x_n \equiv P_n[a(t) x_n(t)+S_{\gamma}(P_{n}^{\tau}[h(t,\tau)]
x_n(\tau))]=P_n[f(t)]. \eqno (2.9)
$$

{\bf Theorem 2.6}  \cite{Boy16},  \cite{Boy25}. Let the operator $ L $ is continuously invertible
in the Holder spaces $ X =H_\beta (0<\beta < \alpha \leq 1).$
  Then for $ n $ such that
$ q=A(n^{\beta} \ln^2n(E_n(a)+E_n(b)+ E_n(\psi) + E_n^t(h)+E_n^{\tau}(h))<1, $
 the system of the equations (2.9) has a
unique solution $ x_n^* $ and the inequality $ \Vert x^*-x_n^* \Vert
\le A(n^{\beta} \ln^2n(E_n(a)+E_n(b)+E_n(\psi) +E_n^t(h)+E_n^{\tau}(h)) +
n^\beta E_n(f)) $ is valid, where $ x^* $ is a solution of the
equation (2.2), 
\[\psi(z) = \exp \left\{\frac{1}{2\pi}\int\limits_{\gamma}(\ln
G(\tau))(\tau-z)^{-1}d\tau \right\},
\]
 $G(t) = (a(t)-b(t))/(a(t)+b(t)),$ $b(t)
= h(t,t).$

Let \ $\overline P_n $ is \ the \  projector for interpolation onto the set of
trigonometrical polynomials, constructed on the knots
$ \overline t_k=e^{i\overline s_k}$,
$\overline s_k=2k\pi/(2n+1), k=0,1,\dots,2n.$

Let an approximate solution of the equation (2.2) is sought in the
form of the polynomial (2.3), the coefficients ${\alpha_k}$
$(k=-n,\ldots,-1,0,1,\ldots,n)$  of which
are defined from the system of the linear algebraic equations
\[
\tilde L_n x_n \equiv \overline P_n \left[ a(t) x_n(t)+
\frac {1}{\pi i} \int\limits_{\gamma} P_n^\tau \left [\frac{h(t,\tau)
x_n(\tau)} {\tau-t} \right] d\tau \right]=
\]
$$
=\overline P_n [f(t)]. \eqno (2.10)
$$

{\bf Theorem 2.7} \cite{Boy5},  \cite{Boy16},  \cite{Boy25}. Let the operator $ L $ is continuously invertible.
Then for $ n $ such that $ q=A(E_n(a)+ E_n(\psi) + E_n^t(h)+E_n^{\tau}(h))\ln^2n<1, $
the system of the equations (2.10) has a unique solution $ x_n^* $ and
the inequality $ {\Vert x^*-x_n^* \Vert}_\beta \le A(q+E_n(f)) $ holds,
where $ x^* $ is a solution of the equation (2.2),
$$\psi(z) = \exp \left\{\frac{1}{2\pi}\int\limits_{\gamma}(\ln
G(\tau))(\tau-t)^{-1}d\tau \right\},$$ $G(t) = (a(t)-b(t))/(a(t)+b(t)),$ $b(t)
= h(t,t).$

     Let us build the optimal with respect to order calculating scheme
for approximate solution of the singular integral equations  (2.2).
We shall seek an approximate solution in the form of the polynomial (2.3), the
coefficients of which $ {\alpha_k} $ $(k=-n,\ldots,n)$ are defined from the system of
linear algebraic equations
$$
\overline L_n x_n \equiv \overline P_n \left[a(t)x_n(t) +
\frac{h(t,t)}{\pi i}\int\limits_{\gamma}\frac{x_n(\tau)}{\tau-t}d\tau +
\right.
$$
$$
+\left. \frac {1}{\pi i}\int\limits_{\gamma} P_n^\tau
\left[\frac{h(t,\tau)-h(t,t)}{\tau-t} x_n(\tau) \right] d\tau \right]=
\overline P_n[f(t)]. \eqno (2.11)
$$

{\bf Theorem 2.8} \cite{Boy5},  \cite{Boy7},  \cite{Boy25}.  Let the operator $ L $ is continuously invertible.
Among all possible algorithms for approximate solution of the singular
integral equation (2.2), using $ n $ values of the functions
$ a(t), f(t) $ and $ n^2 $ values of the function $ h(t, \tau) $, the
calculating scheme (2.3), (2.11) is optimal with respect to
 order and has the error
\[ 
\Vert x^*-x_n^* \Vert_C \le
\]
\[
\le A(E_n(a) + E_n(\psi^{\pm})+E_n^t(h(t,\tau)+E_n^{\tau}(h(t,\tau))+
E_n(f))\ln n, 
\]
where \ $ x^* $ \ and \ $ x_n^*  $ \  are \  solutions of the
equations \ (2.2) and (2.9) \ \  respectively, \ 
$ \psi(z)=\exp \{ \frac {1}{2\pi i}
\int\limits_{\gamma}(\ln G(\tau))(\tau-z)^{-1}d\tau\},$\\
$G(t) = (a(t)-b(t))/(a(t)+b(t)),$ $b(t) = h(t,t).$
\vskip 10 pt
{\bf 2.1.2. Proofs of Theorems}

We drop proofs of the Theorem 2.1 and the Theorem 2.2, because its are special
cases of the Theorem 2.3.

{\bf Proof of Theorem 2.3.}
Here we give the full proof of the Theorem 2.3, because it had not been
printed before.

At first let us investigate the collocation method
$$
K_n x_n \equiv P_n[a(t)x_n(t)+b(t)S_{\gamma} x_n+U_{\gamma}(h(t,\tau)
{\mid \tau-t \mid}^{-\eta}x_n(\tau)]=
$$
$$
= P_n[f(t)]. \eqno (2.12)
$$

Fix  function $ x(t) \in H_{\beta} $ and associate with it
the functions $ x^+(t) $ and $ x^-(t) $, which are analytical inside
and outside of $ \gamma $ respectively \ and \ connected\   with $ x(t) $
by Sohotzky - Plemel formulas $x^+(t)-x^-(t)=x(t)$,
$ S_{\gamma}x=x^+(t)+x^-(t). $
Since the function $ G(t) $ is represented \cite{Gakh} as
\[ G(t)=\frac{\psi^+(t)}{\psi^-(t)}, \psi^{\pm}(z)=\exp(\Theta^{\pm}(z)),
\Theta(z)=\frac{1}{2\pi i}\int\limits_{\gamma}\frac{\ln G(\tau)d\tau}{\tau-z}, 
\]
then  the equations (2.1), (2.5) and (2.12)
 are equivalent to the following equations
respectively
\[ Fx \equiv \psi^-x^+ - \psi^+ x^- + \psi^-U_{\gamma}(h(t,\tau)
{\mid \tau-t \mid}^ {-\eta}x(\tau))=\psi^- f; 
\]
\[ \widetilde F_n x \equiv P_n[\psi^ - x_n^+ -\psi^+ x_n^- + \psi^-
U_{\gamma} [P_{n}^{\tau} [h(t,\tau)d(t,\tau) x_n(\tau)]]]=
\]
\[
=P_n[\psi^-f];
\]
\[ F_n x_n=P_n[\psi^- x_n^+ -\psi^+ x_n^- + \psi^-U_{\gamma}
(h(t,\tau){\mid \tau - t\mid}^{-\eta} x_n(\tau))]=
\]
\[
=P_n[\psi^-f].
\]

     The operator $ K $ is continuously invertible. So the operator $F$ make one-to-one mapping of the space $X$ onto inself. It is follows from  Banach
Theorem that operator $ F $ has  linear inverse operator in the space
$ X $ also.

Let us introduce the operator
\[ Vx= \psi_n^ - x^+ - \psi_n^+ x^- +T_n[\psi^ -U_{\gamma}(h(t,\tau)
{\mid \tau-t \mid}^{-\eta}x(\tau))], \]
where
$ \psi_n=\psi_n^{+} - \psi_n^{-},\psi_n^{+} ( \psi_n^{-}) $ is the polynomial
of best uniform approximation of degree $ n $ for the function
$ \psi_n^{+} ( \psi_n^{-}): \psi^{\pm }_n = T_n \psi^{\pm }. $
Here $T_n$ is the projector, which put to according each continuous
function $f(x)$ its polynomial of best uniform approximation.
It was marked in the first section of this chapter, that the operator
$ U_{\gamma}(h(t,\tau){\mid \tau-t \mid}^{-\eta} x(\tau)) $ \ maps any
continuous \ function into  function belonging to the Holder class
$ H_{\delta}(\delta=\min(\alpha,1-\eta)). $

Therefore
$ \Vert Fx-Vx \Vert_\beta \le A n^{-(\delta- \beta)} \Vert x \Vert_\beta; $

$ \Vert P_n[Fx-Vx] \Vert_\beta \le \Vert P_n \Vert_\beta \Vert Fx-Vx \Vert_\beta \le
A n^{-(\delta- \beta)} \ln n \Vert x \Vert_\beta. $

In obtaining the last relation was used the known \cite{Nat}, \cite{Nat1} inequality
$ \Vert P_n \Vert_C \le 8+\frac{4}{\pi}\ln n. $ It follows from the equality
$ P_n Vx_n=Vx_n, $ that
$ \Vert Fx_n-F_n x_n \Vert_\beta= \Vert F x_n-
Vx_n + P_n Vx_n - F_n x_n \Vert_\beta \le $
$ \Vert Fx_n - Vx_n \Vert_\beta + \Vert P_n Vx_n -F_n x_n \Vert_\beta \le
A n^{-(\delta- \beta)}\ln n \Vert x_n \Vert_\beta. $

It is easy to see that
\[\|F_nx_n\|=\|Fx_n-(Fx_n-F_nx_n)\|\ge \|Fx_n\|-\|(Fx_n-F_nx_n)\|\ge 
\]
\[
\ge\frac{1}{\|F^{-1}\|}\|x_n\| -An^{-(\delta-\beta)}\ln n\|x_n\|\ge A\|x_n\|.
\]

From this inequality and the Theorem about left inverse operator \cite{Kan}, p. 449, it follows that,  for $n$ such that
 $A n^{-(\alpha-\beta)}\ln n < 1$, the operator $F_n$ has left inverse operator.
Since $F_n$ is finite-dimensional operator, then the operator $F_n$ has continuously\\ invertable operator $F^{-1}_n.$

Let solutions of the equations (2.1) and (2.12) are denoted by $ x^* $
and $ x_n^* $ respectively. Then,
as follows from the general theory of the
approximate methods of the analysis, the estimate $ \Vert x^*-x_n^* \Vert_\beta
\le A n^{-(\delta-\beta)} \ln n $
is valed.
The collocation method is justified.

Introduce the operator
\[ V_n x_n \equiv P_n[\psi^-x_n^+ - \psi^+x_n^-+\psi^-U_{\gamma}
(h(t,\tau)d(t,\tau)x_n(\tau))].
\]
 
 Let us  estimate the 
$\Vert F_n x_n-V_n x_n \Vert_\beta.$

It is easy to see that
\[ \Vert F_n x_n-V_n x_n \Vert_\beta=\Vert P_n[\psi^-U_{\gamma}
(h(t,\tau)({\mid \tau-t \mid}^{-\eta}- d(t,\tau))x_n(\tau))] \Vert_\beta.
\]

In the space $ \tilde C[0,2\pi] $
\[\Vert P_n[\psi^-U_{\gamma}(h(t,\tau)({\mid \tau-t \mid}^{-\eta}-
d(t,\tau))x_n(\tau))] \Vert_C \le 
\]
\[ \le \Vert P_n \Vert \Vert \psi^-U_{\gamma}(h(t,\tau)({\mid \tau-
t \mid}^{-\eta}-d(t,\tau))x_n(\tau)) \Vert_C \le 
\]
\[ \le A\ln n \Vert \psi^-U_{\gamma}(h(t,\tau)({\mid \tau-
t \mid}^{-\eta}-d(t,\tau))x_n(\tau)) \Vert_C \le 
\]
\[ \le An^{-(1-\eta)} \ln n\Vert x_n \Vert_C \le An^{-(1-\eta)}\ln n
\Vert x_n \Vert_\beta 
\]
and hence
$$
\Vert F_n x_n -V_n x_n \Vert_\beta \le An^{-(1-\eta-\beta)}\ln n \Vert x_n
\Vert_\beta.
\eqno (2.13)
$$

Let us estimate the norm of the difference
\[ \Vert V_n x_n - \widetilde F_n x_n \Vert_\beta = \Vert P_n [\psi^-U_{\gamma}
(R_n^\tau[h(t,\tau)d(t,\tau)x_n(\tau)])] \Vert_\beta \le
\]
\[ \le \left \Vert P_n \left[\frac{\psi^{-*}(s)}{2\pi}\int\limits_0^{2\pi}
h^*(s,\sigma)d^*(s,\sigma)e^{i\sigma}(x_n(e^{i\sigma})-\overline
x_n(e^{i\sigma}))d\sigma \right] \right \Vert_\beta + 
\]
\[ + \left\Vert P_n \left[ \frac{\psi^{-*}(s)}{2\pi}
\sum_{k=0}^{2n}
\int\limits_{s_k}^{s_{k+1}}\left [ h^*(s,\sigma)d^*(s,\sigma)e^{i\sigma}-h^*(s,s_k)
d^*(s,s_k)e^{is_k} \right] \right.\right.\times
\]
\[
\times \left.\left.
\overline {x_n(e^{i\sigma })}d\sigma \right ] \right\Vert_\beta = \Vert I_1(s) \Vert_\beta +\Vert I_2(s) \Vert_\beta, 
\]
where $ \overline x_n (e^{i\sigma})$ is the step-function equal to
$ x_n(e^{is_k}) $ on the interval $ [s_k, s_{k+1}), \psi^{-*}(s)=
\psi^-(e^{is}), h^*(s,\sigma)=h(e^{is}, e^{i\sigma}),$  $ d^*(s,\sigma)=
d(e^{is}, e^{i\sigma}), R_n^\tau=I-P_n^\tau. $

Since $ x_n(e^{is}) \in H_{\beta} $ and $ \overline x_n(e^{i\sigma}) $
is a step-function coinciding with $ x_n(e^{i\sigma}) $ at the points
$ s_k $, then
\[ \Vert I_1(s) \Vert_\beta \le A \ln n \hskip 11 pt
\max \mid x_n(e^{i\sigma}) - \overline
x_n(e^{i\sigma}) \mid \le A \ln n \Vert x_n \Vert_\beta / n^{\beta}. 
\]

Using the fact  that $ P_n f=\sum_{k=0}^{2n}f(s_k)\psi_k(s), $
where $ \psi_k(s) $ are fundamental polynomials constructed on the points $ s_k $, and
 that $ \max \mid P_n f \mid \le (A+B \ln n) \max \mid f \mid, $ one has
\[ \max_{s} \mid I_2(s) \mid \le 
\]
\[
\le A \ln n \max_{0 \le j \le 2n}
\big | \frac {\psi^{-*}(s_j)}{2\pi} \left \{ \sum_{k=0}^{2n}
\int\limits_{s_k}^{s_{k+1}}[h^*(s_j, \sigma)
d^*(s_j, \sigma)e^{i\sigma}\right.-
\]
\[-\left. h^*(s_j,s_k)d^*(s_j,s_k)e^{is_k}]
\overline x_n(e^{i\sigma})d \sigma \right \} \big | \le 
\]
\[ \le A \ln n \max_{0 \le j \le 2n}
 \left \{ \big | \frac{\psi^{-*}(s_j)}{2\pi}
\sum_{k=0}^{2n} \int\limits_{s_k}^{s_{k+1}}[h^*(s_j, \sigma)
e^{i\sigma}-\right.
\]
\[
\left.
-h^*(s_j,s_k)e^{is_k}]
d^*(s_j, \sigma ) \overline x_n (e^{i\sigma})d \sigma \big |\right . +
\]
\[
+\left .\big|\frac{\psi^{-*}(s_j)}{2\pi}
 \sum_{k=0}^{2n} \int\limits_{s_k}^{s_{k+1}}
h^*(s_j,s_k)e^{is_k} [d^*(s_j, \sigma)-d^*(s_j,s_k)]
\overline x_n(e^{i\sigma}) \big | \right\}= 
\]
\[ =A \ln n \max_{0 \le j
\le 2\pi}( {\mid I_3(s_j)\mid + \mid I_4(s_j)\mid}). 
\]

It is easy to see that
\[ \max_j \mid I_3(s_j) \mid \le A \Vert x_n \Vert_\beta /n^{\alpha},
\]
\[ \max_j \mid I_4(s_j) \mid \le
A \Vert x_n \Vert_C \max_j \sum_{k=0}^
{2n} \int\limits_{s_k}^{s_{k+1}} \mid d^*(s_j, \sigma)-d^*(s_j,s_k) \mid
d\sigma. 
\]

Below the prime in the summation indicate that
$k \ne j-1,j.$ Without loss of
generality we assume $j=0.$  Then
\[ \sum_{k=0}^{2n}{}' \int\limits_{s_k}^{s_{k+1}}   (*)
d\sigma= \sum_{k=1}^{n-1} \int\limits_{s_k}^{s_{k+1}} (*)d\sigma +
 \int\limits_{s_n}^
{s_{n+1}} (*) d\sigma + \sum_{k=n+1}^{2n-2} \int\limits_{s_k}^{s_{k+1}} (*)
d\sigma=
\]
\[
=I_5+I_6+I_7, 
\]
where $ (*)=\mid d^*(0,\sigma)-d^*(0,s_k) \mid . $

For values $ \sigma $, belonging to the segment
$ [2\pi/(2n+1), 2\pi-2\pi/(2n+1)] $, the formula $ d^*(0,\sigma)=
d(0,\sigma)=\frac{1}{2}{\mid \csc \frac {\sigma}{2} \mid}^{\eta}$ is
valid. Therefore
\[ I_5=\frac{1}{2} \sum_{k=1}^{n-1} \int\limits_{s_k}^{s_{k+1}}(\mid \csc
\frac{\sigma}{2} \mid^{\eta} - \mid \csc \frac {s_k}{2}
\mid^{\eta})
d\sigma \le 
\]
\[ \le A \sum_{k=1}^{n-1} \int\limits_{s_k}^{s_{k+1}} {(s_k + \theta
(\sigma-s_k))}^{-1-\eta}(\sigma-s_k)d\sigma \le 
\]
\[
 \le A \sum_{k=1}^{n-1}{\left(\frac{2n+1}{2k\pi}\right)}^{1+\eta}
\int\limits_{s_k}^{s_{k+1}}(\sigma-s_k)d\sigma \le A n^{-(1-\eta)}. 
\]

The sum  $I_7$ can be estimated similarly. It is not difficult to see that
\[ I_6 \le A n^{-(1-\eta)}. 
\]

It follows from the estimations for $ I_3-I_7 $ that $ \max \mid I_2(s_j)
\mid \le A n^{-\delta} \ln n \Vert x_n \Vert_\beta. $
Since $ I_2(s) $ is a polynomial of degree $ n $, then $ \Vert I_2
\Vert_\beta \le A n^{-(\delta-\beta)}\ln n \Vert \tilde x_n \Vert_\beta. $
From the last inequality and values of the norm $ \Vert I_1 \Vert_\beta $ we have
$ \Vert \widetilde F_n x_n-V_n x_n \Vert_\beta \le A n^{-\xi} \ln n
\Vert x_n \Vert_\beta. $
Hence, it follows from the last inequality and the inequality (2.13) and Kantorovich
theory of approximate methods of analyses, that for
$ n $ such that $ q=A n^{-\xi} \ln n<1, $ the operator
$ \widetilde F_n $ is continuous invertible and the estimation
$ \Vert \overline x^*_n - x^*_n \Vert_\beta \le A n^{-\xi}\ln n $ is valid.

Theorem is proved.

{\bf Proof of Theorem 2.4.}

The equation (2.1) is equivalent to the equation
$$
K_1 x \equiv x^+(t) - G(t) x^-(t) + U_\gamma(h(t,\tau) |\tau-t|^{-\eta}x(\tau)) = f(t).
\eqno (2.14)
$$
The last equation is uniquely solvable for any right-hand side. Using well known Banach Theorem \cite{Lus}, we see that the operator $K_1$ is continuously
 invertible and $\|K^{-1}_1\| = C.$

The equation (2.14) is equivalent to the equation
$$
K_2 x \equiv \psi^-(t) x^+(t) - \psi^+(t) x^-(t) + \psi^-(t) U_\gamma(h(t,\tau)) |\tau-t|^{-\eta}x(t)) =
$$
$$= \psi^-(t) f(t),
\eqno (2.15)
$$
where
\[
\psi(z) = \exp \left\{\frac{1}{2\pi i}\int\limits_\gamma
\frac{\ln G(\tau)}{\tau-t}d\tau\right\}.
\]

Repeated the previons arguments we see that the operator $K_2$ is continuously invertible and $\|K^{-1}_2\| = C.$

 Let $\psi^+_n(\psi^-_n)$ 
is the polynomial of best uniform approximation of degree $n$ for the function $\psi^+(\psi^-);$ $T_n$ is projector onto the polynomial of 
best uniform approximation of degree $n$. 

Let us consider the equation
\[
K_3 x \equiv \psi^-_n(t) x^+(t) - \psi^+_n(t) x^-(t) +
\]
$$ 
+T_n[\psi^-(t) U_\gamma(h(t,\tau) |\tau-t|^{-\eta}x(\tau))]
= \psi^-(t) f(t).
\eqno (2.16)
$$

From conditions $a(t),$ $b(t) \in H_\alpha$ it follows that $G(t) \in H_\alpha.$
As $G(t) \ne 0$ then $\ln G(t) \in H_\alpha,$ $t \in \gamma.$ Using the Sohotzky-Plemel formulas, we have
\[
\psi^+(t) = \exp\left\{\frac{1}{2} \ln G(t) + \frac{1}{2\pi i}\int\limits_\gamma
\frac{\ln G(\tau)}{\tau-t}d\tau\right\},
\]
\[
\psi^-(t) = \exp\left\{-\frac{1}{2} \ln G(t) + \frac{1}{2\pi i}\int\limits_\gamma
\frac{\ln G(\tau)}{\tau-t}d\tau\right\}.
\]

So, using the Privalov Theorem, we have $\psi^{\pm}(t) \in H_\alpha.$

Let us estimate the norm
$$
\|K_2 x - K_3 x\| \leq \|(\psi^- - \psi^-_n)x^+\| + \|(\psi^+ - \psi^+_n)x^-\| +
$$
$$
+ \|(I-T_n)[\psi^-(t) U_\gamma(h(t,\tau) |\tau-t|^{-\eta}x(\tau))]\| = J_1+J_2+J_3.
\eqno (2.17) 
$$

It is easy to see that
\[
\|(\psi^- - \psi^-_n)x^+\| \leq \|(\psi^- - \psi^-_n)\| \|x^+\|;
\]
\[
\|(\psi^+ - \psi^+_n)x^-\| \leq \|(\psi^+ - \psi^+_n)\| \|x^-\|;
\]
\[
\|\psi^{\pm} - \psi^{\pm}_n\|_{C(\gamma)} \leq A n^{-\alpha}.
\]

Repeating the proof of Bernstein inverse theorem, one can see that 
\[
\|\psi^{\pm} - \psi^{\pm}_n\| \leq A n^{-\alpha + \beta}.
\]

From Sohotzky $-$ Plemel formulas we receive the inequality 
\[
\|x^{\pm}\| = \|\pm\frac{1}{2}x + \frac{1}{2} S(x)\| \leq \frac{1}{2}\|x\| + \frac{1}{2}\|S(x)\| \leq A\|x\|.
\]

From two last inequalities we have 
$$
J_1 + J_2 \leq A n^{-(\alpha-\beta)}\|x\|.
\eqno (2.18)
$$

Let $x(t) \in H_\beta,$ $0<\beta<\alpha \leq 1.$ According the Lemma 1.2 the function
\[
v(t) = \int\limits_0^{2\pi}x(\tau)\left|\cot\frac{\tau-t}{2}\right|^{\eta}d\tau
\]
belongs to Zygmund class of functions if $\beta = \eta,$ or to Holder class $H_1$ if $\beta > \eta,$ or to Holder class $\beta+1-\eta$ if $\beta < \eta.$ 
 
It is easy to see that
\[
\int\limits_\gamma \frac{x(\tau)d\tau}{|\tau-t|^\eta} = \int\limits_0^{2\pi}
\frac{e^{i \sigma } ix(e^{i\sigma})d\sigma}{|e^{i\sigma}-e^{is}|^{\eta}}
= i\int\limits_0^{2\pi}
\frac{e^{i \sigma } x(e^{i\sigma})d\sigma}
{|e^{+i\frac{\sigma+s}{2}}(e^{-i\frac{\sigma+s}{2}}(e^{i\sigma}-e^{is}))|^{\eta}}=
\]
\[
=  i\int\limits_0^{2\pi}
\frac{e^{i \sigma } x(e^{i\sigma})d\sigma}
{|e^{i\frac{\sigma-s}{2}}-e^{-i\frac{\sigma-s}{2}}|^{\eta}}
=  i\int\limits_0^{2\pi}
\frac{e^{i \sigma } x(e^{i\sigma})d\sigma}
{|2\sin\frac{\sigma-s}{2}|^{\eta}}.
\]

Repeating the proof of the Lemma 1.2, we see that the function 
$g(t) = \int\limits_\gamma \frac{x(\tau)d\tau}{|\tau-t|^{\eta}}$ belongs to Zygmund class $Z$ if $\beta=\eta$ or Holder class $H_\mu,$ 
where $\mu = 1$ if $\beta > \eta$ and $\mu = \beta + 1 -\eta$ if $\beta < \eta.$

Let us put $\beta < \eta.$

 In this case function 
$g_1(t)=\psi^-(t) U_\gamma(h(t,\tau)|\tau-t|^{-\eta}x(\tau))$
belongs to the Holder class $H_\chi,$ $\chi = \min(\alpha, \beta+1-\eta).$ 

Using the smoothness of the function $g_1(t),$ we have
\[
\|(I-T_n)g_1(t)\|_{C(\gamma)} \leq A\frac{1}{n^{\chi}}\|x\|.
\]

Repeating the proof of Bernstein inverse theorem, one can see that 
\[
\|(I-T_n)g_1(t)\| \leq A \frac{1}{n^{\chi-\beta}}\|x\|.
\]

Let us put $\beta=\eta.$

 In this case function $g_1(t)$ belongs to the Zygmund class $Z.$ From this it follows that $g(t) \in H_\alpha$ for $\beta < \alpha < 1$ 
and 
\[
\|(I-T_n)g_1(t)\| \leq A \frac{1}{n^{\alpha-\beta}}\|x\|.
\]

Let us put $\beta > \eta.$ 

In this case function $g_1(t) \in H_\alpha$ and 
\[
\|(I-T_n)g_1(t)\| \leq A \frac{1}{n^{\alpha-\beta}}\|x\|.
\]

The cases $\beta < \eta,$ $\beta = \eta,$ $\beta > \eta$ are considered identicaly. For simplicity we consider the case when $\eta < \beta <\min(\alpha,1-\eta).$

In this case
$$
J_3 \leq A \frac{1}{n^{\alpha-\beta}} \|x\|.
\eqno (2.19)
$$

Collecting the inequalities (2.16) - (2.19), we receive the estimate
\[
\|K_2 x -K_3 x\| \leq A \frac{1}{n^{\alpha-\beta}} \|x\|.
\]

From this inequality and Banach Theorem follows that, for $n$ such that 
$q = A\|K_2^{-1}\| n^{-(\alpha-\beta)}<1$, the operator $K_3$ is continuously invertible and the estimate $\|K_3^{-1}\| \leq \|K_2^{-1}\| / (1-q)$ is valid. 

The operator $K_3$ has a left inverse operator on subspace $X_n$ and 
$$
\|K_3 x_n\| \geq \frac{1}{\|K_3^{-1}\|} \|x_n\|, \, x_n \in X_n.
$$

Let us consider the equation
$$
K_{3,n} \equiv P_n[\psi_n^- x_n^+ + \psi_n^+ x_n^- + l(t) U_\gamma(h(t,\tau)|\tau-t|^{-\eta} x(\tau))] = 
$$
$$
= P_n[l(t) f(t)], \ \ \ l(t)=\psi^-(t)/(a(t)+b(t)).
\eqno (2.20)
$$

It is easy to see that
\[
P_n[\psi_n^- x_n^+  + \psi_n^+ x_n^-] \equiv \psi_n^- x_n^+  + \psi_n^+ x_n^- 
\]
and
\[
\|P_n[D_n[l(t) U_\gamma(h(t,\tau)|\tau-t|^{-\eta} x_n(\tau))]] \leq 
A n^{-(\alpha-\beta)}\ln n \|x_n\|.
\]

So,
\[
\|K_3 x_n - K_{3,n} x_n\| \leq A n^{-(\alpha-\beta)}\ln n \|x_n\|
\]
and
\[
\|K_{3,n}, x_n\| = \|K_3 x_n -(K_3 x_n - K_{3,n} x_n)\| \geq
\]
\[
\geq \left(\frac{1}{\|K_3^{-1}\|} - A n^{-(\alpha-\beta)}\ln n\right) \|x_n\| \geq C\|x_n\|.
\]

From the Theorem about left inverse operator \cite{Kan} follows that, for $n$ such 
$q_1 = A\|K_3^{-1}\| n^{-(\alpha-\beta)}\ln n <1$, the operator $K_{3n}$ has left inverse operator. As the $X_n$ is $n$-dimensional space,
 the operator $K_{3,n}$ has continuously inverse operator $K_{3,n}^{-1}.$

Now we introduce the equation
\[
K_{4,n} x_n \equiv P_n[\psi^- x_n^+ + \psi^+ x_N^- + l(t) U(h(t,\tau)|\tau-t|^{\eta} x_n(\tau))] =
\]
\[
=  P_n[l(t) f(t)]
\]
and estimate the norm
\[
\|K_{3,n} x_n - K_{4,n} x_n\| \leq \|P_n[(\psi^- - \psi_n^-)x_n^+]\| +
\]
\[
+ \|P_n[(\psi^+ - \psi_n^+)x_n^-]\| = J_4 + J_5.
\]

It is easy to see that
\[
\|P_n[(\psi^- - \psi_n^-)x_n^+]\|_{C(\gamma)} \leq A n^{-\alpha}\ln n \|x_n\|,
\]
\[
\|P_n[(\psi^+ - \psi_n^+)x_n^-]\| \leq A n^{-\alpha}\ln n \|x_n\|.
\]

So,
\[
\|P_n[(\psi^{\mp} - \psi_n^{\mp})x_n^{\pm}]\| \leq A n^{-\alpha+\beta}\ln n \|x_n\|
\]
and
\[
\|K_{3,n} x_n - K_{4,n} x_n\| \leq A n^{-\alpha+\beta}\ln n \|x_n\|.
\]

From this inequality and Banach Theorem follows that, for $n$ such 
$q_2=A\|K_{3,n}^{-1}\|n^{-(\alpha-\beta)}\ln n < 1$, the operator $K_{4,n}$ has continuously invertable operator $K_{4,n}^{-1}$ and the estimate
$\|K_{4,n}^{-1}\| \leq \|K_{3,n}^{-1}\|/(1-q_2)$ is valid.

Let us consider the equations
$$
K_{5,n} x_n \equiv P_n[x_n^+(t) + G(t)x_n^-(t) + l(t) U(h(t,\tau)|\tau-t|^{-\eta} x(\tau))] =f_n(t),
$$
$$
K_{6,n} x_n \equiv P_n[a(t) x_n(t) + b(t) S_\gamma(x_n(\tau)) +
$$
$$
+ U_\gamma(h(t,\tau)|\tau-t|^{-\eta} x_n(\tau))] =f_n(t).
\eqno (2.21)
$$
Here  $ f_n(t)=  P_n[f(t)]$.

It is obviously that the equation $K_{5,n} x_n = f_n$ is equivalent to the equation 
$K_{4,n} x_n = P_n[l(t) f]$ and the equation $K_{6,n} x_n = f_n$ is equivalent to the equation $K_{5,n} x_n = f_n.$ So, the operator $K_{6,n}$ 
has the continuously invertable operator $K_{6,n}^{-1}$ and $\|K_{6,n}^{-1}\| \leq A.$ Here $f_n = P_n[f].$

The equation (2.21) is the operator form of collocation method for solution of the equation (2.1). Using  methods of functional analysis 
we receive the estimate 
$\|x^* - \bar x^*_n\| \leq A n^{-(\alpha-\beta)}\ln n,$ where $x^*$ and $\bar x^*_n$ are solutions of equations (2.1) and (2.21) respectively.

Let us estimate the norm
\[
\|K_{n} x_n - K_{6,n} x_n\| = 
\]
$$
= \left\|P_n\left[\sum\limits^{2n}_{k=0}\int\limits_{t'_k}^{t'_{k+1}}
\left(\frac{h(t,\tau) x_n(\tau)}{|\tau-t|^\eta} - \frac{h(t,t_k) x_n(t_k)}{|\tau-t|^\eta}\right)d\tau\right]\right\|.
\eqno (2.22)
$$

It is easy to see that 
\[
\|K_{n} x_n - K_{6,n} x_n\| \leq 
\]
\[
\leq \left\|P_n\left[\sum\limits^{2n}_{k=0}\int\limits_{t'_k}^{t'_{k+1}}
\frac{(h(t,\tau) - h(t,t_k)) x_n(\tau)}{|\tau-t|^\eta}d\tau\right]\right\| +
\]
$$
+ A \ln n \left\|\sum\limits^{2n}_{k=0}\int\limits_{t'_k}^{t'_{k+1}}
\frac{h(t,t_k) (x_n(\tau) - x_n(t_k))}{|\tau-t|^\eta}d\tau\right\| = J_6 + J_7.
\eqno (2.23)
$$

Obviously,
$$
\left\|P_n\left[\sum\limits^{2n}_{k=0}\int\limits_{t'_k}^{t'_{k+1}}
\frac{h(t,\tau) - h(t,t_k)}{|\tau-t|^\eta}x(\tau)d\tau\right]\right\|_{C(\gamma)} \leq
A\frac{\ln n}{n^\alpha}\|x_n\|.
\eqno (2.24)
$$

Then
$$
J_6 \leq \frac{A \ln n}{n^{\alpha-\beta}}\|x_n\|.
\eqno (2.25)
$$

Let
\[
g_2(t) = \sum\limits^{2n}_{k=0}\int\limits_{t'_k}^{t'_{k+1}}
\frac{h(t,t_k) (x_n(\tau) - x_n(t_k))}{|\tau-t|^\eta}d\tau.
\]

For simplisity we will restrict youself with the function
\[
g^*_2(s) = \int\limits^{2\pi}_0
\frac{(x_n(e^{i\sigma}) - \bar x_n(e^{i\sigma}))i e^{i\sigma}}
{\left|\sin \frac{\sigma-s}{2}\right|^\eta}d\sigma,
\]
where $\bar x_n(e^{i\sigma}) = x_n(e^{i s_k}),$ $\sigma \in [s'_{k-1}, s'_k),$ $k=0,1,\ldots,2n.$

It is easy to see that
\[
\|g^*_2(s)\|_{\tilde C[0,2\pi]} \leq A \frac{1}{n^\beta}\|x_n\|.
\]
The difference $g^*_2(s+\delta) - g^*_2(s)$ is estimated more hard:
\[
|g^*_2(s+\delta) - g^*_2(s)| =
\]
\[
= \left|\int\limits^{2\pi}_0
(x_n(e^{i\sigma}) - \bar x_n(e^{i\sigma}))i e^{i\sigma}
\left(\frac{1}{\left|\sin \frac{\sigma-s-\delta}{2}\right|^\eta} -
\frac{1}{\left|\sin \frac{\sigma-s}{2}\right|^\eta}\right)d\sigma\right| \leq
\]
\[
\leq \sum\limits^4_{j=1} \int\limits_{\Delta_j} |\psi(\sigma)| = J_{7,1} + \cdots + J_{7,4},
\]
where $\Delta_1 = [s-2\delta, s+ 2\delta],$ $\Delta_2 = [s+2\delta, s+ \pi - 2\delta],$
$\Delta_3 = [s+\pi -2\delta, s+ \pi + 2\delta],$ $\Delta_4 = [s+\pi +2\delta, s+ 2\pi - 2\delta].$

One can see that
\[
|J_{7,1}| \leq A \frac{1}{n^\beta} \delta^{1-\eta} \|x_n\|,
\]
\[
|J_{7,2}| \leq A \frac{1}{n^\beta} \delta^{1-\eta} \|x_n\|,
\]
\[
|J_{7,3}| \leq A \frac{1}{n^\beta} \delta^{1-\eta} \|x_n\|,
\]
\[
|J_{7,4}| \leq A \frac{1}{n^\beta} \delta^{1-\eta} \|x_n\|,
\]

So,
\[
|g^*_2(s+\delta) - g^*_2(s)| \leq A \frac{1}{n^\beta} \delta^{1-\eta} \|x_n\|
\]
and
$$
J_7 \leq A \frac{\ln n}{n^\beta} \|x_n\|.
\eqno (2.26)
$$

From (2.22) - (2.26) it follows that 
\[
\|K_{n} x_n - K_{6,n} x_n\| \leq A\left(\frac{1}{n^\beta} + \frac{1}{n^{\alpha-\beta}}\right)\ln n \|x_n\|.
\]

So, the operator $K_n$ has continuously invertible operator $K_n^{-1}$ with the norm 
$\|K_n^{-1}\| \leq A.$

Let $x^*(t)$ is a unique solution of the equation (2.1). So,
$$
P_n[a(t) x^*(t) + b(t) S_\gamma(x^*(\tau)) +U_\gamma(h(t,\tau)|\tau-t|^{-\eta} x^*(\tau))] = P_n[f(t)].
\eqno (2.27)
$$

The identity (2.27) we can rewrite as 
\[
P_n[a(t) x^* + b(t) S_\gamma(x^*(\tau)) + \sum\limits^{2n}_{k=0} h(t,t_k) x^*(t_k)
\int\limits^{t'_{k+1}}_{t'_k}|\tau-t|^{-\eta}d\tau] =
\]
\[
= P_n[f(t) - U_\gamma(h(t,\tau)|\tau-t|^{-\eta} x^*(\tau))+
\]
$$
+\sum\limits^{2n}_{k=0} h(t,t_k) x^*(t_k)
\int\limits^{t'_{k+1}}_{t'_k}|\tau-t|^{-\eta}d\tau].
\eqno (2.28)
$$

Let $x^*_n(t)$ is a unique solution of the equation (2.6).

Substracted the identity 
\[
P_n[a(t) x^*_n(t) + b(t) S_{\gamma}(x^*_n(\tau)) + \sum\limits^{2n}_{k=0} h(t,t_k) x^*_n(t_k)\int\limits^{t'_{k+1}}_{t'_k}|\tau-t|^{-\eta}d\tau] =
\]
\[
=  P_n[f(t)]
\]
from the identity (2.28) we have
\[
P_n[a(t)(x^*(t) - x^*_n(t)) + b(t) S_{\gamma}(x^*(\tau) - x^*_n(\tau)) +
\]
\[
+\sum\limits^{2n}_{k=0} (h(t,t_k)(x^*(t_k) - x^*_n(t_k))\int\limits^{t'_{k+1}}_{t'_k}|\tau-t|^{-\eta}d\tau] =
\]
$$
= P_n\left[\int\limits_\gamma \frac{h(t,\tau) x^*(\tau)}{|\tau-t|^{-\eta}}d\tau -
\sum\limits^{2n}_{k=0} h(t,t_k)x^*(t_k) 
\int\limits^{t'_{k+1}}_{t'_k}|\tau-t|^{-\eta}d\tau\right].
\eqno (2.29)
$$

From the Theorem 1.1 it follows that a unique solution $x^*(t)$ of equation (2.1) belongs to the Holder class $H_\alpha.$

In this case
\[
\left\|P_n\left[\int\limits_\gamma 
\frac{h(t,\tau) x^*(\tau)}{|\tau-t|^{-\eta}}d\tau -
\sum\limits^{2n}_{k=0} h(t,t_k)x^*(t_k) 
\int\limits^{t'_{k+1}}_{t'_k}|\tau-t|^{-\eta}d\tau\right]\right\| \leq
\]
\[
\leq A \frac{\ln n}{n^{\alpha-\beta}}.
\]

Applied the inverse operator $K_n^{-1}$ to both side of the equation (2.29) we have
$$
\|x^* - x^*_n\| \leq A n^{-(\alpha-\beta)}\ln n.
\eqno (2.30)
$$

Let us prove that this estimate can not be improved. For this aim we consider
the equation
$$
Kx \equiv \frac{1}{2\pi} \int\limits_0^{2\pi}
\varphi(\tau) \cot\frac{\sigma -s}{2}d\tau = f(s).
\eqno (2.31)
$$

Well known \cite{Gakh} that the solution of this equation is
\[
K^{-1}f \equiv x(t) = \frac{1}{2\pi} \int\limits_0^{2\pi}
f(\sigma) \cot\frac{\sigma -s}{2}d\sigma.
\]

Let $n$ be a integer. Let us solve the equation (2.31) with collocation
method, using arbitrary knots $t_k,$ $k=1,2,\ldots,n.$ Without loss of
generality we can put $t_1 =0.$ Let us introduce the function
$f^*(s) = \min\limits_{k}|s-t_k|^\alpha,$ $0 \leq s \leq 2\pi.$ It easy to see
that the approximate solution $x^*_n(s)$ of the equation (2.31) received with
collocation method based on knots $t_k,$ $k=1,2,\ldots,n,$ is equal to zero
$(x^*_n(s) \equiv 0).$

From the Theorem of inverse operator \cite{Kan}  follows that\\
$\|x\| \geq \frac{1}{m}\|f\|,$ where $m =\|K\|.$ So
$\|x^*(s)\| \geq m^{-1}\|f^*(s)\|,$ where $x^*(s)$ is a solution of the
equation (2.31).

It is easy to see that the function $f^*(s),$ $0 \leq \sigma \leq 2\pi,$
belongs to the Holder class $H_\alpha$ and $\|f^*(s)\|
\geq An^{-\alpha+\beta}.$

So,
$$
\|x^*(s) -x^*_{n}(s)\| = \|x^*(s)\| \geq An^{-\alpha+\beta}.
\eqno (2.32)
$$

Correctness of  Theorem follows from the comparison
of the  inequalities (2.30) and (2.32).

{\bf Proof of Theorem 2.5.}
In proving of  Theorem 2.4 it was noticed that the
equations (2.1) and (2.14) are equivalent.
The equation (2.14) is equivalent to the  Riemann boundary value problem
$$
K_7 x\equiv \psi^-(t)x^+(t) - \psi^+(t)x^-(t)+
\psi^-(t)U_{\gamma}(h(t,\tau){\mid\tau-t\mid}^{-\eta}
x(\tau))=$$
$$=\psi^-(t)f(t),
\eqno (2.33)
$$
where
\[
\psi(z)=\exp\left( \frac{1}{2\pi i}
\int\limits_{\gamma}[\ln G(\tau)](\tau-z)^{-1}d\tau\right).
\]

Hence, the operator $K_7$ has the continuous inverse operator $K_7^{-1}$ with the norm $\|K_7^{-1}\|\le A.$

Let us approximate  functions $\psi^+(t)$ and $\psi^-(t)$
by interpolation polynomials $\psi^+_n(t)$ and $\psi^-_n(t),$
constructed on the knots $t_k=\exp\{i s_k \},$ $s_k = 2k\pi/(2n+1),$
$k=0,1,\ldots,n.$ It is easy to see that the functions
$\psi^+_n(t)$ and $\psi^-_n(t)$ can be represented in the forms
$\psi^+_n(t) = \sum\limits^n_{k=0} \gamma_k t^k$ and
$\psi^-_n(t) = \sum\limits^0_{k=-n} \gamma_k t^k.$

Let us introduce the equation
\[
K_{8} x \equiv \psi^-_n(t) x^+(t) - \psi^+_n(t)x^-(t)+
\]
$$
+\psi^-(t)U_{\gamma}[P_{n}^{\tau}[h(t,\tau)]{\mid\tau-t\mid}^{-\eta}
x(\tau)]
=\psi^-(t)f(t).
\eqno (2.34)
$$

Easy to see that
\[
\|K_7 x - K_8 x\| \leq \|\psi^- - \psi^-_n\| \|x^+\| + 
\|\psi^+ - \psi^+_n\| \|x^-\| \leq A\frac{1}{n^{r+\alpha-\beta}} \|x\|.
\]

From this inequality and Banach Theorem it follows that, for $n$ such that $q = A\|K_7^{-1}\| n^{-(r+\alpha-\beta)} < 1$, the operator $K_8$ 
has continuous invertible operator  $K_8^{-1}$ with the norm  $\|K_8^{-1}\| \leq \|K_7^{-1}\|/(1-q).$

Method of mechanical quadrature for solution of the equation (2.34) can be written as
\[
K_{8,n} x_n \equiv P_n [\psi^-_n(t) x^+_n(t) - \psi^+_n(t)x^-_n(t)+
\]
$$
+\psi^-(t)U_{\gamma}[P_n^\tau [h(t,\tau)]{\mid\tau-t\mid}^{-\eta}
x_n(\tau)]] = P_n[\psi^-(t)f(t)].
\eqno (2.35)
$$

Using Kantorovich theory of approximate methods of analysis, one can prove that,
under condition $q=An^{-(r+\alpha-\beta)}\ln n <1,$ the equation (2.35) has a
unique solution $x^*_{4,n}$ and the estimate $\|x^*-x^*_{8,n}\| \le An^{-(r+
\alpha-\beta)}\ln n$ is valid.

Let us consider the sequence of equations 
$$
K_{9,n} x_n \equiv P_n [\psi^-(t) x^+_n(t) - \psi^+(t)x^-_n(t)+
$$
$$
+\psi^-(t)U_{\gamma}[P_n^\tau [h(t,\tau)]{\mid\tau-t\mid}^{-\eta}
x_n(\tau)]] =P_n[\psi^-(t)f(t)];
\eqno (2.36)
$$
$$
K_{10,n} x_n \equiv P_n [x^+_n(t) G(t)x^-_n(t)+
$$
$$
+U_{\gamma}[P_n^\tau [h(t,\tau)]{\mid\tau-t\mid}^{-\eta}
x_n(\tau)]] =P_n[f(t)];
\eqno (2.37)
$$
$$
K_{11,n} x_n \equiv P_n [a(t)x_n(t) + b(t) S_{\gamma}(x_n(\tau)+
$$
$$
+ U_{\gamma}[P_n^\tau [h(t,\tau)]{\mid\tau-t\mid}^{-\eta}
x_n(\tau)]] =P_n[f(t)].
\eqno (2.38)
$$

Using the arguments repeatedly cited above, we prove the unique solvability of the equations (2.36) - (2.38).

So, the equation $K_{11,n}x_n = f_n, \ f_n=P_n[f(t)],$ has a unique solution $x^*_{11,n}$ and the
estimate
$$
\|x^* - x^*_{11,n}\| \le An^{-(r+\alpha-\beta)}\ln n
\eqno (2.39)
$$
is valid.

The equation $K_{11,n}x_n = f_n$ is the equation (2.7).

Theorem is proved.

Now we will prove that this estimation can not be improved. Let $n$ is a
integer. Let $t_k,$ $k=1,2,\ldots,n,$ are arbitrary knots,
$0=t_1<t_2<\cdots <t_n<2\pi.$ Let us introduce the function $f^*(s)$ which is
equal to $\frac{((s-t_k)(t_{k+1}-s))^{r+\alpha}}{h_k^{r+\alpha}},$ $h_k = t_{k+1}-t_k,$ on
each segment $[t_k, t_{k+1}],$ $k=1,2,\ldots,n-1,$ where $t_{n+1}=2\pi.$ It is
easy to see that $f^*(s) \in W^rH_\alpha(M)$ and
$\|f^*(s)\|_\beta = An^{-(r+\alpha-\beta)}.$

Let us consider the equation
$$
\frac{1}{2\pi} \int\limits_0^{2\pi} x(\tau) \cot\frac{\sigma-s}{2} d\sigma
=f^*(s).
\eqno (2.40)
$$

Well known that a solution of this equation is
\[
x^*(s)=\frac{1}{2\pi} \int\limits_0^{2\pi} f^*(\sigma) \cot\frac{\sigma-s}{2}
d\sigma.
\]

Solving the equation (2.40) by collocation methods with knots $t_k,$
$k=1,2,\ldots,n,$ we receive approximate solution $x^*_n(t) \equiv 0.$
Repeating arguments given behind, we see that
$\|x^* - x_n^*\|_\beta = \|x^* \|_\beta \geq An^{-(r+\alpha -\beta)}\ln n.$

Comparing this estimation with (2.39) we finish the proof of the Theorem.

Theorem is proved.

{\bf Proof of Theorem 2.6.}
Let us consider  the collocation method for the
equation (2.2). In the operator form the collocation method is
written by the expression
$$
L_n^1 x_n\equiv P_n[a(t)x_n(t)+S_{\gamma}(h(t,\tau)x_n(\tau))]=
P_n[f(t)], \eqno (2.41)
$$
where $x_n(t)$ is the polynomial defined in (2.3).

The equations (2.2) and (2.41) are equivalent to the Riemann boundary
value problems
\[L^{(2)} x\equiv Vx+Wx=y
\]
and
\[
L^{(2)}_n x_n\equiv \tilde V_n x_n+
\tilde W_n\tilde x_n=y_n,
\]
where
\[ Vx=\psi^-x^+-\psi^+x^-,
Wx=\psi^-S_{\gamma}((h(t,\tau)-h(t,t))x(\tau)),
\]
\[
\tilde V_n x_n=P_n[Vx_n],
 \tilde W_n x_n=P_n[Wx_n],
y=\psi^-f, y_n=P_n[y],
\]
\[
 G(t)=(a(t)-b(t))/(a(t)+b(t)), \quad
b(t)=h(t,t),
\]
\[
 \psi(z)=\exp\left\{\frac{1}{2\pi i}\int\limits_{\gamma}
(\ln G(\tau)){(\tau-z)}^{-1}d\tau\right\} .
\]

Let us introduce the polynomial
\[
\varphi_n(t)=V_n x_n+T^t_{[n/3]}[\psi^-(t)]S_{\gamma}
((\tilde h(t,\tau)-\tilde h(t,t))x_n(\tau)),
\]
where
\[
V_n x_n=\psi_n^{-}x_n^+-\psi_n^+x_n^-,
\psi_n=T_n(\psi), \tilde h(t,\tau)=T_{[n/3]}^t T_{[n/3]}^{\tau}
[h(t,\tau)].
\]

It is obvious, that
\[\Vert V_n x_n-V x_n\Vert_\beta\le A\Vert x_n\Vert_\beta
\max(E_n(a), E_n(\psi), E_n(b))n^{\beta}.
\]

Since
\[\Vert h(t,\tau)-
\tilde h(t,\tau) \Vert_C \le A\ln n \max({E_n^t[h(t,\tau)],
E_n^{\tau}[h(t,\tau]}),
\]
then the function
$\eta (t,\tau)=
h(t,\tau)-\tilde h(t,\tau)$
belongs to the Holder class with
respect  to $\tau$ with index $1/\ln n$ and with the coefficient
$A \ln n E_n^*(h),$ where
$E_n^*(h)= \ln n \max(E_n^t(h(t,\tau)),E_n^{\tau}(h(t,\tau))).$

Let us prove this statement. Consider the function $\eta (t,\tau) = h(t,\tau) - \tilde h(t,\tau),$ where $t$ is fixed variable. It is easy to see that
\[
E_m(\eta (t,\tau)) \leq A \ln n \max(E^t_m[h(t,\tau)], E^\tau_m[h(t,\tau)])
\]
for $1 \leq m \leq n$ and
\[
E_m(\eta (t,\tau)) \leq A \ln m \max(E^t_m[h(t,\tau)], E^\tau_m[h(t,\tau)])
\]
for $n < m < \infty.$

These inequalities we can rewrite as  
\[
E_m(\eta (t,\tau)) \leq A \ln m \max(E^t_m[h(t,\tau)], E^\tau_m[h(t,\tau)])
\]
for  $n < m < \infty$  and
\[
E_m(\eta (t,\tau)) \leq \frac{A}{m^\gamma}(m^\gamma \ln n \max(E^t_m[h(t,\tau)], 
E^\tau_m[h(t,\tau)])) \leq
\]
\[
\leq \frac{A}{m^\gamma}(n^\gamma \ln n \max(E^t_m[h(t,\tau)], 
E^\tau_m[h(t,\tau)]))
\]
for $1 \leq m \leq n$.

Here $\gamma$ is a arbitrary number, $0 < \gamma < 1.$ 

Repeating the proof of Bernstein inverse theorem we see that
\[
\omega(\eta (t,\tau),\delta) 
\leq \left(\frac{2^{1-\gamma}}{2^{1-\gamma}-1}(1+2^{\gamma}) + 2\frac{1+2^\gamma)}{1-2^{-\gamma}}\right)\times
\]
\[
\times A n^\gamma \ln n
\max(E^t_m[h(t,\tau)], E^\tau_m[h(t,\tau)])\delta^\gamma.
\]

Let us put $\gamma = 1/\ln n.$ In result we have
\[
\omega(\eta (t,\tau),\delta) \leq A n^\gamma \ln^2 n
\max(E^t_m[h(t,\tau)], E^\tau_m[h(t,\tau)])\delta^{1/\ln n}.
\]

Therefore
\[\Vert S_{\gamma}((h,\tau)-\tilde h(t,\tau)-(h(t,t)-
\tilde h(t,t)))x_n(\tau))\Vert_\beta \le
\]
\[
\le  A\Vert x_n\Vert_\beta n^{\beta} \ln^2 n
E_n^*(h).
\]

It follows from this inequality and recieved above estimation of the
norm $\Vert V_n x_n-V x_n\Vert$ that the operator
$L^{(2)}_n$ is continuously invertible for $n$ such that
\[q_1=
\]
\[
=
A \ln^3 n^{\beta} \max(E_n(a),E_n(b), E_n(\psi), E_n^t(h(t,\tau)),
E_n^{\tau}(h(t,\tau)))<1.$$

In addition
$\Vert x^*-x_{1n}^*\Vert\le A_2(q_1+E_n(f)),$
 where
$x_{1n}^*$ is a solution of the equation (2.41).

As marked above the equations $L_n^{(2)} x_n = y_n$ and $L_n^{(1)} x_n = y_n$ are equivalent. So, the operator $L_n^{(1)}$ has the continuous invertable operator $(L_n^{(1)})^{-1}$ and
$\|(L_n^{(1)})^{-1}\| \leq A.$

It is easy to see that
\[
\|L_n^{(1)} x_n - L_n x_n\| \leq A E^\tau_n[h(t,\tau)] \ln^2 n \|x_n\|.
\]

From this inequality and Banach Theorem follows that the operator $L_n$ has the continuous invertable operator $L_n^{-1}$ with the 
norm $\|(L_n^{(1)})^{-1}\| \leq A$ and the following estimate $\|x^*_{1,n} - x^*_n\| \leq A E^\tau_n[h(t,\tau)] \ln^2 n$ is valid. Here $x^*_n$ is 
a unique solution of the equation (2.7).

Using the estimations for $\Vert x^*-x_{1n}^*\Vert_\beta$ and for
$\Vert x_{1n}^*-x_n^*\Vert_\beta$, we verify in correctness of
Theorem.

{\bf Proof of  Theorem 2.7.}

It was shown in proving of  preceding Theorem that the
operator $L^{(2)}_n$ is continuously invertible and
$\Vert [L^{(2)}_n]^{-1} \Vert \le A$ for $n$ such that $q_1<1.$
Repeating arguments of proof of the Theorem 2.5, one can see that the
operator
\[
L^{(3)}_n x_n \equiv \bar P_n \left[a(t)x_n(t)+\frac{1}{\pi i}
\int\limits_{\gamma} \frac{h(t,\tau)x_n(\tau)}{\tau-t}d\tau\right]
\]
is
continuously invertible and
$\|L^{(3)}_n\|^{-1} \leq A$ for $n$ such that\\
 $q_1 < 1.$
 
The correctness of the identity
\[\overline P_n\left[\frac{1}{\pi i}\int_{\gamma}P_n^\tau \left
[\frac{h(t,\tau)x_n(\tau)}{\tau-t}\right]d\tau\right]\equiv
\]
\[
\equiv
\overline P_n\left[\frac{1}{\pi i}\int_{\gamma}\frac{P_n^\tau
[h(t,\tau)x_n(\tau)]}{\tau-t}d\tau\right]\equiv
\]
$$
\equiv\overline P_n\left[\frac{1}{\pi i}\int_{\gamma}
\frac{P_n^\tau[h(t,\tau)]x_n(\tau)}{\tau-t}d\tau\right] 
\eqno (2.42)
$$
follows from the results on the quadrature rules of the highest
algebraic exactness. 

Using this identity we can show that
\[\Vert L^{(3)}_n - L_n\Vert_\beta \le An^{\beta}\ln^2 n(E_n^t(h(t,\tau))+
E_n^{\tau}(h(t,\tau))).
\]

Therefore the operator $L_n$ is continuously invertible for $n$
such that $q=An^{\beta} \ln^2 n(E_n(a(t))+E_n(h(t,t))+E_n(\psi) +
E_n^t (h(t,\tau)) + E_n^{\tau}(h(t,\tau))))<1.$
The estimation of error $\Vert x^*-x_n^*\Vert_\beta$ is evaluated
just as in preceding Theorems. Theorem is proved.

{\bf Proof of  Theorem 2.8.}
Joined proofs of the Theorem 2.7 and the Theorem 2.5 we are verified in the
correctness of the Theorem 2.8.

\begin{center}
{\bf 3.    Approximate Solution of Singular Integral Equations on
Closed Paths of Integration (Basis in Space  $ L_2 $)}
\end{center}

Let us extend investigation of approximate methods for solution of
singular integral equations as
\[
Kx \equiv a(t)x(t)+b(t)S_{\gamma}(x)+U_{\gamma}(h(t,\tau){\mid
\tau-t \mid}^{-\eta} x(\tau))=
\]
\[
=f(t) \eqno (3.1)
$$
and
$$
Lx \equiv a(t)x(t)+S_{\gamma}(h(t, \tau)x(\tau))=f(t). \eqno (3.2)
$$

According to accepted beyond designations
\[ S_{\gamma}x=\frac{1}{\pi i} \int\limits_{\gamma}\frac{x(\tau)d\tau}
{\tau-t},\hskip 10 pt  U_{\gamma}(h(t,\tau) x(\tau))=\frac{1}{2\pi i}
\int\limits_{\gamma}h(t,\tau)x(\tau)d\tau, 
\]
where
$ \gamma $ is a unit circle with the center in  origin of
coordinates.

 The verification of
the suggested below calculating schemes is carried out in the
functions space $ X=L_2(\gamma) $ with scalar product
\[ (f_1, f_2)= \frac{1}{2\pi}\int\limits_0^{2\pi}f_1(e^{is}) \overline
{f_2(e^{is})}ds 
\]
and its subspace $ X_n $ consisting of the
polynomials 
$$
x_n(t)=\sum_{k=-n}^n \alpha_k t^k. \eqno (3.3)
$$

In this paragraph we sake all designations, which are introduced
in the previous paragraph.

{\bf 3.1. Basis Statements}

     The approximate solution of the equations (3.1) we will seek
in the form of the polynomial (3.3), the coefficients $ {\alpha_k} $
of which are determined from the system of linear algebraic equations
written in the operator form  as
\[ K_n x_n \equiv P_n [a(t) x_n(t)+b(t) S_{\gamma} (x_n)+U_{\gamma}
(P_n^\tau[h(t,\tau) d(t,\tau) x_n(\tau)])]= 
\]
$$
= P_n[f(t)], \eqno (3.4)
$$
where
\[
d(t,\tau) =
\left \{ \begin{array}{ccc}
\mid \tau-t \mid^{-\eta} \quad {\rm for} \quad \mid \sigma-s \mid \ge 2\pi/(2n+1),\\
\mid e^{i2\pi/(2n+1)}-1 \mid^{-\eta} \quad {\rm for} \quad
 \mid \sigma-s \mid < 2\pi/(2n+1),\\
\end{array} \right.
\]
$  \tau=e^{i\sigma}, t=e^{is}; P_n $ is the projector of interpolation onto the set of
trigonometrical interpolated polynomials  constructed on the knots
$ t_k=e^{is_k}, s_k=2k\pi/(2n+1), k=0,1,\dots,2n.$

{\bf Theorem 3.1} \cite{Boy4}, \cite{Boy6}, \cite{Boy16}, \cite{Boy25}.  Let the operator $ K \in [X,X] $  has the
linear invertible operator and  functions $ a,b \in C[0,2\pi], h \in C([0,2\pi]^2).$ Then for
$ n $ such that $ q=A[\omega(a; n^{-1/2})+\omega(b; n^{-1/2})+n^{-1/2+\varepsilon}+
[\omega(h;n^{-1})]^{(1-\eta)/(1+\eta)}]<1, $ the equation (3.4) is uniquely
solvable for any right-hand side and the estimate
$ \Vert x^*-x_n^* \Vert \le A[q+\omega(f; n^{-1})] $ is valid, where $ x^* $
and $ x_n^* $  are solutions of equations (3.1) and (3.4), $\varepsilon$ is
arbitrary small number $(\varepsilon >0).$

     The approximate solution of the  equation (3.2) we will seek
in the form of the polynomial (3.3),  coefficients $ {\alpha_k} $
of which are determined from the system of linear algebraic
equations
$$
L_n x_n \equiv \overline P_n \left[ a(t) x_n(t)+
\frac {1}{\pi i} \int_{\gamma} P_n^\tau \left [\frac{h(t,\tau)
x_n(\tau)} {\tau-t} \right] d\tau \right]=\overline P_n [f(t)].
\eqno (3.5)
$$

{\bf Theorem 3.2} \cite{Boy6}, \cite{Boy7}, \cite{Boy16}, \cite{Boy25}.  Let the operator $ L $ has a linear inverse one
in the space $ X. $ Then for $ n $ such that
$ q=A(E_n(a)+E_n(b)+n^{-1/2+\varepsilon}+E_n^t(h(t,\tau))+E_n^{\tau}(h(t,\tau)))\ln^3n<1, $
the system of equations (3.5) is uniquely solvable for any
right-hand side and the estimate
$ \Vert x^*-x_n^* \Vert \le A(q+E_n(f)) $ is valid, where $ x^* $ and
$ x_n^* $ are solutions of the equations (3.2) and (3.5), $\varepsilon$ is
arbitrary small number $(\varepsilon >0).$

     The approximate solution of the equation (3.2) is sought in the
form of the polynomial (3.3),  coefficients $ {\alpha_k} $ of which
are determined from the system of equations
$$
\overline L_n x_n \equiv \overline P_n [a(t)x_n(t)+S_{\gamma}(P_n^{\tau}
[h(t, \tau)] x_n(\tau))]=\overline P_n[f(t)]. \eqno (3.6)
$$

{\bf Theorem 3.3}  \cite{Boy16}, \cite{Boy25}. Let the operator $ L $ is continuously
invertible and the functions $ a, f \in W^r, h \in W^{r,r}
(r=1,2, \dots.) $ Then for $ n $, such that $ q=An^{-r} \ln^3n<1 $, the
system of equations (3.6) is uniquely solvable for any
right-hand side and the estimate $ \Vert x^*-x_n^* \Vert \le A n^{-r}
\ln n $ is valid, where $ x^* $ and $ x_n^* $ are solutions of equations (3.2)
and (3.6).

{\it Note.} These results are diffused to the space $L_p(\gamma),$
$1 \leq p < \infty,$ in the papers I.V. Boykov and  Zhechev \cite{Boy32} - \cite{Boy36}.

{\bf 3.2.  Proofs of Theorems}

{\bf  Proof of Theorem 3.1.} Let us introduce the equation
\[
K_1 x \equiv \tilde a_m(t) x(t)+\tilde b_m S_{\gamma}(x)+
U_{\gamma}(h(t, \tau)d^*(t, \tau) x(\tau))=
\]
$$
=f(t), \eqno (3.7)
$$
where $ a_m(t)=T_m[a(t)], b_m(t)=T_m[b(t)], $
$ d^*(t,\tau)=|\tau-t|^{-\eta} $ for $ |\tau-t| \ge
\rho, d^*(t,\tau)=\rho^{-\eta} $ for $ |\tau-t|<\rho, \rho  $
is  positive number, $ \rho \ge 2 /\pi (2n+1). $
The numbers $ \rho $ and $ m $ are fixed below. It follows from Banach
Theorem that for $ m, \rho $ such that $ q_1=A(\rho^{1-\eta}+\omega(a;
m^{-1})+\omega(b;m^{-1}))<1, $ the operator $ K_1 $ has inverse one with the
norm $ \Vert K_1^{-1} \Vert \le \Vert K^{-1} \Vert/(1-q_1). $

The method of mechanics quadratures in the operator form for the equation (3.7) is
\[  K_{1,n} x_n \equiv P_n[\tilde a_m(t) x_n(t)+\tilde b_m(t)
S_{\gamma} (x_n)+ 
\]
$$
 + U_\gamma(P_{n}^{\tau}[h(t, \tau)d^*(t, \tau) x_n(\tau)])]=
P_n[f(t)]. \eqno (3.8)
$$

     The equations (3.7) and (3.8) are equivalent to the Riemann boundary value
problems
$ K_2 x=Vx+W^*x=y $ and $  K_{2,n} x_n= V_nx_n + \tilde W_n^*
x_n=y_n, $
where
\[ Vx=\psi^-x^+ -\psi^+x^-,
\tilde V_n x_n=P_n[V x_n], 
\]
\[
W^*x=\psi^-(t) U_{\gamma}(h(t,\tau) d^*(t,\tau) x(\tau)), 
\]
\[
\tilde W_n^* x_n=P_n[\psi^-(t) U_{\gamma}(P_{n}^{\tau}[h(t,\tau)d^*(t,\tau)
x_n(\tau)])], 
\]
\[
 \psi(z)=\exp \left\{\frac{1}{2\pi i}
\int\limits_{\gamma}\frac{ \ln((a_m(\tau)-b_m(\tau))/(a_m(\tau)+b_m(\tau)))}{
\tau-z}d\tau\right\}, 
\]
\[
y=\psi^-f, y_n=P_n[\psi^-f].
\]

Input the polynomial
\[ \tilde \varphi_n(t)=\psi_n^-(t)x_n^+(t) - \psi_n^+(t)x_n^-(t) +
T_n [U_{\gamma} (\psi^-(t) h(t,\tau)d^*(t,\tau) x_n(\tau))], 
\]
where $ \psi_n=T_n \psi. $

It is easy to see that
$$
\Vert K_2 x_n - \tilde \varphi_n \Vert + \Vert P_n K_2 x_n - \tilde
\varphi_n \Vert=A\left(\frac{m^{1-\varepsilon}}{n^{1-\varepsilon}}+\frac{\omega(h; n^{-1})}{\rho^{2\eta}}\right) \Vert x_n \Vert, \eqno (3.9)
$$
where $ \varepsilon $ is an arbitrary small number $ (0<\varepsilon<1). $
Necessity of introduction the $ \varepsilon $ follows from Privalov
Theorem \cite{Gakh}, as $ \tilde a_m, \tilde b_m $ are belong to the class
$ H_1. $ 

It is evident,
$$
\Vert  K_{2,n} x_n - P_n K_{2,n} x_n \Vert \le A\omega(h; n^{-1})
\Vert x_n \Vert /\rho^{2\eta}. \eqno (3.10)
$$

It follows from (3.9)-(3.10) (assuming $ \rho=(\omega(h,n^{-1}))^{1/(1+
\eta)},$ $\varepsilon=1/\ln n,$ $m=n^{1/2}) $ and Kantorovich  theory of approximate
methods of analysis, that, for $ n $ such that
\[
q_2=
\]
\[=A[\omega(a;n^{-1/2})+\omega(b, n^{-1/2})+E_n(\psi^{\pm})+
[\omega(h; n^{-1})]^{(1-\eta)/(1+\eta)}]<1,
\]
the linear operator $ K_{2,n}^{-1} $ with the norm $ \Vert 
K_{2,n}^{-1} \Vert \le \Vert K_2^{-1} \Vert/(1-q_2) $ exists.

In addition
\[ \Vert x^*-\tilde x^* \Vert \le A(\omega(a, \frac{1}{n^{1/2}})+
\]
\[+\omega(b; \frac{1}{n^{1/2}})+
E_n(\psi^{\pm})+[\omega(h; \frac{1}{n})]^{(1-\eta)/(1+\eta)}+\omega(f; n^{-1})), $$
where $ x^* $ and $  x^*_n $ are solutions  of equations (3.7)
and $  K_{2,n} x_n=y_n. $

Existence of the linear operator $  K_{1,n}^{-1} $
follows from equivalence of equations
$  K_{2,n} x_n=y_n $ and (3.8). 

Let us estimate its norm. Since
\[ \Vert x_n^* \Vert \le  \Vert
 K_{2,n}^{-1} \Vert \mid \psi^-\mid \Vert P_n[f] \Vert, 
 \]
then
$ \Vert  K_{1,n}^{-1} \Vert \le \Vert  K_{2,n}^{-1} \Vert
\max\limits_t |\psi^{-}(t)|.$

It is easy to see that
$$
\Vert  K_n x_n- K_{1,n} x_n \Vert \le A[[\mid a-\tilde a_m
\mid + \mid b-\tilde b_m \mid] \Vert x_n \Vert +I_1], \eqno (3.11)
$$
where
\[
 I_1=\Vert P_n[\frac {1} {2\pi i} \int\limits_{\gamma} P_n^\tau [h(t,\tau)
x_n(\tau) b(t, \tau)] d\tau] \Vert,
\]
$ b(t, \tau)=|d^*(t, \tau)-d(t,\tau)|. 
$

$ I_1 $ may be estimated as
\[I_1 \le \max\limits_s \left \{ \frac{1}{2\pi} \int\limits_0^{2\pi}\mid P_n^\sigma
[h(s,\sigma)\rm{sgn}[b(s,\sigma)] \mid b(s, \sigma) \mid^{1/2} \mid^2
d\sigma \right\}^{1/2}\times 
\]
\[
\times \left\{\frac{1}{2\pi}\int\limits_0^{2\pi} \left\{\frac{1}{2\pi}
\int\limits_0^
{2\pi} \mid P_n^\sigma[\tilde x_n(\sigma) e^{i\sigma} P_n^s[\mid
b(s,\sigma) \mid^{1/2}]] \mid^2 \right\}d\tau \right\}^{1/2}=
\]
$$
=I_2 I_3,
\eqno (3.12)
$$
where
$ h(s, \sigma)=h(e^{is}, e^{i \sigma}),
b(s, \sigma)=b(e^{is}, e^{i\sigma}),
x_n(\sigma)=x_n(e^{i\sigma}). $

Having designated the value $ [\rho(2n+1)/2\pi]+1 $ by $ v $ we
estimate $ I_2: $
$$
I_2^2 \le A \frac{1}{2n+1} \sum_{k=1}^v \left(\frac{1}{s_k}
\right)^{\eta} \le A \left[\frac{1}{(2n+1)^{1-\eta}}+ \rho^{1-\eta}
\right]. \eqno (3.13)
$$

It is easy to see that
\[
 I_3^2=\frac{1}{(2n+1)^2} \sum_{k=0}^{2n} \sum_{i=0}^{2n} \mid
x_n(s_i) \mid^2 \mid b(s_k, s_i) \mid \le 
\]
$$
\le \frac{1}{(2n+1)^2} \sum_{i=0}^{2n} \mid x_n(s_i) \mid^2
\sum_{k=0}^{2n} \mid b(s_k, s_i) \mid \le A \Vert x_n \Vert^2
\rho^{1-\eta}. \eqno (3.14)
$$

The validity of Theorem follows from the estimations (3.11) - (3.14)
for reduced beyond values $ m, \rho $ and Banach Theorem.

{\bf Proof of Theorem 3.2.}

The collocation method for equation (3.2) may be written in the form
of the following expression
$$
L_{1,n} x_n \equiv \overline P_n[a(t) x_n(t)+
$$
$$
+b(t) S_{\gamma}
(x_n(\tau))+S_{\gamma}((h(t,\tau)-h(t,t))x_n(\tau))]=
$$
$$
=\overline P_n
[f(t)], \eqno (3.15)
$$
where $ b(t)=h(t,t). $

Let us represent the equations (3.2) and (3.15) in the form of
Riemann boundary value problems $ L_2 x\equiv Vx+Wx=y $ and $ L_{2,n} x_n\equiv\tilde
V_n x_n+\tilde W_n x_n=y_n,$ the operators $V,W ,\tilde V_n ,\tilde W_n $ and the
elements $ y $ and $ y_n $ of which are defined beyond in proving  of the Theorem 3.1.

Introduce the polynomial 
\[\varphi_n(t)=\psi_n^-x_n^+ - \psi_n^+ x_n^- +
T_{[n/3]}[\psi^-(t)]S_{\gamma} ((\tilde h(t,\tau)-\tilde h(t,t))
x_n(\tau)),
\]
where
$ \tilde h(t,\tau)=T_{[n/3]}^t T_{[n/3]}^{\tau} h(t,\tau).$

It is obvious
$$
\Vert Vx_n-V_n x_n \Vert \le A \Vert x_n \Vert(E_n(a)+E_n(b)+E_n(\psi)). \eqno (3.16)
$$

Since
\[\mid h(t,\tau)-\tilde h(t,\tau) \mid < A E_n^*(h)=A \ln n (E_n^t(h)+
E_n^{\tau} (h)),
\]
then the function $ \eta(t,\tau)=h(t,\tau)-\tilde h(t,\tau) $ belongs
to the Holder class with respect to variable value $ \tau $ with the
exponent $ 1/\ln n $ and with the coefficient $ A \ln n E_n^*(h). $

Therefore
\[
 \Vert S_{\gamma}(((h(t,\tau)-\tilde h(t,\tau))-(h(t,t)-\tilde
h(t,t)) x_n(\tau) \Vert \le  
\]
\[
\le A \Vert x_n \Vert \ln n (E_n(b)+
E_n^*(h)).
\]

An analogous estimate is valid for $ \Vert \overline P_n [L_2 x_n-
\tilde \varphi_n] \Vert.$
Indeed\\
 $ \Vert \overline P_n[[V_n-V]x_n] \Vert \le A \Vert x_n \Vert
(E_n(a)+E_n(b)+ E_n(\psi)).$

It remains to carry out the following proof
\[ \Vert\overline P_n[S_{\gamma}[((h(t,\tau)-\tilde h(t,\tau))-(h(t,t)-
\tilde h(t,t))) x_n(\tau)]] \Vert \le 
\]
\[
 \le \Vert \overline P_n[S_{\gamma}(\upsilon(\tau, \tau)(x_n(\tau)-
x_n(t)))] \Vert + \Vert \overline P_n[x_n(t) S_{\gamma}(\upsilon
(\tau, \tau))] \Vert + 
\]
\[
+ \Vert \overline P_n \mid[S_{\gamma}(\upsilon(t,\tau)
-\upsilon(\tau,\tau)) x_n(\tau)] \mid +
\]
\[
+ \Vert \overline
P_n[(h(t,t)-\tilde h(t,t))S_{\gamma}(x_n(\tau))] \Vert \le 
\]
\[ \le A \Vert x_n \Vert E_n^*(h) \ln^2n, 
\]
where $ \upsilon(t, \tau)=\overline P_n^t[h(t,\tau)-
\tilde h(t,\tau)]. $

Since the operator $ L_2 $ has the linear inverse one, then,
for $ n $ such that $ q_1=A \ln^2n [E_n(a)+E_n(b)+E_n(\psi)+
E_n^*(h)]<1 $, the operator $ L_{2,n} $ has the linear inverse operator  $ L_{2,n}^{-1} $ with the norm
$ \Vert L_{2,n}^{-1} \Vert \le \Vert L_2^{-1} \Vert/(1-q_1). $

Since the equation $ L_{2,n} x_n=y_n $ is equivalent to (3.15), then
the linear operator $ L_{1,n}^{-1} $ exists. Let us estimate its
norm. It is evident
\[ \Vert \tilde x_1^* \Vert \le \Vert L_{2,n}^{-1} \Vert \Vert
\tilde y \Vert=\Vert \tilde K_2^{-1} \Vert \Vert \overline P_n
[\psi^-f] \Vert \le 
\]
\[
\le \Vert L_{2,n}^{-1} \Vert \max\limits_t \mid \psi^-(t)\mid \Vert
\overline P_n[f] \Vert, 
\]
 i.e. $ \Vert L_{1,n}^{-1} \Vert
\le A. $

In addition $ \Vert x^*- x_{1,n}^* \Vert \le A \ln^2 (n) [E_n^*(h)+
E_n(a)+E_n(b)+E_n(\psi)], $ where $ x^* $ and $ x_1^* $ are solutions of
equations (3.2) and (3.15).

Let us now show that the linear operator $ L_n^{-1} $
exists. To this end let us estimate
\[
 \Vert L_n  x_n - L_{1,n}  x_n \Vert \le \Vert
\frac {1}{ \pi i} \int\limits_{\gamma} \frac {\omega(t, \tau)-\omega(\tau, \tau)}
{\tau-t} x_n(\tau)d\tau \Vert + 
\]
\[
+ \Vert \overline P_n\left[\frac {1} {\pi i} \int\limits_{\gamma}
\frac{\omega(\tau, \tau)
x_n(\tau)}{\tau-t}d\tau\right] \Vert, 
\]
where $w(t,\tau)=\overline
P_n^t[h(t,\tau)-\tilde h^*(t,\tau)], \tilde h^*(t,\tau)=P_n^{\tau}[h(t,\tau)]. $

Since the function $ w(t,\tau)-w(\tau, \tau) $ belongs to the Holder class
with respect to variable $ t $ with index $ 1/\ln n $ and with
coefficient $ A \ln^2n E_n^{\tau} (h) $, then, by repeating preceding
arguments, we make sure of correctness of the estimate
\[
 \Vert L_n x- L_{1,n} x_n \Vert \le A \ln^2n E_n^{\tau}
(h) \Vert x_n \Vert. 
\] 

So, for $n$ such that $q_2=\max(q_1, AE_n^\tau(h)\ln^2 n)<1$, the operator $L^{-1}_n$ exists and $\|L_n^{-1}\|\le \|L_{1,n}^{-1}\|/(1-q_2), $ i.e. $ \Vert L_n^{-1} \Vert \le A. $ 

To this end
$ \Vert x_n^* -   x_{1,n}^* \Vert \le A \ln^2n E_n^{\tau}(h), $
where $ x_n^* $ is a solution of the equation (3.5). The Theorem
is proved.

Proof of the Theorem 3.3 is similar to the proof of the Theorem 3.2.

\begin{center}
{\bf 4. An Approximate Solution of Singular Integral Equations with 
Discontinuous Coefficients and on Open Contours of Integration }
\end{center}

Let us investigate the projective method for solution of
equation
\[
Kx \equiv a(t) x(t) + b(t) S_{\gamma} (x(\tau)) + U_{\gamma}
(h(t, \tau){\mid t-\tau \mid}^{-\eta} x(\tau))=
\]
$$
=f(t), \eqno (4.1)
$$
$$
Lx \equiv e(t) x(t) + S_L(k(t, \tau) x(\tau))=g(t). \eqno (4.2)
$$

Here $ \gamma $ is a unit circle with the center at origin,
$ L=(c_1, c_2) $ is a segment of $ \gamma. $ We will consider that
the functions $ a(t), b(t),$ $ f(t) \in H_\alpha $,  $ h(t, \tau) \in H_{\alpha,\alpha}$  $(0<\alpha \le 1) $ everywhere on the circle
$ \gamma $  except the point $ t=1, $ where the functions
$ a(t), b(t) $ have discontinuity of the first kind. Let $ e(t), 
g(t) \in H_{\alpha}$, $ k(t, \tau)\in H_{\alpha,\alpha} $ $ (0<\alpha \le 1). $ In the complex plane
the cut from the origin across the point $ c  \ (t=1) $ to
infinity is carried out. On the complex plane with thus slit the functions
$ (t-1)^{ \delta} $ and $ t^{\delta} $, which are used below, are analytical.

{\bf 4.1. Fundamental Statements}

The singular integral equation (4.2) on the open contours of integration are
reduced to the singular integral equations with discontinuous
coefficients. Because we will pay main attention to singular
integral equations (4.1). It follows from theory of singular integral
equations \cite{Gakh}, \cite{Mus} that the singularities of a solution $ x^*(t) $ of equation
(4.1) coincides with singularities of canonical function of  characteristic  operator
\[ K^0x \equiv a(t) x(t) + b(t) S_{\gamma} (x(\tau)). 
\]

Let the
solution $ x^*(t) $ in the vicinity of the point $ c $ has the form
$ (t-c)^{\delta} \varphi(t), $ where
\[ 
\delta=\frac{1}{2\pi i}
\ln (G(c-0) /G(C+0))=\xi+i \zeta, -1<\xi<1, 
\]
 $ \varphi \in H_{\alpha},
G(t)=d(t) s(t), d(t)=a(t)-b(t), s(t)= (a(t)+b(t))^{-1}. $
Depending on the value $ \xi $ (it is  assumed that $ -\eta+
\xi>-1 $) it is necessary to distinguish two cases:
a) $ 0 < \xi <1, $ b) $ -1 < \xi \le 0. $ In each case it is necessary
to consider a separate calculating scheme.

At first we turn  our attention to the case when $ 0< \xi <1. $

The approximate solution of equation (4.1) is sought as the polynomial
$$
x_n(t)=\sum_{k=-n}^n\alpha_kt^k, \eqno (4.3)
$$
which
coefficients  are defined from the system of equations 
\[
K_n x_n \equiv \overline P_n [a(t) x_n(t) + b(t) S_{\gamma}
(x_n(\tau)) + 
\]
$$
+ U_{\gamma}( \overline P_{n}^{\tau} [h(t,\tau) d(t, \tau)
x_n(\tau)])]= \overline P_n[f(t)]. \eqno (4.4)
$$

The justification of this calculating scheme is carried out in
space $ X=H_{\beta} \  (\beta< \lambda_1=\min(\alpha, \xi, 1-\eta)) $
and in its subspace $ X_n $ consisting of $n$-order  polynomials.

{\bf Theorem 4.1}  \cite{Boy14}, \cite{Boy16}, \cite{Boy25}.  Let the following conditions are safisfied:
the functions $ a,b,f \in H_{\alpha}, h \in H_{\alpha \alpha}(0<\alpha \le 1) $ everywhere
except the point $ t=1; $ in the point $t=1$ functions $ a(t) $ and $ b(t) $ have the
discontinuity of the first  kind; the Riemann boundary value problem $ \psi^+(t)=G(t)
\psi^-(t) $ has a solution of the form $ \psi=(t-1)^{\delta} \varphi(t),
\delta = \xi + i\zeta; \xi >0; $ the operator $ K \in [X,X] $ is
continuously invertible. Then for $ n $ such that $ q=A(n^{-(\lambda_1-
\beta)} + n^{-\beta}) \ln n<1, $ the system (4.4) is uniquely
solvable and the estimate $ \Vert x^* - x_n^* \Vert_\beta \le
A(n^{-(\lambda_1-\beta)}+n^{-\beta}) \ln n $ is valid.
Here $ x^* $ and $ x_n^* $ are  solutions of equations (4.1) and (4.4).

Let us consider another numerical scheme.

The approximate solution of equation (4.1) is sought in the form
of the polynomial (4.3), which coefficients  are defined from the
system
\[
 K_n x_n \equiv \overline P_n [a(t) x_n(t) + b(t) S_{\gamma}
(x_n(\tau)) + 
\]
$$
+ \frac{1}{2\pi i} \sum_{k=0}^{2n} h(t,t_k) x_n(t_k) \int\limits_{t'_{k-1}}^
{t'_{k}} \mid \tau-t \mid^{-\eta} d\tau]=\overline P_n[f(t)], \eqno (4.5)
$$
where $ t'_k=e^{is'_k}, s'_k=(2k+1)\pi/(2n+1),$
$k=0,1. \dots. 2n. $

{\bf Theorem 4.2} \cite{Boy14}, \cite{Boy16}, \cite{Boy25}. Let all hypothesis of preceding Theorem are
realized and besides a solution $ x^*(t) $ of equation (4.1) has the
form $ x^*(t)=(t-1)^{\delta}\varphi^*(t), $ where $ \varphi^*(t) \in
H_{\alpha}. $ Then for $ n $ such that $ q=A(n^{-(\alpha-\beta)}+
n^{-(\xi-\beta)}) \ln n <1 $, the system (4.5) is uniquely solvable
for any right-hand side and the estimate $ \Vert x^*-x_n^* \Vert_\beta \le
A(n^{-(\alpha-\beta)} +n^{-(\xi-\beta)}) \ln n  $ is valid ,
where $ x^* $ and
$ x_n^* $ are  solutions of equations (4.1) and (4.5).

     Let us get over to the case when $ -1<\xi \le 0. $

     We denote by $ X^* $ \  the space of functions 
$ x(t)=(t-1)^{\delta}\varphi(t)$, $  \varphi\in H_\beta, $ with the norm
\[
 \Vert x(t) \Vert=\max_{t\in \gamma}\mid \varphi(t) \mid + \sup_{t_1
\ne t_2, 1 \notin (t_1,t_2)}[\mid \varphi(t_1)-\varphi(t_2) \mid / \mid
t_1-t_2 \mid^{\beta}], 
\]
where $ \beta<\lambda_2=\min(\alpha, 1-\eta-\xi). $ We denote by $ Y $
the space of functions satisfying the Holder condition $ H_{\beta} $
with the norm
\[ \Vert y \Vert= \max_{t\in \gamma}\mid y(t) \mid + \sup_{t_1
\ne t_2, 1 \notin (t_1,t_2)}[\mid y(t_1)-y(t_2) \mid/ \mid
t_1-t_2 \mid^{\beta}]. 
\]

We denote by $ Y_n $ subspace of space $ Y $ consisting of  polynomials
 $ \sum_{k=-n}^n \alpha_kt^k. $

The approximate solution of the equation (4.1) is sought in the space $X^*_n$, consisted from functions
$ x_n(t)=x_n^+(t)+x_n^-(t), $ where
$$
\begin{array}{ccc}
 x_n^+(t)=(t-1)^{\delta}
\varphi_n^+(t)=(t-1)^{\delta} \sum_{k=0}^n \alpha_k t^k, \\
 x_n^-(t)= \left( \frac{t-1} {t} \right)^{\delta} \varphi_n^-(t)=
\left(\frac{t-1}{t}\right)^{\delta} \sum_{k=-n}^{-1} \alpha_k t^k,\\
\end{array}
\eqno (4.6)
$$
coefficients $\{ \alpha_k \}$ of which are defined from the system
\[
 \overline K_n x_n \equiv \overline P_n[x_n^+(t)-G(t)x_n^-(t)+s(t)
U_{\gamma}(\overline P_n^\tau [h(t, \tau) d(t, \tau)
x_n(\tau)])]=
\]
$$
=\overline P_n[s(t)f(t)]. \eqno (4.7)
$$

{\bf Theorem 4.3}  \cite{Boy16}, \cite{Boy25}. Let the following conditions are realized:\\
1) the functions $ a,  b, f \in H_{\alpha}, \  h \in H_{\alpha \alpha}$ \ $(0<\alpha<1) $ \  everywhere
except the point $ t=1 $, in which  functions $ a(t)$ and $ b(t) $ have a discontinuity
of the first kind;\\
 2) the operator $ K $ acting from the space $ X^* $ into
the space $ Y $ has continuously invertible one;\\
 3) the Riemann boundary value
problem $ \psi^+(t)=G(t) \psi^-(t) $ has a solution of the form
$ (t-1)^{\delta}\varphi(t), \delta= \xi+i \zeta, -1<\xi<0. $\\
Then for $ n $ such that $ q=A(n^{-(\lambda_2-\beta)}+n^{-\beta})
\ln n<1,  \lambda_2=\min(\alpha, 1- \eta+ \xi, 1+\xi)  $,
the system of equations  (4.7) has a unique solution $ x_n^* $ and the
estimate $ \Vert x^*-x_n^* \Vert \le A(n^{-(\lambda_2-\beta)}+
n^{-\beta}) \ln n $ is valid. Here $ x^* $ is a solution of the
equation (4.1).

Let us note the changes which appear in constructing
numerical scheme if one supposes that a free
term of equation (4.1) has a singularity as $ (t-1)^ \upsilon $ at
the point $ c=1. $ Let us input the space $ Y^* $ of the functions 
$ y(t)=(t-1)^{\upsilon} \varphi(t), \varphi \in H_{\beta}, $ with the norm
\[ \Vert y \Vert=\max_{t \in \gamma} \mid \varphi(t) \mid+ \sup_{t_1 \ne
t_2, 1 \notin (t_1,t_2)} \mid \varphi(t_1)- \varphi(t_2) \mid / \mid
t_1-t_2 \mid^{\beta} 
\]
and \  its \ subspace \ $ Y_n^* \in Y^* $\   consisting \ of \  functions of form
$ y_n^*=$\\ $(t-1)^{\upsilon} \sum_{k=-n}^n \alpha_k t^k. $ We denote by
$ \overline P_n^*\in [Y^*,Y_n^*] $ the projector $ \overline P_n^* y(t)=\overline P_n^*[(t-1)^{\upsilon}
\varphi(t)]=(t-1)^{\upsilon} \overline P_n[\varphi(t)]. $

We propose that the operator $ K $ acts from $ X^* $ to
$ Y^* $ and has the continuous inverse one. An approximate solution of
equation (4.1) we will seek in the form of a function $ x_n(t) $, defined by the
expression (4.6), which coefficients  are defined from the system
of algebraic equations:
\[
 K_n^* x_n \equiv
 \]
 \[
\equiv  \overline P_n^* [x_n^+(t)-G(t) x_n^-(t)+s(t)
U_{\gamma} [\overline P_{n\tau}[h(t, \tau) d(t,\tau) x_n(\tau)] d\tau]]=
\]
$$
=\overline P_n^*[s(t)f(t)]. \eqno (4.8)
$$

{\bf Theorem 4.4}  \cite{Boy14}, \cite{Boy16}, \cite{Boy25}. Let the hypothesises of   Theorem 4.3
are realized. Then for $ n $ such that $ q=A(n^{-(\lambda_2-\beta)}+
n^{-\beta}) \ln n<1, $ the system (4.8) is uniquely solvable for any
right-hand side and the estimate $ \Vert x^*-x_n^* \Vert_{X^*}  \le A
(n^{-(\lambda_2-\beta)}+n^{-\beta}) \ln n $ is
valid, where $ x^* $ and $ x_n^* $ are  solutions of equations (4.1)
and (4.8).

     In a certain case it appears to consider the equations (4.1) and (4.2) is
more preferable  as operator equations in the space $ L_p
(1<p<\infty) $ and to lead the proof of calculating schemes in
spaces $ L_p $ and their subspaces $ L_{n,p} $ consisting of
 polynomials as (4.3).

{\bf Theorem 4.5}  \cite{Boy14}, \cite{Boy16}, \cite{Boy25}. Let the operator $ K $ is continuously invertible
in the space $ L_2, $ the coefficients $ a(t), b(t), f(t) \in H_{\alpha}, $
$h(t,\tau) \in H_{\alpha \alpha}$
are continuous everywhere over $ \gamma $
except the point $ t=1 $, in which functions $ a(t), b(t) $ have a discontinuity
of the first kind. Then for $ n $ such that $ q=A(n^{-\alpha}n^{\theta}+
n^{-\eta(1-\eta)/(1+\eta)}) <1 $, the system (4.8) has a unique solution
$ x_n^* $ and an estimate $ \Vert x^*-x_n^* \Vert_{L_2} \le A
(n^{-\lambda_2-\beta)}+n^{-\beta}) \ln n $ is
valid, where $ x^* $ is a solution of equation (4.1). Here $ \Theta=
-(1-\mid \xi \mid) $ for $ \xi \le 0, \Theta=-\xi $ for $ \xi>0; $ the
function $ (t-1)^ {\delta} \varphi_0(t), $ where $ \delta=\xi+i \zeta,
\varphi_0 \in H_{\alpha}$ is a solution of the Riemann boundary value problem
$ \varphi^+(t)=G(t) \varphi^-(t). $

{\bf Theorem 4.6}  \cite{Boy14}, \cite{Boy16}, \cite{Boy25}. Let the operator $ K $ is continuously invertible
in the space $ L_p \ (1<p \le 2), $ functions $ a,\  b,\  h, \  f $
everywhere, except the point $ t=1, $ satisfy the Holder conditions
with the exponent $\alpha$
$(0 < \alpha <1), $ in the point $ t=1 $  functions
$ a,\ b $ have a discontinuity of the first kind, and $ (t-1)^{\delta}
\varphi_0(t)$ $(\delta=\xi+i \zeta) $ is a solution of Riemann boundary value problem
$ \varphi^+(t)=G(t) \varphi^-(t). $ Then for $ n $ such that $ q=A(n^
{-\alpha}+n^{-\Theta}+n^{-\eta(1-\eta) /(1+\eta)}) <1 \ (\Theta=
\xi $ for $ \xi >0, \Theta=1-\mid \xi \mid $ for $ \xi \le 0) $
the system of equations  (4.8) has a unique solution $ x_n^* $ and the estimate
$ \vert x^*-x_n^* \Vert_{L_p} \le A(n^{-\alpha}+n^{-\Theta}+n^{-\eta
(1-\eta)/(1+\eta)}) $ is valid. Here $ x^* $ is  a solution of
the equation (4.1).

\begin{center}
{\bf 5. Singular Integral Equations with Constant Coefficients}
\end{center}

     Let us consider the singular integral equation
\[
Kx \equiv ax(t)+\frac{b}{\pi} \int\limits_{-1}^1 \frac{x(\tau)}{\tau-t}
d\tau +
\frac{1}{\pi} \int\limits_{-1}^1 h(t,\tau)x(\tau)d\tau \equiv
\]
$$ 
\equiv ax(t)+bSx +Hx
=f(t)
\eqno (5.1)
$$
with constant coefficients $a$ and $b$.

The index of the equation (5.1) is define as $\xi = -(\alpha+\beta),$
where
$$
\alpha = \frac{1}{2\pi i}\ln \left(\frac{a-ib}{a+ib}\right)+N,
$$
$$
\beta = -\frac{1}{2\pi i}\ln \left(\frac{a-ib}{a+ib}\right)+M,
$$
where $N$ and $M$ are integer, which we choose as follows:

1) $\xi =1,$ $-1 < \alpha,\beta <0;$

2) $\xi =-1,$ $0 < \alpha,\beta <1;$

3) $\xi =0,$ $\alpha=-\beta,$ $0<|\alpha|<1.$

These cases cover the well-known problems of mechnics.

Under these indexes a
solution  of the equation (5.1) has the form $ x(t)=\omega_{\xi}(t) z(t), $ where
$$ \omega_{\xi}(t)=(1-t)^{\alpha}(1+t)^{\beta}, 0<\mid \alpha \mid,
\mid \beta \mid<1, \xi=-(\alpha+\beta), $$
$ z(t) $ is a smooth function.

Also, the number $ \alpha $ can be determined by the formula
$$
a+b \cot \pi \alpha=0. 
$$

L. Gori and E. Santi \cite{Gori} used  a projective-splines methods for  solution of the
equation  (5.1).

A.V. Dzhishkariani are devoted  papers \cite{Dzi}, \cite{Dzi1} to numerical methods for 
solution of singular integral equations (5.1). 

His results are summed in the review \cite{Dzi2}.

Here, following papers  \cite{ABoy1}, \cite{ABoy2}, we will give reviews of some numerical algorithms for solution the equation
based on other approach.

We will consider the particle case of the equation (5.1) - the equation
$$
Kx \equiv \frac{1}{\pi} \int\limits_{-1}^1 \frac{x(\tau)}{\tau-t}
d\tau
\equiv Sx 
=f(t).
\eqno (5.2)
$$

This equation is important in the aerodynamics.

Results, which  will be given for the equation (5.2), one can easily diffuse to the equation  
\[
Kx \equiv \frac{b}{\pi} \int\limits_{-1}^1 \frac{x(\tau)}{\tau-t}
d\tau +
\frac{1}{\pi} \int\limits_{-1}^1 h(t,\tau)x(\tau)d\tau \equiv
\]
$$ 
\equiv bSx +Hx
=f(t),
$$
where $Hx$ is a compact  operator.

It is known, that the index of the equation (5.2) can take three values:

$$
\xi =\left \{
\begin{array}{ccc}
1,  \ \ \ \ \alpha=\beta=-1/2 \\
-1, \ \ \ \ \ \alpha=\beta=1/2 \\
0,  \ \ \ \ \alpha=-1/2,\beta=1/2\ (or \  \alpha=1/2, \beta=-1/2)\\
\end{array} 
\right.
$$

{\bf The first case. Index $\xi=0$.}

Let $\xi=0, \alpha=-\beta, |\alpha|=|\beta|=1/2.$

There are two possibilities:\\
1) $\alpha=-1/2, \beta=1/2;$\\ 2) $\alpha=1/2, \beta=1/2.$

Consider the first case. The second case is similar.

The solution of equation (5.2) has the form
$$
x(t)=(1-t)^{-1/2}(1+t)^{1/2}\varphi(t).
$$ 

We introduce the function space
$L_{2,\rho_1}[-1,1]$ with the weight $\rho_1=[(1-t)(1+t)]^{1/2}$.

The integral operator $ Hx$ is  completely continuous in the space $L_{2,\rho_1}[-1,1].$

Known \cite{Sege}, that the system of functions
$$
z_k(t)=c_k(1-t)^{-1/2} (1+t)^{1/2} P_k^{(-1/2,1/2)}(t),
$$
where $P_k^{(-1/2,1/2)}(t)$ is the 
Jacobi polynomial, $k=0,1,...,$
$$
c_0=\pi; c_k=(h_k^{(-1/2,1/2)})^{-1/2}, k=1,2,\ldots,
$$
$$
h_k^{(-1/2,1/2)}=h_k^{()1/2,-1/2)}=\left(\frac{2\Gamma(k+1/2)\Gamma(k+3/2)}{(2k+1)(k!)^2}  \right);
$$ 
complete and orthonormal in the space $L_{2,\rho_1}[-1,1].$
Besides,
$$
Sz_k=c_k  P_k^{(1/2,-1/2)}(t).
$$

Let $t_k (k=0,1,...,n)$ are the knots of the polynomial $P_{n+1}^{(1/2,-1/2)}(t).$

Let 
$$
\gamma_k^{(1/2,-1/2)}= \sum\limits^n_{l=0}(P_l^{(1/2,-1/2)}(t_k))^2, k=0,1,...,n.
$$

The polynomial 
$$
f_n(t)=\sum\limits^n_{k=0}\frac{1}{\gamma_k^{(1/2,-1/2)}}\sum\limits^n_{l=0}P_l^{(1/2,-1/2)}(t_k)P_l^{(1/2,-1/2)}(t) f(t_k)
$$
interpolate a function $f(t)$ on the knots 
$t_k,k=0,1,...,n.$

An approximate solution of equation (5.1) will be sought in the form of function
$$
x_n(t)=\sum\limits^n_{k=0}\frac{1}{\gamma_k^{(1/2,-1/2)}}\sum\limits^n_{l=0}P_l^{(1/2,-1/2)}(t_k)P_l^{(-1/2,1/2)}(t)x_k.
$$

The right-hand part  $ f(x) $ of the equation (5.1) is approximated by a polynomial
$$
f_n(t)=\sum\limits^n_{k=0}\frac{1}{\gamma_k^{(1/2,-1/2)}}\sum\limits^n_{l=0}P^{(1/2,-1/2)}(t_k)P_l^{(1/2,-1/2)}(t)f_k.
$$
where $f_k=f(t_k).$

Substituting in the equation (5.2) $ x_n (t) $ instead of $ x (t) $ and $ f_n (t) $ instead of $ f (t) $ and setting the left and right side of this expression in the points  $t_i, i=0,1,...,n,$ we have  $x_i=f_i, i=0,1,...,n.$

Thus, an approximate solution of equation (5.2) has the form
$$
x^*(t)=\sum\limits^n_{k=0}\frac{1}{\gamma_k^{(1/2,-1/2)}}\sum\limits^n_{l=0}P^{(1/2,-1/2)}(t_k) P_l^{(-1/2,1/2)}(t)  f_k.
$$

Since $ f_n (t) $ is the interpolation polynomial   for a function $ f (t) $ on the nodes of the Jacobi polynomial
$P_{n+1}^{(1/2,-1/2)}(t) $ then $||f(t)-f_n(t)||_{C[-1,1]}\leq (1+\lambda_{n+1}) E_n(f),$ where $E_n(f)$ is the best uniform approximation of $ f (t) $ by $n$-order polynomials, $\lambda_{n+1}$ is the  
Lebesgue constant  over the nodes
of the polynomial
$P_{n+1}^{(1/2,-1/2)}(t).$

Since the operator $ K $ is continuously invertible in the space 
$L_{2,\rho_1}[-1,1]$ \cite{Dzi},  
then 
$$
\|x^*(t)-x^*_n(t)\_{[-1,1]}|\leq C(1+\lambda_{n+1}) E_n(f),
$$
where $x^*(t)$ is a solution of the equation  (5.1).

{\bf The second case. Index $\xi=1$.}

Let $\xi =1,$ $\alpha = \beta = -1/2.$ In this case, the equation (5.1) has a solution with the singularities at both ends:
$x(t) = (1-t^2)^{-1/2}\varphi(t),$ where $\varphi(t) $ is a smooth function.

Consider the Chebyshev polynomials of the first kind
$T_n(t) = \cos n \arccos t,$ $n=0,1,\ldots$ $|t| \leq 1.$ 

The polynomials $ T_n (t) $ are orthogonal in the space $ L_2 [-1,1] $ with weight
$(1-t^2)^{-1/2}.$

Let $ U_n (t) $ is the Chebyshev polynomial of the second kind:
$$
U_n(t) = \frac{1}{(1-t^2)^{1/2}} \sin(n+1) \arccos t, \, 
n=0,1,\ldots, \, |t| \leq 1.
$$

It is known that
$$
\frac{1}{\pi} \int\limits^1_{-1} \frac{T_n(\tau) d\tau}{\sqrt{1-\tau^2}(\tau-t)} = U_{n-1}(t), \, n \geq 1;
$$
$$
\frac{1}{\pi} \int\limits^1_{-1} \frac{d\tau}{\sqrt{1-\tau^2}(\tau-t)} = 0.
$$

Let $ t_k, $ $ k = 0,1, \ldots, n-1, $ are nodes of the Chebyshev polynomial of the second kind $ U_{n} (t). $

On nodes $ t_k, $ $ k = 0,1, \ldots, n-1, $
we construct the polynomial $f_{n-1}(t)$, which  interpolate the function $ f: $
$$
f_{n-1}(t) = \sum\limits^{n-1}_{k=0} \frac{1}{\gamma_k} \sum\limits^{n-1}_{l=0}
U_l(t_k) U_l(t) f(t_k),
$$
где $\gamma_k = \sum\limits^{n-1}_{l=0}(U_l(t_k))^2.$

An approximate solution of equation (5.2) is sought in the form of a polynomial
$$
x_n(t) = \frac{1}{\sqrt{1-t^2}} \sum\limits^{n-1}_{k=0} \frac{1}{\gamma_k} \sum\limits^{n-1}_{l=0} U_l(t_k) T_{l+1}(t)x_{l+1} + \frac{x_0}{\sqrt{1-t^2}}. 
$$

Substituting the functions $ x_n (t) $ and $ f_n (t) $ in the equation (5.2) instead of functions $ x (t) $ and $ f (t), $ we have
$$
\sum\limits^{n-1}_{k=0} \frac{1}{\gamma_k} \sum\limits^{n-1}_{l=0}
U_l(t_k) U_l(t) x_{l+1} = 
\sum\limits^{n-1}_{k=0} \frac{1}{\gamma_k} \sum\limits^{n-1}_{l=0}
U_l(t_k) U_l(t) f(t_l).
$$

Equating both sides of this equality in knots $ t_k = 0,1, \ldots, n-1$, we have 
$x_{l+1} = f_l,$ $l=0,1,\ldots,n-1.$

So,
$$
x_n(t) = \frac{1}{\sqrt{1-t^2}} \sum\limits^{n-1}_{k=0} \frac{1}{\gamma_k} \sum\limits^{n-1}_{l=0} U_l(t_k) T_{l+1}(t) f_l + \frac{x_0}{\sqrt{1-t^2}}. 
$$

To find a unknown value $  x_0$ we need an additional condition.

In problems of mechanics are usually given an additional condition
$$
\int\limits^1_{-1} x(\tau) d\tau = p, \, p=\rm{const}.
\eqno (5.3)
$$

Substituting the function $ x_n (t) $ in the last equation, we find $  x_0. $

Thus, we received a unique  solution of the equation (5.2) under the presence of additional condition (5.3).

\begin{center}
{ \bf 6. Approximate Methods for Solution of Nonlinear Singular
Integral Equations}
\end{center}

Iterative methods for solution of nonlinear singular integral equations with
Hilbert  and Cauchy kernels are devoted many papers. It
thise papers are considered simple methods of iteration, Newton-Kantorovich
method, functional correction method and other. Reviews of these methods and
rich bibliography are given in \cite{Gus}.

Obviously, the first works devoted to projective methods for solution of
nonlinear singular integral equations was the  paper \cite{Boy2}.

\begin{center}
{ \bf 6.1.  Projective Methods for Solution of Nonlinear Equations on Closed
	Paths of Integration}
\end{center}
Let us consider the nonlinear singular integral equation
$$
Kx \equiv a(t,x(t))+\frac{1}{\pi i} \int\limits_{\gamma}\frac{h(t,\tau,
x(\tau))}{\tau-t}d\tau=f(t), \eqno (6.1)
$$
where $ \gamma $ is the unit circle with the center at origin.

{\bf Basis Statements.}

{ \bf Calculating Scheme 1.} An approximate solution of equation
(6.1) is sought in the form of polynomial
$$
x_n(t)=\sum_{k=-n}^n \alpha_k t^k, \eqno (6.2)
$$
the coefficients of which are determined from the system 
$$ a(t_j, x_n(t_j))+\frac{1}{2n+1}\sum_{k=0}^{2n}{}'h(t_j,t_k,x_n(t_k))
(1-i \cot \frac{s_k-s_j}{2}) - $$
$$
- \frac{2i}{2n+1}h_u'(t_j,t_j,x_n(t_j))x_{ns}'(e^{is_j})=f(t_j),
j=0,1, \dots, 2n, \eqno (6.3)
$$
where $ t_j=\exp(is_j), s_j=2\pi j/(2n+1),$
the prime in the summation indicate that $ k\neq j, $ $ x_{ns}' $
means differentiation of the function $ x_{n}(s) $ with respect to $ s$ and $h'_u(t,\tau,u)$ means differentiation  of the function  $h(t,\tau,u)$
with respect to u.

{ \bf Theorem 6.1} \cite{Boy2}, \cite{Boy3}, \cite{Boy16}, \cite{Boy25}. Let the equation (6.1) has a unique solution
$ x^* $ in certain sphere $ S, $ it exists the linear operator $ [K'
(x^*)]^{-1} $ and the conditions $ x^*(t), f(t) \in H_{\alpha},$ $
a(t,u), a_u'(t,u)$, $ a_u''(t,u)\in H_{ \alpha 1},$ $ h(t,\tau,u),
h_u'(t,\tau,u),$ $ h_u''(t,\tau,u) \in H_{\alpha \alpha 1}, $ where
$ 0<\alpha \le 1, \mid u\mid<\infty, $ are fulfilled. Then for $ n $ such
that $ q=An^{-(\alpha/2-\beta)} \ln^5n<1, $ the system (6.3) has
a solution $ x_n^* $ and in the metric of the space $ X=H_{\beta} $ the
inequality $ \Vert x^*-x_n^* \Vert \le An^{-(\alpha/2- \beta)}
\ln^2n$ is valid.

{\bf Calculating Scheme 2.} An approximate solution of equation (6.1)
is sought in the form of polynomial
$$
x_n(s)=\sum_{k=0}^{2n} \alpha_k \psi_k(s), \eqno (6.4)
$$
the coefficients of which are defined from the system of equations
\[
a(t_j,x_n(\overline t_j))+\frac{1}{2n+1} \sum_{k=0}^{2n} h(t_j,t_k,
\alpha_k)(1-i \cot \frac{2(k-j) \pi- \pi}{4n+2})=
\]
$$
=f(t_j), \eqno (6.5)
$$
where
\[ t_j =\exp(is_j), \psi_k(s)= \frac{1}{2n+1} \frac{\sin\frac{2n+1}{2}(s-s_k)} {\sin
\frac{s-s_k}{2}}, s_k=\frac{2k \pi}{2n+1},
\]
$\bar t_k=\exp(i \bar s_k), \, \bar s_k=(2k+1)\pi/(2n+1).$

Let us plan of obtaining the system (6.5). The equation (6.1) may be
written as
\[
a(s,x(s))-\frac{i}{2\pi}\int_0^{2\pi}h(s,\sigma,x(\sigma))\cot\frac{\sigma-
s}{2}d\sigma-
\]
$$
-\frac{i}{2\pi}\int_0^{2\pi}h(s,\sigma,x(\sigma))d\sigma=
f(s), \eqno (6.6)
$$
where
$ a(s,x(s))=a(e^{is},x(e^{is}))$, $h(s,\sigma,x(\sigma))=h(e^{is}$, $
e^{i\sigma},x(e^{i\sigma}))$, $f(s)=f(e^{is}). $

Having subtracted the polynomial $ x_n(s) $ into equation (6.6) for
$ x(s), $ setting equal left-hand sides of this equation to their
right-hand sides at the points $ \overline s_j=(2 \pi j+\pi)/(2n+1)$
$(j=0,1, \dots, 2n)$ and taking as quadrature rule for both integrals
the quadrature rule of left rectangle on the point $ s_j $ we come
to the following algebraic scheme
\[ a(\overline s_j, x_n(\overline s_j))-\frac{1}{2n+1} \sum_{k=0}^{2n}
h(\overline s_j,s_k, \alpha_k)\cot \frac{s_k-\overline s_j}{2}+
\]
\[
+\frac {1}
{2n+1} \sum_{k=0}^{2n}h(\overline s_j,s_k, \alpha_k)=
f(\overline s_j), $$
$ j=0,1, \dots, 2n. $

Changing $ a(\overline s_j,x_n(\overline s_j)), h(\overline s_j, s_k, \alpha_k),
f(\overline s_j) $ into $ a(s_j,x_n(\overline s_j))$,\\ $ h(s_j, s_k, \alpha_k)$, $ f(s_j) $
we arrive at the system (6.5).

{\bf Theorem 6.2}  \cite{Boy2}, \cite{Boy3}, \cite{Boy16}, \cite{Boy25}.  Let the equation (6.1) has a unique solution
$ x^*(t) $ in certain sphere $S$, it exists a right bounded inverse
operator $ [K'(x)]_r'(x \in S)$ and the conditions $ a(t)$, $,  f(t)$, $
x^*(t)\  \in H_{\alpha}$, \ $h(t, \tau, u)$, $ \  h_u'(t, \tau, u) \   \in
H_{\alpha,\alpha, 1}$ $\  (0<\alpha \le 1, \  |u| <\infty ) $ are fulfilled.
Then for $ n $ such that $ q=A n^{\alpha^2/2-\beta}\ln^5n<1 $
the  system (6.5) has a unique solutions $ x_n^* $ and in the metric
of the space $ X=H_{\beta} (0<\beta<\alpha^2/4) $ the estimate $ \Vert
x^*-x_n^* \Vert \le A n^{-(\alpha/2-\beta)} \ln^2n $ is valid.

{\bf Calculating Scheme 3.} An approximate solution of equation (6.1) is
sought in the form of the polynomial (6.4), the coefficients
$ \alpha_k $ of which are defined from the equation
$$
K_n x_n \equiv \overline P_n \left [a(t, x_n(t))+\frac{1}{\pi i}
\int\limits_{\gamma} P_{n}^{\tau} \left [\frac {h(t,\tau, x_n(\tau))}{\tau-t}
\right] d\tau \right]=
$$
$$
=
\overline P_n[f(t)]. \eqno (6.7)
$$

{\bf Theorem 6.3}  \cite{Boy2}, \cite{Boy3}, \cite{Boy16}, \cite{Boy25}. 
Let the equation (6.1) has a unique solution
$ x^*(t) $ in the certain sphere $ S $,  exists a right bounded
inverse operator $ [K'(x)]_r^{-1}(x\in S) $ and it is carried out
one of the following conditions:

a) $ x^*(t)\in H_\alpha , a(t,u) \in H_{\alpha, 1}^{r,r+1},h(t,\tau,u) \in
H_{\alpha, \alpha, 1}^{r,r,r+1}, $

b) the functions $ x^*(t), a(t, x^*(t)), h(t, \tau, x^*(t)), h(\tau,
t, x^*(t)) $ are analytical that in the domains $ R_1<\mid t \mid<R_2,
R_1<\mid t \mid, \mid \tau \mid<R_2, $
where $ R_1<1, R_2>1. $

If the condition "a" was fulfilled, then,
for $ n $ such that 
$ q=An^{-r-\alpha+2\beta}\ln^6n<1, $ the system
(6.7) has a solution $ x_n^* $ and in the metric of the space $ X=H_{\beta}
(0<\beta<(r+\alpha)/2)$ \ \ \   $\Vert x^*-x_n^* \Vert \le An^{-r-\alpha+\beta}
\ln^2n. $ If the condition "b" was fulfilled, then, for $ n $ such that
$ q=A[R_1^{n+1}+R_2^{-n-1}]n^{2\beta}\ln^6n<1, $ the equations system
(6.7) has a solution $ x_n^* $ and the estimate $ \Vert x^*-x_n^* \Vert_\beta<A[R_1^
{n+1}+R_2^{-n-1}]n^{\beta}\ln^2n  $ is valid.

     In series of cases the condition to nonlinear singular integral
equations in the space $ L_2 $ is more preferable. We are looking for
approximate solution of the equation (6.1) by means of the iterative process
$$
x_n^{m+1}=x_n^m-[\|K_n x_n^m \|^2/\|K_n'(x_n^0)K_n x_n^m\|^2][K_n'
(x_n^0)]^* K_n x_n^m, \eqno (6.8)
$$
where
$$ K_n x_n \equiv \overline P_n[a(s, x_n (s))-\frac{i}{2\pi}
\int_0^{2\pi} P^{\sigma}_n[h(s, \sigma, x_n(\sigma)) \cot\frac{\sigma-s}
{2}] d\sigma+ $$
$$ +\frac{1}{2\pi} \int_0^{2\pi} P_{n}^{\sigma}[h(s,\sigma, x_n(\sigma))]
d\sigma-f(s)],\hskip 11 pt  x_n^0(s)=P_n[x^0(s)], $$
$ x^0(s) $ is well enough approximation for the solution $ x^* $ of  equation
(6.1), $ x_n(s) $ is trigonometric polynomial of degree not higher
than $ n, a(s)=a(e^{is}),$ $h(s, \sigma)=h(e^{is}, e^{i\sigma}), x_n(s)=x_n(e^{is}). $

Let the following condition are fulfilled:

a) functions $ a(t, x^0(t)), f(t), x^0(t), h(t,\tau, x^0(\tau)) $ (with respect to each variable)
either are analytical inside the ring $ R_1 \le \mid t \mid \le R_2,
R_1<1, R_2>1 $ or belong to the class $ W^r H_{\alpha};$

b) in some sphere defined below
\[ \max_{t,\tau \in \gamma}\{\mid h_u'(t, \tau, u_1)-h_u'(t,\tau,u_2)
\mid, \mid h^*(t,\tau,u_1)-h^*(t,\tau,u_2)\mid \le F,
\]
 where $(u_1,u_2\in S),$
 \[h^*(t,\tau,u)=[h_u'(t,\tau,u)-h_u'(\tau,\tau,u)]/
\mid \tau-t \mid^
{\beta} \exp(i\theta_1)\}, 
\]
$ \beta $ is an arbitrary value $ 0<\beta<\alpha, \theta_1=\theta_1
(\tau,t)=\arg \mid \tau-t \mid. $

Let   us  introduce the designations: $ \delta_0=\| Kx_0-f \|$,$ \tilde
\delta_0=\| K_n x_n^0 \|, B_0=\| [K_n'(x_n^0)]^{-1} \|,
B_*=\max\{\| K'(x_0) \|, \| K_n'(x_n^0) \|\}. $
The existence of constants $ B_0 $ and $ \| K'(x_n^0) \| $ and
their relation to the values $ \| [K'(x^0)]^{-1}\| $ and $ \| K'
(x^0) \| $ follow from results of section 3.

{\bf Theorem 6.4} \cite{Boy5}, \cite{Boy16}, \cite{Boy25}.   Let in the sphere $ S[x; \| x-x^0 \| \le r],$ $
r=B_* B_0^2 \tilde \delta_0/(1-q)+A \ln n E_n(x^0),$ $ q=AF_1 B_* B_0^2+
\sqrt{(1-B_*^{-2}B_0^{-2)}}<1, $ the condition "b" are carried out and
the operator $ K'(x^0) $ has a linear inverse ( it is enough left inverse). Then for $ n $ such that
$ p=A \ln^2n \max[E_n(a)$, $ E_n(\psi)$, $ E_n(x^0)$, $ E_n(h_u'(t,t, x_n^0(t))), $ $ E_n^{t,
\tau}(h_u'(t,\tau, x_n^0(\tau)))]<1, $ the equation $ K_n x_n=0 $
has in $ S $ a unique solution $ x_n^* $ and to which the iterative
process (6.8) converges with the rate $ \| x_n^*-x_n^m \| \le q^m
\tilde \delta_0 B_*(1-q). $ The distance between $ x_n^* $ and
solution $ x^* $ of equation (6.1) is estimated by inequalities
$$ \| x^*-x_n^* \| \le \ln^2 n \max(E_n(a), E_n(\psi), E_n(t), E_n^{t, \tau}(h(t,
\tau,x_n^*(\tau))), $$
$$ \| x^*-x_n^* \| \le B_{*} \tilde \delta_0 /(1-q)+E_n(x_0). $$

{\bf Proof of Theorem 6.1.} An approximate solution is sought in the subspace
$ X_n \subset X $, consisted of polynomials of the form (6.2). Under the conditions of the theorem the operator $ K $
 has  Frechet derivative
$$
K'(x)z \equiv a'_u(t,x(t)) z(t)+\frac{1}{\pi i}\int\limits_{\gamma}\frac{
h'_u(t,\tau,x(\tau))z(\tau)}{\tau-t}d\tau, \eqno (6.14)
$$
satisfying in a ball $ S (x^*, r) $ with an some radius $ r $ the 
Lipschitz condition $ \| K'(x_1) -K'(x_2) \| \le A (r) \| x_1-x_2 \|. $
The system of equations (6.3) in the operator form is written in the form of expression
$$
K_n(x_n) \equiv P_n[a(s,x_n(s))+ \frac{1}{2\pi}\int\limits_0^{2\pi}
\sum\limits_{k=0}^{2n}{}'h(s,s_k,x_n(s_k))\psi_k(\sigma)d\sigma-
$$
$$
-\frac{i}{2\pi}\int\limits_0^{2\pi}
P_n\left[[h(s,s,x_n(\sigma))-
h(s,s,x_n(s))]
\cot\frac{\sigma-s}{2}\right]d\sigma -
$$
$$
\left.
 -\frac{i}{2\pi}\int\limits_0^{2\pi}
\sum\limits_{k=0}^{2n}{}'\left\{[h(s,s_k,x_n(s_k))-h(s,s,x_n(s_k))]\times\right.\right.
$$
$$
\times \left.\left.
\cot\frac{s_k-s}{2}\psi_k(\sigma)\right\}\right]d\sigma
=P_n[f(s)].
$$

Here $a(s,x_n(s)) = a(e^{is}, x_n(e^{is})).$
Functions $h(s,s_k,x_n(s_k)), f(s)$ are defined similary.
Here $\psi_k(\sigma)$ are basic trigonometric
polynomials, constructed on knots $s_k=2k\pi/(2n+1),$ $k=0,1,\ldots,2n.$

It is easy to see that the Frechet derivative of the operator $ K_n $ has the form
$$
K'_n(x_n)z_n \equiv P_n[a'_u(s,x_n(s))z_n(s)-I_1(x_n)-I_2(x_n)+I_3(x_n)],
$$
where
$$I_1(x_n)
=\frac{i}{2\pi} P_n\{\int\limits_0^{2\pi}
\sum\limits_{k=0}^{2n}{}'\{[h'_u(s,s_k,x_n(s_k))-
$$
$$
-
h'_u(s,s,x_n(s_k))] z_n(s_k)
\cot\frac{s_k-s}{2}\}] \psi_k(\sigma)d\sigma\};
$$
$$
I_2(x_n)
=\frac{i}{2\pi} P_n \{\int\limits_0^{2\pi}
P_n^{\sigma} \{[h'_u(s,s,x_n(\sigma))z_n(\sigma)-
$$
$$
-
h'_u(s,s,x_n(s))z_n(s)] \cot\frac{\sigma-s}{2} \}d\sigma\};
$$
$$I_3(x_n)=$$
$$
=\frac{1}{2\pi} P_n\{\int\limits_0^{2\pi}
\sum\limits_{k=0}^{2n}{}'[h'_u(s,s_k,x_n(s_k))z_n(s_k)]
\psi_k(\sigma)d\sigma\}.
$$

We will show that the Frechet derivative of $ K_n $ 
 belongs to
Lipschitz class of functions:
$$
\|K'_n(x'_n)-K'_n(x''_n)\| \le A n^{\beta} \ln^2 n\|x_n'-x''_n\|.
\eqno (6.15)
$$

We will give a proof only for  $I_2$ (for  $I_1$ and $I_3$ the proof
is more simply than for $I_2$). 

It was proved \cite{Nat}, \cite{Nat1} that
$\sum\limits_{k=0}^{2n}|\psi_k(s) \le A \ln n.$ 

Therefore
\[
|I_2(x_n')-I_2(x_n'')| \le
\]
\[
\le  A \ln n \{\max\limits_{0 \le j \le 2n}
\frac{2\pi}{2n+1}|\sum\limits_{k=0}^{2n}{}'\{[h'_u(s_j,s_j,x'_n(s_k))-
\]
\[
-h'_u(s_j,s_j,x''_n(s_k))]z_n(s_k)-
\]
\[
-[h'_u(s_j,s_j,x'_n(s_j))
-h'_u(s_j,s_j,x''_n(s_j))]z_n(s_j)\}
\cot\frac{s_k-s_j}{2}+
\]
\[
+\frac{4}{2n+1}|[h''_u(s_j,s_j,x'_n(s_j))x'_n(s_j)
-h''_u(s_j,s_j,x''_n(s_j))x''_n(s_j)]z_n(s_j)|+
\]
\[
+\frac{4}{2n+1}|[h'_u(s_j,s_j,x'_n(s_j))
-h'_u(s_j,s_j,x''_n(s_j))]z'_n(s_j)|\} \le 
\]
\[
\le A \ln n[I_4+I_5+I_6].
\]

Let us eatimate $|I_4|:$
$$
|I_4|\le A|x'_n-x''_n| |z_n|\sum\limits_{k=1}^{n}
k^{-1} \le A\|x'_n-x''_n\| \|z_n\| \ln n.
$$

It is  known the following inequality of Riesz \cite{Ris}.
If $f(s)$ is a trigonometric
polynomial of degree $n,$ then
$|f'(s)| \le n \max\limits_{s}|f(s+\pi/2n)-f(s)|/2.$
Using this inequality, we obtain the estimates
$|I_5| \le A\|x'_n-x''_n\| \|z_n\|n^{-\beta},$
$\|I_6\| \le A\|x'_n-x''_n\| \|z_n\|n^{-\beta}.$

As $((I_2(x'_n)-I_2(x''_n))$ is the  trigonometric
polynomial of degree $n,$
then $|(I_2(x'_n)-I_2(x''_n)|\le A n^{\beta}\ln n \|x'_n-x''_n\| \|z_n\|.$
Now it is easy to verify the validity of the estimate (6.15).

According the conditions of the Theorem, $x^*(t) \in H_{\alpha}$. So, it exists such polynomial $x_n^0(t) \in X_n$,
that 
 $\|x^*-x_n^0\| \le A n^{-(\alpha-\beta)}$.

According the conditions of the Theorem 
the operator $[K'(x^*)]^{-1}$ exists.
So, the linear operator  $[K_n'(x^0_n)]^{-1}$ exists too.
Indeed, let $\|[K'(x^*)]^{-1}\|=B_0$. From Banach Theorem follows that, for  $n$
such  $q=An^{-(\alpha-\beta)}<1,$ the operator 
$[K'(x_n^0)]^{-1}$ exists with the norm $\|[K'(x_n^0)]^{-1}\| \le B_0/(1-q)$.

In \S2 was shown that, for $n$ such that $q=A \ln n/n^{\alpha-\beta}<1,$
the operator $K'_{1n}(x_n^0)$  has linear inverse operator  $[K'_{1n}(x_n^0)]^{-1}$  with the norm $\|[K_{1n}'(x_n^0)]^{-1}\| \le A \ln n$.
Here 
\[K'_{1n}(x_n^0)z_n \equiv
P_n[a'_u(t,x_n^0(t))z_n(t)+S_\gamma[h'_u(t,\tau,x_n^0(\tau))z_n(\tau)]]
.\]

 Let us estimate the norm
$\|[K_{n}'(x_n^0)-K'_{1n}(x_n^0)\|$.
 
At first we estimate the norm
$$
\|J\|
= \|P_n\left[\int\limits_0^{2\pi}
R_n^{\sigma}[h'_u(s,s,x_n^0(\sigma))z_n(\sigma)-\right.
$$
$$
-\left.
h'_u(s,s,x_n^0(s))z_n(s)]\cot\frac{\sigma-s}{2}\right]d\sigma\|,
$$
where $R_n=I-P_n.$ 

Let us approximate the function  $h'_u(s,s,x_n^0(\sigma))$ (with respect to variable  $\sigma$) with  interpolation polynomial, which is  constructed on knots 
$s_k=2k\pi (2n+1)^{-1},$ $k=0,1,\ldots,2n$. This polynomial will be
denoted by $h_n(s,s,\sigma).$ It is easy to see that function
$[h_n(s,s,\sigma)z_n(\sigma)-h_n(s,s,s)z_n(s)]\cot\frac{\sigma-s}{2}$ is a polynomial of  $2n$
degree with respect to variable  $\sigma$.

It is easy to see that $\|J\| \le \|J_1\|+\|J_2\|,$ where
\[
J_1=P_n\left[\int\limits_0^{2\pi}
\left[h'_u(s,s,x_n^0(\sigma))-h_n(s,s,\sigma)\right]z_n(\sigma)
\cot\frac{\sigma-s}{2}d\sigma\right],
\]
\[
J_2=\sum\limits_{k=0}^{2n} S_k(s)\psi_k(s),
\]
\[
S_k(s)=\int\limits_0^{2\pi}P_n^{\sigma}[[[h'_u(s_k,s_k,x_n^0(\sigma))-
h_n(s_k,s_k,\sigma)]z_n(\sigma)-
\]
\[
-[h'_u(s_k,s_k,x_n^0(s_k))
-h_n(s_k,s_k,s_k)]z_n(s_k)]\cot\frac{\sigma-s_k}{2}]d\sigma.
\]

Obviously, $\|J_1\| \le A \ln^2 n \|z_n\|/n^{\alpha-\beta},$
\[
|S_k(s)| \le
\]
\[
 \le \frac{4\pi}{2n+1}\{|[h''_u(s_k,s_k,x_n^0(s_k))
(x_n^{0})'(s_k)-h'_{n\sigma}(s_k,s_k,s_k)]z_n(s_k)|+
\]
\[
+|h'_u(s_k,s_k,x^0_n(s_k))-h_n(s_k,s_k,s_k)]z'_n(s_k)|\}=\frac{4\pi}{2n+1}(|J_3|+|J_4|)=
\]
\[
=\frac{4\pi}{2n+1}|J_3|,
\]
since $ J_4 = 0 $ by the definition of the polynomial $h_n(s,s,\sigma).$

Let us estimate $|J_3|.$ The function $h'_u(s,s,x_n^0(\sigma))$
has the first derivative  with respect to variable  $\sigma$. This derivative is belong to the Holder class
 $H_{\alpha}$
with the coefficient $An. $ It follows from conditions of the Theorem and the Riesz inequality.   It is easy to see that $|J_3| \le An^{1-\alpha}\|z_n\| \ln n.$

From two last inequalities follows
 $|S_k(s)|\le A\|z_n\|n^{-\alpha}\ln n$.
So $\|J_2\| \le A \ln^2 n \|z_n\|n^{-(\alpha-\beta)}.$
From this and estimate for  $\|J_1\|$ we have $\|J\| \le A \ln^2 n \|z_n\|n^{-(\alpha-\beta)}.$

Let us estimate the following integral
\[
\|J_5\|=
\]
\[
=\|\frac{i}{2\pi}P_n\int\limits_0^{2\pi}
\left[h'_u(s,\sigma,x_n^0(\sigma))-h'_u(s,s,x_n^0(\sigma))\right]z_n(\sigma)
\cot\frac{\sigma-s}{2}d\sigma\|-
\]
\[
-I_1(x_n^0) \le A \ln n\|\int\limits_0^{2\pi}[[h'_u(s,\sigma,x_n(\sigma))-
h^*_n(s,\sigma,x_n(\sigma))]-
\]
\[
-[h'_u(s,s,x_n^0(\sigma))
-h^*_n(s,s,x_n^0(\sigma))]]z_n(\sigma)\cot\frac{\sigma-s}{2}d\sigma\|+
\]
\[
+\left\{\|P_n[\int\limits_0^{2\pi}\sum\limits_{k=0}^{2n}{}'
[[h'_u(s,s_k,x^0_n(s_k))-h^*_n(s,s_k,x^0_n(s_k))]- \right.
\]
\[
-[h'_u(s,s,x^0_n(s_k))-h^*_n(s,s,x^0_n(s_k))]]z_n(s_k)
\psi_k(\sigma)\cot\frac{s_k-s}{2} d\sigma]\|+
\]
\[
\left.
+\|\sum\limits_{k=0}^{2n}\frac{2i}{2n+1}[h^*_n(s_k,\sigma,x^0_n(s_k))]'_
{\sigma=s_k}z_n(s_k)\psi_k(\sigma)\|\right\}=\|J_6\|+\|J_7\|,
\]
where $h^*_n(s,\sigma,x^0_n(\nu))$ is the interpolation trigonometric polynomial of   $[n/2]$ degree (with respect to variable 
$\sigma$) for   $P^{\nu}_{[n/2]}[h(s,\sigma,x^0_n(\nu))]$.

 Polynomial is constructed on knots $s_k(k=0,1,\ldots,2[n/2])$. So,
$
h^*(s,\sigma,x^0_n(\nu))=P^{\sigma}_{[n/2]}P^{\nu}_{[n/2]}
[h(s,\sigma,x^0_n(\nu))].
$

It is easy to see that $\|J_6\| \le A\ln^2 n\|z_n\|n^{-(\alpha-\beta)}.$

Using the Riesz inequality we have 
$\|J_7\| \le A \ln^2 n\|z_n\|n^{-(\alpha-\beta)}.$ 

From estimates for $\|J_6\|$ and
$\|J_7\|$ it follows that \[
\|J_5\| \le A\ln^2 n\|z_n\|n^{-(\alpha-\beta)}.
\]

It is easy to receive the estimate
\[
\|P_n[\int\limits_0^{2\pi}h'_u(s,\sigma,x_n^0(\sigma))z_n(\sigma)d\sigma]-
I_3(x_n^0)\| \le 
\]
$$
\le A \ln n \|z_n\|n^{-(\alpha-\beta)}. \eqno (6.16)
$$

 From estimates for  $\|J\|,$ $\|J_5\|$ and from the inequality  (6.16)
it follows that
\[
\|K'_{1n}(x_n^0)-K'_n(x_n^0)\| \le A\ln^2 n/n^{\alpha-\beta}
\].

We have proved that the operator $K'_{1n}(x_n^0)$ has continuously invertible operator. Using the Banach \ \  theorem \  we see that, for  $n$ such that  $q=An^{-\alpha+\beta}
\ln^2 n<1,$ the operator $K'_n(x_n^0)$ has continuously invertible operator with the norm
$\|(K'_n(x_n^0))^{-1}\|\le A \ln n.$

Let us note that $\|K_n x_n^0\| \le \|K_n x_n^0-K x_n^0\|+\|K x_n^0-K x^*\|
\le A n^{-(\alpha-\beta)}\ln n.$ In the ball  $S[x: \|x-x_n^0\| \le
A \ln^2 n/n^{\alpha-\beta}]$, for $n$ such that
$q=\frac{A\ln n}{n^{\alpha-\beta}}<1,$
all conditions of the Theorem 2.6 from the  Chapter 1 is valid. Using this Theorem, we prove the existance of a unique solution 
 $x_n^*(t)$ of the equation (6.3) and the correctness of the estimate
$\|x^*-x^*_n\| \le A \ln^2 n/n^{\alpha-\beta}.$ 

The Theorem is proved.

{\bf Proof of the Theorem 6.2.}

Justification for computing algorithm 2 is performed
in space $X=H_{\beta} \quad (\beta< \alpha/4)$ and its subspace $X_n,$
which consists of polynomials (6.2).

Repeating the calculations made in the proof of the previous theorem, we can show
that the operator $ K $ has  Frechet derivative
\[
K'(x)z \equiv a(t) z(t)+\frac{1}{\pi i}\int\limits_{\gamma}
h'_u(t,\tau,x(\tau))z(\tau)(\tau-t)^{-1}d\tau,
\]
satisfying in the sphere  $\|x-x^0\| \le r$ the Holder condition
$\|K'(x_1)-K'(x_2)\| \le A\|x_1-x_2\|^{\alpha-\beta},$ where $x^0$ is a element of the space $X,$ and $A$ is a constant, which is  dependent only from $r$ and $x_0.$

The system of equations (6.7) in the operator form is written in the form of expression
\[
K_n(x_n)\equiv \bar P_n[\tilde a(s) x_n(s)-\frac{i}{2\pi}\int\limits_0^{2\pi}
P_n^{\sigma}[\tilde h(s,\sigma,x_n(\sigma))\cot\frac{\sigma-s}{2}]d\sigma+
\]
\[
+\frac{1}{2\pi}\int\limits_0^{2\pi}
P_n^{\sigma}[\tilde h(s,\sigma,x_n(\sigma))]d\sigma]=\bar P_n[\tilde f],
\]
where
\[
\tilde a(s)=a\left(s-\frac{\pi}{2n+1}\right), \quad \tilde h(s,\sigma,x_n(\sigma))=
h\left(s-\frac{\pi}{2n+1},\sigma,x_n(\sigma)\right),
\]
\[
\tilde f(s)=f\left(s-\frac{\pi}{2n+1}\right).
\]

The Frechet derivative of the operator $K_n(x_n)$ is
\[
K'_n(x_n)z_n \equiv \bar P_n[\tilde a(s)z_n(s) -\frac{i}{2\pi}\int\limits_0^{2\pi}
P_n^{\sigma}[\tilde h'_u(s,\sigma,x_n(\sigma))z_n(\sigma) \cot\frac{\sigma-s}{2}]d\sigma+
\]
\[
+\frac{1}{2\pi}\int\limits_0^{2\pi}
P_n^{\sigma}[\tilde h'_u(s,\sigma,x_n(\sigma))z_n(\sigma)]d\sigma],
\]
where $\tilde h'_u(s,\sigma,x_n(\sigma))=h'_u(s-\pi/(2n+1),\sigma,x_n(\sigma)).$

The arguments that differ little from the proof of Theorem
6.1, suggest that \[\|K'_n(x'_n)-K_n(x''_n)\| \le A n^{\beta}
\ln^2 n\|x'_n-x''_n\|^{\alpha}.\]

Let $x_n^0$ is the polynomial of best approximation of the degree $n$
for function  $x^*(t).$ 

Let the operator $K'(x)$ has   bounded right inverse operator $[K'(x)]^{-1}_r$ $(x \in S_1)$ with the norm $\|[K'(x)]^{-1}_r\|=B_0.$

Let us prove the existence of the bounded right inverse operator $[K'(x_n^0)]^{-1}_r(x \in S_1).$ 

The existence of a right inverse operator $[K'(x_n^0)]^{-1}_r$
with such  $n$  that $q=A n^{-(\alpha-\beta)}<1,$ follows from the results of the Section 2 of Introduction. Note that
 $\|[K'(x_n^0)]^{-1}_r\|=B_0/(1-q).$ 
 
 As in the \S2, we can see that, for such
 $n$  that $q=A n^{-(\alpha-\beta)}\ln^2 n<1,$ the operator
$K'_{1n}(x_n^0),$ where
\[
K'_{1n}(x_n^0)z_n \equiv
\]
\[
 \equiv \bar P_n[a(t)z_n(t) -\frac{i}{2\pi}\int\limits_0^{2\pi}
P_n^{\sigma}[h'_u(s,\sigma,x_n^0(\sigma))z_n(\sigma) \cot\frac{\sigma-s}{2}]d\sigma+
\]
\[
+\frac{1}{2\pi}\int\limits_0^{2\pi}
P_n^{\sigma}[h'_u(s,\sigma,x_n^0(\sigma))z_n(\sigma)]d\sigma],
\]
has the \ linear \ invertable operator  $[K'_{1n}(x_n^0)]^{-1}$ with the norm $\|[K'_{1n}(x_n^0)]^{-1}\|
\le A \ln n.$

We will  show that if in some neighborhood
$ S_1 $ of the point $ x^* $  exists a bounded right inverse operator
$[K'(x)]^{-1}_r(x \in S_1)$ with the norm $\|[K'(x)]^{-1}_r\|=B_0,$
that  there exists a linear operator $[K'_n(x_n^0)]^{-1}_r$ for 
sufficiently large $n$.

To prove the existence of a linear operator $[K'_{n}(x_n^0)]^{-1}$ we
need  in assessment of the norm
\[
\|K'_{1n}(x_n^0)z_n-K'_{n}(x_n^0)z_n\| \le \|\bar P_n[[a(t)-\tilde a(t)]z_n(t)\|+
\]
\[
+\|\frac{1}{2\pi}\bar P_n[\int\limits_0^{2\pi}
P_n^{\sigma}[[h'_u(s,\sigma,x_n^0(\sigma))-\tilde h'_u(s,\sigma,x_n^0(\sigma))]
z_n(\sigma)\cot\frac{\sigma-s}{2}]d\sigma]\|+
\]
\[
+\|\frac{1}{2\pi}\bar P_n[\int\limits_0^{2\pi}
P_n^{\sigma}[[h'_u(s,\sigma,x_n^0(\sigma))-\tilde h'_u(s,\sigma,x_n^0(\sigma))]
z_n(\sigma)]d\sigma]\|=
\]
\[
=\|J_1\|+\|J_2\|+\|J_3\|.
\]

Easy to see that
\[\|J_1\|+\|J_3\| \le A \ln^2 n \|z_n\|n^{-(\alpha-\beta)},
\]
\[
|I_2| \le A \ln n \max\limits_{0\le k \le 2n}\left|\int\limits_0^{2\pi}
P_n^{\sigma}[[h'_u(s'_k,\sigma,x_n^0(\sigma))-
\tilde h_u(s_k,\sigma,x_n^0(\sigma))]\times\right.
\]
\[
\left.\left.\times z_n(\sigma)\cot\frac{\sigma-\bar s_k}{2}\right]d\sigma\right|\le
A \ln^2 n \|z_n\|n^{-\alpha}.
\]

Since $I_2$ is the trigonometrical polynomial of the order $n,$ that
\[
\|I_2\| \le A\ln^2 n\|z_n\|n^{-(\alpha-\beta)}.
\]

From estimates for  $\|J_1\|,\|J_2\|,\|J_3\|$ and Banach Theorem it follows, that,
for such $n$ that $q=An^{-(\alpha-\beta)}\ln^2 n<1,$ the linear operator  $[K'_n(x_n^0)]^{-1}$ exists with the norm $\|[K'_n(x_n^0)]^{-1}\|=B_2 \le B_1/(1-q).$ 

Was proved (see the proof of the  Theorem 6.1) that $\|K'_n(x_n^0)\|\le An^{-(\alpha-\beta)}\ln n.$

Equation (6.7) satisfies all  conditions of Theorem 2.6,
given in the introduction.

So, we proved the existence of f unique solution $x_n^*$ of the equation
(6.7).  The estimate of norm $\|x^*-x^*_n\|$ follows from Theorem 2.6, given in the introduction. 

{\bf Proof of the Theorem 6.3.} 
The proof of the Theorem 6.3 is similar to the proof of the Theorem 6.2.

{\bf Proof of the Theorem 6.4.}

We will show that under the conditions of "b"
the operator $ K $ has the Frechet derivative  in the space $ L_2: $
\[
K'(x)z_n \equiv a'(t,x(t))z_n(t)+S_{\gamma}(h'_u(t,\tau,x(\tau))z_n(\tau)),
\]
and inequality
$$
\|K'(x_1)-K'(x_2)\| \le AF \eqno (6.17)
$$
is valid.

To prove the existence of the derivative $ K'(x) z_n $, let us estimate the difference
\[
\|I_1\|/\|z_n\|=\|K(x+z_n)-K(x)-K'(x)z_n\|/\|z_n\|.
\]

We restrict ourselves to estimate the integral term
in the expression for $ K'(x) z_n. $

Obviously,
\[
|I_1|=\left|\frac{1}{\pi i}\int\limits_{\gamma}\left\{\int\limits_0^1
(1-\nu)h''_u(t,\tau,x(\tau)+\nu z_n(\tau))d\nu\right\}\frac{z^2_n(\tau)
d\tau}{\tau-t}\right|.
\]

Represent  $h''_u(t,\tau,u)$ in the form
 $h''_u(t,\tau,u)=[h''_u(t,\tau,u)-h''_u(\tau,\tau,u)]+
h''_u(\tau,\tau,u).$ 

It is easy to see that
 $\|I_1\| \le A\|z_n\|\max|z_n|$
and
 $\|I_1\|/\|z_n\| \le An^{1/2}\|z_n\|.$

From this estimate follows that
$\lim\limits_{\|z_n\| \to 0} \|I_1\|/\|z_n\|=0.$ The inequality (6.17) is proved.

It is easy to see that the Frechet derivative of the operator $K_n$ is
\[
K'_n(x_n)z_n \equiv
\]
\[ \equiv \bar P_n^s[a'_u(s,x_n(s))z_n(s)+\frac{1}{2\pi}\int\limits_0^{2\pi}P_n^{\sigma}[h'_u(s,\sigma,x_n(\sigma))
z_n(\sigma)]d\sigma]-
\]
\[
-\frac{i}{2\pi}\int\limits_0^{2\pi}P_n^{\sigma}[h'_u(s,\sigma,x_n(\sigma))
z_n(\sigma)\cot\frac{\sigma-s}{2}]d\sigma.
\]

Let us proved  justice  of the inequality 
$$
\|K'_n(x'_n)-K'_n(x''_n)\| \le AF. \eqno (6.18)
$$

Indeed, $\|K'_n(x'_n)z_n-K'_n(x''_n)z_n\|\le I_2+I_3,$ 
where
\[
I_2=\| \bar P_n[S_{\gamma}(P_n^{\tau}[h'_u(\tau,\tau,x'_n(\tau))-
h'_u(\tau,\tau,x''_n(\tau)))z_n(\tau)])]\|,
\]
\[
I_3=\| \bar P_n[\frac{1}{\pi i}\int\limits_{\gamma}
P_n^{\tau}[\frac{h^*(t,\tau,x'_n(\tau))-h^*(t,\tau,x''_n(\tau))}
{|\tau-t|^{1-\beta}}z_n(\tau)]]\|=
\]
\[
=\| \bar P_n[\frac{1}{\pi}\int\limits_0^{2\pi}
P_n^{\sigma}[[h^*(s,\sigma,x'_n(\sigma))-h^*(s,\sigma,x''_n(\sigma))]
z_n(\sigma)e^{is}p(s,\sigma)]d\sigma]\|.
\]
Here
 $p(s,\sigma)=|e^{i\sigma}-e^{is}|^{\beta-1}$ for $|\sigma-s|\ge \pi/(2n+1)$
and $p(s,\sigma)=|e^{i\pi/(2n+1)}-1|^{\beta-1}$ for $|\sigma-s|< \pi/(2n+1),$
$h^*(t,\tau,x(\tau))=(h(t,\tau,x(\tau))-h(\tau,\tau,x(\tau)))(\cot\frac{\tau-t}{2})
|\tau-t|^{1-\beta}.$

It is easy to see that
\[
I_2 \le A\|P[[h'_u(\tau,\tau,x'_n(\tau))-
h'_u(\tau,\tau,x''_n(\tau))]z_n(\tau)]\| \le
\]
\[
\le A\max|h'_u(t,t,x'_n(t))-h'_u(t,t,x''_n(t))\|z_n\|.
\]

We proceed to the evaluation of the expression $I_3.$ 
If  integrand member in $ I_3 $
 is a real function then
 $I_3 \le  I_4  I_5,$ where
\[
I_4 =\max\{\frac{1}{\pi}\int\limits_0^{2\pi}
\{P_n^{\sigma}\{h^*(s,\sigma,x'_n(\sigma))-h^*(s,\sigma,x''_n(\sigma))]
[p(s,\sigma)]^{1/2}\}^2d\tau\}^{1/2},
\]
\[
I_5 =\{\frac{1}{2\pi}\int\limits_0^{2\pi}ds
\frac{1}{\pi}\int\limits_0^{2\pi}
\{P_n^{\sigma}[z_n(\sigma)e^{i\sigma}[\bar P_n^s[p(s,\sigma)]]^{1/2}]\}^2
d\sigma\}^{1/2}.
\]

One can see that
\[
I_4 \le A \max|h^*(s,s,x'_n(s))-h^*(s,s,x''_n(s))|.
\]

Similar,
\[
I_5=\{\frac{2}{(2n+1)^2}\sum\limits^{2n}_{i=0}\sum\limits^{2n}_{k=0}
|z_n(s_k)|^2 p(s_i,s_k)\}^{1/2}=
\]
\[
=\{\frac{2}{(2n+1)^2}\sum\limits^{2n}_{k=0}|z_n(s_k)|^2
\sum\limits^{2n}_{i=0}p(s_i,s_k)\}^{1/2} \le A\|z_n\|.
\]

From these calculations, it follows that
$I_3 \le AF||z_n\|.$  From estimates of $I_2$ and $I_3$
the inequality (6.18) is followed.

Now we can prove the existence of the linear operator $[K'(x_n^0)]^{-1}$ with the norm
$\|[K'(x_0)]^{-1}\| \le \|[K'(x_n^0)]^{-1}\|/(1-q_1).$

It follows, at such $n$ that $q_1=A\ln n E_n(x_0)<1$, from existence of the linear operator $[K'(x_0)]^{-1}$ and Banach Theorem.

In the section  3 was proved, at such  $n$ that\\
$q_2=\ A \ln^2 n \max \ [\ E_n(h'_u(t,t,x^0_n(t))),\ \ \    
E_n^{t,\tau}(h'_u(t,\tau,x^0_n(\tau))), \ \ \ $  $E_n(\psi(t))]<1,
$
 the existence of the linear operator 
$[K_n'(x_n^0)]^{-1}$ with norm $B_0.$

The validity of the theorem follows from the results of the item 2 of Introduction.

\begin{center}
{\bf 6.2. Projection Methods for Solving Nonlinear
	Singular Integral Equations
	on Open Contours of Integration}
\end{center}

Let us consider nonlinear singular integral equation
$$
G(x) \equiv a(t,x(t)) + \frac{1}{\pi i} \int\limits_L
\frac{h(t, \tau, x(\tau))}{\tau-t} d\tau = f(t),
\eqno(6.19)
$$
where  $L = (c_1, c_2)$ is a segment of the unit circle $\gamma$, centered at the origin. 

Approximate solution of the equation  (6.19) we will search as polynomial
\[
x_n(t) = \sum_{k=-n}^{n}\alpha_k t^k. 
\]

Coefficients $\{\alpha_k\}$ are determined from the system
\[
G_n(x_n) \equiv \bar{P}_n\left[a^*(t,x_n(t)) + \frac{1}{\pi i} \int\limits_L
P_{n}^{\tau}\left[\frac{h^*(t,\tau,x_n(\tau))}{\tau-t}\right]d\tau\right] =
\]
$$
=\bar{P}_n[f^*(t)],
\eqno(6.20)
$$
where
\[
a^*(t,x_n(t)) =
\left\{\begin{array}{lll}
a(t,x(t)) & \mbox{if $t \in L$,} \\
x(t)      & \mbox{if $t \not\in L$;}
\end{array}
\right.
\]

\[
h^*(t,\tau,x_n(\tau )) =
\left\{\begin{array}{lll}
h(t,\tau,x_n(\tau )) & \mbox{if $t$ и   $\tau \in L$,} \\
0                  & \mbox{if $t$ или $\tau \not\in L$;}
\end{array}
\right.
\]

\[
f^*(t) =
\left\{\begin{array}{lll}
f(t) & \mbox{if $t \in L$,} \\
0    & \mbox{if $t \not\in L$.}
\end{array}
\right.
\]

Justification of the computational scheme (6.20) is carried out
along the same lines as the study of numerical schemes,
cited in the preceding paragraph.

The  difference is
that, when communication between  reversibility operators
$ G $ and $  G_n$ is researched, we need to use
results of the section \S4 instead of the results of sections \S2 and \S3. 

Omitting the intermediate calculations, we formulate statements
which are similar to statements   of the preceding paragraph.

{\bf Theorem 6.5}  \  \cite{Boy17}, \cite{Boy25}. \ 
Let \ in \  the \ ball \  $S(x^*,\tau),$ \ where \\ $r = B_*B_0^2 / (1-q)$, $ B_*=\|R(x^*)\|
+ O(n^{-\alpha}\ln n)$,\\ $q=O(F_1 B_*B_0^2) + \sqrt{1-1/( B_*^4B_0^4)}<1$,\\
$R(x^*)z \equiv a'_u(t,x^*(t))z(t) + \frac{1}{\pi i}
\int\limits_L h'_u(t,\tau, x^*(\tau))z(\tau)(\tau-t)^{-1}d\tau ,$\\
the following conditions are fulfilled:

1)~$ a'_u(t,u) \in H_{\alpha\alpha} , h'_u(t,\tau,u)
\in H_{\alpha\alpha\alpha} , f(t) \in H_{\alpha} $ $(0<\alpha<1)$;

2)~$\max\limits_{t,\tau \in L} \bigl\{|{h'_u(t,\tau,u_1)-
	h'_u(t,\tau,u_2)}|$, $|{h^*(t,\tau,u_1)- h^*(t,\tau,u_2)}|
\bigr\}\le F_1$;

3)~the Frechet derivative  $R(x^*)z$ of the operator $G(x^*)$ is invertible in the space
 $L_2$ and $\|[ R(x^*)]^{-1}\|
\le B_0$;

4)~characteristic equation
\[
a'_u(t,x^*(t))z(t) + \frac{ h'_u(t,t, x^*(t))}{\pi i} \int\limits_L \frac{z(\tau)d\tau }{\tau-t} = 0
\]
has a solution $z(t) = (t-c_1)^{\delta_1} (t-c_2)^{\delta_2}\varphi(t)$,
where $\delta_1=\zeta_1 + i\xi_1$, $\delta_2=\zeta_2 + i\xi_2$, $\varphi(t) \in H_{\alpha}$.

Than, for such $n$ that $O((n^{-\alpha}+ n^{-\theta})\ln n)<1$
($\theta=1$ for $\zeta>0$, $\theta=1-\zeta$ for $\xi \le 0,$ $\zeta=\min(\zeta_1,
\zeta_2)$), the system (6.20) has a unique solution
$x_n^*$ and the estimate  $\|x^*-x_n^*\| = O((n^{-\alpha}+ n^{-\theta})
\ln n)$ is valid. Here  $x^*$ is a solution of the equation  (6.19).

\begin{center}
{\bf 7. Spline-Collocation Method}
\end{center}

\begin{center}
{\bf 7.1. Linear Singular Integral Equations}
\end{center}

Let us consider the singular integral equation
$$
Kx \equiv a(t)x(t)+b(t)\int\limits_{-1}^1 \frac{x(\tau)d\tau}
{\tau-t}+\int\limits_{-1}^1h(t,\tau)x(\tau)d\tau=f(t), \eqno (7.1)
$$
where $a,b,f \in W^r(1),$ $h(t,\tau) \in W^{r,r}(1),$ $b(t) \ne 0,$
$a^2(t) + b^2(t) \geq C >0.$

Methods are discussed below may be extended to singular integral equations on piecewise continuous loops. It works in exceptional  cases too.

Introduce points $t_k=-1+2k/n,$ $k=0,1,\ldots,n,$ and
$t_{kj}=t_k+jh/(r+1),$ $j=1,2,\ldots,r,$ $k=0,1,\ldots,n-1,$ $h=2/n.$
Let $\Delta_k=[t_k,t_{k+1}],$ $k=0,1,\ldots,n-1.$

Let $f(t)\in C[-1,1].$ Let us construct on each segment $\Delta_k$ $(k=0,1,\ldots,n-1)$
interpolated polynomial $L_r(f,\Delta_k).$ The polynomial $L_r(f,\Delta_k)$
interpolate the function $f(t)$ on the segment $\Delta_k$ on knots
$t_{kj},$ $j=1,2,\ldots,r,$ $k=0,1,\ldots,n-1.$

The polynomial $L_r(f,\Delta_k)$ has the form
$$
L_r(f,\Delta_k)=\sum\limits^r_{j=1}f(t_{kj})\psi_{kj}(t),
$$
where $\psi_{kj}(t)$  are fundamental polynomials on knots
$t_{kj},$ $j=1,2,\ldots,r,$ $k=0,1,\ldots,n-1.$

Spline $f_n(t),$ $-1 \leq t \leq 1,$ consists of polynomials $L_r(f,\Delta_k),$ $k=0,1,\ldots,n-1.$

To each knot $t_{kj}$ we put in correspondence the segment
$\Delta_{kj}=[t_{kj}-qh^*, t_{kj}+h^*],$ where $h^* (0 < h^* < h/(r+1))$
and $q$  are parameters, values of which will be determined below.

Approximate solution of the equation (7.1) we will seek in
the form of the spline $x_n(t),$  which is composed from polynomials $L_r(x,\Delta_k).$
Values $x_{kj}=x(t_{kj}),$ $j=1,2,\ldots,r,$ $k=0,1,\ldots,n-1,$
of these polynomials are determined from the system
\[
a(t_{kl})x_{kl}+b(t_{kl})\int\limits_{\Delta_{kl}}\frac{x_n(\tau)}
{\tau-t_{kl}}d\tau+\sum\limits^{n-1}_{i=0}{}'b(t_{kl})
\int\limits_{\Delta_{i}}\frac{x_n(\tau)}{\tau-t_{kl}}d\tau+
\]
$$
+\sum\limits^{n-1}_{i=0}\int\limits_{\Delta_{i}}
h(t_{kl},\tau) x_n(\tau)d\tau=f(t_{kl}), 
\eqno (7.2)
$$
$k=0,1,\ldots,n-1,$ $l=1,2,\ldots,r,$
where the prime in the summation indicates that $i \ne k-1, k, k+1.$

Let us prove that the parameter $h^*$ can be chosen so way, that the
system (7.2) has a unique solution. This proof based on Hadamard theorem
about invertibility of matrix.

Let us write the system (7.2) as the matrix equation
$$
C X=F, \eqno (7.3)
$$
where $X=(x_1,\ldots,x_N),$ $F=(f_1,\ldots,f_N),$ $C=\{c_{ij}\}, \,
i,j=\overline{1,N},$ $N=nr.$

Here $x_l=x_{ij},$ $f_l=f_{ij},$ where $l=ri+j,$ $i=0,1,\ldots,n-1,$
$j=1,2,\ldots,r,$ $x_{ij}=x(t_{ij}),$ $f_{ij}=f(t_{ij}).$

Let $l=ri+j,$ $i=0,1,\ldots,n-1,$ $j=1,2,\ldots,r,$ $k=rv+w,$
$v=0,1,\ldots,n-1,$ $w=1,2,\ldots,r.$ Then elements $c_{kl}$
of matrix $C$ have the form
\[
c_{ll}=a(t_{ij})+b(t_{ij})\int\limits_{\Delta_{ij}}\frac{\psi_{ij}(\tau)d\tau}
{\tau-t_{ij}}+
\]
\[
+\int\limits_{\Delta_i}\psi_{ij}(\tau) h(t_{ij},\tau) d\tau,
\quad l=1,2,\ldots,N,
\]
\[
c_{lk}=b(t_{ij})\int\limits_{\Delta_{ij}}\frac{\psi_{vw}(\tau)d\tau}
{\tau-t_{ij}}
+\int\limits_{\Delta_i}\psi_{vw}(\tau) h(t_{ij},\tau) d\tau,
\]
if $t_{vw}\in \Delta_i$ and
\[
c_{lk}=\int\limits_{\Delta_v}\psi_{vw}(\tau) h(t_{ij},\tau) d\tau,
\]
if $t_{vw} \bar\in \Delta_i.$

Let us estimate from below  the diagonal  coefficients $c_{ll}.$

 If $b(t_{ij}) \not = 0,$
then
$$
|c_{ll}| \ge |b(t_{ij})|\left|\int\limits_{\Delta_{ij}}
\frac{\psi_{ij}(\tau)d\tau}{\tau-t_{ij}}\right|-|a(t_{ij})|-
\left|\int\limits_{\Delta_i}\psi_{ij}(\tau) h(t_{ij},\tau) d\tau\right|.
\eqno (7.4)
$$

Easy to see that
$|a(t_{ij})| \le A,$
$\left|\int\limits_{\Delta_i}\psi_{ij}(\tau) h(t_{ij},\tau) d\tau\right|=
O(n^{-1}),$ where $A$ is $a$ constant, which is not depend from indexes
$i,j,$ $i=0,1,\ldots,n-1,$
$j=1,2,\ldots,r.$

Well known that $\psi_{ij}(t_{ij})=1.$
Let us show, that, for  big enough number $M$, exist parameters $h^*$ and $q$,
under which  the  inequalities
$$
\int\limits_{\Delta_{ij}}\frac{\psi_{ij}(\tau)d\tau}{\tau-t_{ij}} \ge M, \ i=0,1,\ldots,n-1,\ j=1,2,\ldots,r,
\eqno (7.5)
$$
is valid.

Indeed
\[
\left|\int\limits_{\Delta_{ij}}\frac{\psi_{ij}(\tau)d\tau}{\tau-t_{ij}}\right| \ge
\left|\int\limits_{\Delta_{ij}}\frac{d\tau}{\tau-t_{ij}}\right|-
\]
$$
-\left|\int\limits_{\Delta_{ij}}\frac{\psi_{ij}(\tau)-1}{\tau-t_{ij}}d\tau\right| \ge
|\ln q|-
\left|\int\limits_{\Delta_{ij}}\frac{\psi_{ij}(\tau)-1}{\tau-t_{ij}}d\tau\right|.
\eqno (7.6)
$$

We can chose the parameter $q$ so way that $|\ln q| \ge M+1+a.$
Also we can chose the parameter  $h^*$ so way that
$$
\left|\int\limits_{\Delta_{ij}}\frac{\psi_{ij}(\tau)-1}{\tau-t_{ij}}
d\tau\right| \le \varepsilon, \eqno (7.7)
$$
where $\varepsilon (\varepsilon>0)$  is a arbitrary shall number.

So, we can chose parameters $q$ and $h^*$ so way, that
$|c_{ll}| \ge M,$ $l=1,2,\ldots,N,$ where $M$  is a arbitrary big number.

Let us estimate coefficients $|c_{kl}|$ for $k \ne l,$ $l,k=1,2,\ldots,N.$

The function
$\psi_{vw}(t_{ij})=0$ for $v \ne i$ and $w \ne j.$
Using this fact one can see that
$$
\left|\int\limits_{\Delta_{ij}}\frac{\psi_{vw}(\tau)d\tau}{\tau-t_{ij}}\right|
=O\left(\frac{1}{n}\right) \eqno (7.8)
$$
and
$$
\left|\int\limits_{\Delta_{i}}\psi_{vw}(\tau)h(t_{ij},\tau) d\tau\right|=
o\left(\frac{1}{n}\right). \eqno (7.9)
$$

Collecting inequalities (7.4) $-$ (7.9) we see, that exist such parameters $h^*$
and $q$,  that conditions of Hadamard theorem on
invertibility of matrices are valid. From Hadamard theorem it follows, that the
system (7.2) has a unique solution. Obviously,  the system (7.3) has a unique solution too.

Easy to see that, in the space $R_N$ of vectors $X_N=(x_1,\ldots,x_N)$
with the norm $\|X_N\|=\max\limits_{1 \le k \le N} |x_k|$, the norm
of the matrix $C^{-1}$ is evaluated as $\|C^{-1}\| \le A/M.$

Let $P_N$ is the projector from space $C[-1,1]$ onto set of vector-functions: $P_Nf={f(t_{ij})},$
$i=0,1,\ldots,n-1,$
$j=1,2,\ldots,r.$
 This is defined as $P_N[f]=f_N(t).$

In the operator form the system (7.2) can be written as
\[
K_N x_N \equiv 
\]
\[
 \equiv P_N \left[a(t) x_nt)+b(t) \int\limits^1_{-1} e(t,\tau)\frac{x_n(\tau)}{\tau-t}d\tau+\int\limits_{-1}^1h(t,\tau)x(\tau)d\tau\right]=
 \]
 \[
 = P_N[f(t)].
\]

The function $e(t,\tau)$  is defined by formula
\[
e(t_{kl},\tau) = \left\{
\begin{array}{cc}
1 \, {\rm if} \, \tau \in \Delta_{kl} \, {\rm or} \, \tau \in [-1,1] \backslash
(\Delta_{k-1} \cup \Delta_k \cup \Delta_{k+1}),\\
0 \, {\rm if} \, \tau \in (\Delta_{k-1} \cup \Delta_{k} \cup \Delta_{k+1})
\backslash \Delta_{kl}.\\
\end{array}
\right.
\]

The estimate $\|K_N^{-1}\| \le A/M$ is valid.

Let us estimate the error of the offered method.

Let $x^*(t)$ be a solution of the equation (7.1). Equating left
and right sides of the equation (7.1) in points $t_{kl},$
$k=0,1,\ldots,n-1,$ $l=1,2,\ldots,r,$ we have
\[
a(t_{kl}) x^*_{kl}
+\sum\limits^{n-1}_{i=0} b(t_{kl})
\int\limits_{\Delta_{i}}\frac{x^*(\tau)}{\tau-t_{kl}}d\tau+
\]
$$
+\sum\limits^{n-1}_{i=0}\int\limits_{\Delta_{i}}
h(t_{kl},\tau)x^*(\tau)d\tau=f(t_{kl}), \eqno (7.10)
$$
$k=0,1,\ldots,n-1,$ $l=1,2,\ldots,r,$ where $x^*_{kl}=x^*(t_{kl}).$

Subtracting from the system (7.2) the system (7.10) we have
\[
a(t_{kl})(x_{kl} - \tilde x^*_{kl})+b(t_{kl})
\int\limits_{\Delta_{kl}}\frac{x_n(\tau)-\tilde x^*_n(\tau)}
{\tau-t_{kl}}d\tau+
\]
\[
+\sum\limits^{n-1}_{i=0}{}'b(t_{kl})
\int\limits_{\Delta_{i}}\frac{x_n(\tau)-\tilde x^*_n(\tau)}
{\tau-t_{kl}}d\tau+\sum\limits^{n-1}_{i=0}\int\limits_{\Delta_{i}}
h(t_{kl},\tau)(x_n(\tau)-\tilde x^*_n(\tau) d\tau=
\]
\[
=b(t_{kl})\int\limits_{\Delta_{kl}}\frac{x^*(\tau)-\tilde x^*_n(\tau)}
{\tau-t_{kl}}d\tau+
b(t_{kl})\int\limits_{\Delta^*_{kl}}\frac{x^*(\tau)}{\tau-t_{kl}}d\tau+
\]
$$
+\sum\limits^{n-1}_{i=0}{}'b(t_{kl})
\int\limits_{\Delta_{i}}\frac{x^*(\tau)-\tilde x^*_n(\tau)}
{\tau-t_{kl}}d\tau+
\sum\limits^{n-1}_{i=0}\int\limits_{\Delta_{i}}
h(t_{kl},\tau)(x^*(\tau)-\tilde x^*_n(\tau)) d\tau, \eqno (7.11)
$$
$k=0,1,\ldots,n-1,$ $l=1,2,\ldots,r,$\\
where $\Delta^*_{kl}=\Delta_{k-1} \bigcup (\Delta_{k}\setminus \Delta_{kl})
\bigcup \Delta_{k+1},$ $x^*$ and $x^*_n$ are solutions of equations (7.1) and (7.2).

Here $\tilde x^*_n(\tau)$ is the local spline,  which  approximate the solution
$x^*(\tau)$ on the set of knots  $t_{k,l}, k=0,\ldots,n-1, l=1,2,\ldots,r.$  

Let us estimate the right side of the expression (7.11).

Easy to see that
\[
r_1=\left|b(t_{kl})\int\limits_{\Delta_{kl}}
\frac{x^*(\tau)- \tilde{x}^*_n(\tau)}{\tau-t_{kl}}d\tau\right| \le
\]
\[
\le \left|b(t_{kl})\int\limits_{\Delta_{kl}}
\frac{x^*(\tau)- \tilde x^*_n(\tau)-({x}^*(t_{kl})- \tilde{x}^*_n(t_{kl}))}
{\tau-t_{kl}}d\tau\right| \le An^{-r};
\]
\[
r_3= \left|\sum\limits^{n-1}_{i=0}{}'b(t_{kl})
\int\limits_{\Delta_{i}}\frac{x^*(\tau)- \tilde x^*_n(\tau)}
{\tau-t_{kl}}d\tau \right| \le An^{-r} \ln n,
\]
\[
r_4=\left|\sum\limits^{n-1}_{i=0}\int\limits_{\Delta_{i}}
h(t_{kl},\tau)(x^*(\tau)- \tilde x^*_n(\tau)) d\tau \right|\le An^{-r}.
\]

For evaluating the integral
\[
r_2=\left|b(t_{kl})\int\limits_{\Delta^*_{kl}}
\frac{x^*(\tau)}{\tau-t_{kl}}d\tau\right|
\]
we introduce a function $\psi(\tau),$ which satisfy the following
conditions:

$ 1) \int\limits_{\Delta^*_{kl}}
\frac{\psi(\tau)}{\tau-t_{kl}}d\tau=0;$

2) in the domain $\Delta^*_{kl}$ the function $\psi(t)$ realize the
best approximation for the function $x^*(t)$ in the metric of the
space $C.$

Introduce the following designations:
\[
E^*_{kl}(x^*)=\left|\int\limits_{\Delta^*_{kl}}
\frac{x^*(\tau)-\psi(\tau)}{\tau-t_{kl}}d\tau\right|,
\]
\[
E^*(x^*)=\max\limits_{k,l} E^*_{kl}(x^*).
\]

Here $\psi(t)$ satisfy to conditions 1) and 2).

Obviously, $r_2 \le BE^*(x^*),$ where $B=\|b(t)\|_C.$

Now it is easy to see that
$$
\max |x^*(t_{kl})-x_n(t_{kl})| \le A(n^{-r} \ln n+E^*(x^*)).
$$

{\bf Theorem 7.1} \cite{Boy25}, \cite{Boy28}, \cite{Boy29}.
Let $a,b,f \in W^r,$ $h \in W^{rr},$ $|b(t)| \geq C >0 $
on the segment $[-1,1].$
Let the equation (7.1) has a unique solution $x^*(t).$
Then exist such parameters $h^*$ and $q$, that the system of equations (7.2)
has a unique solution $x_n^*$ and estimate
$\|x^*(t_{kl})-x^*_n(t_{kl})\|_C \le A(n^{-r}\ln n+E^*(x^*))$ is valid.

\begin{center}
{\bf 7.2. Nonlinear Equations}
\end{center}

Let us consider the nonlinear singular integral equation
\[ Kx \equiv a(s)x(s)+S(b(s,\sigma,x(\sigma)))+T(h(s,\sigma,x(\sigma)))
\equiv
\]
\[ 
\equiv a(s)x(s)+
\frac{1}{2\pi}\int\limits_0^{2\pi} b(s,\sigma,
x(\sigma))\cot
\frac{\sigma-s}{2}d\sigma+\int\limits_0^{2\pi}h(s,\sigma,x(\sigma))d\sigma=
\]
$$
= f(s). \eqno (7.12)
$$

Let us assume that $ a(s),f(s) \in H_\alpha $, $ b'_3(s,\sigma,u)
\in H_{\alpha,\alpha,\alpha} $, $ 0<\alpha\le 1 $, where
$ b'_3(s,\sigma,u) \equiv \partial b(s,\sigma,u) / \partial u.$

Easy to see, that Frechet derivative of the nonlinear operator $K(x)$ in the space
$X=H_\beta $ $ (0<\beta<\alpha)$ has the form
\[
K'(x_0)z\equiv a(s)z(s)+
\frac{1}{2\pi}\int\limits_0^{2\pi}b'_3(s,s,x_0(s))z(\sigma)
\cot\frac{\sigma-s}{2}d\sigma+
\]
\[
+\frac{1}{2\pi}\int\limits_0^{2\pi}
h'_3(s,\sigma,x_0(\sigma))z(\sigma)d\sigma.
\]

We assume, that the function $ a^2(s)-(b'_3(s,s,x_0(s)))^2 $
can be equal to zero on sets of points with measure not bigger than zero.
For construction numerical scheme for solution of the equation (7.12),
we choose the sets of knots
\[ s_k=\frac{\pi k}{n}, \hskip 5pt s_k^*=s_k+h, \hskip 5pt
0<h\le\frac{\pi}{2n}, \hskip 5pt k=0,...,2n. 
\]
The value of the parameter $ h $ will be defined later.

Approximate solution of the equation (7.12)  is defined from the following system of nonlinear
equations
\[ a(s_j^*)x(s_j^*)+\frac{1}{2\pi}\sum_{k=0,k\neq j-1,j+1}^{2n-1}
b(s_j^*,s_k^*,x(s_k^*))
\int\limits_{s_k}^{s_{k+1}}\cot\frac{\sigma-s_j^*}{2}d\sigma+ 
\]
$$ +\frac{\pi}{n}\sum_{k=0}^{2n-1}h(s_j^*,s_k^*,x(s_k^*))=f(s_j^*),
\hskip 5pt j=0,...,2n-1.
\eqno (7.13) $$

For solution of the system (7.13) we use the Newton - Kantorowich method.

After found solution  $ x(s_j^*) $, $ j=0,1,...,2n-1, $
of nonlinear system of equations (7.13), we restore the approximate solution of the equation
(7.12) in the form of trigonometric polynomial
$$ x_n(s)=\sum_{k=0}^{2n-1}x(s_k^*)\psi_k(s), \eqno (7.14) $$
where
\[ \psi_k(s)=\frac{1}{2n+1} \frac{\sin \frac{2n+1}{2}(s-s_k^*)}
{\sin \frac{s-s_k^*}{2}}.
\]

One can restore the approximate solution of the equation (7.13)
in the form of polygons, constructed on  knots $ s^*_k, x(s^*_k) $, $ k=0,...,2n-1 $.

Let $ P_n $ is the projector from the space $ X $ onto
the space $ X_{n}\subset X. $ The space $X_n$ consist of trigonometric
polynomial of $n$ order. The projector $P_n$ is introduced by the formula
\[
 P_nx=\sum_{k=0}^{2n-1} x(s_k^*)\psi_k(s). 
 \]

In the space $ X_n $ the system (7.13) can be written as operator
equation $ K_n x_{n} =f_n.$

The Frechet derivative $ K'_n(x_0)z_n, z_n\in X_n, $ of the operator $ K_n $
on a element $ x_0 $ can be written as the vector
\[
 a(s_j^*)z_n(s_j^*) +
 \]
 \[
 +\frac{1}{2\pi}\sum_{k=0,k\neq j-1,j+1}^{2n-1}
b'_3(s_j^*,s_k^*,x_0(s_k^*))z_n(s_k^*)\int\limits_{s_k}^{s_{k+1}}\cot\frac{\sigma-
s_j^*}{2}d\sigma+ 
\]
\[
+\frac{\pi}{n}\sum_{k=0}^{2n-1}h'_3(s_j^*,s_k^*,x_0(s_k^*))z_n(s_k^*)
, j=0,...,2n-1. 
\]

We will show that the parameter $h$ can be chosen so way, that the
Frechet derivative $ K'_n(x_0) $ is invertible. Under this condition
the approximate solution of the equation (7.12) we will seek by
Newton-Kantorovich modified method
$$ \tilde  x_{m+1}=\tilde  x_m-[K_n'(x_0)]^{-1}K_n(\tilde x_m),
\hskip 5pt m=0,1,2,...
\eqno (7.15) $$

Here $ K_n'(x_0) $ is the Frechet derivative of the operator $ K_n $
in the space $X_n $ on a initial element $\tilde x_0= P_n(x_0) $.

Let us transform the iterative method (7.15) to the following form
$$ K'_n(x_0)\tilde x_{m+1}=K'_n(x_0)\tilde x_m-K_n\tilde x_m.
\eqno (7.16) $$

On the each step of the iterative method (7.16) we must decide the
following system of linear algebraic equations
$$
Lx=g_n,  \eqno (7.17)
$$
where $ g_n=K'_n(x_0)\tilde x_m-K_n\tilde x_m $, $ L=K'_n(x_0). $

Let us prove that the system (7.17) has a unique solution on each
step of the iterative process (7.15). For this aim we must prove that the
operator $ K_n'(x_0) $ is invertible.
Let be $a_j=a(s_j^*), $ $ b_{jk}=b'_3(s_j^*,s_k^*,x_0(s_k^*)), $  $ h_{jk}=
h'_3(s_j^*,s_k^*,x_0(s_k^*))$, $ j,k=0, ...,2n-1 $.
Assume that derivatives $b'_3(s_j^*,s_j^*, x_0(s_k^*)) $, $ j,k=0,...,2n-1, $
are not equal to zero on elements $ x_0(s) $ in neighborhood of the solution
$ x^*(s) $ of the equation (7.12).
Diagonal elements of the matrix $L$ of the system (7.17) have the form
\[ |l_{jj}|=\left|a_j+\frac{b_{jj}}{\pi} \ln\left|\frac{\sin\frac{s_{j+1}-
s_j^*}{2}}{\sin\frac{s_j-s_j^*}{2}}\right|+\frac{\pi h_{jj}}{n}\right|=
\]
\[
=
\left|a_j+\frac{b_{jj}}{\pi} \ln\left|\frac{\sin\left(\frac{\pi}
{2n}-\frac{h}{2}\right)}{\sin\frac{h}{2}}\right|+\frac{\pi h_{jj}}{n}
\right|, 
\]
$ j=0,1,...,2n-1. $

From conditions $ b_{jj}\neq 0 $, $ j=0,...,2n-1, $ follow that for sufficient
small $ h $  $ |l_{jj}| > A+B \ln n.$ Let us estimate the sum
$$ \sum_{k=0,k\neq j}^{2n-1}|l_{jk}|\le
\frac{1}{2\pi}
\sum_{k=0,k\neq j,j-1,j+1}^{2n-1}|b_{jk}|\left|\int\limits_{s_k}^{s_{k+1}}
\cot\frac{\sigma-s_j^*}
{2}d\sigma\right|+ $$
$$ +\frac{\pi}{n}\sum_{k=0,k\neq j}^{2n-1}h_{jk} \le
A+B \ln n. $$

Using the Hadamard theorem about invertibility of matrix, we see that the system (7.17) has a
unique solution.

Let us prove the convergence of the iterative process (7.15)
to the exact solution of the equation (7.12).

The proof of the convergence we will give in the space  $ R_{2n} $
of vectors $ v=(v_1,...,v_{2n}) $ with the norm $ \|v\|=\displaystyle
\max_{1\le i\le 2n} |v_i|. $

One can write the matrix of the operator  $ K'_n(x_0)$ in the form
\[ K'_n(x_0)=D+E, 
\]
with elements
\[ d_{jk}=
\left\{
\begin{array} {ccc}

0, \hskip 5 pt j \neq k; \\
a_j+\frac{\displaystyle b_{jj}}{\displaystyle 2\pi}
\displaystyle \int\limits_{s_j}^{s_{j+1}}\cot\frac{\displaystyle
\sigma-s_j^*}{\displaystyle 2}d\sigma+\frac{\pi h_{jj}}{n},
 \hskip 5pt j=k; \\
\end{array}
\right.
\]
and
$$ e_{jk}=
 \left\{
\begin{array} {ccc}

0, \hskip 5 pt j=k; \\
\frac{\displaystyle\pi h_{jk}}{\displaystyle n}, \hskip 5pt j=k-1,k+1;\\
\frac{\displaystyle b_{jk}}{\displaystyle 2\pi}
\displaystyle \int\limits_{s_k}^{s_{k+1}}\cot\frac{\displaystyle
\sigma-s_j^*}{\displaystyle 2}d\sigma+\frac{\pi h_{jk}}{n},  \\
{\rm in \hskip 5pt other \hskip 5pt cases}.\\
\end{array}
\right.
\]

Easy to see that
\[
||D^{-1}||\le \left(\left|C\ln\frac{\displaystyle\sin\left(\frac{\pi}
{2n}-\frac{h}{2}\right)}{\displaystyle\sin\frac{h}{2}}\right|+G
\right)^{-1}, \   ||E||\le A+B\ln n, \]

\[ ||D^{-1}||||E||\le (A+B\ln n)\left|G+C
\ln\frac{\displaystyle\sin\left(\frac{\pi}
{2n}-\frac{h}{2}\right)}{\displaystyle\sin\frac{h}{2}}\right|^{-1}
\le q<1, 
\]
where $C$ and $G$ are constants.

From the Banach theorem follows that
\[ ||[K'_n(x_0)]^{-1}||\le\left(\left|G+C
\ln\frac{\displaystyle\sin\left(\frac{\pi}
{2n}-\frac{h}{2}\right)}{\displaystyle\sin\frac{h}{2}}\right|-A-B\ln n\right)
^{-1}. 
\]

Assume that in some ball $ S(x_0,r) $ the inequality
$$ \|[K'_n(x_0)]^{-1}\|\|K'_n(u_1)-K'_n(u_2)\|\le q<1,
u_1,u_2 \in S(x_0,r) \eqno (7.18) 
$$
is valid.

 Using statements  for convergence of Newton - Kantorovich method , which are given in the section 2 of the Introduction,
we prove that the equation (7.13)
has a unique solution   $ x_n^*(s) $ in ball $ S(x_0,r)$ and iterative
process (7.15) converges  to this solution.

Combining proof of convergence of Newton-Kantorovich method, given   in the section 2 of the Introduction,  
and statements of previous items, we see that
$\|x^*-x^*_n\|_C\le An^{-\alpha}\ln n,$
$||x^*-x^*_n||_X\le An^{-\alpha+\beta}\ln n.$

{\bf  Theorem 7.2} \cite{Boy25}.
Let the equation (7.12) has a unique solution $ x^*(s)
\in H_\alpha. $ Let $ a(s), $ $ f(s) \in H_\alpha, $
$ b'_3(s,\sigma,u), h'_3(s,\sigma,u)\in H_{\alpha,\alpha,\alpha}, $
$ 0<\alpha\le 1.$ Then exist values of the parameter $ h $, for which
conditions (7.18) is valid and the system of equations (7.17) has a unique solution
on each step of the iterative process (7.15).
Under these values of the parameter $h$ and under some additional
conditions estimates
$\|x^*(s)-x_n^*(s)\|_C\le A n^{-\alpha}\ln n,
\ \   ||x^*(s)-x_n^*(s)||_X\le A n^{-\alpha+\beta}\ln n$
is valid. Here $x_n^*(s)$ is a solution of the system (7.13).

\vskip 25 pt
\begin{center}

{\bf 8. Singular Integral Equations in Exceptional Cases}

\end{center}

\vskip 25 pt

     Numerical methods for solution of singular integral equations
in exceptional cases are investigated since 70 years of the past century.
As a rule, they was investigated in cases, when conditions of normality of singular
operators was broken at separate points. The results in this direction was
printed in the books \cite{Mich}, \cite{Pr}. \cite{Pr3}.

     There \ are \  many \  physical and technical problems, where conditions of
normality are violated on whole segments or on all domains of
definition of singular operators. M.M. Lavrentyev formulated \cite{MLavr} some
problems for  investigation of the existence and stability of
solution of  singular integral equations  of the kinds
$$
x(t)+ \frac{1}{\pi i} \int\limits_{-\infty}^ \infty \frac{x(\tau)d\tau}
{\tau-t}=f(t), \eqno (8.1)
$$
$$
\frac{1}{\pi i} \int\limits_{-\infty}^ \infty \frac{x(\tau_1,t_2)d\tau_1}
{\tau_1-t_1} + \frac{1}{\pi i} \int\limits_{-\infty}^ \infty \frac{x(t_1,
\tau_2)d\tau_2}{\tau_2-t_2}=f(t_1,t_2), \eqno (8.2)
$$
$$
\int\limits_{-\infty}^ \infty \int
\limits_{-\infty}^ \infty x(\tau_1,\tau_2) (\frac{1} {\tau_1-t_1}-
\frac{1}{\tau_2-t_2})d\tau_1 d\tau_2=f(t_1,t_2). \eqno (8.3)
$$

     These problems was decided in \cite{Boy12}.

     We need in numerical methods for solution of the singular integral
equations of the types (8.1) - (8.3). In this section we describe some of
these methods.

\begin{center}
{\bf 8.1. Equations on Closed Circles}
\end{center}

Let us consider singular integral equation
$$
a(t)x(t)+\frac{b(t)}{\pi}\int\limits_{\gamma}\frac{x(\tau)d\tau}{\tau-t}+
\int\limits_{\gamma}h(t,\tau)x(\tau)d\tau=f(t), \eqno (8.4)
$$
where $\gamma=\{\gamma:|z|=1\},$ $b(t) \ne 0,$ $a,b,f \in H_{\alpha},$
$h(t,\tau) \in H_{\alpha \alpha},$ $0< \alpha \le 1.$

Assume that the function $a^2(t)-b^2(t)$ can be equal to zero on
arbitrary sets, which are belong to $\gamma,$ or $a^2(t)-b^2(t) \equiv 0$
on whole circle $\gamma.$

Using Hilbert transform we  pass to the equation
\[
a(e^{is})x(e^{is})-\frac{i b(e^{is})}{2\pi } \int\limits_0^{2\pi}
x(e^{i\sigma}){\mbox {ctg}} \frac {\sigma-s}{2}d\sigma+
\]
$$
+i\int\limits_0^{2\pi}h(e^{is},e^{i\sigma})x(e^{i\sigma})e^{i\sigma}d\sigma+
\frac{b(e^{is})}{2\pi}\int\limits_0^{2\pi}x(e^{i\sigma})d\sigma=f(e^{is}),
\eqno (8.5)
$$
where $ 0\le s < 2\pi$.

For simplicity, instead of the equation (8.5), we will consider the equation
$$
a(s)x(s)+\frac{b(s)}{2\pi}\int\limits_{0}^{2\pi}
x(\sigma){\mbox {ctg}}\frac{\sigma-s}{2}
d\sigma+\frac1{2\pi}\int\limits_0^{2\pi } h(s,\sigma)x(\sigma)d\sigma=
$$
$$
=f(s).
\eqno (8.6)
$$

Let us introduce knots
$$ s_k=\frac{\pi k}{n}, \hskip 5 pt s_k^{*}=\frac{\pi k}{n}+h,
\hskip 5pt 0< h \leq \frac{\pi}{2n}, \hskip 5pt k=0,\dots,2n, $$
where parameter $h$ will be defined later.

Solution of the equation (8.6) we will find in the form
\[ x_n(s)=\sum_{k=0}^{2n-1}a_k\psi_k(s), \]
where
\[
\psi_k(s)=\frac{1}{2n+1}(\sin \frac{2n+1}{2}(s-s^*_k))/(\sin \frac{s-s^*_k}{2}).
\]

Unknown coefficients $a_k$ we will find from the following system of linear
algebraic equations
\[
a(s_j^*)x(s_j^*)+\frac{b(s_j^*)}{2\pi }
\sum_{k=0,k\neq j-1,j+1}^{2n-1}x(s_k^*)
\int\limits_{s_k}^{s_{k+1}}{\mbox {ctg}}\frac{\sigma-s_j^{*}}{2}d\sigma+
\]
$$
+\frac{\pi}{n}\sum_{k=0}^{2n-1}h(s_j^*,s_k^*)
x(s_k^*)=f(s_j^*),
\hskip 5 pt j=0,\dots, 2n-1.
\eqno (8.7)
$$

 Under some conditions
 for $h$,
we will prove invertibility of the system (8.7), using the well-known Hadamard Theorem.

The system (8.7) can be rewritten as the algebraic equations system
\[
CX = F,
\]
where $C = \{c_{kl}\},$ $k,l=1,2,\ldots,2n+1,$ $X=(x_1,\ldots,x_{2n+1}),$
$F=(f_1,\ldots,f_{2n+1}).$

Let be $ a_j=a(s_j^*) $, $ b_j=b(s_j^*) $, $ h_{jk}=
h(s_j^*,s_k^*) $.
It is easy to see that
\[ |c_{jj}|=
\left|a_j+\frac{b_j}{\pi}\ln\left|\frac{\sin\left(
\frac{\pi}{2n}-\frac{h}{2}\right)}
{\sin\frac{h}{2}}\right|+\frac{\pi h_{jj}}{n}\right| 
\]
and
\[ \sum_{k=0,k\neq j}^{2n-1}|c_{jk}|\le\frac{\pi}{n} \sum_{k=0,
k\neq j}^{2n-1}|h_{jk}|+\frac{|b_j|}{2\pi}
\sum_{k=0,k\neq j,j-1,j+1}^{2n-1}\left|\int\limits_{s_k}^{s_{k+1}}
{\mbox {ctg}}\frac{\sigma-s_j^*}
{2}d\sigma\right| \le 
\]
\[
 \le A+B \ln n. 
 \]

From two last inequality follows, that parameter $h$ can be chosen so way, that
$|c_{jj}|>\sum\limits^{2n-1}_{k=0, k \ne j}|c_{jk}|.$ So from Hadamard
Theorem follows, that the system (8.7) has a unique solution.

For each knot $ s_j^*,\hskip 5 pt
j=0,1,\ldots,2n-1, $ exists a function $ \psi(s) $ for which
\[
 \psi(s_j^*)=x(s_j^*) 
 \]
and
$$ \int\limits_{s_{j-1}}^{s_j}\psi(\sigma)
{\mbox {ctg}}\frac{\sigma-s_j^*}{2}d\sigma+\int\limits_{s_{j+1}}^{s_{j+2}}\psi(\sigma)
{\mbox {ctg}}\frac{\sigma-s_j^*}{2}d\sigma=0. \eqno (8.8) 
$$

As a example  we can  to conside  the straight line
$$
\psi(\sigma)=x(s_j^*)+k_j(\sigma-s_j^*),
\eqno (8.9)
$$
$k_j=-\frac{nx(s_j^*)}{2\pi}
\ln\left(\frac{\sin\frac{h}{2}}{\sin\left(\frac{\pi}{n}-
\frac{h}{2}\right)}\frac{\sin\left(\frac{2\pi}{n}-\frac{h}{2}\right)}{\sin\left(\frac{\pi}{n}+
\frac{h}{2}\right)}\right).$

Assume that the equation (8.6) with right side $f(s)$ has a unique
solution $ x^* \in H_\alpha.$ We proved that the system (8.7), under
corresponding parameter $h$, has a unique solution $ x^*_n(t).$
Let $P_n$ be the projector from space
 $ X=H_\beta $ $ (0<\beta<\alpha) $ onto interpolated polynomials on
knots $ s_k^* $, $ k=0,1\dots,2n-1 $.
We can write the system (8.7) in operator form as $K_n x_n =f_n.$
Then
 \[
  x_n^*-P_n x^* =K_n^{-1}(K_n(x_n^*-P_nx^*))= 
  \]
 \[ =K_n^{-1}(P_nf-K_nP_nx^*)=K_n^{-1}(P_nKx^*-K_nP_nx^*)= 
 \]
\[
= K_n^{-1}(P_n K x^*-P_n K_n x^*)+ K_n^{-1}(P_n K_n x^*-P_nK_nP_nx^*). \]

Let us estimate in the metrix $C(\gamma)$ the difference
$P_n K x^* - P_n K_n x^*.$

One can see that
\[
\left| \sum_{l=0,l\neq j-1,j+1}^{2n-1}
\int\limits_{s_l}^{s_{l+1}}\left[(x^*(\sigma)-
(x^*(s_l^*)\right]{\mbox {ctg}}\frac{\sigma-
s_j^{*}}{2} d\sigma+\right.
\]
\[ 
+
\int\limits_{s_{j-1}}^{s_{j}}(x^*(\sigma)-
\psi(\sigma)){\mbox {ctg}}\frac{\sigma-
s_j^{*}}{2} d\sigma+
\]
\[
\left.
+\int\limits_{s_{j+1}}^{s_{j+2}}(x^*(\sigma)-
\psi(\sigma)){\mbox {ctg}}\frac{\sigma-
s_j^{*}}{2} d\sigma\right| \le
\]
\[
\le\left|\int\limits_{s_j}^{s_{j+1}}(x^*(\sigma)-x^*(s_j^*)){\mbox {ctg}}\frac
{\sigma-s_j^{*}}{2} d\sigma\right| +
\]
\[+ \left|\sum_{k=0,k\neq j,j-1,j+1}^{2n-1}
\int\limits_{s_k}^{s_{k+1}}(x^*(\sigma)-x^*(s_j^*)){\mbox {ctg}}\frac
{\sigma-s_j^{*}}{2} d\sigma \right|+
\]
\[
+\left| \int\limits_{s_{j-1}}^{s_{j}}(x^*(\sigma)-
x^*(s_j^*)){\mbox {ctg}}\frac{\sigma-
s_j^{*}}{2} d\sigma\right|+
\]
\[ +\left|\int\limits_{s_{j-1}}^{s_{j}}(\psi_j(s_j^*)-
\psi_j(\sigma)){\mbox {ctg}}\frac{\sigma-
s_j^{*}}{2} d\sigma\right|+
\]
\[
+\left|\int\limits_{s_{j+1}}^{s_{j+2}}(x^*(\sigma)-
x(s_j^*)){\mbox {ctg}}\frac{\sigma-
s_j^{*}}{2} d\sigma\right|+
\]
\[
+\left|\int\limits_{s_{j+1}}^{s_{j+2}}(\psi(\sigma)-
\psi(s_j^*)){\mbox {ctg}}\frac{\sigma-
s_j^{*}}{2} d\sigma\right|=I_1 \div I_6.
\]

One can see that:
$I_1 \le A n^{-\alpha};$ $I_2 \le A n^{-\alpha}\ln n;$
$I_3+I_4 \le A n^{-\alpha};$ $I_5+I_6 \le A k n^{-1},$
where $k=\max\limits_{1 \le j \le 2n}|k_j|.$

If coefficients $k_j$ in the equations (8.9) satisfy the inequality
$k=\max\limits_{1 \le j \le 2n}|k_j| \le n^{1-\alpha},$ then
$R_n \le An^{-\alpha} \ln n.$

In the metric of space $ R_{2n+1} $ we have
$$ \|x_n^*-P_nx^*\|\le A\frac{\ln^3 n}{n^\alpha} $$

So,
$$ \|x^*-x_n^*\|_C\le \|x^*-P_nx^*\|_C+\|P_nx^*-x_n^*\|_C\le A\frac{\ln^3 n}
{n^\alpha}. $$

{\bf  Theorem 8.1} \cite{Boy25}, \cite{Boy28}. Let the equation (8.6) has a unique solution
$ x^*(t)\in H_\alpha, $ $ 0<\alpha\le 1. $
Then, for sufficient small value $h$, the system (8.7) has a unique solution
$ x^*_n(t)$ and for $h$ so, that the coefficient $k$ of straight line (8.9) is
smaller than $n^{1-\alpha}$, the estimate
$ \|x^*-x_n^*\|_C\le A\frac{\displaystyle \ln^3 n}{\displaystyle n^\alpha} $
is valid.

\newpage

\begin{center}
{\bf 8.2. Equations on Segments}
\end{center}

Let us consider a singular integral equation
$$
a(t)x(t)+\frac{b(t)}{\pi } \int\limits_{-1}^1\frac{x(\tau)}{\tau-t}d\tau
+\int\limits_{-1}^1 h(t,\tau)d\tau=f(t),
\hskip 5 pt -1 < t < 1, \eqno (8.10)
$$
where $ f(t),$ $ a(t),$ $  b(t) \in H_{\alpha},$ $h(t,\tau) \in H_{\alpha,
\alpha},$ $0<\alpha \leq 1.$

We assume that $b(t)\neq 0 $,
but $ a^2(t)-b^2(t) $ can be equal to zero on sets with measure bigger
than zero.

We will use knots
$ t_k=-1+\frac{k}{n}, \hskip 5pt  k=0,\dots,2n, \hskip 5 pt
t_j^{*}=t_j+h, \hskip 5 pt j=1,2,\dots,n-1, \hskip 5 pt t^*_j=t_{j+1}-h,
$
$\hskip 5 pt j=n,n+1,\dots,2n-2, \hskip 5 pt 0<h< \frac{1}{2n}.
$

The parameter $ h $ will be chosen later.

The integral
$$
I_2x=\int\limits_{-1}^1\frac{x(\tau)}{\tau-t}d\tau \quad \eqno (8.11)
$$
we approximate by the following quadrature rule.

This quadrature rule for $t \in [t_j,t_{j+1}),j=1,2,\dots,2n-2,$ is
$$
I_2x=R_N+
\left\{
\begin{array} {ccc}
\int\limits_{t_0}^{t_2}
\frac{x(t_1^{*})}{\tau-t}d\tau+
\sum\limits_{k=2,k\neq j-1,j+1}^{2n-3}\int\limits_{t_k}^{t_{k+1}}
\frac{x(t_k^{*})}{\tau-t}d\tau+\\
+\int\limits_{t_{2n-2}}^{1}
\frac{x(t_{2n-1}^{*})}{\tau-t}d\tau,
\hskip 5pt j \neq 1,2n-2; \\
\int\limits_{t_1}^{t_2}
\frac{x(t_1^{*})}{\tau-t}d\tau+
\sum\limits_{k=3}^{2n-3}\int\limits_{t_k}^{t_{k+1}}
\frac{x(t_k^{*})}{\tau-t}d\tau+\\
+\int\limits_{t_{2n-2}}^{1}
\frac{x(t_{2n-1}^{*})}{\tau-t}d\tau,
\hskip 5pt j = 1; \\
\int\limits_{t_0}^{t_2}
\frac{x(t_1^{*})}{\tau-t}d\tau+
\sum\limits_{k=2}^{2n-4}\int\limits_{t_k}^{t_{k+1}}
\frac{x(t_k^{*})}{\tau-t}d\tau+\\
+\int\limits_{t_{2n-2}}^{t_{2n-1}}
\frac{x(t_{2n-1}^{*})}{\tau-t}d\tau,
\hskip 5pt j = 2n-2; \\
\end{array}
\right. \eqno (8.12)
$$

Using this quadrature rule we receive the system of algebraic
equations for approximate solution of the singular integral
equation (8.10):
\[
a(t_j^{*})x(t_j^{*})+\frac{b(t_j^*)}{\pi }
\left (\int\limits_{t_0}^{t_2}\frac{x(t_1^{*})}{\tau-t^*_j}d\tau+
\sum\limits_{k=2,k\neq j-1,j+1}^{2n-2}\int\limits_{t_k}^{t_{k+1}}
\frac{x(t_k^{*})}{\tau-t^*_j}d\tau\right.+
\]
\[
\left.+\int\limits_{t_{2n-1}}^{1}
\frac{x(t_{2n-1}^{*})}{\tau-t^*_j}d\tau\right)+\frac{\pi}{n}
\sum_{k=0}^{2n-1}h(t_j^*,t_k^*)x(t_k^*)=
\]
\[
=f(t_j^{*}),\hskip 5pt j \neq 1,2n-2;
\]
\[
a(t_j^{*})x(t_j^{*})+\frac{b(t_j^*)}{\pi }
\left (
\int\limits_{t_1}^{t_2}
\frac{x(t_1^{*})}{\tau-t^*_j}d\tau+
\sum\limits_{k=3}^{2n-3}\int\limits_{t_k}^{t_{k+1}}
\frac{x(t_k^{*})}{\tau-t^*_j}d\tau+\right.
\]
\[
\left.+\int\limits_{t_{2n-2}}^{1}
\frac{x(t_{2n-1}^{*})}{\tau-t^*_j}d\tau\right)
+\frac{\pi}{n}
\sum_{k=0}^{2n-1}h(t_j^*,t_k^*)x(t_k^*)=
\]
\[
=f(t_j^{*}),
\hskip 5pt j = 1;
\]
\[
a(t_j^{*})x(t_j^{*})+\frac{b(t_j^*)}{\pi }
\left (\int\limits_{t_0}^{t_2}\frac{x(t_1^{*})}{\tau-t^*_j}d\tau+
\sum\limits_{k=2}^{2n-4}\int\limits_{t_k}^{t_{k+1}}
\frac{x(t_k^{*})}{\tau-t_j}d\tau+\right.
\]
\[
\left.+\int\limits_{t_{2n-2}}^{t_{2n-1}}
\frac{x(t_{2n-1}^{*})}{\tau-t^*_j}d\tau
\right)+\frac{\pi}{n}
\sum_{k=0}^{2n-1}h(t_j^*,t_k^*)x(t_k^*)=
\]
$$
=f(t_j^{*}),
\hskip 5pt j = 2n-2.   \eqno (8.13)
$$

Using Hadamard Theorem, we  prove, that the system (8.13) has a unique solution.

We use following designations:
$ a_j=a(t_j^*) $, $ b_j=b(t_j^*) $, $ h_{jk}=h(t_j^*,t_k^*) $.
The system (8.13) can be written as $Cx=b,$ where $C$
is the matrix with elements $\{c_{ij}\},$ $i,j=1,2,\ldots,2n-2,$ $x=(x_1,\ldots,x_{2n-2}),$
$f=(f_1,\ldots,f_{2n-2}).$ 

Easy to see that
\[ |c_{jj}|=\left|a_j+\frac{b_j}{\pi}\int\limits_{t_j}^{t_{j+1}}\frac{d\tau}
{\tau-t_j^*}+\frac{h_{jj}}{2n}\right| \ge \left|a_j+\frac{b_j}{\pi}\ln\frac
{t_{j+1}-t_j^*}{t_j^*-t_j}\right|-\left|\frac{h_{jj}}{n}\right|= 
\]
\[
 = \left |a_j+\frac{b_j}{\pi}\ln\frac{\frac{1}{n}-h}{h}
\right |-\left |\frac{h_{jj}}{n}\right |. 
\]

We assumed that $ b_j \neq 0, j=1\dots, 2n-2. $
So we can select the parameter $h$ such way, that $|c_{jj}|$ will become bigger than any arbitrary integer $M.$

From the second hand, for $ j\neq 1,2n-2, $ we have
\[ \sum_{k=1,k\neq j}^{2n-1}|c_{jk}|\le A\left (\left |\int\limits_0^{t_2}
\frac{d\tau}{\tau-t_j^*}\right|+\sum_{k=2,k\neq j,j-1,j+1}^{2n-2}
\left|\int\limits_{t_k}^{t_{k+1}}\frac{d\tau}{\tau-t_j^*}\right|+\right.
\]
\[
\left.+\left|\int\limits_{t_{2n-2}}^1 \frac{d\tau}{\tau-t_j^*}\right|\right)+
\frac{\pi}{n}
\sum_{k=0,k\neq j}^{2n-1}|h_{jk}|\le
\]
\[
\le A\left( \ln\frac{t_j^*-t_0}{t_j^*-t_{j-1}}+\ln
\frac{t_{2n}-t_j^*}{t_{j+2}-t_j^*}\right)+D\le E\ln n+D. 
\]

For other cases we have similar estimations.
From Hadamard Theorem follows that system (8.13) has a unique solution.

Let us estimate the error of offered method.

For each point $ t_j^*,\hskip 5 pt j=1,\dots,2n-2, $
we find a function $ \psi(t) $ with properties
$$
\psi(t_j^*)=x(t_j^*),\hskip 5 pt
  \int\limits_{t_{j-1}}^{t_j}\psi(\tau)
\frac{1}{\tau-t_j^*}d\tau+\int\limits_{t_{j+1}}^{t_{j+2}}\psi(\tau)
\frac{1}{\tau-t_j^*}d\tau=0. \eqno (8.14)
$$

As example of such function we can take the straight line
$ \psi_j(t)=x(t_j^*)+k(t-t_j^*), $
where
$$
k=-\frac{x(t_j)n}{2} \ln \frac {h(2/n-h)}{(1/n+h)(1/n-h)}. \eqno (8.15)
$$

For case when $j \neq1,2n-2,$ we have the following estimate of error of the
quadrature formula (8.12)
\[ \left| R_n\right| \le \left|\int\limits_{-1}^{t_2}
\frac{x(\tau)-x(t_1^*)}{\tau-t_j^*}d\tau+
\int\limits_{t_{2n-2}}^{1}
\frac{x(\tau)-x(t_{2n-1}^*)}{\tau-t_j^*}d\tau+\right.
\]
\[
+\sum_{k=2,k\neq j-1,j+1}^{2n-2}
\int\limits_{t_k}^{t_{k+1}}
(x(\tau)-x(t_k^{*}))\frac{d\tau}{\tau-t_j^{*}}
+ 
\]
\[
 + \left.\int\limits_{t_{j-1}}^{t_j}\frac{(x(\tau)-\psi(\tau))d\tau}
{\tau-t_j^{*}} +\int\limits_{t_{j+1}}^{t_{j+2}}\frac{(x(\tau)-\psi(\tau))
d\tau}
{\tau-t_j^{*}} \right|\le 
\]
\[
\le \left|\int\limits_{-1}^{t_2}
\frac{x(\tau)-x(t_1^*)}{\tau-t_j^*}d\tau\right|+
\left|\int\limits_{t_{2n-2}}^{1}
\frac{x(\tau)-x(t_{2n-1}^*)}{\tau-t_j^*}d\tau\right|+
\]
\[
+\left|\int\limits_{t_j}^{t_{j+1}}(x(\tau)-x(t_j^{*}))\frac{d\tau}
{\tau-t_j^{*}}\right| +
\]
\[
+\sum\limits_{k=3,k\neq j,j-1,j+1}^{2n-1}
\left|\int\limits_{t_k}^{t_{k+1}}(x(\tau)-x(t_k^{*}))\frac{d\tau}
{\tau-t_j^{*}} \right|+ 
\]
\[
+\left|\int\limits^{t_j}_{t_{j-1}} \frac{(x(\tau)-\psi(\tau))d\tau}{\tau-t^*_j}
\right|+\left|\int\limits^{t_{j+2}}_{t_{j+1}} \frac{(x(\tau)-\psi(\tau))d\tau}
{\tau-t^*_j}\right|=r_1+ \cdots r_6.
\]

Estimating each sum, we see that
$$
R_N \le An^{-\alpha} +Akn^{-1} + \frac A{n^\alpha} \ln \left(
\frac{(t_{2n}-t_j^*)\cdot (t^*_j-t_0))}{n^2}\right).
$$

Similar estimations we have in other cases.
Repeated arguments of previous item, we formulate the following
statement.

{\bf Theorem 8.2}  \cite{Boy25}, \cite{Boy28}. Let $a(t), b(t), f(t) \in H_{\alpha},$
$h(t,\tau) \in H_{\alpha,\alpha},$ $0 < \alpha \le 1;$ $b(t) \ne 0$ for $t \in [-1,1].$ Let the equation
(8.10) has a unique solution $x^*(t) \in H_{\alpha}(0<\alpha \le 1).$
Then it  exists a small parameter $h$ such  that the system (8.13) has
a unique solution $x_n^*(t).$ If coefficients $k_j$ in the equations
$\psi_j(t)=x(t^*_j)+k_j(t-t^*_j)$ satisfy the inequality
$\max |k_j| \le n^{1-\alpha},$
then the estimate $\|x^*-x^*_n\|_{H_{\beta}} \le An^{-(\alpha-\beta)} \ln^3 n$ is valid.

\begin{center}
{\bf 8.3. Approximate Solution of the Equation (8.1)}
\end{center}

For approximate solution of the equation (8.1) we use two sets of knots:
$t_k=-A+kA/N,$ $k=0,1,\ldots,2N,$ $t^*_k=t_{k}+h,$ $k=1,\ldots,N-1,$
$t^*_k=t_{k+1}-h,$ $k=N,\ldots,2N-2,$
$0<h\le A/2N.$ We assume that $A$ is a sufficient large constant and $h$
is a small constant.

Repeating arguments of previous item we can construct the numerical scheme
for solution of the equation (8.1) and formulate conditions of solvability.

\begin{center}
{\bf 8.4. Approximate Solution of the Equation (8.2)}
\end{center}

     Let us consider the approximate method for solution of the
equation (8.2). We will use  two sets of knots:
$ t_k=-A+kA/N, k=0,1, \dots, 2N $ and $ t_k^*=t_k+h,$ $k=1, \dots,
N-1, $ $t^*_k=t_k-h,$ $k=N+1,\ldots,2N-2,$
where $ A $ and $ h $ are constants. We assume that $ A $ is a
sufficient large constant and $ h $ is a small constant, $ 0<h<A/N. $

We assume that the right side of the equation (8.2) has the form of function
\[ f(t_1,t_2)=\frac{g(t_1,t_2)}{(1+t_1^2+t_2^2)^\lambda} , 
\]
$ \lambda>3/2 , $ where  $ g(t_1,t_2) $   is a function which is belongs to the Holder class of
functions. 

An approximate solution of  equation (8.2) we will find in the form of
\[ x_N(t_1,t_2)=\sum_{k=1}^{2N-2}\sum_{l=1}^{2N-2}a_{kl}\bar\psi_{kl}(t_1,t_2), 
\]
where
\[
\bar\psi_{kl}(t_1,t_2)=
\]
\[
=\left\{
\begin{array} {ccc}
\psi_{kl}(t_1,t_2), \, 2 \le k,l \le 2N-2,\\
\psi_{00}(t_1,t_2) \cup \psi_{01}(t_1,t_2) \cup \psi_{10}(t_1,t_2) \cup \psi_{11}(t_1,t_2), \\
k=l=1;\\
\psi_{2n-1,0}(t_1,t_2) \cup \psi_{2n-2,0}(t_1,t_2)  \cup \psi_{2n-1,1}(t_1,t_2) \cup \psi_{2n-2,1}(t_1,t_2), \\
k=2n-2, l=1,\\
\psi_{2n-1,2n-1}(t_1,t_2) \cup \psi_{2n-1,2n-2}(t_1,t_2) \cup \psi_{2n-2,2n-1}(t_1,t_2) \cup\\
\cup \psi_{2n-1,2n-2}(t_1,t_2), \, k=l=2n-2,\\
\psi_{0,2n-1}(t_1,t_2) \cup \psi_{1,2n-1}(t_1,t_2) \cup \psi_{0,2n-2}(t_1,t_2) \cup \psi_{1,2n-2}(t_1,t_2), \\
k=1, l=2n-2,\\
\psi_{0,l}(t_1,t_2) \cup \psi_{1,l}(t_1,t_2), \, k=1, 1 \le l \le 2n-2,\\
\psi_{2n-1,l}(t_1,t_2) \cup \psi_{2n-2,l}(t_1,t_2), \, k=2n-2, 1 \le l \le 2n-2;\\
\psi_{k,0}(t_1,t_2) \cup \psi_{k,1}(t_1,t_2), \, 1 \le k \le 2n-2, l=1,\\
\psi_{k,2n-2}(t_1,t_2) \cup \psi_{k,2n-1}(t_1,t_2), \, 1 \le k \le 2n-2, l=2n-2;\\
\end{array}
\right.
\]
\[
\psi_{kl}(t_1,t_2)=
\left\{
\begin{array} {ccc}
1,(t_1,t_2) \in \Delta_{kl},\\
0, (t_1,t_2) \notin \Delta_{kl};\\
\end{array}
\right.
\]
$ \Delta_{kl}=[t_k,t_{k+1})\times [t_l,t_{l+1}). $

     Unknown coefficients $ \{a_{k,l}\},$ $k,l=1, \dots, 2N-2, $ are defined from
the system
$$
\frac{A}{\pi N i}\sum_{k=1}^{2N-2}{}'\frac{x_N(t_k^*,t_{l_1}^*)}{t_k-t_{k_1}^*}+
\frac{A}{\pi N i}\sum_{l=1}^{2N-2}{}''\frac{x_N(t_{k_1}^*,t_l^*)}{t_l-t_{l_1}^*}=
f(t_{k_1}^*,t_{l_1}^*),
\eqno (8.16)
$$
$ k_1,l_1=1,\dots,2N-2, $
where one prime in the first sum indicates that $k \ne k_1-1, k_1+1,$
two primes in the second sum indicate that $l \ne l_1-1, l_1+1.$

     If the constant $ h $ is a sufficient small number, the system (8.16)
has a unique solution $ x^*. $

Error of numerical scheme (8.16) is estimated as in previous items.

\begin{center}
{\bf 8.5. Approximate Solution of the Equation (8.3)}
\end{center}

     Let us consider the approximate method for solution of the
equation (8.3). We will use  two sets of knots:
$ t_k=-A+kA/N, k=0,1, \dots, 2N $ and $ t_k^*=t_k+h_1, k=1, \dots, N-1,$
$t_k^{**}=t_k-h_2, $ $k=N+1,N+2,\ldots,2N-2,$ where $ A $ and $ h_1, $ $ h_2 $ are constants. We assume that $ A $ is a
sufficient large constant and $ h_i $ are small constants, $ 0<h_i \leq A/N, $ $i=1,2.$

The approximate solution of  equation (8.3) we will find in the form of
the function
\[
 x_N(t_1,t_2)=\sum_{k=1}^{2N-2}\sum_{l=1}^{2N-2}a_{kl}\bar\psi_{kl}(t_1,t_2),
\]
where $ \bar\psi_{kl}(t_1,t_2)$ was defined in the previous item.

     Unknown coefficients $ \{a_{kl}\},$ $k,l=1, \dots, 2N-2, $ are defined from
the system
$$
\frac{A}{N}\sum_{k=1}^{2N-2}\sum_{l=1}^{2N-2}{}'
\left(\frac{1}{t_k-t_{k_1}^*}-\frac{1}{t_l-t_{l_1}^*}\right)x_N(t_k,t_l)=
f(t_{k_1}^*,t_{l_1}^*), \eqno (8.17)
$$
$ k_1,l_1=1,\dots,2N-2, $ where the prime in the sum indicates that
$(k,l) \ne (k_1-1,l_1-1), (k_1-1,l_1), (k_1-1,l_1+1),
(k_1,l_1-1),  (k_1,l_1+1), (k_1+1,l_1-1), (k_1+1,l_1), (k_1+1;l_1+1).$

Constants $h_1$ and $h_2$ can be chosen so way that the system (8.17) has a unique solution.

	\begin{center}
		{\bf 9. Approximate Solution of Singular Integro-Differential Equations}

	\end{center}
	
	{\bf 9.1. Linear Equations}
	
	{\bf Introduction. } In the section 9.1.1 proposed a calculation scheme of the mechanical quadratures method for solution of singular integro - differential equations
	(SIDE)
	
	\[
	Kx \equiv  \sum\limits_{k=0}^m \left[ a_k(t) x^{(k)}(t)+
	\frac{ b_k(t)}{ \pi i}
	\int\limits_L \frac{ x^{(k)}(\tau)} { \tau - t} d \tau\right.
	+
	\]
	$$
+\left. 	\frac{ 1}{ 2 \pi i}
	\int\limits_L
	\frac{ h_k(t,\tau)x^{(k)}(\tau)} { \left| \tau - t \right|^\gamma} d \tau
	\right]
	=f(t)
	\eqno (9.1)
	$$
	under conditions
	$$
	\int\limits_L x(t) t^{-t-1} dt =0, k=0,1,\dots,m-1,
	\eqno (9.2)
	$$
	where $a_k(t), b_k(t),h_k(t,\tau),f(t) \in C_{2 \pi}, L$ is the unit circle centered at the origin, $0 \le \gamma < 1,$ and given the justification of this scheme.
	
	In the section 9.1.2 proposed calculation scheme of mechanical quadrature method for approximate solution of SIDE
	$$
	Fx \equiv \sum\limits_{k=0}^m \left[ a_k(t) x^{(k)}(t)+
	\frac{ 1}{  \pi i}
	\int\limits_L
	\frac{ h_k(t,\tau)x^{(k)}(\tau)} { \tau - t } d \tau
	\right]
	=f(t)
	\eqno (9.3)
	$$
	under conditions (9.2),  where functions $f(t)\in H_{\alpha}, a_k(t)\in H_{\alpha}, k=\overline{0,m},h_k(t,\tau)\in H_{\alpha,\alpha},k=\overline{0,m},$
	$ (0 < \alpha <1)$, and provided its justification.
	
	\begin{center}
		{\bf 9.1.1. Approximate Solution of the Boundary Value Problem (9.1), (9.2)
			}
	\end{center}
	
	$1^o.$ {\bf Computational scheme.} 
	
	An approximate solution of the boundary value problem (9.1), (9.2)
	is sought in the form of a polynomial
	$$
	\tilde x_n(t) = \sum\limits_{k=0}^n \alpha_k t^{k+m}+\sum\limits_{k=-n}^{-1} \alpha_k t^k
	\eqno (9.4)
	$$
	which coefficients $ \{\alpha_k \} $ are determined from the system of equations
	\[
	\tilde K_n \tilde x_n =
	P_n \left[\sum\limits_{k=0}^m \left[ a_k(t) \tilde x_n^{(k)}(t)+
	\frac{ b_k(t)}{ \pi i} \int\limits_L\frac{ \tilde x_n^{(k)}(\tau) d \tau}{ \tau - t}+\right. \right. 
	\]
	$$
	\left. \left.
	+\frac{ 1}{ 2 \pi i} \int\limits_L P^{\tau}_n
	\left[ h_k (t, \tau) \tilde x_n^{(k)}(\tau) d(t,\tau) \right] d \tau \right] \right] =P_n \left[ f(t) \right],
	\eqno (9.5)
	$$
	where  $P_n$  is the projector onto the set of interpolating
	trigonometric polynomials of degree $ n $ built over the knots 
	$t_k = e^{is_k}$, $s_k=2k \pi/(2n+1)$ $ (k=\overline{0,2n})$,
	$d(t, \tau)=\left| t -\tau \right|^{-\gamma}$ for
	$\left| \sigma - s\right| \ge 2 \pi/(2n+1)$,
	$d(t, \tau)= \left| e^{i s_1} -1 \right|^{-\gamma}$ for
	$\left| \sigma - s\right| \le 2 \pi/(2n+1)$,
	$\tau = e^{i \sigma}$, $t = e^{i s}$.
	
	{\bf Justification of the method.} 
	First of all, the boundary value problem
	(9.1), (9.2) and the equation (9.5) are reduced to the 
	equivalent Riemann boundary value problems. For this we introduce the function
	$\Phi (z) = \frac{ 1}{ 2 \pi i} \int\limits_L \frac{x(\tau) d \tau}{\tau - z}$
	
	We will need in  Sohotzky - Plemel formulas
	
	$\Phi^+(t) - \Phi^-(t) = x(t)$,
	$\Phi^+(t) + \Phi^-(t) = \frac{ 1}{ \pi i} \int\limits_L \frac{x(\tau) d \tau}{\tau - t}, \cdots,$
	
	$\Phi^{(m)+}(t) - \Phi^{(m)-}(t) = x^{(m)}(t)$,
	$\Phi^{(m)+}(t) + \Phi^{(m)}(t) = \frac{ 1}{ \pi i} \int\limits_L \frac{x^{(m)}(\tau) d \tau}{\tau - t}$.
	
	Substituting the previous formulas in equations (9.1) and (9.5), we arrive to the following Riemann boundary value problems:
	
	\[
	Kx \equiv  \sum\limits_{k=0}^m \left[ \left( a_k(t)+b_k(t) \right) \Phi^{(k)+}(t)-
	\left( a_k(t)-b_k(t) \right) \Phi^{(k)-}(t)+\right.
	\]
	$$
	+
	\left.
	\frac{ 1}{2  \pi i}
	\int\limits_L \frac{h_k(t,\tau) \left( \Phi^{(k)+}(\tau) - \Phi^{(k)-}(\tau) \right)}{ \left| \tau- t \right|^{\gamma}} d \tau \right]= f(t)
	\eqno (9.6)
	$$
	under conditions (9.2) and 
	\[
	\tilde K_n \tilde x_n \equiv
	\]
	\[ \equiv P_n \left[ \sum\limits_{k=0}^m \left[
	\left( a_k(t)+b_k(t) \right) \tilde \Phi_n^{(k)+}(t)-
	\left( a_k(t)-b_k(t) \right) \tilde \Phi_n^{(k)-}(t)+\right. \right.
	\]
	$$
	+\left.\left. \frac{ 1}{ 2 \pi} \int\limits_L
	P_n^{\tau} \left[h_k(t,\tau) d(t,\tau)
	\left(
	\tilde \Phi_n^{(k)+}(\tau) - \tilde \Phi_n^{(k)-}(\tau)
	\right)
	\right] d\tau 
	\right] \right]
	=
	$$
	$$
	= P_n \left[f(t) \right].
	\eqno (9.7)
	$$
	
	Easy to see that $\tilde \Phi_n(t) = \sum\limits_{k=0}^n \alpha_k t^{k+m}$,
	$\tilde \Phi_n(t) = \sum\limits_{k=-n}^{-1} \alpha_k t^{k}$.
	
	We reduce the problem (9.6), (9.2)
	and (9.7) to the equivalent singular integral equations. For this 
	we  use the well known integral representation of  Krikunov \cite{Krik}.
	Functions $ \frac {d ^ m \Phi ^ + (z)} {dz ^ m} $ and
	$ \frac {d ^ m \Phi ^ - (z)} {dz ^ m} $
	are represented  into the integral of Cauchy type with the same density
	$$
	\frac{ d^m \Phi^+(z)}{ dz^m}=
	\frac{ 1}{ 2 \pi i}
	\int\limits_L
	\frac{ v(\tau) d \tau}{\tau -z},
	\frac{ d^m \Phi^-(z)}{ dz^m}=
	\frac{ z^{-m}}{ 2 \pi i}
	\int\limits_L
	\frac{ v(\tau) d \tau}{\tau -z}.
	\eqno (9.8)
	$$
	
	To transform the approximate system (9.7) for the  Riemann boundary value problem (9.6) into 
	the equivalent approximate system for singular integral equation,  we  will use the integral representation
	$$
	\frac{ d^m \tilde \Phi_n^+(z)}{ dz^m}=
	\frac{ 1}{ 2 \pi i}
	\int\limits_L
	\frac{ \tilde v_n(\tau) d \tau}{\tau -z},
	\frac{ d^m \tilde \Phi_n^-(z)}{ dz^m}=
	\frac{ z^{-m}}{ 2 \pi i}
	\int\limits_L
	\frac{ \tilde v_n(\tau) d \tau}{\tau -z}.
	\eqno (9.9)
	$$
	where $ \tilde v_n(t)= \sum\limits_{k=0}^n \frac{ (m+k)!}{k!} \alpha_k t^k
	-\sum\limits_{k=1}^n(-1)^m \alpha_{-k}
	\frac{ (m+k-1)!}{(k-1)!} t^{-k}$.
	
	Repeating arguments given in the work \cite{Krik}, and using integral representations
	(9.8) and (9.9), we reduce the boundary value problems (9.6),
	(9.2) and (9.7) to equivalent singular integral equations
	\[
	K_1 v \equiv a_1^* v(t)+\frac{ b_1^*(t)}{ \pi i}
	\int\limits_L \frac{ v(\tau) d\tau}{ \tau - t}+
	\]
	\[
	+	\frac{ 1}{ (2 \pi i)^2}
	\int\limits_L \frac{ h_0(t,\tau)}{ |\tau - t|^{\gamma}} 
	(\int\limits_L k_0(\tau,\sigma) \ln(\tau - \sigma) v(\sigma) d \sigma)d\tau +
	\]
	\[
		+\dots+ \frac{ 1}{ 2 \pi i}
	\int\limits_L \frac{  h_m(t,\tau)}{|\tau - t|^{\gamma}} 
	\left(\frac{1 }{ 2 \pi i} \int\limits_L
	\frac{ v(\sigma)+ \tau^{-m} v(\sigma)}{\sigma - \tau} d \sigma\right)d\tau =
	\]
	$$
	= f_1(t)
	\eqno (9.10)
	$$
	and
	\[
	\tilde K_{1,n} \tilde v_n = P_n \left[
	a_1^* \tilde v(t)+\frac{ b_1^*(t)}{ \pi i}
	\int\limits_L \frac{ \tilde v_n(\tau) d\tau}{ \tau - t}+\right.
	\]
	\[
	\frac{ 1}{ (2 \pi i)^2}
	\int\limits_L P_n^{\tau} \left[
	h_0(t,\tau) d(t,\tau)
	\left[\int\limits_L k_0(\tau,\sigma) \ln(\tau - \sigma) \tilde v(\sigma) d \sigma\right]\right] d \tau
	+ 
	\]
	\[
	+\dots+ \left.\frac{ 1}{ (2 \pi i)^2}
	\int\limits_L P_{\tau}\left[ h_m(t,\tau)
	d( t,\tau) \left[
	\int\limits_L
	\frac{ \tilde v(\sigma)+ \tau^{-m} \tilde  v(\sigma)}{\sigma - \tau} d \sigma \right]d \tau \right] \right] =
	\]
	$$
	= P[f_1(t)],
	\eqno (9.11)
	$$
	where $ k_0 (t, \tau) $ $, \dots, $ $ k_{m-1} (t, \tau) $ are Fredholm kernels;
	$a_1^*(t)$,$b_1^*(t)$, $f_1^*(t) \in C_{2 \pi}$; explicit form of these
	functions can be discharged on the basis of presentation
	Yu. M. Krikunov \cite{Krik}. \footnote{The explicit form of these functions is not issued,  since below
		we	use only their characteristics.}
	
	Let $ X $ is the  space of
	square-integrable functions with scalar product
	$$
	(g_1,g_2)=\frac{ 1}{ 2 \pi} \int\limits_0^{2 \pi}
	g_1(t) \overline{g_2 (t)}  ds, t = e^{i s}.
	$$
	
	Approximate solution  $\tilde v_n(t)$ we will seek in the subspace 
	$X_n$ of the space $X.$ Subspace $ X_n$ consists  of $n$-order polynomials $ x_n=\left\{ \sum\limits_{k=-n}^n \alpha_k t^k \right\}.$ Obviously,
	$ \tilde K_{1,n} \in \left[ X_n \to X_n\right]$.

	We assume that the operator $K_1  $ has the linear inverse.
	(So, the boundary value problem (9.1), (9.2) has
	a unique solution for any $ f $).
	
	As
	$
	\left\| \int\limits_L
	k_0(\tau,\sigma) \ln (\tau - \sigma) v(\sigma) d \sigma
	\right\| \le A \| v \|,
	\left\| \int\limits_L
	\frac{ v(\sigma)} { (\sigma- \tau )}  d \sigma
	\right\| \le A \| v \|,
	$
	where (as everywhere else) $A$ are  well-defined constants independent
	from $ v $ and $ n $, instead of equations
	(9.10), (9.11)
	we can restrict ourself by the equations
	$$
	K_1 v \equiv a(t) v(t)  +
	\frac{ b(t)}{ \pi i}
	\int\limits_L \frac{v(\tau) d \tau}{ \tau - t}  +
	\frac{ 1}{ 2 \pi i}
	\int\limits_L \frac{ h(t, \tau) v(\tau)}{ \left| \tau -t\right|^{\gamma}} d \tau =f(t),
	\eqno (9.12)
	$$
	\[
	K_{1,n} \tilde v_n \equiv P_n \left[ a(t) \tilde v_n(t)  +
	\frac{ b(t)}{ \pi i}
	\int\limits_L \frac{\tilde v_n(\tau) d \tau}{ \tau - t}
	\right.
	+
	\]
	$$
	+\left.
	\frac{ 1}{ 2 \pi i}
	\int\limits_L P_n^{\tau} \left[
	h(t, \tau) d (t,\tau) \tilde v_n(\tau) \right] d \tau \right] =P_n[f(t)].
	\eqno (9.13)
	$$
	
	Consider the equation
	$$
	K_2  v\equiv \tilde a(t) v(t) + \tilde b(t) S(v) + H(v) =f(t),
	\eqno (9.14)
	$$
	where
	$
	S(v)=
	\frac{ 1}{ \pi i}
	\int\limits_L \frac{v(\tau) d \tau}{ \tau - t}
	$,
	$
	H(v) =
	\frac{ 1}{ 2 \pi i}
	\int\limits_L \frac{ h(t, \tau) v(\tau)}{ \left| \tau -t\right|^{\gamma}} d \tau;
	$
	$\tilde a(t)$, $\tilde b(t)$ are polynomials of best uniform
	approximation degree at most $ m $ for functions $ a (t) $, $ b (t), t\in \gamma, $
	respectively. The value of $ m $ will be fixed below.
	
	From Banach theorem it implies that, for such $ n $  that
	$q_1 = A \left[ \omega(a;m^{-1})+\right.$ $\left.\omega(b;m^{-1})\right] <1$,
	operator $ K_2 $ has a linear inverse operator with the norm
	$\| K_2^{-1}  \| \le \| K_1^{-1} \| /(1 - q_1)$.
	
	Collocation method for the equation (9.14)  is written as follows
	$$
	K_{2,n} \tilde v_n \equiv P_n[\tilde a(t) \tilde v_n(t) + \tilde b(t) S(\tilde v_n) + H(\tilde v_n)] =P_n[f].
	\eqno (9.15)
	$$
	
	Following the method proposed in the paragraph 3, equations (9.14) and (9.15)
	can be written in the equivalent form:
	$$
	K_3 v \equiv \Psi v +Wv=y, K_3 \in [X \to X]
	\eqno (9.16)
	$$
	and
	$$
	K_{3,n} \tilde v_n \equiv \tilde \Psi_n \tilde v_n +\tilde W_n \tilde v_n= \tilde y_n, \tilde K_{3,n} \in [\tilde X_n \to \tilde X_n],
	\eqno (9.17)
	$$
	where $\Psi v = \psi^- v^+ - \psi^+ v^-$,
	$Wv=lH(v) $, $l =\psi^-/(\tilde a  +\tilde b) $,
	\[\psi(z) = \exp \left\{
	\frac{ 1}{ 2 \pi i} 
	\int\limits_L
	\left[
	\frac{\ln\left[
		\left( \tilde a(\tau) - \tilde b(\tau) \right)/
		\left( \tilde a(\tau) + \tilde b(\tau) \right)
		\right]}{\tau - z}
	\right]
	d \tau
	\right\}
	,\]
	$ \tilde \Psi_n \tilde v_n = P_n[\Psi \tilde v_n]$,
	$ \tilde W_n \tilde v_n =P_n[lH(\tilde v_n)]$,
	$y = lf$, $ \tilde y_n = P_n[y]$.
	
	 Let us
	introduce the polynomial
	\[
	\tilde \varphi(t) =
	\Psi_n \tilde v_n +
	\left(
	T^{[n/2]} [l]
	\right)
	\frac{ 1}{ 2 \pi i}
	\int\limits_L
	T_t^{[n/2]}
	\left[
	h(t, \tau) d^*(t, \tau)
	\right]
	\tilde v_n(\tau) d \tau,
	\]
	where
	$\Psi_n \tilde v_n = \psi_n^- \tilde v_n^+ - \psi_n^+ \tilde v_n^-$, $[n/2]$ is
	antje $n/2$; $T^{[n/2]}[f]$ and
	$\psi_n$  are polynomials of the best uniform approximation
	degree to $ [n / 2] $ and $ n $ for functions
	$ f $ and $ \psi $, respectively; $d^*(t,\tau) = |\tau - t|^{-\gamma}$ for
	$|\tau - t|>\rho$, $d^*(t,\tau) = \rho^{-\gamma}$ for $|\tau - t| \le \rho$;
	$\rho$ is a fixed positive number, $\rho \ge |e^{is_1} - 1|$.
	
	Let us estimate
	$$
	\left\|
	K_3 \tilde v_n - \tilde \varphi
	\right\|
	\le\left\|
	I_1
	\right\|
	+
	\left\|
	I_2
	\right\|
	+
	\left\|
	I_3
	\right\|
	+
	\left\|
	I_4
	\right\|,
	\eqno (9.18)
	$$
	where\\
	$I_1  =\Psi_n \tilde v_n - \Psi \tilde v_n$, \ \ 
	$I_2 = D^{[n/2]}[l]H(\tilde v_n)$,\\
	$I_3 = \left( T^{[n/2]}[l]\right) \frac{ 1}{ \pi i} \int\limits_L h(t,\tau) \tilde v_n(\tau)
	\left[ |\tau - t|^{-\gamma} - d^*(t, \tau) \right] d \tau$,\\
	$I_4 = \left( T^{[n/2]}[l]\right) \frac{ 1}{ \pi i} \int\limits_L
	D_t^{[n/2]}[h(t,\tau) d^*(t,\tau)] \tilde v_n(\tau) d \tau$,\\
	$D=E-T,E$ is the identity operator.
	
	Since $ \tilde a, \tilde b $ are included in the Holder class $ H_1 $
	with the  factor $A m$, than 
	$\psi \in H_{\delta}, \delta = 1- \epsilon$,  where $\epsilon$ is a arbitrary number, 
	$0 < \epsilon < 1$. It follows from Privalov theorem  \cite{Gakh}.
	
	So
	$$
	\|I_1 \| \le A m \| \tilde v_n \| / n^{1-\epsilon}
	\eqno (9.19)
	$$

	Easy to see that
	\[
	\|I_2 \| \le A \frac{m}{n^{1-\epsilon}} \| \tilde v_n \|,
	\]
	\[
	\|I_3 \| \le A  \| \tilde v_n \| \rho^{1-\gamma},
	\]
	$$
	\|I_4 \| \le A  \| \tilde v_n \| \frac{\omega (h;n^{-1})}{\rho^{2 \gamma}}.
	\eqno (9.20)
	$$
	
	Let us estimate 
	$$
	\| P_n K_3 \tilde v_n - \tilde \varphi \| \le
	\| I_5 \|+
	\| I_6 \|+
	\| I_7 \|+
	\| I_8 \|,
	\eqno (9.21)
	$$
	where $I_5=P_n I_1$, $I_6=P_n I_2$, $I_7=P_n I_3$, $I_8=P_n I_4$. 
	
	Obviously,
	\[
	\| I_5 \| =
	\left[ \frac{ 2 \pi}{2n+1}
	\sum\limits_{k=0}^{2n}
	| \psi_n(t_k) - \psi(t_k)|^2
	| \tilde v_n (t_k) |^2
	\right]^{\frac{ 1}{ 2}}
	\le
	\]
	$$
	\le A m \| \tilde v_n\| / n^{1-\epsilon}.
	\eqno (9.22)
	$$
	
	We proceed to the estimation of $ \| I_7 \| $. For this we represent $ I_7 $
	as follows:
	$$
	I_7=I_9+I_{10},
	\eqno (9.23)
	$$
	where
	\[
	I_9=P_n \left[ \frac{ 1}{ \pi i}
	\int\limits_L \frac{ c(t, \tau) \tilde v_n(\tau)}{ \tau - t} d \tau \right]=
	P_n \left[ \frac{ 1}{ \pi i}
	\int\limits_L \frac{ P_n^t[c(t, \tau)] \tilde v_n(\tau)}{ \tau - t} d \tau \right]=
	\]
	\[
	=  \frac{ 1}{ \pi i}
	\int\limits_L \frac{ c(t, \tau) - c(\tau,\tau) }{ \tau - t} \tilde v_n(\tau) d \tau+
	\frac{ 1}{ \pi i}
	\int\limits_L \frac{ c(\tau,\tau)[\tilde v_n(\tau) - \tilde v_n(t)] }{ \tau - t} d \tau+
	\]
	\[
	+P_n \left[
	\frac{ \tilde v_n(t)}{ \pi i}
	\int\limits_L
	\frac{ c(\tau,\tau) d \tau}{ \tau - t}
	\right]=
	I_{11}+I_{12}+I_{13},
	\eqno (9.24)
	$$
	\[
	I_{10} = P_n \left[
	\frac{ 1}{ \pi i}
	\left( T^{n/2} [l] \right)
	\int\limits_L h(t,\tau) \tilde v_n(\tau)
	[d(t,\tau) - d^*(t,\tau) ]
	d \tau
	\right],
	\]
	$c(t,\tau) = \  \left( T^{n/2} [l] \right)
	h(t,\tau)(1-|\tau - t|^{\gamma} d(t,\tau))
	|\tau - t|^{1-\gamma} e^{i \theta}$,	
	$\theta = \theta(\tau,t) = {\rm arg}(\tau - t)$,
	$c_1(t,\tau) = P^t_n[c(t,\tau)]$.
	
	From the definition of function $ d (t, \tau) $ follows that $ | c_1 (t, \tau) | \le A n ^ {\gamma-1} $.
	Since $ c (t, \tau) $ with respect to the variable $ t $ is a trigonometric polynomial of degree
	not higher than the $n,$ it is possible to show that $ c_1 (t, \tau) $ belongs to the class
	Holder with degree $ 1 /\ln n $ and with the factor
	$ A \left (n^{\gamma - 1}n^{1 /\ln n} \right) = A n^{\gamma - 1}. $

	Therefore
	$$
	\| I_{11} \| \le A n^{\gamma - 1} \ln n \| \tilde v_n \|.
	\eqno (9.25)
	$$
	
	For  $I_{12}$, $ I_{13}$ we have:
	$$
	\| I_{12} \| + \| I_{13}\| \le A  n^{\gamma - 1} \ln n \| \tilde v_n \|.
	\eqno (9.26)
	$$

	Now we will estimate $ \| I_9 \| $. Let us  introduce the notations
	$$
	b(t,\tau) = d(t, \tau) - d^*(t,\tau), \ \ 
	g(t,\tau) = \left(T^{[n/2]}[l] \right) h(t,\tau).
	$$
	
	Obviously
	\[
	\| I_9 \| =
	\]
	\[
	= \left[ \frac{ 2 \pi }{ 2 n +1} \sum\limits_{k=0}^{2n}
	\left[ \frac{ 1}{ \pi i }
	\int\limits_L g(t_k, \tau) \tilde v_n (\tau)  [b(t_k,\tau)]^{\frac{1}{2}}
	P_n^t \left[ [b(t, \tau)]^{\frac{1}{2}} \right]
	d \tau \right]^2
	\right]^{\frac{1}{2}}\le
	\]
	\[
	\le\left[
	\frac{ 2 \pi}{ 2n+1} \sum\limits_{k=0}^{2n}
	\frac{1}{\pi i}
	\left[ \int\limits_L \left[|g(t_k, \tau) |
	| [b(t_k,\tau)]^{\frac{1}{2}} | \right]^2 d \tau \right]\times\right.
	\]
	\[
	\times \left.
	\left[ \int\limits_L |P_n^t \left[[b(t, \tau)]^{\frac{1}{2}} \tilde v_n(\tau)\right]|^2
	d \tau \right] \right]^{\frac{1}{2}} \le 
	\]
	\[
	\le
	\left[ \max_t \int\limits_L\left[
	|g(t,\tau)| |[b(t,\tau)]^{\frac{1}{2}}| \right]^2 d \tau \right]^{\frac{1}{2}} \times
	\]
	\[
	\times
	\left[
	\int\limits_0^{2\pi}
	ds \frac{ 1}{ \pi}
	\int\limits_0^{2 \pi}
	\left| \left[  P_n^s \left[ [b(e^{i s},e^{i \sigma})]^{\frac{1}{2}} \right]
	\tilde v_n(e^{i \sigma})
	\right]^2
	e^{i \sigma}
	\right| d \sigma
	\right]^{\frac{1}{2}} \le
	\]
	$$
	\le A \rho^{1-\gamma} \| \tilde v_n\|.
	\eqno (9.27)
	$$
	
	Repeating the previous discussion, we obtain
	$$
	\| I_6 \| \le A m \ln n \| \tilde v_n \| /n^{1-\epsilon}.
	\eqno (9.28)
	$$
	
	Norm of $ \| I_8 \| $ is estimated in the same way as the $ \| I_4 \| $:
	$$
	\| I_8 \| \le A \| \tilde v\| \omega(h;n^{-1})/ \rho^{2 \gamma}.
	\eqno (9.29)
	$$
	
	From the estimates (9.18) - (9.29) it follows that
	\[
	\|
	(K_3 \tilde v_n - \tilde \varphi) +
	(P_n \tilde \varphi - P_n K_3 \tilde v_n)
	\|
	\le
	\]
	\[
	\le 	A
	[
	\frac{m}{n^{1-\epsilon}}+\frac{\omega(h;n^{-1})}{\rho^{2 \gamma}}
	+\rho^{1-\epsilon} + n^{\gamma-1}
	]
	\| \tilde v_n\|.
	\]
	
	{\it Note.} It is assumed that $ \gamma> 0 $. Otherwise it is possible to put
	$ \rho = 0 $,  $ m / n^{1 - \epsilon} <\min [w (a; n^{- 1}), \omega (b; n^{- 1})]. $

	Assume 
	$\varepsilon = \gamma/2$, $ m = n^{\gamma/2}$,
	$\rho= [\omega(h;n^{-1})^{1/(1+\gamma)}]$.
	Using the results of the paragraph 3,
	we obtain,
 for such  $n$ 	that 
 \[ q_2 =   A \left [\omega (h; n^{- 1}) \right]^{\frac {1 - \gamma} {1+ \gamma}} <  1,
 \]
  that
	the equation (9.17), and hence (9.15), is uniquely solvable for any right-hand side, and the estimate
	$ \| \tilde K_3^{- 1} \| \le A $ is valid.

	Let us estimate
	$\| K_2^{-1}\|$. Let $\tilde v_1^*$ is the solution of the equation  (9.14),
	and hence (9.17).
	
	Then
	\[
	\| \tilde v_1^*\| \le \| \tilde K_{3,n}^{-1} \tilde y \| \le
	\| \tilde K_{3,n}^{-1}\| \|P_n[lP[f]] \| \le \| \tilde K_{3,n}^{-1}\||l| \|P_n[f]\|.
	\]
	
	Hence $\| \tilde K_2^{-1}\| \le A$.
	
	Repeating the above discussion, we see that	
	$$
	\| \tilde K_{2,n} - \tilde K_{1,n} \| \le A
	\left[
	[w(h;n^{-1})]^{\frac{ 1- \gamma}{ 1+ \gamma}} +
	w(a;n^{\gamma - 1}) +\omega(b;n^{\gamma - 1})
	\right].
	$$
	
	Hence, from the Banach theorem we have
	
	{\bf Theorem 9.1} \cite{Boy32}, \cite{Boy34}.
 Let the problem (9.1), (9.2)
		has a unique solution for any right-hand side, and let the functions
		$ a_i $, $ b_i $, $ h_i $, $ f \in C_{2 \pi} $, $ i = \overline {0, m} $.
		Then for $ n $ such that
		$$
		q=A \sum\limits_{k=0}^m
		\left[
		\omega(a_k;n^{\gamma-1}) +
		\omega(b_k;n^{\gamma-1}) +
		[\omega(h_k;n^{-1})]^{\frac{ 1 -\gamma}{ 1+\gamma}}
		\right] < 1,
		$$
		the system of equations (9.5) has a unique solution
		$\tilde x_n^*$ and the estimate 
		$\| x^* - \tilde x_n^* \| \le
		A[q+\omega(f;n^{-1})]$ is valid. Here $x^*$ is the solution of the
		boundary value problem (9.1), (9.2).

	\begin{center}
	{\bf
		9.1.2. Approximation  Solution of the Boundary Value Problem (9.3), (9.2)}
	\end{center} 
	
	{\bf Computational scheme.}
	
	An approximate solution of the boundary value problem
	(9.3), (9.2) is sought in the form of the polynomial (9.4)
	with coefficients $ \{\alpha_k \} $, which are determined from the system of
	algebraic equations
	\[
	\tilde F_n \tilde x_n \equiv
	\]
	\[
	\equiv
	\overline{P_n} \left[ \sum\limits_{k=0}^m \left[
	a_k(t) \tilde x_n^{(k)}(t) + \frac{ 1}{ \pi}
	\int\limits_L \left[\frac{P_n^{\tau}[h_k(t,\tau) \tilde x_n^{(k)}(\tau)]}{\tau - t}
	\right] d \tau \right] \right] =
	\]
	$$
	= \overline{P}[f],
	\eqno (9.30)
	$$
	where the operator $ \overline {P_n} $ is the  projection operator onto the set of
	interpolating polynomials of degree at most $ n $, constructed on nodes  $\overline t_k = e^{i \overline {s_k}}$,
	$\overline {s_k} = (2 k \pi + \pi)/(2n+1), k=0,1,\ldots,2n$.
	
	{\bf Justification of the method.} As in the previous item
	boundary value problem (9.3), (9.2) and the equation (9.30)
	is reduced to an equivalent singular integral equations.
	
	To do this, we will use the identity
	$$
	\begin{array}{c}
	\overline P_n \left[ \int\limits_0^{2 \pi}
	P_n^{\sigma} \left[ h(s,\sigma) \tilde g(\sigma)
	\cot \frac{ \sigma - s}{ 2}
	\right] d \sigma \right] \equiv\\
	\equiv
	\overline P_n \left[ \int\limits_0^{2 \pi}
	P_n^{\sigma} \left[ h(s,\sigma) \tilde g(\sigma)
	\right]
	\cot \frac{ \sigma - s}{ 2}
	d \sigma \right] \equiv\\
	\equiv
	\overline P_n \left[ \int\limits_0^{2 \pi}
	P_n^{\sigma} \left[ h(s,\sigma)\right] \tilde g(\sigma)
	\cot \frac{ \sigma - s}{ 2}
	d \sigma \right], (\tilde g \in \tilde X_n).
	\end{array}
	\eqno (9.31)
	$$
	
	The validity of this identity follows from the fact, that if
	$ \tilde g (s) $ is a polynomial of degree $ n $, then
	$ [\tilde g (\sigma) - \tilde g(s)] \cot \frac {\sigma - s} {2} $
	is a polynomial of degree $ n $, and
	$ \sum \limits_{k = 0}^{2n} \cot \frac {2 k \pi + \pi} {4 n +2} = 0 $.
	
	Using the identity (9.31), the equations (9.3)
	(9.30) can be represented as follows:
	\[
	F x \equiv \sum\limits_{k=0}^m
	\left[a_k(t) x^{(k)} (t)  + \frac{ b_k(t)}{\pi i}
	\int\limits_L \frac{ x^{(k)}(\tau)}{ \tau - t} d \tau +\right.
	\]
	$$
	+\left.
	\frac{ 1}{ 2 \pi i } \int\limits_L
	\frac{ h_k(t,\tau ) - h_k(t,t)}{ \tau -t } x^{(k)}(\tau)
	d \tau \right] = f(t)
	\eqno (9.32)
	$$
	and 
	\[
	\tilde F_n \tilde x_n \equiv \overline P_n \left[ \sum\limits_{k=0}^m
	\left[ a_k(t) \tilde x_n^{(k)} (t)  + \frac{ b_k(t)}{\pi i}
	\int\limits_L \frac{ \tilde x_n^{(k)}(\tau)}{ \tau - t} d \tau +\right. \right.
	\]
	$$
	+\left. \left.
	\frac{ 1}{ 2 \pi i } \int\limits_L
	\frac{ P_n^{\tau} \left[ \left [ h_k(t,\tau ) - h_k(t,t) \right]\tilde x_n^{(k)}(\tau) \right]}{ \tau -t }
	d \tau \right] \right] = \overline P_n[f],
	\eqno (9.33)
	$$
	where $b_k(t) = h_k(t,t), n \ge m$.
	
	Repeating the arguments of the previous section, the boundary value problem (9.32)
	(9.2) and the equation (9.33) can be reduced to  equivalent singular integral equations.
	
	Therefore, the justification of the computational scheme (9.30) for the boundary value problem
	(9.3), (9.2) can be reduced to the justification of the computing scheme
	\[
	\tilde F_{1,n} \tilde v_n \equiv \overline P_n
	\left[ a(t) \tilde v_n(t)  + \frac{ b(t)}{ \pi i}
	\int\limits_L \frac{ \tilde v_n(\tau) d \tau}{ \tau - t} + \right.
	\]
	$$
	+ \left.
	\frac{ 1}{ 2 \pi i } \int\limits_L
	\frac{ P_n^{\tau} \left[ \left [ h(t,\tau ) - h(t,t) \right] \tilde v_n(\tau) \right]}{ \tau -t }
	d \tau \right] = \overline P_n[f]
	\eqno (9.34)
	$$
	for the  singular integral equation
	$$
	F_1 v \equiv
	a(t) v(t)  + \frac{ b(t)}{ \pi i}
	\int\limits_L \frac{ v(\tau) d \tau}{ \tau - t}
	+
	\frac{ 1}{ 2 \pi i } \int\limits_L
	\frac{ h(t,\tau ) - h(t,t) }{ \tau -t }
	v(\tau) d \tau =
	\]
	 $$
	= f(t),
		\eqno (9.35)
	$$
	where
	$$
	a(t),b(t),f(t)\in H_{\alpha} ,h(t,\tau) \in H_{\alpha,\alpha},(0<\alpha<1).
	\eqno (9.36)
	$$
	
	Justification held in the spaces $ X $ and $ \tilde X_n $, introduced in the previous section 
	$ (v \in X, \tilde v_n \in \tilde X_n) $.
	
	As in the previous section, equations (9.35) and (9.34)
	can be written in the following equivalent forms:
	$$
	F_2 v \equiv \Psi v +W y = y, F_2 \in [X \to X]
	\eqno (9.37)
	$$
	and
	$$
	\tilde F_{2,n}  \tilde v_n \equiv \tilde \Psi_n \tilde v_n +\tilde W_n \tilde v_n = \tilde y_n,
	\tilde F_{2,n} \in [\tilde X_n \to \tilde X_n],
	\eqno (9.38)
	$$
	where 
	\[
	W v = \frac{ l}{ 2 \pi i}
	\int\limits_L \frac{ h(t, \tau) - h(t,t)}{ \tau - t}
	v(\tau ) d \tau
	, 
	\]
	\[
		\tilde W_n \tilde v_n = \overline{P_n}
	\left[
	\frac { l}{ 2 \pi i}
	\int\limits_L
	\frac{ P_n^{\tau}
		\left[ \left[
		h(t, \tau) - h(t,t)\right] \tilde v_n(\tau) \right]}{ \tau - t}
	d \tau
	\right],
	\]
	$\tilde {\Psi} = \overline{P_n}[\Psi]$, $\tilde y_n = \overline{P_n} [y].$
	The functions  $\Psi,l,y$ are introduced in the section 9.1.1.
	
	First of all, we justify the method of collocation
	for the equation (9.37). The collocation method for this
	equation can be written as
	$$
	\tilde F_{3,n}  \tilde v_n\equiv \overline{P_n}
	\left[
	\Psi \tilde v_n + W \tilde v_n \right]  = \tilde y_n,
	\tilde F_{3,n} \in [\tilde X_n \to \tilde X_n].
	\eqno (9.38)
	$$
	
	To justify the collocation method we introduce the polynomial 
	\[
	\tilde \varphi_n(t) = \Psi_n \tilde v_n + \frac{ 1}{ 2 \pi i}
	\int\limits_L \frac{[
	\tilde h(t,\tau) - \tilde h(t,t)]
	\tilde v_n(\tau)}{\tau-t} d \tau,
	\] 
	where
	\[
	\tilde h(t,\tau) =  T^t_{[n/2]} T^{\tau}_{[n/2]}[h(t,\tau)l(t)]].
	\]
	Symbols 
	$
	\Psi_n
	$,
	$ T$, $[n]$ was introduced in the section 9.1.1.
	
	It is easy to see that  
	\[|l(t) h(t,\tau) - \tilde h(t,\tau)| \le A \ln n [E_n^t[h(t,\tau)l(t)] +E_n^{\tau}[h(t,\tau)l(t)]],
	\]
	where  $E_n[f]$ is the the best approximation of $ f $ by trigonometric
	polynomials of degree at most $ n.$ 
	
	Hence it follows that
	$ \{l (t) h (t, \tau) - \tilde h (t, \tau) \} $ belongs to the class
	Holder with exponent $ 1 /  \ln n $ and
	coefficient $ A \ln n [E_n^t (h (t, \tau) l (t)) + E_n^{\tau} (h (t, \tau) l (t))]. $
	
	Then  $\|F_2 \tilde v - \varphi \| \le  A	\ln^2 n [E_n(\psi)+E_n^t(hl)+E_n^{\tau}(hl)]$. 
	
	From the previous
	arguments follow that a similar estimate is valid for
	$\| \overline P_n F_2 \tilde v_n - \tilde \varphi_n \|$.
	Therefore  the operator $ \tilde F_{3,n} $ has a linear inverse with the norm
	$\| \tilde F_{3,n}^{-1}\| \le \|F_{2,n}^{-1} \| / (1-q_1)$, when $ n $ such that
	$q_1 = A_4 \ln^2 n [E_n (\psi) + E_n^t (hl) + E_n^{\tau} (hl)] < 1.$
	Using the equality (9.31) we obtain the estimate
	$\|\tilde F_{3,n} - \tilde F_{2,n} \| \le A_5 \ln^2 n [E_n^t(hl)+E_n^{\tau}(hl)]$.
	From this estimate and Banach theorem implies the following theorem.
	
	{\bf Theorem 9.2} \cite{Boy32}, \cite{Boy34}.
 Let the boundary value problem (9.3), (9.2) has 
		a unique solution for any right-hand side.
		Then for such $ n $  that
		$
		q = A n^{-\alpha} \ln^2 n <1,
		$
		the system of equations (9.30) has a unique solution $ \tilde x_n^* $
		and the estimate
		$
		\| x^* - \tilde x_n^* \| \le A[q+\ln^2nE_n(f)]
		$
		is valid, where $x^*$ is a unique solution of the  boundary value problem (9.3), (9.2).

	\begin{center}
	{\bf 9.3. Approximate Solution of Nonlinear
		Singular Integro-differential 
		Equations on Closed
		Contours of Integration	}
	\end{center}
	
	In the items 9.1.1 and 9.1.2 we investigated the numerical methods for solution of linear singular integro-differential equations on closed
	contours of integration.
	
	Now we examine the application of these methods to nonlinear
	singular integro-differential equations of the form
	\[
	Kx\equiv a(t,x(t),\ldots,x^{(m)}(t))+
	S_{\gamma}(h(t,\tau,x(\tau),\ldots,x^{(m)}(\tau)))=
	\]
	$$
	=f(t)
	\eqno (9.39)
	$$
	under conditions
	$$
	\int\limits_{\gamma}x(t)t^{-k-1}dt=0,
	\;k=0,1,\ldots,m-1.
	\eqno (9.40)
	$$
	
	Here $ \gamma $ is the unit circle
	centered at the origin, 
	\[
	S_{\gamma}(x(\tau))=\frac{1}{\pi i}
	\int\limits_{\gamma}\frac{x(\tau)}{\tau-t}d\tau,
	\]
	\[
	a_{u_i}'(t,u_0,u_1,\ldots,u_m)\in H_{\alpha,\ldots,\alpha},
	\]
	$$
	h_{u_i}'(t,\tau,u_0,u_1,\ldots,u_m)\in H_{\alpha,\ldots,\alpha},\;f(t)\in H_{\alpha}.
	\eqno (9.41)
	$$
	
	{\bf Computational scheme.} An approximate solution of the boundary value problem
	(9.39) (9.40) is sought in the form of the polynomial
	$$
	\tilde x_n(t)=\sum_{k=0}^{n}\alpha_k t^{k+m}+
	\sum_{k=-n}^{-1}\alpha_k t^{k},
	\eqno (9.42)
	$$
	the coefficients of which are determined from the system of nonlinear algebraic equations
	\[
	\tilde{K_n}\tilde{x_n}\equiv
	\]
	$$
	 \equiv \bar{P_n^t}[a(t,\tilde{x}_n(t),
	\ldots,\tilde{x}^{(m)}_n(t))+
	S_{\gamma}(P_n^{\tau}[h(t,\tau,
	\tilde{x_n}(\tau),
	\ldots,\tilde{x}^{(m)}_n(\tau))])]=
	$$
	$$
	=\bar{P}_n^t[f(t)],
	\eqno (9.43)
	$$
	where $P_n(\bar{P_n})$ is  the  projector of interpolation onto the set of  $ n $ order trigonometric polynomials
	on knots
	 $$s_k=2k\pi/(2n+1)\;(\bar{s}_k=(2k\pi+\pi)/(2n+1)),
	\;k=0,1,\ldots,2n.$$
	
	{\bf Justification of the method.} We transform the boundary value problem (9.39),
	(9.40) and approximating its equation (9.43) to the
	equivalent  nonlinear equations. 
	
	Let us introduce the function
	$\Phi(z)=\frac{1}{2\pi
		i}\int\limits_{\gamma}\frac{x(\tau)}{\tau-z}d\tau. $ Using the   Sohotzky - Plemel  formulas, we receive the 
	equations
	\[ Kx\equiv
	a(t,\Phi^{+}(t)-\Phi^{-}(t),\ldots,
	\Phi^{(m)+}(t)-\Phi^{(m)-}(t))+
	\]
	$$
	+S_{\gamma}(h(t,\tau,\Phi^{+}(\tau)-\Phi^{-}(\tau),
	\ldots,\Phi^{(m)+}(\tau)-\Phi^{(m)-}(\tau)))=
	$$
	$$
	=f(t),
	\eqno (9.44)
	$$
	\[
	\tilde{K_n}\tilde{x_n}\equiv
	\bar P_n^t[a(t,\tilde{\Phi}_{n}^{+}(t)-\tilde{\Phi}_{n}^{-}(t),
	\ldots,\tilde{\Phi}_{n}^{(m)+}(t)-\tilde{\Phi}_{n}^{(m)-}(t))+
	\]
	$$
	+S_{\gamma}(P_n^{\tau}[h(t,\tau,\tilde{\Phi}_{n}^{+}(\tau)-
	\tilde{\Phi}_{n}^{-}(\tau),\ldots,
	\tilde{\Phi}_{n}^{(m)+}(\tau)-\tilde{\Phi}_{n}^{(m)-}(\tau))])]=
	$$
	$$
	=\bar{P}_n^t[f(t)],
	\eqno (9.45)
	$$
	where $\tilde{\Phi}_{n}(z)=\frac{1}{2\pi
		i}\int\limits_{\gamma}\frac{\tilde{x}_{n}(\tau)}{\tau-z}d\tau.$
	
	Using the integral representation of Y. M.  Krikunov, the
	equation (9.44) under the conditions (9.40) and the equation (9.45)   can be transformed to 
	nonlinear singular integral equations
	\[
	K_1x\equiv a(t,\eta_0(v(t)),\ldots,\eta_m(v(t)))+
	\]
	$$
	+S_{\gamma}(h(t,\tau,\eta_0(v(\tau)),
	\ldots,\eta_m(v(\tau))))=f_1(t),
	\eqno (9.46)
	$$
	\[
	\tilde{K}_{1,n}\tilde{x}_{n}\equiv
	\bar{P}_n^t[a(t,\eta_0(\tilde{v}_{n}(t)),
	\ldots,\eta_m(\tilde{v}_{n}(t)))]+
	\]
	$$
	+S_{\gamma}(P_n^{\tau}[h(t,\tau,\eta_0(\tilde{v}_{n}(\tau)),
	\ldots,\eta_m(\tilde{v}_{n}(\tau)))])]=\bar{P}_n [f_{1,n}(t)],
	\eqno (9.47)
	$$
	where 
	\[
	\eta_0(v(t))=\frac{1}{2\pi
		i}\int\limits_{\gamma}k_0(t,\tau)v(\tau)d\tau,
	\]
	\[\ldots
	\]
	\[
	\eta_{m-1}(v(t))=\frac{1}{2\pi
		i}\int\limits_{\gamma}k_{m-1}(t,\tau)v(\tau)d\tau,
	\]
	\[
	\eta_m(v(t))=\frac{1}{2}
	v(t)(1+t^{-m})+S_{\gamma}(v(\tau)-t^{-m}S_{\gamma}(v(\tau))),
	\]
	$k_0(t,\tau),\ldots,k_m(t,\tau)$ are Fredholm kernels (explicit
	form of these functions \  discharged \ in \ \cite{Krik}), $ v $ \ is  \ the integral \  density of  Y.M. Krikunov
	presentation,
	$\tilde{v}_{n}(t)=\sum\limits_{k=-n}^{n}\beta_k t^k,$ $\beta_k$ are 
	constants, uniquely expressed in terms of $ \alpha_k. $
	
	Since \  $ \eta_i (v) \  (i = 0,1, \ldots, m) $ are linear operators, then
	for  simplicity of calculations, instead of the equations (9.46) and (9.47), we can 
	restrict ourself with equations
	$$
	K_2(v)\equiv a(t,v(t))+S_{\gamma}(h(t,\tau,v(\tau)))=f(t)
	\eqno (9.48)
	$$
	and
	$$
	\tilde{K}_{2,n}(\tilde{v})\equiv \bar{P}_{n}[a(t,\tilde{v}_{n}(t))
	+S_{\gamma}(P_n^{\tau}[h(t,\tau,\tilde{v}_{n}(\tau))])]=
	$$
	$$
	=\bar{P}_{n}[f(t)],
	\eqno (9.49)
	$$
	where \[
	a_{u}'(t,u)\in
	H_{\alpha,\alpha},\; h_{u}'(t,\tau,u)\in
	H_{\alpha,\ldots,\alpha},\; f(t)\in
	H_{\alpha}\;(0<\alpha<1).
	\]
	
	We can choice  the space, in which we will justificate
	the  collocation method (9.49) for the equation
	(9.48). At first, we will work in the space  $L_2$ with the
	scalar product 
	\[
	(v_1,v_2)=\frac{1}{2\pi}\int\limits_{0}^{2\pi}
	v_1(t)\overline{v_2(t)}dt,\;t=e^{is}.
	\] 
	
	The approximate  solution is sought in the subspace $ \tilde {L}_{2,n} \subset L_2, $
	consisting of  $ n $ order polynomials.	
	In these conditions, method
	collocation (9.49) for the equation (9.48) was justified in the item 3.
	Using the results given there, we receive the following statement.
	
	{\bf Theorem 9.2} \cite{Boy32}, \cite{Boy36}.  Suppose that in a  ball $ S $ boundary problem
		(9.39) (9.40) has a unique solution $ x^* \in W^r H_{\alpha}^{(m)} $ and the equation $ K'(x^*) z = f $ for arbitrary $ f $ has
		a unique  solution that satisfies (9.40). Then,
		for such $ n $  that $ q = A \ln^2n / n^{\alpha} <1$, the system of equations 
		(9.43) has a unique solution $ \tilde {x}_{n}^*, $, and the
		estimation
		\[
		\|(x^*)^{(m)}-(\tilde{x}_{n}^*)^{(m)}\|\leq A\ln n/n^{\alpha}
		\]
		is valid.

	
	Let us now take as the space, in which 
	the collocation method (9.49) for the equation
	(9.48) will be studied, the space $ X=H_\beta$ $(0<\beta<\alpha) $ with the norm
	\[
	\|x\|=\max|x(t)|+\sup\limits_{t_2\ne
		t_1}[|x(t_1)-x(t_2)|/|t_1-t_2|^{\beta}].
	\] 
	
	Approximate solution of the equation (9.48) we  seek
	 in the space $\tilde{X}_{n}\subset X$ consisting of $ n $ order  polynomials. The justification of collocation method for the equation (9.48) is given in the paragraph 2.
	 
	Using the results given there, we receive the following statement.
	
	{\bf Theorem 9.3} \cite{Boy36}  Suppose that in a  ball $ S $ the boundary problem
		(9.39) (9.40) has a unique solution $x^*\in
		H_{\alpha}^{(m)},$ and the equation $ K '(x ^ *) z = f $ with an arbitrary right side
		has a unique solution satisfying the conditions (9.40).
		Then, for $ n $ such that $ q = A \ln^6 n / n^{\alpha} <1, $ the system of 
		equations (9.43) has a unique solution $ \tilde {x}_{n}^* (t), $
		and the estimate
		$$
		\|(x^*)^{(m)}-(\tilde{x}_{n}^*)^{(m)}\|\le A\ln^2 n/n^{\alpha}
		$$
		is valid.
	
	\begin{center}
		{\bf 9.4. Approximate Solution of Linear
			Singular Integro-Differential Equations with Discontinuous
			Coefficients and on Open Contours of Integration}
	\end{center}
	
	In this and the following items are studied directly (without regularization) \ 
	methods \ for \  approximate \ solution of singular integro-differential
	equations with discontinuous coefficients and on open contours
	of integration. For the first time, approximate methods for solving singular
	integral equations with discontinuous coefficients have been studied in
	monograph \cite{Ivan}, in which an efficient algorithm for transformation singular integral equations with \  discontinuous \ 
	coefficients \  to \  equivalent singular integral equations
	with continuous coefficients was proposed.

	Let us consider  the singular integro-differential equation 
	$$
	Kx\equiv \sum_{k=0}^{m}\Bigl[a_k(t)x^{(k)}(t)+
	b_k(t)S_{\gamma}(x^{(k)}(\tau))+
	U_{\gamma}\left(\frac{h_k(t,\tau)}{|\tau-t|^{\eta}}x^{(k)}(\tau)\right)\Bigr]=
	$$
	$$
	=f(t)
	\eqno (9.50)
	$$
	under conditions
	$$
	\int\limits_{\gamma}x(t)t^{-k-1}dt=0,\;\:k=0,1,\ldots,m-1,
	\eqno (9.51)
	$$
	where  $a_k(t)$, $ b_k(t)$, $h_k(t,\tau)\;(k=\overline{0,m}),f(t)$ are
	functions, which are continuous in the metric  $ C $ everywhere on $ \gamma $, 
	except of the point $ t = 1, $ in which $ a_k (t) $ and
	$ b_k (t) \; (k = \overline {0, m}) $ have discontinuity of the first kind.
	
	An approximate solution of the boundary value problem (9.50), (9.51) is
	sought in the form of  polynomial (9.42), whose coefficients
	determined from the system of algebraic equations
	$$
	\tilde{K}_{n}\tilde{x}_{n}\equiv P_n\Bigg[\sum_{k=0}^{m}\Bigl[a_k(t)\tilde{x}_{n}^{(k)}(t)+
	b_k(t)S_{\gamma}(\tilde{x}_{n}^{(k)}(\tau))+
	$$
	$$
	+U_{\gamma}(P_n^{\tau}[h_k(t,\tau)d(t,\tau)\tilde{x}_{n}^{(k)}(\tau)])\Bigr]
	\Bigg]=P_n[f(t)],
	\eqno (9.52)
	$$
	where $d(t,\tau)=|\tau-t|^{-\eta}$ for $|\sigma-s|\ge
	2\pi/(2n+1),$ $d(t,\tau)=|e^{i2\pi/(2n+1)}-1|^{-\eta}$ for
	$|\sigma-s|\le 2\pi/(2n+1),$ $t=e^{is},\;\tau=e^{i\sigma}.$

	Just as in the case of singular integro-differential
	equations with continuous coefficients, boundary value problem
	(9.50), (9.51) and the system (9.52) are reduced to
	equivalent singular integral equation and
	approximating its algebraic system. They are written in
	the item 9.1 (formulas (9.6) and (9.7)). Instead of considering these equations, we confine ourselves to equations 
	$$
	Kx\equiv a(t)x(t)+b(t)S_{\gamma}(x(\tau))+
	U_{\gamma}\left(\frac{h(t,\tau)}{|\tau-t|^{\eta}}x(\tau)\right)=f(t)
	\eqno (9.53)
	$$
	and
	\[
	\tilde{K}_{n}\tilde{x}_{n}\equiv
	\]
	\[
	 \equiv P_n\Bigl[a(t)\tilde{x}_{n}(t)+
	b(t)S_{\gamma}(\tilde{x}_{n}(\tau))+
	U_{\gamma}(P_n^{\tau}[h(t,\tau)d(t,\tau)\tilde{x}_{n}(\tau)])\Bigr]\]
	$$
	=P_n[f(t)],
	\eqno (9.54)
	$$
	where  $a(t)$ and $b(t)$ have discontinuity of the first kind in the point $t=1.$
	
	Justification of the computational scheme (9.54) will be carried out in the space $L_2.$
	
	Equations (9.53) and (9.54) are transformed  to the equivalent Riemann 
	boundary value problems
	\[
	K_1x\equiv x^+(t)+G(t)x^-(t)+
	D(t)U_{\gamma}(h(t,\tau)|\tau-t|^{-\eta}x(\tau))=
	\]
	$$
	=D(t)f(t),
	\eqno (9.55)
	$$
	\[
	\tilde{K}_{1,n}\tilde{x}_{n}\equiv
	\]
	\[
	\equiv P_n\Bigl[\tilde{x}_{n}^+(t)+
	G(t)\tilde{x}_{n}^-(t)+
	D(t)U_{\gamma}(P_n^{\tau}[h(t,\tau)d(t,\tau)\tilde{x}_{n}(\tau)])\Bigr]=
	\]
	$$
	=P_n[D(t)f(t)],
	\eqno (9.56)
	$$
	where $G(t)=S(t)D(t),\;S(t)=a(t)-b(t),\;D(t)=(a(t)+b(t))^{-1}.$
	
	Let a solution of the  Riemann boundary value problem  $\varphi^+(t)=G(t)\varphi^-(t)$ is
	\[
	\varphi(t)=(t-1)^{\delta}\varphi_0(t), \rm{where } \, 
	\delta=\zeta+i\xi,\;\varphi_0\in H_{\alpha}, 0<\alpha<1.
	\]
	
	As in the case of singular integral equations with continuous
	coefficients, equations (9.55) and (9.56) can be
	represented in the following equivalent forms:
	$$
	K_2x\equiv Vx+Wx=y,
	\eqno (9.57)
	$$
	$$
	\tilde{K}_{2,n}\tilde{x}_{n}\equiv
	\tilde{V}_{n}\tilde{x}_{n}+\tilde{W}_{n}\tilde{x}_{n}=\tilde{y}_{n},
	\eqno (9.58)
	$$
	where $V x=\psi^-x^+-\psi^+x^-,\;
	Wx=lU_{\gamma}(h(t,\tau)|\tau-t|^{-\eta}x(\tau)),$
	$l=\psi^-/(a+b),$ $y=lf,$ $\tilde{y}=P[y],$
	$\tilde{V}_{n}\tilde{x}_{n}=P_n[V\tilde{x}_{n}],$
	$\tilde{W}_{n}\tilde{x}_{n}=P_n[lU_{\gamma}(P_n^{\tau}[h(t,\tau)d(t,\tau)\tilde{x}_{n}(\tau)])],$
	\[
	\psi(z)= \left \{ 
	\begin{array}{cc}
	(z-1)\exp\left\{\frac{1}{2\pi i}
	\int\limits_{\gamma}\frac{\ln[S(\tau)D(\tau)]}{\tau-z}d\tau\right\}
	\rm{ for } \  \zeta\le 0, \\
	\exp\left\{\frac{1}{2\pi i}
	\int\limits_{\gamma}\frac{\ln[S(\tau)D(\tau)]}{\tau-z}d\tau\right\}
	\rm { for }\  \zeta>0.
	\end{array}
	\right.
	\]
	
	To justify the proposed computational scheme we introduce a polynomial
	$$\tilde{\varphi}_{n}(t)=V_n\tilde{x}_{n}+U_{\gamma}
	(T_t^{(n)}[l(t)h(t,\tau)d^*(t,\tau)]\tilde{x}_{n}(\tau)),$$ where
	$V_n\tilde{x}_{n}=\psi_n^-\tilde{x}_{n}^+-\psi_n^+\tilde{x}_{n}^-,$
	$\psi_n=T^{(n)}(\psi),$ $T^{(t)}_n$ is projector of interpolation with respect to   variable $ t  $ onto
	the set of  the $ n $ order trigonometric polynomials of the best uniform approximation;
	$d^*(t,\tau)=|\tau-t|^{-\eta}$ for $|\tau-t|\ge\rho,$
	$d^*(t,\tau)=\rho^{-\eta}$ for $|\tau-t|\le\rho,$ a constant $\rho$ is fixed below.
	It can be shown that\footnote{modules
		continuity of functions $a(e^{is})$ and $b(e^{is})$ are determined in
		open interval   $0<s<2\pi.$}
	\[
	\|K_{2,n}\tilde{x}_n-\tilde{\varphi}_{n}\|
	\le   
	A\Bigl[w(a;n^{-1})+w(b;n^{-1})+
	\]
	$$
+	w(h;n^{-1})\rho^{1-\eta}+
	\rho^{-2\eta}n^{-\eta}+n^{\theta}\Bigr]\|\tilde{x}_{n}\|,
	\eqno (9.59)
	$$
	\[
	\|P_n K_{2,n}\tilde{x}_{n} - \tilde{\varphi}_{n}\|\le
	A\Bigl[w(a;n^{-1})+w(b;n^{-1})+
	\]
	$$  
	+w(h;n^{-1})\rho^{1-\eta}+
	\rho^{-2\eta}n^{-\eta}+n^{\theta}\Bigr]\|\tilde{x}_{n}\|,
	\eqno (9.60)
	$$
	where $\theta=-(1-|\zeta|)$ for $\zeta\le 0$ and  $\theta=-\zeta$
	for $\zeta>0.$
	
	Let us estimate $\|P_nK_{2,n}\tilde{x}_{n}-\tilde{K}_{2,n}\tilde{x}_{n}\|.$
	Repeating the arguments, given in the paragraph 3, we obtain the estimate
	$$
	\|P_n K_{2,n}\tilde{x}_{n}-\tilde{K}_{2,n}\tilde{x}_{n}\| \le
	A\Bigl[w(h;n^{-1})+\rho^{-2\eta}n^{-\eta}\Bigr]\|\tilde{x}_{n}\|.
	\eqno (9.61)
	$$
	
	Let $\rho=n^{-\eta/(1+\eta)}.$ Collecting the estimates (9.59) - (9.61), we receive the following statement.
	
	{\bf Теорема 9.4 } \cite{Boy36}.  Let the operator $K$ has the linear inverse operator $K^{-1}$ in the space $L_2;$ functions $a(t),\;b(t),\;h(t,\tau),\;f(t)\in C$ everywhere, except the point $t=1,$ where the functions $a(t), b(t)$ have a first-order gap. 
		Then, for $n$ such that
		$$
		q=A\Bigl[w(a;n^{-1})+w(b;n^{-1})+w(h;n^{-1})+
		n^{\theta}+n^{-\eta(1-\eta)/(1+\eta)}\Bigr]<1,
		$$ 
		the system of equation (9.54) has a unique solution
		$\tilde{x}^*$. The estimate $\|x^*-\tilde{x}^*\|\le
		A[q+w(f;n^{-1})],$ where $x^*$ is a solution of the equation (9.53), is valid. Here $\theta=-(1-|\zeta|)$ for $\zeta\le 0,$
		$\theta=-\zeta$ for $\zeta>0,$ $\delta=\zeta+i\xi,$
		$(t-1)^{\delta}\varphi_0(t)$ is a solution of Riemann boundary task
		$\psi^+(t)=[(a(t)-b(t))/(a(t)+b(t))]\psi^-(t).$ 
	
	A corollary of the Theorem 9.4 is the following statement.
	
	{\bf Theorem $\bf 9.4'$. } Let the value problem (9.50),
		(9.51) has a unique solution for any right-hand side, and
		 conditions $ a_k(t)$, $ b_k(t)$, $ h_k(t,\tau)$
		$(k=\overline{0,m}),$ $f(t)\in C$ are performed everywhere, except the point
		$t=1,$ in which the functions  $a_k(t)$ и $b_k(t)$ have a first-order gap.
		Then for such $n$  that 
		\[
		q=A\max_{0\le k\le m}
		\Bigl[w(a_k;\frac{1}{n})+w(b_k;\frac{1}{n})+w(h_k;\frac{1}{n})+
		n^{\theta}+\frac{1}{n^{\eta(1-\eta)/(1+\eta)}}\Bigr]<1,
		\]
		the system of equations (9.52) has a unique solution
		$ \tilde {x}^*_n $ and the estimate $\|x^*-\tilde{x}^*_n\|\le
		A[q+w(f;n^{-1})]$ is valid. Here $x^*$ is a solution of the boundary value problem
		(9.50), (9.51);  $\theta=-(1-|\zeta|)$ for
		$\zeta\le 0,$ $\theta=-\zeta$ for $\zeta>0;$
		$\delta=\zeta+i\xi,$ $(t-1)^{\delta}\varphi_0(t)$ is a solution of Riemann boundary value problem
		$\Phi^+(t)=[(a_m(t)-b_m(t))/(a_m(t)+b_m(t))]\Phi^-(t).$ 
	
	Consider the singular integro-differential equation
	$$
	Fx\equiv \sum_{k=0}^{m}\Bigl[a_k(t)x^{(k)}(t)+
	S_L(h_k(t,\tau)x^{(k)}(\tau))\Bigr]=f(t), \, t \in [c_1, c_2],
	\eqno (9.62)
	$$
	with conditions
	$$
	x(1)=x'(1)=\cdots=x^{(m-1)}(1)=0,
	\eqno (9.63)
	$$
	where  $L$ is a segment of the counter  $\gamma$, wherein one end of segment
	$ L = (c_1, c_2) $ (say, $ c_1 $) coincides with the point $ t = 1. $
	
	An approximate solution of the boundary value problem (9.62), (9.63) is
	sought in the form of the polynomial
	$$
	\tilde{x}_{n}(t)=(t-1)^{m-1}\sum_{k=-n}^{n}\alpha_kt^k,
	\eqno (9.64)
	$$
	coefficients $ \{\alpha_k \} $ of which are determined from the system of 
	equations
	$$
	\tilde{F}_{n}\tilde{x}_{n}\equiv \sum_{k=0}^{m}\Bigl[\bar{P}_{n}
	[\bar{a}_k(t)\tilde{x}_{n}^{(k)}(t)+
	S_{\gamma}(P_n^{\tau}[\bar{h}_k(t,\tau)\tilde{x}_{n}^{(k)}(\tau)])]\Bigr]
	=\bar{P}_{n}[\bar{f}(t)],
	\eqno (9.65)
	$$
	where $\bar{a}_k(t)=a_k(t)\;(k=\overline{0,m})$ for $t\in L,\;
	\bar{a}_k(t)=0 \;(k=\overline{0,m-1}),\;\bar{a}_m(t)=1$ for
	$t\notin L,$ $\bar{h}_k(t,\tau)=h_k(t,\tau)$ for $t\in L$ and $\tau\in L,$
	$\bar{h}_k(t,\tau)=0$ for $t \notin L$ and  $\tau \in \gamma$ or $t \in \gamma$ and
	$\tau\notin L\; (k=\overline{0,m}),$ 
	$\bar{f}(t)=f(t)$ for $t\in L,$ $\bar{f}(t)=0$ for $t\notin L.$
	
	Equation (9.62) can be written in the following equivalent
	forms:
	\[
	Fx\equiv
	\]
	\[ \equiv\sum_{k=0}^{m}\Bigl[a_k(t)x^{(k)}(t)+
	b_k(t)S_L(x^{(k)}(\tau))+
	S_L((h_k(t,\tau)-b_k(t))x^{(k)}(\tau))\Bigr]=
	\]
	\[
	=   f(t), \, t \in L,
	\]
	\[
	Fx\equiv
	\]
	\[ \equiv\sum_{k=0}^{m}\Bigl[\bar{a}_k(t)x^{(k)}(t)+
	\bar{b}_k(t)S_{\gamma}(x^{(k)}(\tau)+
	U_{\gamma}(g_k(t,\tau)x^{(k)}(\tau))\Bigr]=
	\]
	$$
	+\bar{f}(t), \, t \in \gamma,
	\eqno (9.66)
	$$
	where
	$b_k(t)=h_k(t,t),\;g_k(t,\tau)=\bigl(h_k(t,\tau)-h_k(t,t)\bigr)/(\tau-t)$
	for $t\in L$ and $\tau\in L,$ $g_k(t,\tau)=0$ for $t\in L$ and $\tau\in \gamma $ or  $t\in \gamma$ and $\tau\in L. $ 
	
	Taking advantage of the following identity 
	\[
	\bar{P}_{n}\left[\int\limits_{\gamma}P_n^{\sigma}\left[
	h(s,\sigma)\tilde{\varphi}(\sigma)\cot\frac{\sigma-s}{2}
	\right]d\sigma\right]=
	\]
	\[
	=
	\bar{P}_{n}\left[\int\limits_{\gamma}P_n^{\sigma}[h(s,\sigma)]
	\tilde{\varphi}(\sigma)\cot\frac{\sigma-s}{2}d\sigma\right],
	\]
	true if $ \tilde {\varphi} (\sigma) $ is a polynomial of degree
	not exceeding $ n $, the system (9.66) reduces to
	\[
	\tilde{F}_{n}\tilde{x}_{n}\equiv \bar{P}_{n}\Bigl[\sum_{k=0}^{m}
	\Bigl[\bar{a}_k(t)\tilde{x}_{n}^{(k)}(t)+
	\bar{h}_k(t,t)S_{\gamma}(\tilde{x}_{n}^{(k)}(\tau))+
	\]
	$$
	+U_{\gamma}(P_n^{\tau}[g_k(t,\tau)
	\tilde{x}_{n}^{(k)}(\tau)])\Bigr]
	\Bigr]=\bar{P}_{n}[\bar{f}(t)].
	\eqno (9.67)
	$$
	
	The relationship between the equations (9.66) and (9.67) is studied above and
	formulated in the Theorem 9.4. Using this theorem and
	the results of the item 9.1, we receive the  following statement.
	
	{\bf Theorem 9.5 } \cite{Boy36}.  Let the problem (9.62),
		(9.63) has a unique solution for any right-hand side, and
		the conditions   $h_k(t,\tau)\in H_{\alpha,\alpha},\  k=\overline{0,m},$ for $ t$ and $\tau \in L,;  a_k(t) \in H_{\alpha}, f(t)\in H_{\alpha},$ for $ \  t \in L,$ are implemented.			
		 Then for $ n $ such that
		$ q = A \ln^2n (n^{- \alpha} + n^{- \theta}), $ the system of equations
		(9.65) has a unique solution $\tilde{x}_{n}^*,$ and
		the estimate $\|x^{*(m)}-\tilde{x}^{*(m)}_n\|\le
		Aq$ is valid.  Here $x^*$ is s solution of the boundary value task (9.62)
		(9.63); $\theta=\min\{\theta_1,\theta_2\},$
		$\theta_1=1-|\zeta_1|$ for $\zeta_1\le 0,$
		$\theta_1=\zeta_1$ for $\zeta_1>0,$ $\theta_2=1-|\zeta_2|$ for
		$\zeta_2\le 0,$ $\theta_2=\zeta_2$ for $\zeta_2>0;$
		$(t-c_1)^{\delta_1}(t-c_2)^{\delta_2}\varphi_0$ ($\varphi_0\in H$)
		is a solution of  Riemann boundary value task
		$\Phi^+(t)=[(a(t)-b(t))/(a(t)+b(t))]\Phi^-(t),$
		$\delta_1=\zeta_1+i\xi_1,$ $\delta_2=\zeta_2+i\xi_2$.

	\begin{center}
		{\bf {9.5. Approximate Solution of Nonlinear
				Singular Integro-Differential Equations on the Open
				Contour of Integration}}
	\end{center}
	
	Consider the nonlinear singular integro-differential
	equation
	\[
	Gx\equiv
	 a(t,x(t),\ldots,x^{(m)}(t))+
	S_{L}(h(t,\tau,x(\tau),\ldots,x^{(m)}(\tau)))=
	\]
	$$
	=f(t)
	\eqno (9.68)
	$$
	with conditions
	$$
	x(1)=x'(1)=\cdots=x^{(m-1)}(1)=0.
	\eqno (9.69)
	$$
	
	An approximate solution of the boundary value problem (9.68), (9.69)
	is sought in the form of the polynomial (9.64), the coefficients $ \{\alpha_k \} $ of
	which are determined from the system of nonlinear algebraic
	equations
	\[
	\tilde{G}_{n}\tilde{x}_{n}\equiv
	\]
	\[
	\equiv \bar{P}_{n}\Bigl[
	\bar{a}(t,\tilde{x}_{n}(t),\ldots,\tilde{x}_{n}^{(m)}(t))+
	S_{\gamma}(P_n^{\tau}[\bar{h}(t,\tau,\tilde{x}_{n}(\tau),
	\ldots,\tilde{x}_{n}^{(m)}(\tau))])\Bigr]=
	\]
	$$
	=  \bar{P}_{n}[\bar{f}(t)],
	\eqno (9.70)
	$$
	where \[
	\bar{a}(t,\tilde{x}_{n}(t),\ldots,\tilde{x}_{n}^{(m)}(t))=
	a(t,\tilde{x}_{n}(t),\ldots,\tilde{x}_{n}^{(m)}(t))
	\]
	for $t\in L,$
	\[
	\bar{a}(t,\tilde{x}_{n}(t),\ldots,\tilde{x}_{n}^{(m)}(t))=\tilde{x}_{n}^{(m)}(t)
	\]
	for $t\notin L,$
	\[
	\bar{h}(t,\tau,\tilde{x}_{n}(\tau),
	\ldots,\tilde{x}_{n}^{(m)}(\tau))=h(t,\tau,\tilde{x}_{n}(\tau),
	\ldots,\tilde{x}_{n}^{(m)}(\tau))\]
	 for $t\in L$ and $\tau \in L,$
	\[
	\bar{h}(t,\tau,\tilde{x}(\tau),\ldots,\tilde{x}^{(m)}(\tau))=0
	\]
	for  $t \notin L,$ $\tau \in \gamma$ or $t \in \gamma,$ $\tau \notin L,$ 
	\[
	\bar f(t) = f(t)
	\]
	 for $t\in L,$ 
	 \[
	 \bar{f}(t)=0
	 \] 
	 for $t\notin L.$
	
	Justification of the computational scheme (9.64) (9.70) for 
	approximate solution of the boundary value problem (9.68), (9.69) is
	conducted considering the results of the preceding item, as well as 
	the results of the item 9.1.
	
	{\bf Theorem 9.6 } \cite{Boy33}, \cite{Boy36}.  Let the boundary value problem (9.68), \ 
		(9.69) \  has \  a \  unique \  solution
		$x^*$ in a ball $ S=B(x^*,r)$ \ $ r>0; $ \  the \ conditions  
		$h'_{u_i}(t,\tau,u_0,u_1,\ldots,u_m)\in H_{\alpha,\ldots,\alpha}$ $(i=\overline{0,m}),$
		$a'_{u_i}(t,u_0,u_1,\ldots,u_m)\in H_{\alpha,\ldots,\alpha},$
		$f(t)\in H_{\alpha}$ are fulfilled  on  $L,$ and equation  $G'(x^*)z=f$ has
		a unique  solution, that satisfies (9.69) . Then, for such $n$  that $q=A\ln^5n(n^{-\alpha}+n^{-\theta})<1,$ the
		system of equations (9.70) has in a ball $S_1\subset S$
		a unique solution  $\tilde{x}_{n}^*$ and the estimate
		$$\|x^{*(m)}-\tilde{x}_{n}^{*(m)}\|_{L_2}\le
		A\ln^2n(n^{-\alpha}+n^{-\theta})$$ is valid. Here
		$\theta=\min\{\theta_1,\theta_2\},$ $\theta_1=1-|\zeta_1|$ for
		$\zeta_1\le 0,$ $\theta_1=\zeta_1$ for $\zeta_1>0,$
		$\theta_2=1-|\zeta_2|$ for $\zeta_2\le 0,$
		$\theta_2=\zeta_2$ for $\zeta_2>0;$
		$(t-c_1)^{\delta_1}(t-c_2)^{\delta_2}\varphi_0$ is a 
		solution of the Riemann boundary value problem
		\[
		\Phi^+(t)=
		\]
		\[
		=\frac{a'_{u_m}(t,x^*(t),\ldots,x^{*(m)}(t))-
			h'_{u_m}(t,t,x^*(t),\ldots,x^{*(m)}(t))}
		{a'_{u_m}(t,x^*(t),\ldots,x^{*(m)}(t))+
			h'_{u_m}(t,t,x^*(t),\ldots,x^{*(m)}(t))}\Phi^-(t),\]
		$\delta_1=\zeta_1+i\xi_1,\;\delta_2=\zeta_2+i\xi_2.$

	{\it Note 1. } The above  we have studied the computing
	schemes under the assumption that the solutions of the corresponding boundary value problems 
	owned to $L_2.$ If we assume that the solutions of boundary value problems
	belong to the class $L_p\;(1<p\le 2),$ then, after similar
	in nature, but more cumbersome calculations, we will see that
	the results, which are set forth above,  are valid in the 
	space $L_p\;(1<p\le 2).$
	
	{\it Note 2. } Justification of applicability of the method
	Newton-Kantorovich to the approximate solution of nonlinear
	algebraic systems, approximating non-linear singular
	integral-differential equations, carried out in the same manner as  in the case of nonlinear singular integral equations.

\begin{center}

{\bf Chapter 2}

{\bf Approximate Solution of Multi-Dimensional Singular Integral
Equations}
\end{center}

The chapter consists of the  three sections.
In the first section we will consider polysingular integral equations.
In the second section we will consider
multi-dimensional Riemann boundary value problems.
The third section is devoted to numerical methods for solution of multi-dimensional singular integral
equations with Zygmund-Calderon type of kernels.
Also in this section we will investigate parallel  iterative-projection
methods for
solution of multi-dimensional singular integral equations.

The main attention we will spare to collocation method, because the proofs of convergence for moment and the Galerkin methods
are similar to collocation method.

There \ are \  three \  principle different methods for proof of
the convergence of
projector methods for polysingular integral equations. Historical the first
method is based on transform the polysingular integral equations into
multi-dimensional Riemann boundary \  value \  problems. \  The \ second \ method \  is based on Simonenko local principle \cite{Sim}, \cite{Dud}. The third method is based on the
theory of commutative rings. The main results of the second direction and the
rich bibliography are given in \cite{Pr}. In the introduction to the
work it was noticed, that we will not touch  two last directions.

\begin{center}

{\bf 1. Bisingular integral equations}

\end{center}

     Let us consider bisingular integral equations of
the following kind
\[
Kx \equiv a(t_1,t_2)x(t_1,t_2)+b(t_1,t_2)S_{1}(x(\tau_1,t_2))+
\]
\[
+c(t_1,t_2)S_{2}(x(t_1,\tau_2))+d(t_1,t_2)S_{12}(x(\tau_1,\tau_2))+
\]
$$
+U_{12}[h(t_1,t_2,\tau_1,\tau_2)x(\tau_1,\tau_2)]=f(t_1,t_2),
\eqno (1.1)
$$
where
\[
S_{1}(x(\tau_1,t_2))=\frac{1}{\pi i}\int \limits_{\gamma_1}
\frac{x(\tau_1,t_2)}{\tau_1-t_1}d\tau_1 ,
\]
\[
S_{2}(x(t_1,\tau_2))=\frac{1}{\pi i}
\int \limits_{\gamma_2}\frac{x(t_1,\tau_2)}{\tau_2-
t_2} d\tau_2,
\]
\[
S_{12}(x(\tau_1,\tau_2))=-\frac{1}{\pi^2}\int \limits_{\gamma_1}
\int \limits_{\gamma_2}\frac{x(\tau_1,\tau_2)}{(\tau_1-t_1)(\tau_2-
t_2)}d\tau_1 d\tau_2,
\]
\[
U_{12}[h(t_1,t_2,\tau_1,\tau_2)x(\tau_1,\tau_2)]=\int \limits_{\gamma_1}
\int \limits_{\gamma_2}h(t_1,t_2,\tau_1,\tau_2)x(\tau_1,\tau_2)
d\tau_1 d\tau_2,
\]
$ \gamma_i=\{z_i: \ |z_i|=1, \  i=1,2 \}, \  \gamma_{12}= \gamma_1\times \gamma_2. $ 

The equation (1.1) we will consider under conditions
$a^2(t_1,t_2)-b^2(t_1,t_2) \ne 0.$

Numerical methods for solution of the equation (1.1) was devoted many works.

Convergence of one common projective method was proved by
A.V. Kozak and I.B.
Simonenko \cite{Koz}.

In  cases when $b=c=0$ or $a(t_1,t_2)=a(t_1)a(t_2),$ $b(t_1,t_2)=b(t_1)b(t_2),$
$c(t_1,t_2)=c(t_1)c(t_2),$ $d(t_1,t_2)=d(t_1)d(t_2)$ convergence of the method
of mechanical quadrature was proved by I.V. Boykov \cite{Boy13}, \cite{Boy16}, \cite{Boy25}.

 I.V. Boykov proved convergence of projector method in common case when functions $a(t_1,t_2),$ $b(t_1,t_2),$ $c(t_1,t_2),$ $d(t_1,t_2)$
are belong to Holder class and factorable.

I.K. Lifanov \cite{Lif} used the discrete vortex method for solution of the equation
(1.1) under conditions
$a(t_1,t_2) \equiv b(t_1,t_2) \equiv c(t_1,t_2) \equiv h(t_1,t_2,\tau_1,
\tau_2) \equiv 0.$

At first we consider the equation
$$
Kx \equiv a(t_1,t_2)x(t_1,t_2)
+d(t_1,t_2)S_{12}(x(\tau_1,\tau_2))
=f(t_1,t_2).
\eqno (1.2)
$$

Let us introduce the function
\[
X(z_1,z_2)=-\frac1{\pi^2}\int \limits_{\gamma_1}\int
\limits_{\gamma_2}\frac{x(\tau_1,\tau_2)}{(\tau_1-z_1)(\tau_2-z_2)}
d\tau_1 d\tau_2.
\]

     Let $ D_i^{+} (D_i^{-}), i=1,2, $ be the set of points $ z_i,
i=1,2, $ on the complex plane $ z_i, i=1,2, $ which  satisfy the following conditions: $ |z_i|<1 $ $(|z_i|>1)$, $ i=1,2. $ Let $ D^{++}=D^{+}_1 D^{+}_2$, $
D^{+-}=D^{+}_1 D^{-}_2,$  $ D^{-+}=D^{-}_1 D^{+}_2,$  $ D^{--}=D^{-}_1
D^{-}_2. $ If a function $ f(z_1,z_2) $ is analytical in the
domain $ D^{++}, $  we write $ f^{++}(z_1,z_2). $
If the function $ f(z_1,z_2) $ is analytical in the domain
$ D^{+-}(D^{-+}  \ {\rm or}\   D^{--}), $ we write $ f^{+-}(z_1,z_2)$ $(f^{-+}
(z_1,z_2)\  {\rm or}\  f^{--}(z_1,z_2)). $

     Let\  $ (z_1,z_2) \ \in \  D^{++}.\  $ Let \ $ t_i \ \in \ \gamma_i, i=1,2. $ \ 
     We \ will \ write $ f^{++}(t_1,t_2)=\lim_{(z_1,z_2)\to(t_1,t_2)}
f^{++}(z_1,z_2). $ Limiting values $ f^{+-}(t_1,t_2)$, $f^{-+}(t_1,t_2)$, $
f^{--}(t_1,t_2) $ of the functions $ f^{+-}(z_1,z_2)$, $f^{-+}(z_1,z_2)$, $
f^{--}(z_1,z_2) $ are defined by the similar way.

     We need in Sohotzky - Plemel  formulas \cite{Gakh}
$$
\begin{array}{ccc}
X^{++}(t_1,t_2)+X^{+-}(t_1,t_2)+X^{-+}(t_1,t_2)+X^{--}(t_1,t_2)=
S_{12}x, \\
X^{++}(t_1,t_2)-X^{+-}(t_1,t_2)-X^{-+}(t_1,t_2)+X^{--}(t_1,t_2)=x. \\
\end{array}
\eqno (1.3)
$$

     Using Sohotsky $-$ Plemel formulas (1.3), we can write  the equation (1.2) as
\[
(a(t_1,t_2)+b(t_1,t_2))X^{++}(t_1,t_2)-(a(t_1,t_2)-b(t_1,t_2))
X^{-+}(t_1,t_2)-
\]
\[
-(a(t_1,t_2)-b(t_1,t_2))X^{+-}(t_1,t_2)+(a(t_1,t_2)+
b(t_1,t_2))X^{--}(t_1,t_2)=
\]
$$
=f(t_1,t_2).
\eqno (1.4)
$$
     Dividing the equation (1.4) on $ a+b $, we receive the equation
\[
X^{++}(t_1,t_2)-G(t_1,t_2)(X^{+-}(t_1,t_2)+X^{-+}(t_1,t_2))+
X^{--}(t_1,t_2)=
\]
$$
= f_1(t_1,t_2),
\eqno (1.5)
$$
where $ G(t_1,t_2)=((a(t_1,t_2)-d(t_1,t_2))/(a(t_1,t_2)+d(t_1,t_2)),$\\$
f_1(t_1,t_2)=f(t_1,t_2)/(a(t_1,t_2)+d(t_1,t_2)). $

     Let $ a(t_1,t_2) \pm b(t_1,t_2) $ be not equal to zero, $ a(t_1,
t_2), b(t_1,t_2) \in H_{\alpha \alpha} (0<\alpha<1). $ Let the partial
indexes \cite{Kak} of the operator $ K $ are equal to zero.

     Well known \cite{Kak}  the representation
\[ G(t_1,t_2)= \frac{\displaystyle \psi^{++}(t_1,t_2) \psi^{--}
(t_1,t_2)}{\displaystyle \psi^{+-}(t_1,t_2) \psi^{-+}(t_1,t_2)}, 
\]
where
\[ \psi^{\pm \pm}(t_1,t_2)=\exp \{P^{\pm \pm }[\ln G(t_1,t_2)] \}, 
\]
\[ P^{\pm \pm}=P^{\pm} P^{\pm}, P^{\pm}=\frac12(I\pm S), Sf=\frac{1}
{\pi i} \int \limits_\gamma \frac{f(\tau)d\tau}{\tau-t}, 
\]
$ \gamma =\{z:|z|=1 \}.$

     An approximate solution of the equation (1.5) we will look for in the form
$$
\begin{array}{ccc}
X^{++}_n(t_1,t_2)=\psi^{++}(t_1,t_2) \sum\limits_{k=0}^n \sum\limits_{l=0}^n
\alpha_{kl}t_1^k t_2^l,  \\
X^{+-}_n(t_1,t_2)=\psi^{+-}(t_1,t_2) \sum\limits_{k=0}^n \sum\limits_{l=-n}^{-1}
\alpha_{kl}t_1^k t_2^l,  \\
X^{-+}_n(t_1,t_2)=\psi^{-+}(t_1,t_2) \sum\limits_{k=-n}^{-1} \sum\limits_{l=0}^n
\alpha_{kl} t_1^k t_2^l,   \\
X^{--}_n(t_1,t_2)=\psi^{--}(t_1,t_2) \sum\limits_{k=-n}^{-1} \sum\limits_{l=-n}^{-1}
\alpha_{kl} t_1^k t_2^l. \\
\end{array}
\eqno (1.6)
$$

So, approximate solution of the equation (1.2) we will seek in the form
\[
X_n(t_1,t_2) = X_n^{++}(t_1,t_2) - X_n^{+-}(t_1,t_2) -
X^{-+}(t_1,t_2) + X^{++}(t_1,t_2).
\]

The coefficients $d_{kl},$ $k,l=-n,\ldots,-1,0,1, \ldots,n,$ are defined from
the system of linear algebraic equations
\[
K_nX_n \equiv P_{nn}[X_n^{++}(t_1,t_2)-G(t_1,t_2)
()X_n^{+-}(t_1,t_2)+X_n^{-+}(t_1,t_2))+
\]
$$
+X_n^{--}(t_1,t_2)]
= P_{nn}[f_1(t_1,t_2)].
\eqno (1.7)
$$

The convergence of collocation method (1.7) we investigate
in the Banach space $ E $ of functions $ f(t_1,t_2), $ which are  belong to
the Holder class $ H_{\alpha,\alpha}, 0<\alpha<1, $ with the norm
\[
 \| f(t_1,t_2)\|=
 \]
 $$
 = \max_{\displaystyle(t_1,t_2)\in \gamma_{12}} |f(t_1,t_2)| +
\max_{t_2}\sup_{\displaystyle t_1'' \neq t_1'}
\frac{|f(t_1',t_2)-f(t_1'',t_2)|}{|t_1'-t_1''|^\beta}+$$
$$+\max_{t_1}\sup_{\displaystyle t_2' \neq t_2''}\frac{|f(t_1,t_2')-
f(t_1,t_2'')|}{|t_2'-t_2''|^\beta}, $$ $ 0<\beta<\alpha. $

{\bf Theorem 1.1} \cite{Boy13}, \cite{Boy16}, \cite{Boy25}. Let the following conditions are satisfied:\\
1) the operator $ K $ is continuously invertible in the space $ E; $\\
2) $ a \pm b \neq 0; $\\
3) the partial indexes of the operator $ K $ are
equal to zero;\\
4) $a, b, f \in H_{\alpha, \alpha}, 0<\alpha<1. $

Then, for $ n $ such that $ q=An^{-\alpha+\beta} \ln^2n<1, $
the system (1.7) has a unique solution $ x^*_n(t_1,t_2) $ and the estimate
$ \|x^*-x^*_n\|_E \le An^{-\alpha+\beta} \ln^2n $ is valid, where  $ x^*(t_1,t_2) $ is a unique solution of the equation (1.1).

{\bf Theorem 1.2} \  \cite{Boy16}, \ \cite{Boy25}.\  Let the following \ conditions are satisfied:\\
1) the operator $ K $ is continuously invertible in the space
$ L_2(\gamma_{12}); $\\
2) $ a \pm b \neq 0; $\\
3) the partial indexes
of the operator $ K $ is equal to zero;\\
4)$ a, b, f \in H_{\alpha,
\alpha}, 0<\alpha<1. $

Then, for $ n $ such that
$ q=A \max (\omega(a; n^{-1}), n^{-\alpha} \ln n, \omega(b; n^{-1}))<1, $
the system (1.7) has a unique solution $ x^*_n(t_1,t_2) $ and the estimate
$ \|x^*-x^*_n\| \le A \max (\omega(a; n^{-1}), n^{-\alpha} \ln n, \omega(b; n^{-1}),
\omega(f; n^{-1})) $ is valid in the space $ L_2[\gamma_{12}]. $ Here
$ x^* $ is a unique solution of the equation (1.1).

     The similar statements is valid for singular integral equations
\[
a(t_1,t_2)x(t_1,t_2)-\frac{1}{ \pi^2} \int \limits_{\gamma_1} \int
\limits_{\gamma_2} \frac{x(\tau_1,\tau_2)d\tau_1 d\tau_2}{(\tau_1-
t_1)(\tau_2-t_2)}+ 
\]
$$
+
\int \limits_{\gamma_1} \int \limits_{\gamma_2}
h(t_1,t_2,\tau_1,\tau_2)x(\tau_1, \tau_2) d\tau_1 d\tau_2=f(t_1,t_2) 
$$
and, under special conditions \cite{Boy16}, \cite{Boy25} for coefficients $ a(t_1,t_2),$ $
b(t_1,t_2),$ $ c(t_1,t_2),$ $ d(t_1,t_2), $ for singular integral equations
\[
a(t_1,t_2)x(t_1,t_2)+\frac {b(t_1,t_2)}{\pi} \int \limits_{\gamma_1}
\frac{x(\tau_1,t_2)}{\tau_1-t_1} d\tau_1 + 
\]
\[
 + \frac{c(t_1,t_2)}{\pi} \int \limits_{\gamma_2} \frac{x(t_1,
\tau_2)}{\tau_2-t_2} d\tau_2 + \frac {d(t_1,t_2)}{\pi^2} \int
\limits_{\gamma_1} \int \limits_{\gamma_2} \frac{x(\tau_1,\tau_2)
d\tau_1 d\tau_2} {(\tau_1-t_1)(\tau_2-t_2)} +
\]
\[ + \int \limits_{\gamma_1} \int \limits_{\gamma_2} h(t_1,t_2,
\tau_1,\tau_2)x(\tau_1, \tau_2) d\tau_1 d\tau_2= f(t_1,t_2). 
\]

\begin{center}
{\bf 2. Riemann Boundary Value Problem}
\end{center}

Let us consider iterative methods for solution of
characteristical bisingular equations 
$$
a(t_1,t_2)\varphi(t_1,t_2)+b(t_1,t_2)S_{12}(\varphi(\tau_1,\tau_2))=
f(t_1,t_2).
\eqno (2.1)
$$

Conditions on functions $ a,b,f $ will be imposed below.
We now  assume that they are continuous and $ a \pm b \ne 0 $ on
$ \gamma_{12}. $ In the previous section was shown that
the preceding equation reduces to the Riemann boundary value problem
$$
\varphi^{++}-\frac{a-b}{a+b} \varphi^{-+}-\frac{a-b}{a+b} \varphi^{-+}-
\varphi^{--}=\frac{1}{a+b}f
\eqno (2.2)
$$

Let us suppose that for any poits $ t_2 \in \gamma_2 $ the
values of the function
\[ G(t_1,t_2)=\frac{a(t_1,t_2)-b(t_1,t_2)}{a(t_1,t_2)+b(t_1,t_2)}
\hskip 20 pt (G(t_1,t_2) \ne 0) 
\]
lie inside of an angle $ \Gamma_1 $, which has the vertex  in origin of the plane
of complex variable $ Z_1.$ Also we suppose that    values of the angle $ \Gamma_1 $ are  not exceed $ \pi-\delta_1(\delta_1>0). $
Then one can find  such constant $ \alpha_1 $ that $ \max_{(t_1,t_2)
\in \gamma_{12}} \mid \alpha_1 G(t_1,t_2)-1 \mid \le q_1<1. $

Let for any fixed poits $ t_1 \in \gamma_1 $ values of functions
$ G(t_1,t_2) $ lie inside of an angle $ \Gamma_2 $, which has the center at
origin. Also we suppose that  the values of angle are  not exceed 
$ \pi- \delta_2(\delta_2>0). $ Then it exists a constant $ \alpha_2 $
such that $ \max_{(t_1,t_2) \in \gamma_{12}} \mid \alpha_2 G_1(t_1,
t_2)-1 \mid \le q_2<1. $

Let be $ q=\min(q_1,q_2). $ For certainty we put $q=q_1$ and  $ \alpha =\alpha_1.$ 
We will change the functions $\varphi^{\pm \pm}$ using the formulas
$$\psi^{++}=\varphi^{++}; \psi^{+-}=\alpha^{-1}\varphi^{+-}; \psi^{-+}=
\alpha^{-1}\varphi^{-+}; \psi^{--}=\varphi^{--}$$
and receive the equation
\[ \psi^{++}- \alpha G(t_1,t_2) \psi^{+-}- \alpha G(t_1,t_2)
\psi^{-+}+ \psi^{--}=(a+b)^{-1}f. 
\]

By adding the function $ \psi $ to the both sides of the preceding
equation and by  using the Sohotzky $-$ Plemel formulas, we obtain the
equation
\[ \psi(t_1,t_2)=(\alpha G(t_1,t_2)-1)\psi^{+-}(t_1,t_2)+(\alpha
G(t_1,t_2)-1) \psi^{-+}(t_1,t_2)+ 
\]
\[ +(a(t_1,t_2)+b(t_1,t_2))^{-1}f(t_1,t_2). 
\]

For solving this equation the following iterative process may be used
$$
\psi_{n+1}=(\alpha G-1) \psi_n^{+-}+(\alpha G-1)
\psi_n^{-+}+(a+b)^{-1}f.
\eqno (2.3)
$$

We will investigate its convergence in  the space
$ L_2(\gamma_{12}). $ From  Sohotzky $-$ Plemel formulas  follows
that
$$ (\alpha G-1) \psi^{+-}+(\alpha G-1) \psi^{-+}=(\alpha G-1)(- \frac
{1}{2} \psi+ \frac{1}{2}S_{12} \psi). $$

Since $ \| S_{12} \|_{L_2}=1, $ then
$$ \|(\alpha G-1) \psi^{+-}+
(\alpha G-1) \psi^{-+} \|_{L_2} \le $$
$$ \le \frac{1}{2} \max_{(t_1,t_2) \in \gamma_{12}} \mid \alpha
G(t_1,t_2)-1 \mid(\| \psi \|+ \|S_{12} \psi \|) \le $$
$$ \le \max_{(t_1,t_2) \in \gamma_{12}} \mid \alpha G(t_1,t_2)-1 \mid
\| \psi \| \le q\| \psi \|. $$

The iterative process (2.3) thus converges.
The estimate $ \| \varphi^* - \varphi_n^* \| \le Aq^n $ is valid, where
$\varphi^*$ is a solution of the equation (2.1),
$ \varphi_n^{*++}= \psi_n^{++}, \varphi^{*+-}= \alpha \psi^{+-}, \varphi^{*-+}
=\alpha \psi_n^{-+}, \varphi^{*--}= \varphi_n^{--}; \psi_n $ is result of 
$ n-$ iteration  which is obtained by  the formula (2.3).

{\bf Theorem 2.1} \cite{Boy16}, \cite{Boy25}. Let the following conditions  are fulfilled:\\ 1) $ a(t_1,t_2) \pm b(t_1,t_2)
\ne 0 $ for $ (t_1,t_2) \in \gamma_{12};$\\
2) $a(t_1,t_2), b(t_1,t_2) \in C[\gamma_{12}]; $\\
 3)  for any fixed $ t_1 \in
\gamma_1(t_2 \in \gamma_2) $ the values of  function $ G(t_1,t_2) $
lie inside of an angle  of the solution
$ \pi- \alpha_1(\pi- \alpha_2), \alpha_i>0 (i=1,2) $ with the center in origin. 

Then the equation
(2.1) has a unique solution $ \varphi^* $ to which the iterative process
(2.3) (after returning from the function $\psi$ to the function $\varphi$)
converges with the rate of geometric progression.

It is not difficult to see \ 
(using Sohotzky $-$ Plemel
formulas), that an iterative processes may be built analogously  for
equations 
\[ a(t_1,t_2) \varphi(t_1,t_2)+b(t_1,t_2)S_1 \varphi
(\tau_1,t_2)=f(t_1,t_2), 
\]
\[
 a(t_1,t_2) \varphi(t_1,t_2)+b(t_1,t_2)S_2 \varphi(t_1, \tau_2)=f(t_1,
t_2), 
\]
\[
 a(t_1,t_2)S_1 \varphi(\tau_1,t_2)+b(t_1,t_2)S_2 \varphi(t_1, \tau_2)=f
(t_1,t_2), 
\]
\[ a(t_1,t_2)S_1 \varphi(\tau_1,t_2)+b(t_1,t_2)S_{12} \varphi(\tau_1,
\tau_2)=f(t_1,t_2), 
\]
\[ a(t_1,t_2)S_2 \varphi(t_1, \tau_2)+b(t_1,t_2)S_{12} \varphi(\tau_1,
\tau_2)=f(t_1,t_2). 
\]

We turn our attention to one moment.

The equation 
\[ a(t_1,t_2)x(t_1,t_2)+b(t_1,t_2)S_{12}(x(\tau_1, \tau_2))+ 
\]
$$
+U_{12}(h(t_1,t_2, \tau_1, \tau_2)x(\tau_1, \tau_2))=f(t_1,t_2),
\eqno (2.4)
$$
where $ U_{12}$ is  totally continuous operator, can be
studied as  the equation (2.1) and as above mentioned 
equations.

Let us now consider the more general singular integral equation
\[ a(t_1,t_2)x(t_1,t_2)+b(t_1,t_2)S_1(x(\tau_1,t_2))+ 
\]
$$
+c(t_1,t_2)S_2(x(t, \tau_2))+d(t_1,t_2)S_{12}
(x(\tau_1, \tau_2))=f(t_1,t_2).
\eqno (2.5)
$$

This equation is equivalent to the Riemann boundary value problem
\[ a_1(t_1,t_2) \psi^{++}(t_1,t_2)+b_1(t_1,t_2) \psi^{+-}
(t_1,t_2)+ 
\]
\[+c_1(t_1,t_2) \psi^{-+}(t_1,t_2)+d_1(t_1,t_2) \psi^{--}(t_1,t_2)=f(t_1,t_2), 
\]
in which the coefficients $ a_1,b_1,c_1,d_1 $ are determined, using  Sohotzky $-$ Plemel
formulas  for functions $ x(t_1,t_2),$ $ S_1(x(\tau_1,t_2)),$ $ S_2(x(t_1,\tau_2)),$ $ S_{12}(x(\tau_1,\tau_2)). $ 

     Assume that the conditions, described in   this section
     above, are fulfilled.  These conditions  guarantee the existence of constants  $ \alpha, \beta,
\gamma, \delta $ such, that
\[ \mid \alpha a_1(t_1,t_2)+1 \mid<q_1, \ \ \ 
\mid \beta b_1(t_1,t_2)-1 \mid<q_2, 
\]
\[ \mid \gamma c_1(t_1,t_2)-1 \mid<q_3, \ \ \ 
\mid \delta d_1(t_1,t_2)+1 \mid<q_4. \]

Then, provided $ q_1+q_2+q_3+q_4<1 $, the iterative process
\[ v_{n+1}=(\alpha a_1(t_1,t_2)+1) v_n^{++}+(\beta
b_1(t_1,t_2)-1) v_n^{+-}+ 
\]
\[ +(\gamma c_1(t_1,t_2)-1) v_n^{-+}+(\delta d_1(t_1,t_2)+1)
v_n^{--}+f(t_1,t_2), 
\]
 where $ v^{++}= \alpha^{-1}
\psi^{++}, v^{+-}= \beta^{-1} \psi^{+-}, v^{-+}=
\gamma^{-1} \psi^{-+}, v^{--}= v^{-1} \psi^{--} $,
converges to a solution of equation (2.5) with the rate of
geometric progression, of course, after  replacing the function $ v $ on the function $ \psi $.

    The proof of this statement is conducted in the same way, as in
case of simplest bisingular integral  equation (2.1).

\begin{center}
{\bf 3. Approximate Solution  of  Multi-Dimensional Singular Integral Equations }
\end{center}

In this section we will consider numerical methods for solution of
multi-dimensional singular integral equation 
\[
Kx\equiv a(t)x(t) + b(t)\int\limits_{E_m}
\frac{\varphi(t,\Theta)x(\tau)}{(r(t,\tau))^m}d\tau+
\]
$$
+ \int\limits_{E_m} h(t,\tau)x(\tau)d\tau =f(t),
\eqno (3.1)
$$
where $t=(t_1,\ldots,t_m),$ $\tau=(\tau_1,\ldots,\tau_m),$ $m=2,3,\ldots,$
$r(t,\tau)=(\sum\limits^{m}_{k=1}(t_k-\tau_k)^2)^{1/2},$
$\Theta=((t-\tau)/r(t,\tau)),$ $a(t), b(t), f(t)$ are smooth functions of $m$
variables, $h(t,\tau)$ is a smooth function of $2m$ variables.

These equations  have many applications in the elasticity theory,  in the
oscillating theory  and other. In spite of this there are only   a few works
devoted to approximate methods for solution of multi-dimensional singular
integral equations. At first we notice the papers \cite{Mich1}, \cite{Rad}, in which the method of moments and the method of least squares
are used for solution of the equation $Kx=f$.

Several iterative methods for solution of equations  (3.1) was offered in
\cite{Musa}, \cite{Par}.

In this section we give a review of works \cite{Boy30}, \cite{Boy31} and
\cite{Boy25}.

\begin{center}
{\bf 3.1.  Approximate Solution of Multi-Dimensional Singular
	Integral Equations on Holder Classes of Functions}
\end{center}

Let us consider the equation
$$
Kx \equiv a(t) x(t) + b(t)\int\limits_G {\varphi(\Theta) x(\tau)
\over (r(t, \tau))^2} d\tau = f(t) ,
\eqno (3.2)
$$
where  $G$  is a simply connected domain on $E_2, t=(t_1,t_2),
\tau=(\tau_1,\tau_2),$ $ \Theta=((\tau - t)/r(t,\tau)),
r(t,\tau)=((t_1-\tau_1)^2 + (t_2-\tau_2)^2)^{1/2}.$

Assume that functions $a(t)$, $b(t)$, $f(t),$ $\varphi(\Theta) \in H_{\alpha \alpha}(1),$
 $\alpha$ ($0 < \alpha \le 1$).

All results of this item can be diffused on singular integral equation 
$$
a(t) x(t) + \int\limits_G {\varphi(t, \Theta) x(\tau)
\over (r(t, \tau))^2} d\tau = f(t),
$$
where  $G=[-A,A]^l, l=2,3,\dots,$  $ t=(t_1,\dots,t_l),
\tau=(\tau_1,\dots,\tau_l),$ $ \Theta=((\tau - t)/r(t,\tau)),
r(t,\tau)=((t_1-\tau_1)^2 + \dots+(t_l-\tau_l)^2)^{1/2}.$

For simplicity we consider the case when $G=[-A,A; -A,A]$.

Let us cover the domain $G$ by squares
 $\Delta_{k l} = [t_k,t_{k+1};t_l,t_{l+1}]$, $k,l = \overline{0,N+1}$, where
 $t_k = \displaystyle -A + 2A k/( N+2)$, $k = \overline{0,N+2}$.

Together with squares $\Delta_{k l}$, $k,l = \overline{0,N+1}$,
we will use   rectangles
$\bar \Delta_{k l}$, $k,l = \overline{1,N}$, which are defined by following way:
\\ $\bar \Delta_{k l} = \Delta_{k l}$ for $k = \overline{2,N-1}$ and $l = \overline{2,N-1}$;
\\ $\bar \Delta_{1 1} = \Delta_{0 0} \cup \Delta_{0 1} \cup \Delta_{1 0} \cup \Delta_{1 1}$;
\\ $\bar \Delta_{1, N} = \Delta_{0, N} \cup \Delta_{1, N} \cup \Delta_{0, N+1} \cup \Delta_{1, N+1}$;
\\ $\bar \Delta_{N, 1} = \Delta_{N, 0} \cup \Delta_{N+1, 0} \cup \Delta_{N, 1} \cup \Delta_{N+1, 1}$;
\\ $\bar \Delta_{N, N} = \Delta_{N, N} \cup \Delta_{N, N+1} \cup \Delta_{N+1, N} \cup \Delta_{N+1, N+1}$;
\\ $\bar \Delta_{1, l} = \Delta_{0, l} \cup \Delta_{1, l}$ for $l = \overline{2,N-1}$;
\\ $\bar \Delta_{N, l} = \Delta_{N, l} \cup \Delta_{N+1, l}$ for $l = \overline{2,N-1}$;
\\ $\bar \Delta_{k, 1} = \Delta_{k, 0} \cup \Delta_{k, 1}$ for $k = \overline{2,N-1}$;
\\ $\bar \Delta_{k, N} = \Delta_{k, N} \cup \Delta_{k, N+1}$ for $k = \overline{2,N-1}$.

Approximate solution of the equation (3.1) we will seek as piecewise constant
function $x_N (t_1,t_2).$ The function $x_N (t_1,t_2)$ is equal to unknown constant
$x_{k l}$ on each rectangular $\bar\Delta_{k l}$, $k,l=1,2,\ldots,N.$

Unknown values $x_{k l},$ $k,l=1,2,\ldots,N,$ are defined from the
system of linear algebraic equations
$$
a(\bar t_{k l}) x_{k l} +b(\bar t_{k l}) \displaystyle \sum\limits^N_{i=1} \sum\limits^N_{j=1} {}' \bar d_{i j}
    (\bar t_{k l}) x_{i j} + x_{k l} d_{k l} (\bar t_{k l}) = f (\bar t_{k l}),
\eqno (3.3)
$$
$    k,l =1,2,\ldots,N,$\\
where $\sum \sum'$ denote that the sum is taken over  rectangules
$\bar\Delta_{ij},$ which intersections with squares $\Delta_{kl}$ are empty,
$t_{kl}=(t_k,t_l),$ $\tau_{kl}=(\tau_k,\tau_l),\ \   \bar t_{k l} = (t_k + \bar h_1, t_l + \bar h_2)$,
\[ 
d_{i j} (\bar t_{k l}) = \displaystyle \int\limits_{\Delta_{i j}} \varphi
    \left({\bar t_{k l} - \tau \over r(\bar t_{k l} , \tau)} \right) {d \tau \over r^2 (\bar
    t_{k l} , \tau)}, 
   \]
   \[ 
\bar d_{i j} (\bar t_{k l}) = \displaystyle \int\limits_{\bar \Delta_{i j}} \varphi
    \left({\bar t_{k l} - \tau \over r(\bar t_{k l} , \tau)} \right) {d \tau \over r^2 (\bar
    t_{k l} , \tau)},
  \] 
values of $\bar h_1$ and $\bar h_2$ will be given later.

Let us estimate the expression
\[
\displaystyle \sum\limits_{i=1}^N \sum\limits_{j=1}^N {}' \bar d_{i j} (\bar t_{k l}) \le
    B h^2 \sum\limits_{i=1}^N \sum\limits_{j=1}^N {}' {1 \over (t_k - t_j)^2 + (t_l - t_j)^2} \le
\]
\[
    \le 8 B h^2 \sum\limits_{i=1}^{[{N \over 2}] + 1} \sum \limits_{j=1}^{[{N \over 2}] + 1}
    {1 \over t^2_i + t_j^2} \le
\]
\[
 \le \displaystyle 8 B h^2 \sum\limits_{i=1}^{[{N \over 2}] + 1} \sum\limits_{j=1}^
     {[{N \over 2}] + 1} {N^2 \over (i^2 + j^2) 4 A^2} =
\]
$$
     =8 B \sum\limits_{i=1}^{[{N \over 2}] +
     1} \sum\limits_{j=1}^{[{N \over 2}] + 1} {1 \over i^2 + j^2} \le
 D \ln N,
\eqno (3.4)
$$
where $h=2A/N.$

Let us denote the matrix of the system (3.3) as $C=\{c_{ij}\},$
$i,j=1,2,\ldots,N^*, N^*=N^2.$ So, for arbitrary $j,$ $1 \le j \le N^2,$
we have
$$
\displaystyle\sum\limits^{N^*}_{k=1, k \ne j} |c_{j k}| \le D \ln N.
\eqno    (3.5)
$$

Let us show that we can choose parameters $\bar h_1$,  $\bar h_2$,
so way, that the matrix $C$ will be invertible.

Well known \cite{Mich2}, that the necessary and sufficient condition of the
existence of multi-dimensional singular integrals \\
$ \int\int\limits_{E_2}\frac{\varphi(\Theta)x(\tau)}{r^2(t,\tau)}d\tau,$ 
$t=\bar t_{kl},$
 is the equality
$ \int\limits_S \varphi \left({\bar t_{k l} - \tau \over r(\bar t_{k l} , \tau)}
    \right) d s =0,$
\\ where $S$ is the unit circle with the center in the point $\bar t_{k l}$.

From this condition follows existence at least of two rays from point
$\bar t_{kl},$ where the function $\varphi(\Theta)=0.$ So, varying
values $\bar h_1$ and $\bar h_2$, we  can make the integral
\[
\int\int\limits_{\Delta_{kl}}\varphi\left(\frac{\tau_1-\bar t_k}
{r(\tau,\bar t_{kl})},
\frac{\tau_2-\bar t_l}{r(\tau,\bar t_{kl})}\right)\frac{d\tau}
{r^2(\tau,\bar t_{kl})}
\]
bigger than any fixed number $M.$ We can select $M$ bigger than
$D \ln N.$ This inequality   is the sign, that conditions of Hadamard Theorem about invertibility of matrix
is fulfilled. So, we proved the existence of a unique solution $x^*_N(t_1,t_2)$
of the system (3.3).

We will estimate the error of approximation of  exact solution $x^*(t_1, t_2)$ of the equation (3.2)  with respect to    solution $x^*_N (t_1, t_2)$
of the equation (3.3).

Let $x^* (t_1 , t_2)$ be a solution of the equation (3.2).
Equating left and right sides of the equation (3.3) in points $\bar t_{k l}$,
we have:
\[
a(\bar t_{k l}) x^*(\bar t_{k l}) + \displaystyle \sum\limits^N_{i=1} \sum\limits^N_
    {j=1} \int\limits_{\bar\Delta_{i j}} {\varphi \left( {\tau_1 - \bar t_k \over r (\tau, \bar
    t_{k l})} , {\tau_2 - \bar t_l \over r (\tau, \bar t_{k l})} \right) x^*(\tau) \over
    r^2 (\tau, \bar t_{k l})} d \tau =
\]
$$
=f(\bar t_{k l}),\ \  k,l =\overline{1,N}. \eqno (3.6)
$$

Let $P_N$ is the projector from space $X=C(G)$
onto the set of piece-constant functions. This projector is defined as
\[
P_N f(t)=\left\{
\begin{array}{cc}
f(\bar t_k,\bar t_l), \, {\rm if} \, t \in \bar \Delta_{kl},\\
0, \, {\rm if} \, t \in G \backslash \bar\Delta_{kl}.\\
\end{array}
\right.
\]

The system (3.3) in operator form can be written as
\[
K_Nx_N\equiv
P_N\left[a(t_1,t_2) x_N( t_1,t_2) +\right.
\]
\[
+b(t) \left.\int\limits_{G}P(t,\tau)\left[ {\varphi \left(
{\tau_1 -  t_1 \over r (\tau, t)} , {\tau_2 -  t_2 \over r (\tau, \bar t)}
\right) x_N(\tau) \over
    r^2 (\tau, \bar t)} d \tau\right]\right]=
\]
\[
=P_N[ f (t)],
\]
where
\[
P(\bar t_{kl},\tau)=\left\{
\begin{array}{cc}
1, \, {\rm if} \,  \tau \in G\setminus g_{kl},\\
0, \, {\rm if} \, \tau \in g_{kl} ,\\
\end{array}
\right.
\]
 $g_{k l} = [t_{k-1}, t_{k+2}; t_{l-1}, t_{l+2}] \backslash \Delta_{k l}$.

Then

\[x^*_N - P_N x^*  = K^{-1}_N (K_N(x^*_N - P_N x^*)) =
\]
\[
= K^{-1}_N (P_N f-
K_N P_N x^*) =
K_N^{-1} (P_N K x^* - K_N P_N x^*) =
\]
\[
= K_N^{-1} (P_N K x^* - P_N K P_N x^*) +
K_N^{-1} (P_NK P_N x^* - P_N K_N P_N x^*).
\]

The system (3.3) we can consider in the space $R_{N^*}$ with the norm
$\|u\|=\max_{1\le k\le N^*} |u_k|,$ where  $u=(u_1,\dots,u_{N^*}).$
In the space $R_{N^*}$ 
$\max_{1\le k,l \le N-1} |x_{k,l}|\le B \max_{1\le k,l \le N-1} |f_{k,l}|,$ where 
$f_{kl}=f(\bar t_{k,l}).$
The solution $x_N^*(t_1,t_2)$ is a piece-constant function with the norm
 $\|x_N^*(t_1,t_2)\|_C=\max_{1\le k,l \le N-1} |x_{k,l}|.$ Easy to see that
 $\|P_Nf(t_1,t_2)\|_C=\max_{1\le k,l \le N-1} |f_{k,l}|.$
 From these  equalities follow that  $\|K^{-1}_N\|\le B.$

In the space $C$ we have
$$
||x^*_N - P_N x^*|| \le B N^{-\alpha} + B ||P_N [K P_N x^* - K_N P_N x^*]||.
\eqno (3.7)
$$

Well known, that $||P_N||= 1.$
So, we need to estimate the integral
\\ $I_1 = \Bigg\vert \displaystyle \int\limits_{g_{k l}} \varphi \left({\tau_1 - \bar t_k \over
    r (\tau, \bar t_{k l})} , {\tau_2 - \bar t_l \over r (\tau, \bar t_{k l})} \right)
    {1 \over r^2 (\tau, \bar t_{k l})} P_N [x^*] d \tau \Bigg\vert.$

Let us construct a function $\psi (x_1, x_2) \in H_{\alpha \alpha}$ $(0 < \alpha \le 1),$ which has properties: 1)
$\psi(\bar t_{k l}) = (P_N x^*) (\bar t_{k l})$ and
$$
2) \displaystyle\int\limits_{g_{k l}} \varphi \left({\tau_1 - \bar t_k \over r(\tau, \bar t_{k l})},
    {\tau_2 - \bar t_l \over r(\tau, \bar t_{k l})} \right) {\psi (\tau_1, \tau_2) \over r^2 (\tau,
    \bar t_{k l})} d \tau_1 d \tau_2 = 0. \eqno (3.8)
$$

Existence of such function follows from the fact, that function $\varphi (\Theta)$
is equal to zero at least on two rays outbound from point $\bar t_{k l}$.

Using (3.8), we can prove that
$$
\Bigg\vert \displaystyle \int\limits_{g_{k l}} \varphi \left({\tau_1 - \bar t_k \over
    r (\tau, \bar t_{k l})} , {\tau_2 - \bar t_l \over r (\tau, \bar t_{k l})} \right)
    {1 \over r^2 (\tau, \bar t_{k l})} (P_N x^* -
$$
$$
-\psi(\tau_1, \tau_2)) d \tau_1 d \tau_2
    \Bigg\vert \le
$$
$$
\le \displaystyle \int\limits_{g_{k l}} \Bigg\vert \varphi \left({\tau_1 - \bar t_k \over
    r(\tau, \bar t_{k l})}, {\tau_2 - \bar t_l \over r(\tau, \bar t_{k l})} \right)\Bigg\vert  {1 \over
    r^2 (\tau, \bar t_{k l})} (|(P_N x^*)(\tau_1, \tau_2) -
$$
$$
-(P_N x^*)(\bar t_{k l})| 
    + |\psi (\tau_1, \tau_2) - \psi (\bar t_{k l})|) d \tau_1 d \tau_2 \le
    $$
    $$
    \le AB \displaystyle
    \int\limits_{g_{k l}} {d \tau_1 d \tau_2 \over r^{2-\alpha/2} (\tau, \bar t_{k l})}
    \le     {AB \over N^{\alpha/2}}.  \eqno (3.9)
    $$

So, from (3.7) $-$ (3.9) we have
\\ $||x^* - x_N^*|| \le ||x^* - P_N x^*|| + ||P_N x^* - x^*_N|| \le A N^{- \alpha/2}$.

{\bf Theorem 3.1}  \cite{Boy30}, \cite{Boy31}, \cite{Boy25}. Let the equation (3.1) has a unique solution $x^*(t).$
Let $\varphi (\Theta),$ $f(t)$ are belong to Holder class of functions with
exponent $\alpha.$
Then exist parameters $\bar h_1,$ $\bar h_2,$ such, that the system (3.3) has
a unique solution
$x^*_N(t)$ and the estimate $\|x^*-x^*_N\|_{C(G)} \le AN^{-\alpha/2}$ is valid. 

\begin{center}
{\bf 3.2. Approximate Solution of Linear Multi-Dimensional Singular Integral Equations
	on Sobolev Classes of Functions
	}
\end{center}

Let us consider the two-dimensional singular integral equation
$$
a(t) x(t) + b(t)\int\limits_G \frac{\varphi(\Theta) x(\tau)}
{r^2(t,\tau)}d\tau=f(t),
\eqno (3.10)
$$
where $G=[-A,A]^2.$ $t=(t_1,t_2),$ $\tau=(\tau_1,\tau_2),$
$r(t,\tau)=((t_1-\tau_1)^2+(t_2-\tau_2)^2)^{1/2},$
$\Theta=((t-\tau)/r(t,\tau));$
functions $a, b, f,\varphi$ have derivatives up to $r$ order.

We consider two-dimensional singular integral equations for simplicity.
All received below results can be diffused to multi-dimensional singular
integral equations.

Let us construct the numerical method for solution of the equation (3.1).

Let us cover the domain $G$ by squares
\\ $\Delta_{k l} = [t_k,t_{k+1};t_l,t_{l+1}]$, $k,l = \overline{0,N+1}$,
$t_k = -A + 2A k /(N+2),$ $k = \overline{0,N+2}$.

Together with squares $\Delta_{k l}$, $k,l = \overline{0,N+1}$, we need  
 rectangles
$\bar \Delta_{k l}$, $k,l = \overline{1,N}$, which are defined by following way:
\\ $\bar \Delta_{k l} = \Delta_{k l}$ for $k = \overline{2,N-1}$ and $l = \overline{2,N-1}$;
\\ $\bar \Delta_{1 1} = \Delta_{0 0} \cup \Delta_{0 1} \cup \Delta_{1 0} \cup \Delta_{1 1}$;
\\ $\bar \Delta_{1, N} = \Delta_{0, N} \cup \Delta_{1, N} \cup \Delta_{0, N+1} \cup \Delta_{1, N+1}$;
\\ $\bar \Delta_{N, 1} = \Delta_{N, 0} \cup \Delta_{N+1, 0} \cup \Delta_{N, 1} \cup \Delta_{N+1, 1}$;
\\ $\bar \Delta_{N, N} = \Delta_{N, N} \cup \Delta_{N, N+1} \cup \Delta_{N+1, N} \cup \Delta_{N+1, N+1}$;
\\ $\bar \Delta_{1, l} = \Delta_{0, l} \cup \Delta_{1, l}$ for $l = \overline{2,N-1}$;
\\ $\bar \Delta_{N, l} = \Delta_{N, l} \cup \Delta_{N+1, l}$ for $l = \overline{2,N-1}$;
\\ $\bar \Delta_{k, 1} = \Delta_{k, 0} \cup \Delta_{k, 1}$ for $k = \overline{2,N-1}$;
\\ $\bar \Delta_{k, N} = \Delta_{k, N} \cup \Delta_{k, N+1}$ for $k = \overline{2,N-1}$.

Approximate solution of the equation (3.1) we will seek in the form of local spline
$x_{nn}(t_1,t_2).$

Let us describe the construction of the spline $x_{nn}(t_1,t_2).$

Let us consider a function $f(x_1,x_2),$ which is given on the square $G.$
In the each rectangular $ \bar\Delta_{kl}$ the
function $f(x_1,x_2)$ is approximated by interpolation polynomial,
constructed by following way. Let
$\Delta_{kl} = [t_k, t_{k+1}; t_l, t_{l+1}].$
Let us introduce knots $(t_k^i,t_l^j),$ where
$t_k^i=t_k+(t_{k+1}-t_k)i/(r+1),$ $i=1,2,\ldots,r,$ $k,l=1,\ldots,N.$

The polynomial $L_r(f,\Delta_{kl})$ is introduced by formula
$$
L_r(f,\Delta_{kl})=\sum\limits_{i=1}^{r} \sum\limits_{j=1}^{r}
f(t_k^i,t_l^j)\psi_{ki}(t_1)\psi_{lj}(t_2), \eqno (3.11)
$$
where $\psi_{ki}(t_1)$ and $\psi_{lj}(t_2)$ are fundamental polynomials
constructed on  knots $t_k^i$ and $t_l^j.$

The interpolated polynomial (3.11) is diffused to domains $\bar\Delta_{kl},$
$k,l=1,2,\ldots,N,$ in following way. If $k,l=2,\ldots,N-1,$ then
$L_r(f,\bar\Delta_{kl})=L_r(f,\Delta_{kl}).$ Let us consider the square
$\bar\Delta_{11}.$ We will construct the polynomial $L_r(f,\Delta_{11})$
and diffuse this polynomial on the square $\bar \Delta_{11}.$ In result we have the polynomial
$L_r(f,\bar\Delta_{11}).$ Polynomials $L_s(f,\bar\Delta_{kl}),$ $k,l=1,N,$
we will receive by similar way. Let us consider the rectangular $\bar\Delta_{1,l},$
$l=2,\ldots,N-1.$ In this case we construct the polynomial $L_r(f,\Delta_{1l})$
and diffuse this polynomial on domain $\bar \Delta_{1l}.$ In result we have polynomial
$L_r(f,\bar\Delta_{1l}), l=2,\ldots,N-1.$ Polynomials $L_r(f,\bar\Delta_{N,l}),$
$l=2,\ldots$, $N-1;$ $L_r(f,\bar\Delta_{k,1}),$ $L_r(f,\bar\Delta_{k,N})$ $k=2,\ldots,N-1,$
are constructed by similar way.

Spline $f_{nn}(t_1,t_2)$ consists of polynomials $L_r(f,\bar\Delta_{kl}),$
where $k$, $l=1,\ldots,N.$

Approximate solution of the equation (3.1) we seek as local spline $x_{nn}(t_1,t_2)$
with unknown values $x_{nn}(t_k^i,t_l^j),$ $k,l=1,\ldots,N,$ $i,j=1,\ldots,r.$

To each knots $t^{ij}_{kl} ,\ (t^{ij}_{kl}=(t^i_k,t^j_l) ,\ k,l=1,2,\ldots,N, ,\ i,j=1,2,\ldots,r,$
we put in correspondence the rectangle
$\Delta^{ij}_{kl}=[t^i_k-q_1h,t^i_k+h;
t^j_l-q_2h, t^l_l+h],$ where $h \leq 2A/(rN),$ $q_1,$ $q_2$ are parameters,
values of which will be defined later.

Values $x^{ij}_{kl}=x_{nn}(t^i_k,t^j_l),$ $k,l=1,2,\ldots,N,$ $i,j=1,2,\ldots,r,$
we define from the system of linear algebraic equations
$$
a_{kl}^{ij} x_{kl}^{ij}+b_{kl}^{ij}\int\limits_{\Delta_{kl}^{ij}}
\varphi \left(\frac{\tau_1-t_k^i}{r(\tau,M_{kl}^{ij})},
\frac{\tau_2-t_l^j}{r(\tau,M_{kl}^{ij})}\right)
\frac{x_{nn}(\tau)d\tau}{r^2(\tau,M_{kl}^{ij})}+
$$
$$
+b_{kl}^{ij} \sum^N_{k_1=1}\sum^N_{l_1=1}{}' \int\limits_{\bar \Delta_{kl}}
\varphi \left(\frac{\tau_1-t_k^i}{r(\tau,M_{kl}^{ij})},
\frac{\tau_2-t_l^j}{r(\tau,M_{kl}^{ij})}\right)
\frac{x_{nn}(\tau)d\tau}{r^2(\tau,M_{kl}^{ij})}=
$$
$$
=f_{kl}^{ij}, \eqno (3.12)
$$
$i,j=1,2,\ldots,r,$ $k,l=1,2,\ldots,N,$\\
where
$M_{kl}^{ij}=(t_k^i, t_l^j),$ $a_{kl}^{ij}=a(M_{kl}^{ij}),$
$b_{kl}^{ij}=b(M_{kl}^{ij}),$ $f_{kl}^{ij}=f(M_{kl}^{ij}),$
$i,j
=1,2,\ldots,r,$ $k,l=1,2,\ldots,N,$
the prime in the summation indicate, that $(k_1,l_1) \ne (k+v, l+w),$ 
$v,w=-1,0,1.$

From necessary and sufficient conditions for the existence of singular
integral follow, that the function $\varphi(\Theta),
\Theta=(\frac{\tau_1}{r(0,\tau)},\frac{\tau_2}{r(0,\tau)})$
must be equal to
zero on rays which are radiated from the origin. So, one can select parameters
$q_1,q_2,h$ so way that the integral
$$
\left|\int\int\limits_{\Delta_{kl}^{ij}}
\varphi \left(\frac{\tau_1-t_k^i}{r(\tau,M_{kl}^{ij})},
\frac{\tau_2-t_l^j}{r(\tau,M_{kl}^{ij})}\right)
\frac{\psi_{ki}(\tau_1)\psi_{ij}(\tau_2)}{r^2(\tau,M_{kl}^{ij})}
d\tau_1 d\tau_2\right|
$$
will be arbitrarily large.

Repeating arguments given in the previous item, one can see that conditions of the
Hadamard theorem are fulfilled. So, the system (3.12) has a unique solution.
The convergence of solution $x_{nn}(t_1,t_2)$ of the system (3.12)
to exact solution of
the equation (3.1) is proved by the method similar to the method, which we
have used for polysingular integral equations.

In result we can formulate the following statement.

{\bf Theorem 3.2}  \cite{Boy30}, \cite{Boy31}, \cite{Boy25}.
Let $a(t), b(t)$, $\varphi (t) \in W^{r,r}(M)$, $b(t) \neq 0.$
Let the equation (3.1) has a unique solution $x^*(t)$. Then exist  such parameters
$h, q_1, q_2$, that the system (3.12) has a unique solution
$x^*_{nn}$ and under some additional conditions the estimate $\|x^* - x^*_{nn}\|_C \le
An^{-r} \ln^2n $ is valid.

\begin{center}
{\bf  3.3. Parallel Method for Solution of Multi-Dimensional
	Singular Integral Equations
	}
\end{center}

In this section we investigate numerical methods for solution of
multi-dimensional singular integral equations on parallel computers
with P processors. For determination we assume that offered method will
be realized on computer with MIMD architecture (Multiple \ \ Instruction \  stream/
Multiple \  Date \ \  stream) \ in M.J. Flynn classification.

For simplicity we will describe this method only for 
the equation (3.2).

Approximate solution of the equation (3.2) we seek as piece-constant
function $x_N(t_1,t_2),$ unknown values $x_{kl},$ $k,l=1,2,\ldots,N,$
of  which we define from the system (3.3). The system (3.3)  can be written
as matrix equation
$$
CX=F, \eqno (3.13)
$$
where $C=\{c_{kl}\},$ $k,l=1,2,\ldots,M,$ $M=N^2,$ $X=(x_1,\ldots,x_M)^T,$
$F=(f_1,\ldots,f_M)^T.$

It is easy to see that for the system (3.13) conditions of the Hadamard theorem is valid.

Let $M=PL.$ Let us decompose the matrix $C$ of the system (3.13) on $P^2$
blocks, vectors $X$ and $F$ on $P$ blocks. In this case the
system (3.13) can be written as
$$
\bar C \bar X=\bar F, \eqno (3.14)
$$
where $\bar C=\{B_{ij}\},$ $i,j=1,2,\ldots,P,$ $\bar X=(\bar x_1,\ldots,\bar x_p),$
$\bar F=(\bar f_1,\ldots,\bar f_p),$
$$
B_{ij} = \left(
\begin{array}{ccc}
c_{(i-1)L+1,(j-1)L+1}, \ldots, c_{(i-1)L+1,(j-1)L+L}\\
\cdots\\
c_{(i-1)L+L,(j-1)L+1}, \ldots, c_{(i-1)L+L,(j-1)L+L}\\
\end{array}
\right),
$$
$i,j=1,2,\ldots,P,$
$$
\bar x_i=(x_{(i-1)L+1}, \ldots, x_{(i-1)L+L}), \quad i=1,2,\ldots,P,\\
$$
$$
\bar f_i=(f_{(i-1)L+1}, \ldots, f_{(i-1)L+L}), \quad i=1,2,\ldots,P.
$$

For solving the system (3.14) we use block-iterative method
$$
B_{kk} \bar x_k^{n+1}=f_k-\sum\limits^{P}_{l=1, l\ne k} B_{kl}
\bar x_l^n, \eqno (3.15)
$$
$k=1,2,\ldots,P.$

From conditions, imposed on the matrix C,  one can see that matrixes
$B_{kk},$ $k=1,2,\ldots,P,$ are invertible and norms
$\|B_{kk}^{-1}\|,$ $k=1,2,\ldots,P,$ can be made as small as it is necessary for
the convergence of the iterative process (3.15).

\begin{center}
{\bf Appendix 1}
\end{center}
1. Stability of numerical schemes

In this section we will investigate a criteria of stability of numerical
schemes for solution of singular integral equations.

We investigate stability of the three types of singular integral equations. One
can notice that stability of solutions of numerical schemes for other types of singular integral equations which
was considered in this paper can be investigated by similar way.

\begin{center}
{\bf 1.1. Stability of Solutions of Bisingular Integral Equation}
\end{center}

 Let us consider the bisingular integral equation 
$$
a(t_1,t_2)x(t_1,t_2) + \frac{b(t_1,t_2)}{\pi^2}
\int\limits_{\gamma_1}\int\limits_{\gamma_2}
\frac{x(\tau_1,\tau_2)d\tau_1 d\tau_2}{(\tau_1-t_1)(\tau_2-t_2)}=
f(t_1,t_2).
\eqno (1.1)
$$

Approximate method for solution of the equation (1.1) was given in the section
1 of the chapter 2. Below we use definitions introduced in the section 1 of
the chapter 2. It was shown that one can transform the equation (1.1) into
Riemann  value boundary problem
$$
L\varphi \equiv
$$
$$
\equiv \psi^{++}(t_1,t_2)-G(t_1,t_2)( \psi^{+-}(t_1,t_2)+
\psi^{-+}(t_1,t_2)) + \psi^{--}(t_1,t_2) =
$$
$$
= f^*(t_1,t_2), \quad  f^*(t_1,t_2) = \frac{f(t_1,t_2)}{a(t_1,t_2)+b(t_1,
t_2)}.
\eqno (1.2)
$$

Approximate solution of the equation (1.2) we seek in the form of the
functions
$$
\varphi^{++}_{nn}(t_1,t_2) = \psi^{++}(t_1,t_2)
\sum\limits^n_{k=0}\sum\limits^n_{l=0}\alpha_{kl}t_1^k t_2^l,
$$
$$
\varphi^{+-}_{nn}(t_1,t_2) = \psi^{+-}(t_1,t_2)
\sum\limits^n_{k=0}\sum\limits^{-1}_{l=-n}\alpha_{kl}t_1^k t_2^l,
$$
$$
\varphi^{-+}_{nn}(t_1,t_2) = \psi^{-+}(t_1,t_2)
\sum\limits^{-1}_{k=-n}\sum\limits^n_{l=0}\alpha_{kl}t_1^k t_2^l,
$$
$$
\varphi^{--}_{nn}(t_1,t_2) = \psi^{--}(t_1,t_2)
\sum\limits^{-1}_{k=-n}\sum\limits^{-1}_{l=-n}\alpha_{kl}t_1^k t_2^l.
\eqno (1.3)
$$

Coefficients $\{\alpha_{kl}\},$ $k,l=-n,\ldots,-1,0,1,\ldots,n,$ are defined
from the system of linear algebraic equations
$$
L_n\varphi_{nn} \equiv P_{nn}[L\varphi_{nn}] = P_{nn}[f^*(t_1,t_2)],
\eqno (1.4)
$$
where $P_{nn}$ is the operator from the space $C[\gamma_1 \times \gamma_2]$
onto the set of interpolation polynomials, constructed on knots $v_{kl}=(v_k, v_l),$ $k,l=0,1,\ldots,2n,$  $v_k = \exp\{is_k\},$ $s_k = 2k\pi/(2n+1),$
$k,l=0,1,\ldots,2n.$

Let the operator $L$ is invertible in Holder space $H_\beta$
and the inverse operator $\|L^{-1}\|_\beta $ has the norm $\|L^{-1}\|_\beta = B.$

In the chapter 2 was proved that the operator $L_n$ is invertible under
condition $q=An^\beta(E_{nn}(\psi^{++}) + E_{nn}(\psi^{+-}) + E_{nn}(\psi^{-+})+
E_{nn}(\psi^{--}))\ln^2 n <1.$

 Let the values $a(v_i,v_j),$ $d(v_i,v_j)$
are determined with error $\varepsilon.$ Let $\min_{i,j}|a(v_i,v_j)+d(v_i,v_j)|
\ge d.$ Then the error of determination of values $G(t_i,t_j)$ is not bigger
than $\varepsilon_1 = 2\varepsilon/(d-\varepsilon) \approx 2\varepsilon/d.$

Let us denote by $\tilde a(v_i,v_j),$ $\tilde d(t_i,t_j),$ $\tilde f(v_i,v_j)$
perturbed values of $a(v_i,v_j),  d(v_i,v_j), f(v_i,v_j):
|\tilde a(v_i,v_j) - a(v_i,v_j)| \le \varepsilon,$
$|\tilde d(v_i,v_j) - d(v_i,v_j)| \le \varepsilon,$
$|\tilde f(v_i,v_j) - f(v_i,v_j)| \le \varepsilon,$
$i,j =0,1,\ldots,2n.$

With perturbed values of $a(t_i,t_j),$ $d(v_i,v_j),$ $f(v_i,v_j)$ collocation
method for equation (1.2) has the view
\[
\tilde L_n\varphi_{nn} \equiv
P_{nn}[\varphi^{++}_{nn}(t_1,t_2)-\tilde G(t_1,t_2)( \varphi^{+-}_{nn}(t_1,
t_2)+
\varphi_{nn}^{-+}(t_1,t_2)) +
\]
$$
+ \varphi_{nn}^{--}(t_1,t_2)] 
=P_{nn}[\tilde f^*(t_1,t_2)],
\eqno (1.5)
$$
where $\tilde G(v_i,v_j)=(\tilde a(v_i,v_j)-\tilde b(v_i,v_j))/
(\tilde a(v_i,v_j)+\tilde b(v_i,v_j)],$\\
$\tilde f^*(v_i,v_j)=\tilde f(v_i,v_j)/(\tilde a(v_i,v_j)+
\tilde b(v_i,v_j)),$ $i,j =0,1,\ldots,2n.$

It is easy to see that
$$
|G(v_i,v_j) - \tilde G(v_i,v_j)| \le \varepsilon_1 = 2\varepsilon/(d-\varepsilon),
$$
$$
|f(v_i,v_j) - \tilde f(v_i,v_j)| \le \varepsilon_2 = \varepsilon/(d-\varepsilon).
$$

Using inverse theorems of approximation theory,  one can see that
$$
\|P_{nn}[(G(t_1,t_2)-\tilde G(t_1,t_2))\varphi_{nn}^{\pm \mp}]\|_\beta \le
A\varepsilon_1 n^\beta \ln^2 n \|\varphi_{nn}\|_\beta;
$$
$$
\|P_{nn}[f^*(t_1,t_2)-\tilde f^*(t_1,t_2)]\|_\beta \le
An^\beta\varepsilon_2 \ln^2 n.
$$

Here $\tilde G(t_1,t_2)$ and $\tilde f_1(t_1,t_2)$ are Holder functions, which
are equal to $\tilde G(v_i,v_j)$ and $\tilde f_1(v_i,v_j)$
in  points $v_i,v_j$, $i,j =0,1,\ldots,2n.$

Using Banach Theorem for inverse operator, we see that for such $\varepsilon_1$ and
$\varepsilon_2$ that $A\varepsilon_1 n^\beta \ln^2 n <1,$ $A\varepsilon_2 n^\beta \ln^2 n <1,$
the operator $\tilde L_n$ is continuously invertible and the equation (1.5)
has a unique solution.

So, under conditions $\varepsilon_1 \le An^{-\beta} (\ln^2 n)^{-1}$, $\varepsilon_2 \le An^{-\beta} (\ln^2 n)^{-1},$ the system (1.5) has a unique
solution and collocation method (1.2) is stable.

Under conditions
$\|a(t) - \tilde a(t)\|_C \le \varepsilon,$ $\|b(t) - \tilde b(t)\|_C \le \varepsilon,$
$\|f(t) - \tilde f(t)\|_C \le \varepsilon,$ where
$\varepsilon \le An^{-\beta} (\ln n)^{-1}$, similar statements we have for
collocation method for equation 
$$
a(t)x(t) + \frac{b(t)}{\pi i}
\int\limits_{\gamma}
\frac{x(\tau)d\tau}{\tau-t}= f(t).
$$

The collocation method for singular integral equations was investigated in the
item 2 of the chapter 1.

{\it Note.} We drop compact operator
$$
Tx \equiv \int\limits_{\gamma}h(t,\tau)x(\tau)d\tau,
$$
because stability of collocation method and method of mechanical quadrature
for Fredholm integral equations of the second kind is well known. By this
reason we have  dropped compact operator in the equation (1.1) too.

\begin{center}
{\bf 1.2.  Stability of Solutions of Multi-Dimensional Singular Integral Equations}
\end{center}

Let us consider multi-dimensional singular integral equation 
$$
Kx \equiv a(t)x(t) + b(t)\int\limits_{G} \frac{\varphi(\Theta)x(\tau)}
{(r(t,\tau))^2}d\tau = f(t).
\eqno (1.6)
$$

Mechanical quadrature method for solution of the equation (1.6) was
investigated in the item 3 of the chapter 2. This method in operator form can
be written in the form of the equation (3.3) from the item 3 of the chapter 2.
Below we will use the definitions introduced in the item 3 of the chapter 2.
In the base of the proof of the solvability the equation (3.3) from the
chapter 2 was put the Hadamard Theorem for invertibility of matrix. Using this
theorem one can make the following conclusion.

Let elements of the system (3.3) are given with error $\varepsilon$:
$$
|a(\bar t_{kl}) -\tilde a(\bar t_{kl})| \le \varepsilon, \quad
|f(\bar t_{kl}) -\tilde f(\bar t_{kl})| \le \varepsilon,
$$
$k,l=1,2,\ldots,N;$
$$
|\bar d_{ij}(\bar t_{kl}) - \tilde{\bar  d_{ij}}(\bar t_{kl})| \le \varepsilon, \quad
|d_{ij}(\bar t_{kl}) -\tilde d_{ij}(t_{kl})| \le \varepsilon,
$$
$k,l=1,2,\ldots,N,$ $i,j = 1,2,\ldots,N.$

With values of $\varepsilon$ such, that
$$
|(a(\bar t_{kl}) \pm \varepsilon) + (d_{kl}(\bar t_{kl}) \pm \varepsilon)| >
\sum\limits^N_{i=1}\sum\limits^N_{j=1}{}^*
|\bar d_{ij}(\bar t_{kl}) \pm \varepsilon|,
$$
the system
$$
\tilde a(\bar t_{kl})x_{kl} +
\sum\limits^N_{i=1}\sum\limits^N_{j=1}{}'
\tilde {\bar d_{ij}}(\bar t_{kl}) x_{ij} +
x_{kl}\tilde d_{kl}(\bar t_{kl})= \tilde f(\bar t_{kl}),
\eqno (1.7)
$$
$k,l=1,2,\ldots,N,$\\
has a unique solution $\{\tilde x^*_{kl}\}$ and for $k,l=1,2,\ldots,N,$  estimations
$\| \tilde x^*_{kl} - x^*_{kl}\|_{R_n} \le A\varepsilon$ are valid,
where $\{x^*_{kl}\}$ is a solution of the system (3.3) from of the chapter 2.
Here $\sum\sum^*$ denote that the $\sum$ is taken by $(i,j) \neq (k,l).$

\begin{center}
{\bf 1.3. Stability of Solutions of Nonlinear Singular Integral Equations}
\end{center}

Let us consider stability of solutions of nonlinear
singular integral equation 
$$
Kx \equiv a(t,x(t)) + \frac{1}{\pi i}
\int\limits_{\gamma}
\frac{h(t,\tau,x(\tau))d\tau}{\tau-t}=f(t),
\eqno (1.8)
$$
where $f(t) \in W^r(M),$ $a(t,u) \in W^{rr}(M),$ $a(t,\tau,u) \in W^{rrr}(M).$

Collocation method and method of mechanical quadrature for the equation (1.8)
was investigated in the item 6 of the Chapter 1. Below we will use
designations which was introduced in that item.

Approximate solution of the equation (1.8) we look for in the form of the polynomial
$$
x_n(t)=\sum\limits^n_{k=-n} \alpha_k t^k
\eqno (1.9)
$$
coefficients $\{\alpha_k\}$ of which are defined from the system
\[
K_n x_n \equiv \bar P_n^t \left[a(t,x_n(t)) + \frac{1}{\pi i}
\int\limits_{\gamma} P_n^{\tau}
\left[\frac{h(t,\tau,x_n(\tau))}{\tau-t}\right]d\tau\right] =
\]
$$
= \bar P_n [f(t)].
\eqno (1.10)
$$

In the Theorem 2.6 of the Introduction were given conditions for convergence of
numerical scheme (1.10).

We will investigate stability of solution of the system (1.10) under
conditions, that functions $a(t,\tau),$ $h(t,\tau,u),$ $f(t)$ were given with
error:
$|a(t,\tau) - \tilde a(t,\tau)| \le \varepsilon,$
$|h(t,\tau,u) - \tilde h(t,\tau,u)| \le \varepsilon,$
$|f(t)-\tilde f(t)| \le \varepsilon.$

So, instead of the system (1.10), we have the system
\[
\tilde K_n x_n \equiv \bar P_n^t \left[\tilde a(t,x_n(t)) + \frac{1}{\pi i}
\int\limits_{\gamma} P_n^{\tau}
\left[\frac{\tilde h(t,\tau,x_n(\tau))}{\tau-t}\right]d\tau\right] =
\]
$$
= \bar P_n^t [\tilde f(t)].
\eqno (1.11)
$$

Let us introduce the functions
$$
a_\varepsilon(t,x_n(t)) = \bar P_n^t [\tilde a(t,x_n(t)) - a(t,x_n(t))],
$$
$$
h_\varepsilon(t,\tau,x_n(\tau)) = \bar P_n^t P_n^{\tau}
[\tilde h(t,\tau,x_n(\tau)) - h(t,\tau,x_n(\tau))],
$$
$$
f_\varepsilon(t) = \bar P_n^t[\tilde f(t) - f(t)].
$$

From Bernstein inequality and from Bernstein inverse theorems one can see, that
$$
\|a_\varepsilon(t,x_n(t))\|_\beta \le A\varepsilon n^{\beta}\ln n,
$$
$$
\|h_\varepsilon(t,\tau,x_n(\tau))\|_\beta \le A\varepsilon n^{\beta}\ln^2 n,
$$
$$
\|f_\varepsilon(t)\|_\beta \le A\varepsilon n^{\beta}\ln n,
$$
$$
\|a'_{\varepsilon,2}(t,x_n(t))\|_\beta \le A\varepsilon n^{1+\beta}\ln n,
$$
$$
\|h'_{\varepsilon,3}(t,\tau,x_n(\tau))\|_\beta \le A\varepsilon n^{1+\beta}\ln^2
n.
$$

Frechet derivatives of operators $K_n$ and $\tilde K_n$ have the forms
\[
K'_n(x_n^0)x_n \equiv 
\]
\[
 \equiv  \bar P_n
\left[a'_2(t,x_n^0(t))x_n(t) + \frac{1}{\pi i}
\int\limits_{\gamma} P_n
\left[\frac{h'_3(t,\tau,x_n^0(\tau))x_n(\tau)}{\tau-t}\right]d\tau\right] =
\]
$$
= \bar P_n [f(t)] \eqno (1.12)
$$
and
\[
\tilde K'_n(x_n^0)x_n \equiv
\]
\[
\equiv \bar P_n
\left[\tilde a'_2(t,x_n^0(t))x_n(t) + \frac{1}{\pi i}
\int\limits_{\gamma} P_n
\left[\frac{\tilde h'_3(t,\tau,
x_n^0(\tau))x_n(\tau)}{\tau-t}\right]d\tau\right] =
\]
$$
= \bar P_n [\tilde f(t)] \eqno (1.13)
$$

Repeating the arguments from the item 6 of the Chapter 1 one can see, that
$$
q_1 = \|\tilde K'_n(x_n^0)x_n - K'_n(x_n^0)x_n\|_\beta \le
A \varepsilon n^{1+\beta}\ln n\|x_n\|_\beta.
$$

In the item 6 of the Chapter 1 was proved, that, for such $n$  that
$$
q = A\|[K'(x_0)]^{-1}\|_\beta n^{-r+\beta} \ln n < 1,
$$
the operator $K'_n(x_0)$ is continuously invertible in the space
$H_\beta$ and
$\|[K'_n(x_0)]^{-1}\|_\beta \le \|K'(x_0)\|_\beta  /(1-q).$

From this and Banach Theorem of invertible operators follows, that, if
$\|[\tilde K'_n(x_0)]^{-1}\|_\beta q_1 < 1,$ the operator
$\tilde K'_n(x_0)$ is  invertible  and the estimate
$$
\|[\tilde K'_n(x_0)]^{-1}\|_\beta \le \|[K'_n(x_0)]^{-1}\|_\beta  /(1-q_1)
$$
is valid.

So, the equation (1.11) has a unique solution, which we design as
$\tilde x^*_n.$ Let us estimate the norm $\|x^*_n(t) - \tilde x^*_n(t)\|_\beta,$
where $x^*_n(t)$ is a unique solution of the equation (1.10).
For this aim we use the Newton $-$ Kantorovich method for solution of the equation
(1.11):
$$
\tilde x_n^{m+1}(t) =\tilde x_n^m(t) - [\tilde K'_n(\tilde x_n^0)]^{-1}
(\tilde K_n(\tilde x_n^m) - \tilde
f_n(t)), \quad m=0,1,\ldots,
\eqno (1.14)
$$
where $\tilde f_n(t) = \bar P_n[\tilde f(t)],$
$\tilde x^0_n(t) = x^*_n(t).$

It is easy to see that
\[
\|\tilde x^1_n(t) - \tilde x^0_n\|_\beta \le
\]
$$
\le  \|[\tilde K'(\tilde x_n^0)]^{-1}\|_\beta
\|(\tilde K_n(\tilde x_n^0) - K_n(\tilde x_n^0)) - (\tilde f_n(t) - f_n(t))\|_\beta \le
$$
$$
\le
A\varepsilon n^\beta \ln n.
$$

So, $\eta_0 = A\varepsilon n^\beta \ln n$ and, repeating arguments from the item 6
of the Chapter 1, we see that, under conditions which was described in that item,
the estimate $\|\tilde x_n^*(t) - x_n^*(t)\|_\beta \le A\varepsilon n^\beta \ln^2 n$ is valid.

From given above arguments we can make conclusion that, for so small $\varepsilon$
 that $A\varepsilon n^\beta \ln^2 n < 1$, the
method of mechanical quadrature is
stable in the space $H_\beta.$

\end{document}